\documentclass[dvips,12pt]{article}
\usepackage{amsmath,amssymb,amsfonts,amscd,theorem,pstricks,pst-node,epsfig}
\usepackage[latin1]{inputenc}
\def\simarrow{\mathrel{\raise -0.5mm\hbox{$\sim$}}\hspace{-1.8mm}{\rightarrow} } 

\def\bsimarrow{\leftarrow\hspace{-0.7mm}\mathrel{\raise -0.5mm\hbox{$\backsim$}} }

\def\bt{\begin{tabular}}
\def\te{\end{tabular}}

\def\Abstract{\small
\begin{center}
{\bf Abstract\vspace{-.5em}\vspace{0pt}}
\end{center}
\quotation
}

\def\lettrine#1#2#3{\noindent\hangindent#1\hangafter-#2
\hskip-#1\smash{\hbox to #1{#3\hfill}}\ignorespaces}

\newcommand{\To}[1]{\mathop{\to}\limits_{#1}}

\def\BM{\begin{pmatrix}}
\def\EM{\end{pmatrix}}

\def\ds{\displaystyle}

\def\d=f{\buildrel\hbox{\scriptsize d\'{e}f}\over \Longleftrightarrow}

\def\cit{\text{\it I\hskip -6ptC\/}}

\def\square{\hfill\hbox{\vrule height .9ex width .8ex depth -.1ex}}
\def\rit{\text{\it I\hskip -2pt  R}}

\def\zit{\text{\it Z\hskip -4pt  Z}}

\def\nit{\text{\it I\hskip -2pt  N}}

\def\Bd{{\text B}}

\def\Ed{{\text E}}
\def\Es{{\cal E}}

\def\Fs{{\cal F}}

\def\be{\begin{equation}}
\def\ee{\end{equation}}
\def\beqn{\begin{eqnarray}}
\def\eeqn{\end{eqnarray}}
\def\nobeqn{\begin{eqnarray*}}
\def\noeeqn{\end{eqnarray*}}
\def\ba{\left(\begin{array}}
\def\ea{\end{array} \right) }

\def\bpr{\paragraph{Proof.}}

\def\epr{\square\vskip 6pt}
\def\eop{\hbox{\vrule height .9ex width .8ex depth -.1ex}}

\def\o{\overline}
\def\and{\; \mbox{and} \;}

\newcommand{\half}{\frac{1}{2}}

\def\hfl#1#2{\smash{\mathop{\hbox to 12mm{\rightarrowfill}}
\limits^{\scriptstyle #1}_{\scriptstyle #2}}}

\def\Be{\begin{enumerate}}
\def\Ee{\end{enumerate}}

\def\Bena{\begin{enumerate}
\def\labelenumi{\theenumi)}
\def\theenumi{\arabic{enumi}}
\def\labelenumii{\theenumii)}
\def\theenumii{\alph{enumii}}}

\def\Bean{\begin{enumerate}
\def\labelenumii{\theenumii)}
\def\theenumii{\arabic{enumii}}
\def\labelenumi{\theenumi)}
\def\theenumi{\alph{enumi}}}

\def\Bero{\begin{enumerate}
\def\labelenumii{\theenumii)}
\def\theenumii{\arabic{enumii}}
\def\labelenumi{(\theenumi)}
\def\theenumi{\roman{enumi}}}

\def\BeRo{\begin{enumerate}
\def\labelenumii{\theenumii)}
\def\theenumii{\arabic{enumii}}
\def\labelenumi{(\theenumi)}
\def\theenumi{\Roman{enumi}}}

\def\Bi{\vskip 11pt\begin{itemize}\itemsep=18pt}
\def\bi{\begin{itemize}}
\def\Ei{\end{itemize}\vskip 11pt}
\def\ei{\end{itemize}}

\def\Bd{\begin{description}}
\def\Ed{\end{description}}

\def\R{\right}
\def\L{\left}
\def\F{\frac}

\def\prod{\mathop{\Pi}\limits}
\def\sum{\mathop{\Sigma}\limits}
\def\bigoplus{\mathop{\oplus}\limits}
\def\bigcup{\mathop{\cup}\limits}

\def\vec{\overrightarrow}

\def\Gn{G^{(n)}({F_{\o \omega }}\times {F_\omega })}
\def\gn{G^{(n)}(F_{\o \omega }\times F_\omega )}

\def\Pn{P^{(n)}({F^1_{\o \omega }}\times {F^1_\omega })}
\def\Gnv{G^{(n)}({F^+_{\o v}}\times {F^+_v})}
\def\Pnv{P^{(n)}({F^+_{\o v^1}}\times {F^+_{v^1}})}

\def\GLn{\GL_n({F_{\o \omega }}\times {F_\omega })}
\def\GLnv{\GL_n({F^+_{\o v^1}}\times {F^+_{v^1}})}
\def\AaF{{F_{\o \omega }}\times {F_\omega }}
\def\GLnV{\GL_n({F^+_{\o v}}\times {F^+_{v}})}

\def\des{^{(\rm desing)}}

\def\wF{\widetilde F}
\def\wFL{\widetilde F_L}
\def\wFR{\widetilde F_R}

\def\SMOD{\operatorname{SMOD}}
\def\wt{\widetilde}

\def\Abstract{\small
\begin{center}
{\bf Abstract\vspace{-.5em}\vspace{0pt}}
\end{center}
\quotation}

\newcommand{\myendofpart}{\ifodd\value{page}\vfill\eject
\thispagestyle{empty}
\fi}

\def\Bi{\begin{itemize}}

\def\Aa{{\mathbb{A}\,}}
\def\SS{{\mathbb{S}\,}}
\def\ZZ{{\mathbb{Z}\,}}

\def\cit{{\mathbb{C}\,}}
\def\rit{{\mathbb{R}\,}}
\def\NN{{\mathbb{N}\,}}

\def\labelenumi{\alph{enumi})}
\def\theenumi{\alph{enumi}}
\def\labelenumii{\arabic{enumii}.}
\def\theenumii{\arabic{enumii}}

\def\Aut{\operatorname{Aut}}

\def\hol{\operatorname{hol}}
\def\cusp{\operatorname{cusp}}

\def\Gal{\operatorname{Gal}}

\def\Rep{\operatorname{Rep}}

\def\res{\operatorname{res}}
\def\resp#1{(resp. #1~)}
\def\rresp#1{\qquad \mbox{(resp.} \quad #1\ )}

\def\ELLIP{\operatorname{ELLIP}}
\def\Ellip{\operatorname{ELLIP}}

\def\RL{_{R\times L}}
\def\LR{_{L,R}}

\def\GL{\operatorname{GL}}

\def\ELLIP{\operatorname{ELLIP}}
\def\Irr{\operatorname{Irr}}
\def\Int{\operatorname{Int}}

\def\Diff{\operatorname{Diff}}
\def\Out{\operatorname{Out}}
\def\TAN{\operatorname{TAN}}
\def\EIS{\operatorname{EIS}}

{\theorembodyfont{\rmfamily}
  }
{\theorembodyfont{\rmfamily}
  }

{\theorembodyfont{\rmfamily}
  }
{\theorembodyfont{\rmfamily}
  }

{\theorembodyfont{\rmfamily}
  }

\def\lr{left (resp. right) }

\def\bbf{\boldmath\bf}

\def\To{\begin{CD} @>>>\end{CD}}

\begin{document}

\pagestyle{myheadings}

\def\thepage{\arabic{page}}
\markright{}

{{\pagestyle{empty}

\null
 \vfill
 
\begin{center} 
{\LARGE \bbf $n$-dimensional global correspondences of Langlands over singular schemes (II)\par}%
   \vskip 3em
   {\large C. Pierre
\\[5mm]
Institut de Mathématique pure et appliquée\\
Université de Louvain\\
B-1348 Louvain-la-Neuve\\ Belgium\par}

\end{center}

\vfill

\begin{flushright}
\begin{minipage}{.5\linewidth}
``This paper is dedicated to R. Thom who, by his enthusiasm, convinced me of the importance of the singularities and, by his patience, backed me up along my long research towards the blowups of the versal deformations, the geometries of these processes and the (strange) attractors tied up to these''.
\end{minipage}
\end{flushright}
\vfill
\centerline{MSC (2000): 14B05, 14B07, 14E15, 32S30, 11R39, 11F70}
\eject

\null
 \vfill
\begin{center}
{\LARGE \bbf $n$-dimensional global correspondences of Langlands over singular schemes (II)\par}%
   \vskip 3em
   {\large C. Pierre
\par}

\end{center}

\vfill

{\Abstract{\noindent 
A rather complete phenomenology of the singularities is developed according to a new algebraic point of view in the frame of Langlands global correspondences.

That is to say, a process of:
\Bi
\item[-] singularizations and versal deformations of these,
\item[-] singularizations and monodromies of these,
\Ei
is envisaged on all the sections of sheaves of differentiable (bi)functions on (bi)linear algebraic (semi)groups constituting the $n$-dimensional representations of the global Weil groups.

To get the searched holomorphic and cuspidal representations, it is necessary to consider:
\Bi
\item[-] the resolutions of singularities and the blowups of the versal deformations;
\item[-] the resolutions of the singularities in the monodromy cases.
\Ei

Furthermore, the geometry of the versal deformations and of their blowups is studied, as well as the associated dynamics leading to the consideration of  singular hyperbolic attractors and of singular strange attractors.

}} \vfill\eject

}}


\setcounter{page}{1}
\section*{Introduction}
\addcontentsline{toc}{section}{Introduction}

\thispagestyle{empty}

\vskip 11pt

This paper constitutes the second part of the $n$-dimensional global correspondences of Langlands \cite{Lan}: it points out the cases in which these correspondences can still be stated while the sheaf of differentiable (bi)functions on a bilinear algebraic semigroup, constituting the $n$-dimensional representation of the global Weil group, is affected by all kinds of singularities together with deformations and blowups of these.

The phenomenology of the singularities can be split into two sets based on contracting or dilating morphisms characterized respectively by underlying topological subsets getting closer and closer or farther and farther.

In the set characterized by {\bf dilating morphisms\/}, we find:
\Bean
\item {\bf the desingularizations\/}, or the resolutions of the singularities, of a singular scheme consisting in monomializing polynomial ideals by sequences of blowups \cite{Hau}.
\item {\bf the monodromy\/} transformations of a singular scheme arising in an expanding phase in such a way that non-singular fibres can be generated.
\Ee
while we have in the set characterized by {\bf contracting morphisms\/}:
\Bean
\item {\bf the singularizations\/} introduced as the inverse morphisms of the resolutions of singularities, and which are defined by sequences of contracting surjective morphisms producing singular loci.
\item {\bf the versal deformations\/} of these singularities which can be interpreted as extensions of the sequences of contracting surjective morphisms of singularizations, recovering then the classical definition of the versal deformation as (contracting) fibre bundles of which fibres are the bases of the versal deformations.
\item {\bf the blowups of the versal deformations\/}, introduced in \cite{Pie3}, \cite{Pie4} as the inverse morphisms of the versal deformations: they are based upon Galois antiautomorphisms and are also called spreading-out isomorphisms.
\Ee
\vskip 11pt

{\bf The considered mathematical frame\/}, recalled in chapter 1, is the same as the one which was envisaged in the first part \cite{Pie1} of the $n$-dimensional global correspondences of Langlands, that is to say, considering:
\Bi
\item {\bbf sets $F^+_v$ (resp. $F^+_{\o v}$~) of $r$ packets of \lr real pseudo-ramified equivalent completions\/} associated with the \lr (algebraically closed) extension (semi)field of a number field of characteristic zero and characterized by increasing Galois extension degrees being integers modulo $N$~.

\item {\bbf sets $F_\omega $ (resp. $F_{\o\omega }$~) of $r$ packets of \lr complex pseudo-ramified equivalent completions \/} covered by their real equivalents.

\item {\bf bilinear algebraic semigroups\/} $\GL_n(F_{\o\omega }\times F_\omega )\equiv T_n^t(F_{\o\omega })\times T_n(F_\omega )$ and 
$\GL_n(F^+_{\o v }\times F^+_v )\equiv T_n^t(F^+_{\o v })\times T_n(F^+_v )$ respectively over products of complex and real completions in such a way that their representation spaces $G^{(n)}(F_{\o\omega }\times F_\omega )$ and 
 $G^{(n)}(F^+_{\o v }\times F^+_v )$ be $n^2$-dimensional complex and real bilinear affine semigroups \cite{Pie1} decomposing into sets $\{g^{(n)}\RL[j,m_j]\}^r_{j=1,m_j}$ and
 $\{g^{(n)}\RL[j_\delta ,m_{j_\delta }]\}^r_{j_\delta =1,m_{j_\delta }}$ of $r$ packets of complex and real equivalent conjugacy class representatives:
 
 each conjugacy class representative $g^{(n)}\RL[j,m_j]$ is the product of a right $n$-dimensional complex algebraic semitorus $g^{(n)}_R[j,m_j]$ by its left equivalent
 $g^{(n)}_L[j,m_j]$ verifying\linebreak
 $g^{(n)}_R[j,m_j]\simeq (F_{\o\omega _{j,m_j}})^n$ and
 $g^{(n)}_L[j,m_j]\simeq (F_{\omega _{j,m_j}})^n$ where $F_{\o\omega _{j,m_j}}$ (resp.
 $F_{\omega _{j,m_j}}$~) is the $(j,m_j)$-th corresponding complex completion of $F_{\o\omega }$ (resp. $F_\omega $~) and each conjugacy class representative 
 $g^{(n)}\RL[j_\delta ,m_{j_\delta }]$ is the product of a right $n$-dimensional real algebraic semitorus $g^{(n)}_R[j_\delta ,m_{j_\delta }]$ by its left equivalent
 $g^{(n)}_L[j_\delta ,m_{j_\delta }]$ verifying
 $g^{(n)}_R[j_\delta ,m_{j_\delta }]\simeq (F^+_{\o v _{j_\delta ,m_{j_\delta }}})^n$ and
 $g^{(n)}_L[j_\delta ,m_{j_\delta }]\simeq (F^+_{v _{j_\delta ,m_{j_\delta }}})^n$ where $F^+_{\o v _{j_\delta ,m_{j_\delta }}}$ (resp.
 $F^+_{v _{j_\delta ,m_{j_\delta }}}$~) is the $(j_\delta ,m_{j_\delta })$-th corresponding real completion of $F^+_{\o v }$ (resp. $F^+_v $~).
 
 \item {\bbf a (bisemi)sheaf $\theta _{G^{(n)}_R}\times \theta _{G^{(n)}_L}$ of differentiable bifunctions\/} $\phi ^{(n)}_{G_{j_R}}(x_{g_{j_\delta }})\otimes \phi ^{(n)}_{G_{j_L}}(x_{g_{j_\delta }})$ on the conjugacy class representatives $g^{(n)}\RL[j_\delta ,m_{j_\delta }]$ of the real bilinear algebraic semigroup $G^{(n)}(F^+_{\o v }\times F^+_v )$ and a bisemisheaf 
 $\theta^{(\cit)} _{G^{(n)}_R}\times \theta^{(\cit)} _{G^{(n)}_L}$ of complex-valued differentiable bifunctions $\phi  ^{(\cit)} _{G^{(n)}_R}(x_{g_R})\otimes \phi ^{(\cit)} _{G^{(n)}_L}(x_{g_L})$ on the conjugacy class representatives $g^{(n)}\RL[j,m_j]$ of the complex bilinear algebraic semigroup 
 $G^{(n)}(F_{\o\omega }\times F_\omega )$~.
 \Ei
 \vskip 11pt
 
 The most important part of this paper, i.e. chapters 2, 3 and 4, concerns the study of:
 \Bi
 \item[-] degenerate singularities on the sections
 $\phi ^{(n)}_{G_{j_L}}(x_{g_{j_\delta }})$ (resp. $\phi ^{(n)}_{G_{j_R}}(x_{g_{j_\delta }})$~) of the \lr semisheaf $\theta _{G^{(n)}_L}$ (resp. $\theta _{G^{(n)}_R}$~) as resulting from contracting, generally surjective, morphisms.
 \item[-] the existence of Langlands global correspondences in the ``singular'' context.
 \Ei
 \vskip 11pt
 
 First of all, a {\bf process of singularization\/} is introduced in chapter 2, as being the inverse of monoidal transformations.  It consists in projecting a sequence of normal crossings divisors defined on irreducible completions of rank $N$ onto a singular locus which then becomes the homotopic image of these normal crossings divisors under a sequence of contracting surjective morphisms.
 
 For example, a singularization of type $A_k$~, given by the germs $y_L=x_L^{k+1}$ (resp. $y_R=x_R^{k+1}$~) of differentiable functions 
 $\phi ^{(n)}_{G_{j_L}}(x_{g_{j_\delta }})$ (resp. $\phi ^{(n)}_{G_{j_R}}(x_{g_{j_\delta }})$~) on the $j_\delta $-th conjugacy class of $G^{(n)}(F^+_v)$ (resp. $G^{(n)}(F^+_{\o v})$~), is generated by a sequence of $k+1$ contracting surjective morphisms, being in fact fibre bundles whose contracting fibres are the homotopic images of the normal crossing divisors.
 
 Outside of the singular locus, the contracting morphism of singularization is an isomorphism.  Due to the relative small topological space on which is defined the semisheaf 
 $\theta _{G^{(n)}_L}$ (resp. $\theta _{G^{(n)}_R}$~), it will be assumed that on every function
 $\phi ^{(n)}_{G_{j_L}}(x_{g_{j_\delta }})\in \theta _{G_L^{(n)}}$ (resp. cofunction $\phi ^{(n)}_{G_{j_R}}(x_{g_{j_\delta }})\in \theta _{G_R^{(n)}}$~) a same kind of singularity (or set of singularities) is generated by the singularization process.
 
 So, the singularization of $\theta _{G_L^{(n)}}$ (resp. $\theta _{G_R^{(n)}}$~) is a process transforming it into a singular semisheaf $\theta^* _{G_L^{(n)}}$ (resp. $\theta ^*_{G_R^{(n)}}$~) whose sections $\phi ^{*(n)}_{G_{j_{\delta _L}}}(x_{g_{j_\delta }})$ (resp.
 $\phi ^{*(n)}_{G_{j_{\delta _R}}}(x_{g_{j_\delta }})$~) are the differentiable functions
 $\phi ^{(n)}_{G_{j_{\delta _L}}}(x_{g_{j_\delta }})$ (resp. $\phi ^{(n)}_{G_{j_{\delta _R}}}(x_{g_{j_\delta }})$~) endowed with germs $\phi _{j_\delta }(\omega _L)$ (resp.
 $\phi _{j_\delta }(\omega _R)$~) having (degenerate) singularities assumed to be of corank 1 (to simplify the handling).
 \vskip 11pt
 
 {\bf The versal deformation\/} $\theta^{\rm vers} _{G_L^{(n)}}= \theta^* _{G_L^{(n)}}\times \theta _{S_L}$ (resp. $\theta^{\rm vers} _{G_R^{(n)}}= \theta^* _{G_R^{(n)}}\times \theta _{S_R}$~) of the singular semisheaf $\theta^* _{G_L^{(n)}}$ (resp. $\theta^* _{G_R^{(n)}}$~), whose sections are endowed with germs having degenerate singularities (of corank 1), can be interpreted as the total space of a fibre bundle $D_{S_L}$ (resp. $D_{S_R}$~) of which fibre $\theta _{S_L}$ (resp. $\theta _{S_R}$~) is the family of the (semi)sheaves of the base $S_L$ (resp. $S_R$~) of the versal deformation.
 
 In this context, the versal deformation consists in an extension of the singularization process in the sense that it is generated by a sequence of contracting morphisms extending the sequence of contracting surjective morphisms of singularizations by projecting sets of normal crossing divisors in the neighbourhoods of the singular loci according to the finite determinacies of the considered degenerate singularities on the sections of 
 $\theta^* _{G_L^{(n)}}$ (resp. $\theta ^*_{G_R^{(n)}}$~). If we refer to the degenerate germ $y_L=x_L^{k+1}$ (resp. $y_R=x_k^{k+1}$~) of type $A_k$~, its versal deformation will then result from a sequence of $(k-2)$ contracting morphisms extending the sequence of contracting surjective morphisms of singularization in $(k-2)$ dimensions in such a way that a sequence of $(k-2)$ (sets of) normal crossings divisors be projected in the neighbourhood of the singular locus.
 \vskip 11pt
 
 This constitutes the content of chapter 2, section 2, while section 3 deals with the {\bf geometry of the versal deformation\/}:
 
 It is proved that:
 \Be
 \item the geometry is hyperbolic in the neighbourhood of the singular locus of a not unfolded degenerate singular germ of corank $m\le 3$ and multiplicity $i$~, $1\le i\le n$~, in the sense that:
 \Bi
 \item the limit set of the Kleinian group acting in the neighbourhood of a singular locus corresponds precisely to this singular locus.
 
 \item the ordinary set of the Kleinian group can be associated with the neighbourhood of the singular locus and is characterized by a hyperbolic metric.
 \Ei
 
 \item the neighbourhood of the unfolded germ on the section 
 $\phi ^{(n)}_{G_{j_{\delta _L}}}(x_{g_{j_ \delta }})$ (resp.
 $\phi ^{(n)}_{G_{j_{\delta _R}}}(x_{g_{j_ \delta }})$~) of the unfolded semisheaf 
 $\theta^{\rm vers} _{G_L^{(n)}}$ (resp. $\theta^{\rm vers} _{G_R^{(n)}}$~) is characterized by a spherical geometry except in the neighbourhood of the singular locus where the geometry is  hyperbolic.
 \Ee
 \vskip 11pt
 
 Chapter 3 envisages the blowup of the versal deformation as well as the study of the strange attractors related to the versal deformations of singular germs.
 
 {\bf The blowup of the versal deformation\/}
 \[ \theta^{\rm vers} _{G_L^{(n)}}= \theta^* _{G_L^{(n)}}\times \theta _{S_L}
 \qquad \text{(resp.} \quad
 \theta^{\rm vers} _{G_R^{(n)}}= \theta^* _{G_R^{(n)}}\times \theta _{S_R}\ )\]
 of the singular semisheaf
 $\theta^* _{G_L^{(n)}}$ (resp. $ \theta^* _{G_R^{(n)}}$~) consists essentially in an algebraic endomorphism $\Pi _{S_L}$ (resp. $\Pi _{S_R}$~), based on Galois antiautomorphisms, pulling out partially or completely the sheaves $\theta ^1(\omega ^i_L)$ (resp. $\theta ^1(\omega ^i_R)$~), $1\le i\le s$~, of the fibre (corank 1 case)
 \begin{align*}
 \theta _{S_L} &= \{ \theta ^1(\omega ^1_L),\cdots, \theta ^1(\omega ^i_L),\cdots,\theta ^1(\omega ^s_L)\}\\
 \text{(resp.} \quad 
  \theta _{S_R} &= \{ \theta ^1(\omega ^1_R),\cdots, \theta ^1(\omega ^i_R),\cdots,\theta ^1(\omega ^s_R)\}\ )\end{align*}
  of the versal deformation from the $(n-1)$-dimensional coefficient sheaf
  \begin{align*}
  \theta _L(a) &= \{\theta _L^{n-1}(a_1) , \cdots,\theta _L^{n-1}(a_i),\cdots,
  \theta _L^{n-1}(a_s)\}\in \theta _{G_L^{(n)}}\\
  \text{(resp.} \quad
  \theta _R(a) &= \{\theta _R^{n-1}(a_1) , \cdots,\theta _R^{n-1}(a_i),\cdots,
  \theta _R^{n-1}(a_s)\}\in \theta _{G_R^{(n)}}\ )\end{align*}
  on which $\theta _{S_L}$ (resp. $\theta _{S_R}$~) was projected.
  
  This blowup is maximal when all the base (semi)sheaves of $\theta _{S_L}$ (resp. $\theta _{S_R}$~) have been pulled out from $\theta _L(a)$ (resp. $\theta _R(a)$~).
  
  The blowup is complete if it is given by the composition of maps
  \[ (S\circ T)_L=(\tau _{\vee _{\omega _L}}\circ \Pi _{S_L}) \qquad \text{(resp.} \quad
  (S\circ T)_R=(\tau _{\vee _{\omega _R}}\circ \Pi _{S_R})\ )\]
  where $\tau _{V_{\omega _L}}$ (resp. $\tau _{V_{\omega _R}}$~) is the projective map of the tangent bundle projecting all the disconnected base (semi)sheaves $\theta ^1_I(\omega ^i_L)$ (resp.
  $\theta ^1_I(\omega ^i_R)$~) of $\theta _{S_L}$ (resp. $\theta _{S_R}$~) in the vertical tangent spaces: this blowup then constitutes an extension of the quotient algebra of the versal deformation of the singular semisheaf $\theta^* _{G_L^{(n)}}$ (resp. $\theta^* _{G_R^{(n)}}$~). When it is maximal, it is given by the map:
  \begin{align*}
  (S\circ T)_L^{\max}: \theta^* _{G_L^{(n)}}\times \theta _{S_L}
  & \To \theta^* _{G_L^{(n)}}\cup \theta _{S_L}\\
  \text{(resp.} \quad
    (S\circ T)_R^{\max}: \theta^* _{G_R^{(n)}}\times \theta _{S_R}
  & \To \theta^* _{G_R^{(n)}}\cup \theta _{S_R}\ )\end{align*}
  and corresponds to the inverse of the versal deformation $D_{S_L}$ (resp. $D_{S_R}$~) according to:
  \[ (S\circ T)_L^{\max}=(D_{S_L})^{-1} \qquad \text{(resp.} \quad
   (S\circ T)_R^{\max}=(D_{S_R})^{-1}\ ).\]
   The family of disconnected base semisheaves are then glued together and cover partially the singular semisheaf $\theta^* _{G_L^{(n)}}$ (resp. $\theta^* _{G_R^{(n)}}$~): they are labelled 
   $\theta^* _{(S\circ T)(1)_L}$ (resp. $\theta^* _{(S\circ T)(1)_R}$~) and verify 
   $\theta^* _{(S\circ T)(1)_L} \simeq \theta _{S_L}$ (resp. $\theta^* _{(S\circ T)(1)_R} \simeq \theta _{S_R}$~).  But, $\theta^* _{(S\circ T)(1)_L} $ (resp. $\theta^* _{(S\circ T)(1)_R} $~) can be affected by singularities on its sections involving versal deformations and blowups.
   
   Section 3.2 envisages the versal deformation and its blowup from a differentiable and dynamical point of view.

   {\bf The dynamics\/} is envisaged around singularities on the sections of the tangent bundle on the conjugacy class representatives of the algebraic semigroup $G_L^{(n)}(F^+_v)\simeq T_n(F^+_v)$ (resp. $G_R^{(n)}(F^+_{\o v})\simeq T^t_n(F^+_{\o v})$~).  Then, the neighbourhood 
of the singular germ $\phi _{j_\delta }(\omega _L)$
\resp{$\phi _{j_\delta }(\omega _R)$}
on the $n$-dimensional real-valued differentiable function $\phi ^{\rm TAN}_{G_{j_{\delta _L}}}(x^{\rm TAN}_{g_{j_\delta }})$ (resp. $\phi ^{\rm TAN}_{G_{j_{\delta _R}}}(x^{\rm TAN}_{g_{j_\delta }})$~) of the space of sections $\Gamma (T(G_L^{(n)}(F^+_v)))$   (resp. $\Gamma (T(G_R^{(n)}(F^+_{\o v})))$~) of the tangent bundle on $G_L^{(n)}(F^+_v)$ (resp. $G_R^{(n)}(F^+_{\o v})$~) is a {\bf singular hyperbolic attractor\/} $\Lambda _L^{\rm TAN}$ (resp. $\Lambda _R^{\rm TAN}$~) with respect to the diffeomorphisms $\Diff_L(T(G_L^{(n)}(F^+_v)))$ (resp. $\Diff_R(T(G_R^{(n)}(F^+_{\o v})))$~).
   
   And, the versal unfolding of the germ $\phi _{j_\delta }(\omega _L)$ (resp. $\phi _{j_\delta }(\omega _R)$~) involves the map:
 \begin{align*}  
\vee D_{\Lambda _L} : \quad \Lambda _L^{\rm TAN} &\To \Lambda ^{\rm TAN}_{{\rm str}_L}\\
\text{(resp.} \quad 
\vee D_{\Lambda _R} : \quad \Lambda _R^{\rm TAN} &\To \Lambda ^{\rm TAN}_{{\rm str}_R}\ )\end{align*}
of the singular hyperbolic attractor $\Lambda ^{\rm TAN}_L$ (resp. $\Lambda ^{\rm TAN}_R$~) into the
{\bf singular strange attractor\/} 
\begin{align*}
\Lambda ^{\rm TAN}_{{\rm str}_L}&=\Lambda ^{\rm TAN}_L\times \Lambda ^{\rm TAN}_{{\rm unf}_L}\\
\text{(resp.} \quad
\Lambda ^{\rm TAN}_{{\rm str}_R}&=\Lambda ^{\rm TAN}_R\times \Lambda ^{\rm TAN}_{{\rm unf}_R}\ )
\end{align*}
where $\Lambda ^{\rm TAN}_{{\rm unf}_L}$ (resp. $\Lambda ^{\rm TAN}_{{\rm unf}_R}$~) is an unfolded attractor which can be expressed according to:
\[ \Lambda ^{\rm TAN}_{{\rm unf}_L}=\cup \Lambda ^{\rm TAN}_{\omega ^i_{j_{\delta _L}}}
\qquad \text{(resp.} \quad 
\Lambda ^{\rm TAN}_{{\rm unf}_R}=\cup \Lambda ^{\rm TAN}_{\omega ^i_{j_{\delta _R}}}\ )\]
with $\Lambda ^{\rm TAN}_{\omega ^i_{j_{\delta _L}}}$ (resp. $\Lambda ^{\rm TAN}_{\omega ^i_{j_{\delta _R}}}$~) a singular hyperbolic attractor resulting from a singularity on the generator
$\omega ^i_{j_{\delta _L}}$ (resp. $\omega ^i_{j_{\delta _R}}$~) of the versal deformation of
$\phi _{\delta _j}(\omega _L)$ (resp. $\phi _{\delta _j}(\omega _R)$~).

Finally, a blowup of the singular strange attractor $\Lambda ^{\rm TAN}_{{\rm str}_L}$ (resp.
$\Lambda ^{\rm TAN}_{{\rm str}_R}$~) can disconnect the singular hyperbolic attractors
$\Lambda ^{\rm TAN}_{\omega _{j_{\delta _L}}}$ (resp. $\Lambda ^{\rm TAN}_{\omega _{j_{\delta _R}}}$~) from the basic singular hyperbolic attractor $\Lambda ^{\rm TAN}_L$ (resp. $\Lambda ^{\rm TAN}_r$~).
\vskip 11pt

In chapter 4, it is analysed in what extend it is possible to develop {\bf global correspondences of Langlands\/} for a bisemisheaf of differentiable functions on the real algebraic bilinear semigroup $(G^{(n)}(F^+_{\o v}\times F^+_v)$ affected by degenerate singularities.

Recall that a global correspondence consists in a bijection between the $n$-dimensional irreducible representation $\Irr\Rep^{(n)}_{W _{F^+\RL}}(W^{ab}_{F^+_R}\times W^{ab}_{F^+_L})$ of the product, right by left, of Weil groups and the irreducible cuspidal representation
$\Irr\Ellip(\GL_n(\Aa_{F^{+,T}_{\o v}}\times\Aa_{F^{+,T}_{v}}))$ of
$\GL_n(F^+_{\o v}\times F^+_v)$ as developed in \cite{Pie1}.

Now, $\Irr\Rep^{(n)}_{W _{F^+\RL}}(W^{ab}_{F^+_R}\times W^{ab}_{F^+_L})$ is given by the bilinear affine semigroup  $G^{(n)}(F^+_{\o v}\times F^+_v)$ or by the semisheaf $(\theta _{G^{(n)}_R}\otimes \theta _{G^{(n)}_L})$ on it.

But, {\bbf under singularization, versal deformation and blowup of it, the bisemisheaf
$\theta _{G^{(n)}_R}\otimes \theta _{G^{(n)}_L}$ has been transformed into:\/}
\begin{multline*}
\theta _{G^{(n)}_R} \otimes \theta _{G^{(n)}_L}
\begin{CD} @>{\o\rho _{G_R}\times \o\rho _{G_L}}>>\end{CD}
\theta ^*_{G^{(n)}_R} \otimes \theta^* _{G^{(n)}_L}
\begin{CD} @>{D_{S_R}\times D_{S_L}}>>\end{CD}
(\theta ^*_{G^{(n)}_R} \times \theta _{S_R})
\otimes
(\theta ^*_{G^{(n)}_L} \times \theta _{S_L})\\
\begin{CD} @>{(S\circ T)_R^{\max}\times (S\circ T)_L^{\max}}>>\end{CD}
(\theta ^*_{G^{(n)}_R} \cup \theta ^*_{(S\circ T)(1)_R}) \otimes
(\theta ^*_{G^{(n)}_L} \cup \theta ^*_{(S\circ T)(1)_L})\end{multline*}
where
\Bi
\item[-] $\o\rho _{G_R}\times \o\rho _{G_L}$ is the contracting morphism of singularization.
\item[-] $D_{S_R}\times D_{S_L}$ is the contracting morphism of versal deformation.
\item[-] $(S\circ T)_R^{\max}\times (S\circ T)_L^{\max}$ is the blowup of the versal deformation.
\Ei
So, $(\theta _{G^{(n)}_R}\otimes \theta _{G^{(n)}_L})$ has generated under
$((S\circ T)_R^{\max} \circ D_{S_R}\circ \o\rho _{G_R})\times
((S\circ T)_L^{\max} \circ D_{S_L}\circ \o\rho _{G_L})$ the singular bisemisheaf
$(\theta ^*_{G^{(n)}_R}\otimes \theta ^*_{G^{(n)}_L})$ and the singular compactified base bisemisheaf
$(\theta ^*_{(S\circ T)(1)_R}\otimes \theta ^*_{(S\circ T)(1)_L})$ of the blowup of the versal deformation.

But, these bisemisheaves $(\theta ^*_{G^{(n)}_R} \otimes \theta^* _{G^{(n)}_L})$
and $(\theta ^*_{(S\circ T)(1)_R}\otimes \theta ^*_{(S\circ T)(1)_L})$~, affected by singularities, cannot be endowed with a cuspidal representation.

To reach this objective, it is necessary to:
\Be
\item[1)] desingularize those bisemisheaves.
\item[2)] submit them to a toroidal compactification.
\Ee
The desingularization corresponds to the classical monoidal transformations and is reached by a set of inverse morphisms of those defining   a singularization as developed in section 2.1.

Before considering the cuspidal representations of these bisemisheaves
$(\theta^* _{G^{(n)}_R} \otimes \theta^* _{G^{(n)}_L})$
and $(\theta ^*_{(S\circ T)(1)_R}\otimes \theta ^*_{(S\circ T)(1)_L})$~, we can at this stage envisage holomorphic representations of the corresponding desingularized bisemisheaves
$(\theta _{G^{(n)}_R} \otimes \theta _{G^{(n)}_L})$
and $(\theta _{(S\circ T)(1)_R}\otimes \theta _{(S\circ T)(1)_L})$~.  We shall briefly recall how to get a holomorphic representation for $(\theta _{G^{(n)}_R} \otimes \theta _{G^{(n)}_L})$~, taking into account that the same procedure can be applied to
 $(\theta _{(S\circ T)(1)_R}\otimes \theta _{(S\circ T)(1)_L})$~.
 
 The {\bf global holomorphic representation\/} $\Irr\hol^{(n)}(\theta _{G^{(n)}_R} \otimes \theta _{G^{(n)}_L})$ of the bisemisheaf $(\theta _{G^{(n)}_R} \otimes \theta _{G^{(n)}_L})$ is given by the morphism:
 \[\Irr\hol^{(n)}_{\theta _{G\RL}}: \quad \theta _{G^{(n)}_R} \otimes \theta _{G^{(n)}_L}
 \To f_{\o v}(z^*)\otimes f_v(z)\]
 where $f_{\o v}(z^*)\otimes f_v(z)$ is the holomorphic bifunction (i.e. product of a holomorphic function by the corresponding symmetric cofunction) obtained by gluing together and adding the bisections of the bisemisheaf $\theta _{G^{(n)}_R} \otimes \theta _{G^{(n)}_L}$~.
 
 So, in a few words, a singular bisemisheaf $(\theta ^*_{G^{(n)}_R} \otimes \theta ^*_{G^{(n)}_L})$~, submitted to a versal deformation transforming it into
 $(\theta ^*_{G^{(n)}_R} \times \theta _{S_R}) \otimes (\theta ^*_{G^{(n)}_L} \times \theta _{S_L})$~, can be endowed with a holomorphic representation if a blowup of the versal deformation is considered as well as a desingularization of the resulting singular bisemisheaf 
 $(\theta ^*_{G^{(n)}_R} \otimes \theta^*_{G^{(n)}_L})$~.
 
 To get {\bf a cuspidal representation of the desingularized bisemisheaves\/} $(\theta _{G^{(n)}_R} \otimes \theta _{G^{(n)}_L})$ and $(\theta _{(S\circ T)(1)_R}\otimes \theta _{(S\circ T)(1)_L})$~, a toroidal compactification of the bilinear algebraic semigroups
 $G^{(n)}(F^+_{\o v}\times F^+_v)$ and  $G^{(n)}(F^+_{\o v_{\rm cov}}\times F^+_{v_{\rm cov}})$ on which they are defined must be performed in such a way that the products, right by left, of their corresponding conjugacy class representatives be products, right by left, of $n$-dimensional real semitori.
 
 The bisemisheaves on the toroidal bilinear algebraic semigroups
  $G^{(n)}(F^{+,T}_{\o v}\times F^{+,T}_v)$ and $G^{(n)}(F^{+,T}_{\o v_{\rm cov}}\times F^{+,T}_{v_{\rm cov}})$ will be written $(\theta _{G^{(n)}_{T_R}}\otimes \theta _{G^{(n)}_{T_L}})$ and
  $ (\theta^{\rm cov} _{G^{(n)}_{T_R}}\otimes \theta ^{\rm cov}_{G^{(n)}_{T_L}})$~.
  
  Remark that the toroidal compactifications of the bisemisheaves 
  $(\theta _{G^{(n)}_R} \otimes \theta _{G^{(n)}_L})$ and 
  $(\theta _{(S\circ T)(1)_R}\otimes \theta _{(S\circ T)(1)_L})$ are such that their holomorphic representations are transformed into cuspidal representations according to:
  \[\begin{array}{ccccc}
  \Irr\hol(\theta _{G^{(n)}_R} \otimes \theta _{G^{(n)}_L})&: \quad&
  \theta _{G^{(n)}_R}\otimes \theta _{G^{(n)}_L} &\To&
  f_{\o v}(z^*)\otimes f_v(z)\\
  \downarrow &&\downarrow &&\downarrow \\
    \Irr\Ellip(\theta _{G^{(n)}_{T_R}} \otimes \theta _{G^{(n)}_{T_L}})&: \quad&
  \theta _{G^{(n)}_{T_R}} \otimes \theta _{G^{(n)}_{T_L}} &\To&
\Ellip\RL(n,j_\delta ,m_{j_\delta })\end{array}\]
where $\Ellip\RL(n,j_\delta ,m_{j_\delta })= \Ellip_R(n,j_\delta ,m_{j_\delta })\otimes
\Ellip_L(n,j_\delta ,m_{j_\delta })$~, being the global elliptic representation of
$(\theta _{G^{(n)}_R} \otimes \theta _{G^{(n)}_L})$~, given by the product, right by left, of $n$-dimensional real global elliptic semimodules as introduced in \cite{Pie1},  corresponds to the searched cuspidal representation.

(Note that a similar procedure can be applied to the covering bisemisheaf
$(\theta _{(S\circ T)(1)_R}\otimes S_{(S\circ T)(1)_L})$~.)

We refer to proposition 4.2.10 which states the Langlands global correspondences as resulting from the singularization and the versal deformation of the bisemisheaf 
$(\theta _{G^{(n)}_R} \otimes \theta _{G^{(n)}_L})$~.
\vskip 11pt

In chapter 5, {\bf the monodromy of (isolated) singularities\/} on the (bisemi)sheaf
$(\theta ^{\cit}_{G^{(n)}_R} \otimes \theta ^{\cit}_{G^{(n)}_L})$ of differentiable bifunctions
$\phi ^{(n)}_{G_R^{(\cit)}}(z_{g_R})\otimes \phi ^{(n)}_{G_L^{(\cit)}}(z_{g_L})$ on the complex bilinear algebraic semigroup $G^{(n)}(F_{\o\omega }\times F_\omega )$ is analysed and the Langlands global correspondences on the non singular fibres, generated by monodromy, are developed in the irreducible and reducible cases.

\Bi
\item The monodromy arises in an expanding phase which reflects the expansion of the subvarieties of a given variety with respect to a fixed measure and which is assumed to generate locally surjective morphisms of singularizations.

\item The generated singularities can be non degenerate or degenerate in which case small deformations of these can split them up into simpler ones. So, assume that each section of the semisheaf $\theta ^{\rit}_{G^{(n)}_L} \subset \theta ^{\cit}_{G^{(n)}_L}$ (resp.
$\theta ^{\rit}_{G^{(n)}_R} \subset \theta ^{\cit}_{G^{(n)}_R}$~) is a Morse function affected by an isolated non degenerate singularity on a domain $U_{j_L}$ (resp. $U_{j_R}$~) included into the conjugacy class representative $g_L^{(n)}[j,m_j]$ (resp. $g_R^{(n)}[j,m_j]$~) and described locally by
\[ \phi ^{(2n)}_{G_{j_L}^{(\rit)}}(U_{j_L}) = \sum^{2n}_{i=1} x^2_{i_{L_j}}
\qquad \text{(resp.} \quad
\phi ^{(2n)}_{G_{j_R}^{(\rit)}}(U_{j_R}) = 
\sum^{2n}_{i=1} x^2_{i_{R_j}}\ ).\]

\item The critical level set of $\phi ^{(2n)}_{G_{j_L}^{(\rit)}}(U_{j_L})$ (resp.
$\phi ^{(2n)}_{G_{j_R}^{(\rit)}}(U_{j_R})$~) is the singular fibre $F^{(2n-1)}_{\circ_{j_L}}$
(resp. $F^{(2n-1)}_{\circ_{j_R}}$~) given by
\[ \phi ^{(2n)}_{G_{j_L}^{(\rit)}}(U_{j_L}) = \sum^{2n}_{i=1} x^2_{i_{L_j}}=0
\qquad \text{(resp.} \quad
\phi ^{(2n)}_{G_{j_R}^{(\rit)}}(U_{j_R}) = \sum^{2n}_{i=1} x^2_{i_{R_j}}=0\ )\]
while the non singular fibres $F^{(2n-1)}_{\lambda _{j_L}}$ (resp. $F^{(2n-1)}_{\lambda _{j_R}}$~) are diffeomorphic to the space $TS^{2n-1}_{L_j}$ (resp. $TS^{2n-1}_{R_j}$~) of the tangent bundle to a unit sphere $S^{2n-1}_{L_j}$ (resp. $S^{2n-1}_{R_j}$~), which is diffeomorphic to the vanishing cycle
$\Delta ^{(2n-1)}_{L_j}\subset F^{(2n-1)}_{\lambda _{j_L}}$ (resp. $\Delta ^{(2n-1)}_{R_j}\subset F^{(2n-1)}_{\lambda _{j_R}}$~).

As $\Delta ^{(2n-1)}_{L_j} $ (resp. $\Delta ^{(2n-1)}_{R_j}$~) is diffeomorphic to the unit sphere $S^{2n-1}_{L_j}$ (resp. $S^{2n-1}_{R_j}$~), it must correspond to a function on the corresponding conjugacy class representative of the parabolic subgroup $P^{(2n-1)}(F^+_{v^1})$
(resp. $P^{(2n-1)}(F^+_{\o v^1})$~).

\item So, the mapping:
\begin{align*}
h_{\gamma _{j_L}}: \quad F_{\lambda _{j_L}}^{(2n-1)}\To F^{(2n-1)}_{\lambda _{j_L}}\\
\text{(resp.} \quad
h_{\gamma _{j_R}}: \quad F_{\lambda _{j_R}}^{(2n-1)}\To F^{(2n-1)}_{\lambda _{j_R}}\ )\end{align*}
of the non-singular fibre into itself is the monodromy of the closed loop $\gamma _{j_L}\subset \Delta ^{(2n-1)}_{L_j}$ (resp. $\gamma _{j_R}\subset \Delta ^{(2n-1)}_{R_j}$~) realized by the conjugacy action of the $j$-th conjugacy class representative of the restricted linear algebraic semigroup $G^{(2n-1)}(F^{+(\rm res)}_{v_j})$ (resp. $G^{(2n-1)}(F^{+(\rm res)}_{\o v_j})$~).

\item If a degenerate singularity decomposes by deformation into a set of elementary non degenerate singular points, the single monodromy becomes a monodromy group.  In this context, if every section
$\phi ^{(2n)}_{G^{(\rit)}_{j_L}}$ (resp. $\phi ^{(2n)}_{G^{(\rit)}_{j_R}}$~), $1\le j\le r\le \infty $~, of $\theta ^{(\rit)}_{G_L^{(2n)}}$ (resp. $\theta ^{(\rit)}_{G_R^{(2n)}}$~) is endowed with a set of $k$ non degenerate singularities $\omega _{i_L}$ (resp. $\omega _{i_R}$~), $1\le i\le k$~, on $U_{j_L}$ (resp. $U_{j_R}$~), then the set of bisheaves $\{\Fs_{F^{(2n-1)}_{i_{\lambda _{j_R}}}}
\otimes \Fs_{F^{(2n-1)}_{i_{\lambda _{j_L}}}}\}^k_{i=1}$ of non singular bifibres
$F^{(2n-1)}_{i_{\lambda _{j_R}}}(t)\otimes F^{(2n-1)}_{i_{\lambda _{j_L}}}$ are generated by monodromy above every bisection of 
$\theta ^{(\rit)}_{G_R^{(2n)}}\otimes\theta ^{(\rit)}_{G_L^{(2n)}}$~.

And, if there are $b_i$~, $b_i\in\nit$~, non singular fibres in the sheaf 
$\Fs_{F^{(2n-1)}_{i_{\lambda _{j_l}}}}$ (resp. $ \Fs_{F^{(2n-1)}_{i_{\lambda _{j_R}}}}$~), then we get a set of $k\times b_i$~, $1\le i\le k$~, $1\le \beta _i\le b_i$~, monodromy bi(semi)sheaves above $\theta ^{(\rit)}_{G_R^{(2n)}}\otimes\theta ^{(\rit)}_{G_L^{(2n)}}$~.

\item Let $\{\Fs_{F^{(2n-1)}_{i_{\lambda _R}}}(\beta _i) \otimes \Fs_{F^{(2n-1)}_{i_{\lambda _L}}}(\beta _i)\}_{i,\beta _i}$~, $1\le i\le k$~, $1\le \beta _i\le b_i$~, be the set of $k\times b_i$ monodromy bisemisheaves above the desingularized bisemisheaf $\theta ^{(\rit)}_{G_R^{(2n)}}\otimes\theta ^{(\rit)}_{G_L^{(2n)}}$~.

Then {\bf a global holomorphic correspondence\/} can be stated for the bismisheaf
$\theta ^{(\rit)}_{G_R^{(2n)}}\otimes\theta ^{(\rit)}_{G_L^{(2n)}}$ as developed before and the set of global holomorphic correspondences:
\begin{multline*}
\Irr\Rep^{(2n-1)}_{W^{\rm mon}_{F^+\RL}}(W^{ab}_{F^+_{R_{\rm mon}}}(\beta _i) \times W^{ab}_{F^+_{L_{\rm mon}}}(\beta _i))\\
\To \Irr\hol^{(2n-1)}(\Fs_{F^{(2n-1)}_{i_{\lambda _R}}}(\beta _i)\otimes
\Fs_{F^{(2n-1)}_{i_{\lambda _L}}}(\beta _i)) \qquad \forall\ i,\beta _i\;, \end{multline*}
can be similarly found for the monodromy bisemisheaves
$\Fs_{F^{(2n-1)}_{i_{\lambda _R}}}(\beta _i)\otimes \Fs_{F^{(2n-1)}_{i_{\lambda _L}}}(\beta _i)$~.

After a toroidal compactification of these bisemisheaves, it is proved that:
\Bean
\item a {\bf cuspidal representation\/}, given by the elliptic representation\linebreak
$\Ellip \RL(2n,j,m_j)$~, can be associated with the desingularized bisemisheaf
$\theta ^{(\rit)}_{G_R^{(2n)}}\otimes\theta ^{(\rit)}_{G_L^{(2n)}}$~.

\item {\bf no cuspidal representation can be found for the monodromy bisemisheaves\/}, except if surgeries are performed.
\Ee

\item The {\bf orthogonal reducibility of the bisemisheaf\/} $\theta ^{(\cit)}_{\GL_{2n}}(F_{\o\omega }\times F_\omega )$ leads to the following decomposition:

\[ \theta ^{(\cit)}_{\GL_{2n=2_1+\cdots+2_n}}(F_{\o\omega }\times F_\omega )
= \mathop\boxplus^n_{\ell=1} \theta ^{(\cit)}_{\GL_{2_\ell}}(F_{\o\omega }\times F_\omega )\]
where the irreducible bisemisheaves $\theta ^{(\cit)}_{\GL_{2_\ell}}(F_{\o\omega }\times F_\omega )$
are able to generate monodromy groups.

If, on the domain, $U^{(2)}_{j_R}\times U^{(2)}_{j_L} \subset g^{(2)}\RL[j,m_j]$~, each bifunction of 
$\theta ^{(\cit)}_{\GL_{2_\ell}}(F_{\o\omega }\times F_\omega )$ is locally a Morse (bi)function in such a way that its critical set is the singular bifibre $F^{(1)}_{\circ_{j_R}}\times F^{(1)}_{\circ_{j_L}}$
given by
\[ \phi ^{(2)}_{G^{(\cit)}_{g_{j_R}}}(U^{(2)}_{j_R})\otimes\phi ^{(2)}_{G^{(\cit)}_{g_{j_L}}}(U^{(2)}_{j_L})=z^2_{j_1}+z^2_{j_2}=0\;, \quad (z_1,z_2)\in \cit^2\;,\]
then,
\Be
\item[1)] the corresponding non singular bifibres $F^{(1)}_{\lambda _{j_R}}(t) \times F^{(1)}_{\lambda _{j_L}}(t)$ are diffeomorphic to the product, right by left, $T^2_{\lambda _{j_R}}(t) \times T^2_{\lambda _{j_L}}(t)$ of two semitori.

\item[2)] the homology group $H_1(F^{(1)}_{\lambda _{j_L}};\zit) \simeq \zit$ (resp.
$H_1(F^{(1)}_{\lambda _{j_R}};\zit) \simeq \zit$~) of the semitorus $T^2_{\lambda _{j_L}}$ (resp. $T^2_{\lambda _{j_R}}$~) is generated by the upper (resp. lower) semicircle $\Delta ^{(1)}_{L_j}$
(resp. $\Delta ^{(1)}_{R_j}$~) on $T^2_{\lambda _{j_L}}$ (resp. $T^2_{\lambda _{j_R}}$~) in such a way that $\Delta ^{(1)}_{L_j}$
(resp. $\Delta ^{(1)}_{R_j}$~) shrinking onto the singularity, becomes the vanishing semicycle.
\Ee\Ei
\vskip 11pt

If each bisection of the bisemisheaf $\theta ^{(\cit)}_{\GL_2(F_{\o\omega }\times F_\omega )}$ is endowed with the same singular bifibre $F^{(1)}_{\circ_{j_R}}\times F^{(1)}_{\circ_{j_L}}=z^2_{j_1}+z^2_{j_2}=0$~, then a set of $\beta $ bisemisheaves $\{\theta ^{(\cit)_{\rm mon}}_{\GL_2(F_{\o\omega }\times F_\omega )}(b )\}^\beta _{b=1}$~, isomorphic to (or ``copies of'') the desingularized bisemisheaf $\theta ^{(\cit)}_{\GL_2(F_{\o\omega }\times F_\omega )}$~, can be generated by monodromy if $\beta $ is the number of non singular bifibres above each bisection of $\theta ^{(\cit)}_{\GL_2(F_{\o\omega }\times F_\omega )}$~.

And, a set of {\bbf $\beta $ global holomorphic correspondences\/} can be associated with the $\beta $ monodromy bisemisheaves according to:
\begin{eqnarray*}
&\Irr Rep^{(1)}_{W^{\rm mon}_{F\RL}}(W^{ab}_{F_{R_{\rm mon}}}(b) \times W^{ab}_{F_{L_{\rm mon}}}(b)) &
\To \Irr\hol^{(1)}(\theta ^{(\cit)_{\rm mon}}_{\GL_2(F_{\o\omega }\times F_\omega )}(b))\\
&\| &\\
&\theta ^{(\cit)_{\rm mon}}_{\GL_2(F_{\o\omega }\times F_\omega )}(b) &
\To f_{\o\omega }(z^*_{m_b})\times f_\omega (z_{m_b})\;, \quad 1<b<\beta \;.\end{eqnarray*}

Similarly, on the toroidal compactified monodromy bisemisheaves, the following {\bf Langlands irreducible global correspondences can be stated\/}:
\begin{eqnarray*}
&\Irr Rep^{(1)}_{W^{\rm mon}_{F\RL}}(W^{ab}_{F_{R_{\rm mon}}})(b) \times W^{ab}_{F_{L_{\rm mon}}}(b)) &
\To \Irr\cusp(\theta ^{(\cit)_{\rm mon}}_{\GL_2(F_{\o\omega }\times F_\omega )}(b))\\
&\| &\\
&\theta ^{(\cit)_{\rm mon}}_{\GL_2(F_{\o\omega }\times F_\omega )}(b) &
\To \EIS^{\rm mon}\RL(1,j,m_j)_b \end{eqnarray*}
where $\EIS^{\rm mon}\RL(1,j,m_j)$~, being the product, right by left, of the equivalents of the Eisenstein series, constitutes the cuspidal representation of the $b$-th monodromy bisemisheaf
$\theta ^{(\cit)_{\rm mon}}_{\GL_2(F_{\o\omega }\times F_\omega )}$~.

\section[(Bisemi)sheaf of differentiable (bi)functions on the bilinear algebraic semigroup $\Gn$]{\boldmath (Bisemi)sheaf of differentiable (bi)functions on the bilinear algebraic semigroup $\Gn$}

\addtocontents{toc}{\protect\thispagestyle{empty}}

\subsection{Completions at infinite places of a global number field}

Let $\wF$ denote a finite (algebraically closed) Galois extension of a global number field 
$F^0$ of characteristic zero. In the complex case, the splitting field $\wF = \wFR \cup \wFL$ 
is assumed to be composed of the left and right splitting semifields $\wFL$ and $\wFR$ in one-to-one correspondence in such a way that the \lr algebraic extension semifield $\wFL$ (resp. $\wFR$~) is the set of complex (resp. conjugate complex) simple roots of a polynomial ring over $F^0$~.

In the real case, the symmetric splitting field is noted $\wF^+=\wFR^+\cup \wFL^+$ where $\wFL^+$ (resp. $\wFR^+$~) is the algebraic extension semifield composed of the set of positive (resp. symmetric negative) simple real roots.

The left and right equivalence classes of the local completions $F_L^{(+)}$ and $F_R^{(+)}$ 
respectively of $\wFL^{(+)}$ and $\wFR^{(+)}$ are the left and right complex (resp. real) infinite places of $F_L^{(+)}$ and $F_R^{(+)}$~: they are noted $v=\{v_{1_\delta },\cdots,v_{j_\delta },\cdots,v_{t_\delta }\}$ and $\o v=\{\o v_{1_\delta },\cdots,\o v_{j_\delta },\cdots,\o v_{t_\delta }\}$ in the real case and
$ \omega =\{\omega _1,\cdots,\omega _{j },\cdots,\omega _r\}$ and $\o \omega =\{\o \omega _1,\cdots,\o \omega _{j},\cdots,\o \omega _r\}$ in the complex case and are equal in number.
\vskip 11pt

The \lr complex pseudo-unramified completions $F_{\omega_j}^{nr}$ (resp. $F_{\o\omega_j}^{nr}$~), $1\le j\le r$~, of $\wFL$ (resp. $\wFR$~) are pseudo-unramified $F^0$-semimodules characterized by their ranks, called global residue degrees,
\[ [F_{\omega _j}^{nr}:F^0]=j \qquad \text{(resp.} \quad [F_{\o \omega _j}^{nr}:F^0]=j \; ),\]
and the \lr real pseudo-unramified completions $F_{v_{j_\delta }}^{+,nr}$ (resp. 
$F_{\o v_{j_\delta }}^{+,nr}$~), $1\le j_\delta \le r$~, of $\wFL$ (resp. $\wFR$~) are also characterized by their global residue degrees:
\[ [F_{v _{j_\delta }}^{+,nr}:F^0]=j  \qquad \text{(resp.} \quad [F_{\o v _{j_\delta }}^{+,nr}:F^0]=j  \; ).\]

The \lr complex pseudo-ramified completions $F_{\omega _j}$ (resp. $F_{\o \omega _j}$~) 
of $\wFL$ (resp. $\wFR$~) and the \lr real pseudo-ramified completions $F_{v_{j_\delta }}^{+} $ (resp.  $F_{\o v_{j_\delta }}^{+} $~) are $F^0$-semimodules  generated from irreducible central completions $F_{\omega ^1_j}$ (resp. $F_{\o \omega ^1_j}$~) of rank $N\cdot m^{(j_\delta )}$ in the complex case and from irreducible central completions 
$F^+_{v^1_{j_\delta }}$ (resp. $F^+_{\o v^1_{j_\delta }}$~) of rank $N$ in the real case, where $m^{(j_\delta )}$ is the multiplicity of the $j_\delta $-th real completion $F^+_{v_{j_\delta }}$ covering its complex equivalent $F_{\omega _j}$~.
\vskip 11pt

So, if the irreducible central completions are given by their ranks:
\Bi
\item \quad $[F_{\omega ^1_j}:F^0] = N\centerdot m^{(j_\delta )}$ \qquad (resp. \quad $[F_{\o\omega ^1_j}:F^0] = N\centerdot m^{(j_\delta )}$~),
\item \quad $[F^+_{v^1_{j_\delta }}:F^0] = N\phantom{\;\centerdot m^{(j_\delta )}}$ \qquad (resp. \quad $[F^+_{\o v^1_{j_\delta }}:F^0] = N$~),
\Ei
the pseudo-ramified completions can be expressed from their corresponding pseudo-unramified equivalents as follows:
\Bi
\item
$\begin{array}[t]{rl}
[F_{\omega _j}:F^0] &= [F_{\omega _j}^{nr}: F^0]\times [F_{\omega ^1_j}:F^0]=*_c+j\centerdot N\centerdot m^{(j_\delta )}\\
\text{(resp.} \quad 
[F_{\o\omega _j}:F^0] &= [F_{\o\omega _j}^{nr}: F^0]\times [F_{\o\omega ^1_j}:F^0]=*_c +j\centerdot N\centerdot m^{(j_\delta )}\;).\end{array}$

\item
$\begin{array}[t]{rl}
[F_{v_{j_\delta }}:F^0] &= [F_{v_{j_\delta }}^{nr}: F^0]\times [F_{v^1_{j_\delta }}:F^0]=*+j \centerdot N\\
\text{(resp.} \quad 
[F_{\o v_{j_\delta }}:F^0] &= [F_{\o v_{j_\delta }}^{nr}: F^0]\times [F_{\o v^1_{j_\delta }}:F^0]=*+j \centerdot N\;),\end{array}$

where $*_c$ denotes an integer inferior to $N\centerdot m^{(j_\delta )}$ and $*$ an integer inferior to $N$~.
\Ei
\vskip 11pt

Then, the complex pseudo-ramified completions $F_{\omega _j}$ (resp. $F_{\o\omega _j}$~), 
$1\le j\le r\le\infty $~, can be approximatively cut into a set of $j$ irreducible equivalent completions 
$F_{\omega _j^{j'}}$~, $1\le j'\le j$ (resp. $F_{\o\omega _j^{j'}}$~), of rank $N\centerdot m^{(j_\delta )}$ while the real pseudo-ramified completions $F^+_{v_{j_\delta }}$ (resp. $F^+_{\o v_{j_\delta }}$~), $1\le j_\delta \le r\le\infty $~, can be approximatively cut into a set of $j $ irreducible equivalent completions 
$F^+_{v_{j_\delta }^{j'_\delta }}$~, $1\le j'_\delta \le j_\delta $ (resp. $F^+_{\o v_{j_\delta }^{j'_\delta }}$~), of rank $N$~.
\vskip 11pt

On the other hand, as a place is an equivalence class of completions, we have to take into 
account a set of complex completions $\{F_{\omega _{j,m_j}}\}_{m_j}$~, $m_j\in\NN$~, 
$m_j\ge 1$~, equivalent to the basic completion $F_{\omega _j}$ at the $j$-th complex place $\omega _j$ and characterized by the same rank as $F_{\omega _j}$~.  These complex equivalent completions $F_{\omega _{j,m_j}}$ are generated from the basic completion $F_{\omega _j}$ in a nilpotent way \cite{Pie1}.

Similarly, at a real place $v_{j_\delta }$~, a set of real completions 
$\{F^+_{v_{j_\delta ,m_{j_\delta }}}\}_{m_{j_\delta }}$ equivalent to the basic real completion $F^+_{v_{j_\delta }}$ and characterized by the same rank has to be considered.

As it was indicated before, each complex completion $F_{\omega _j}$ is covered by the set $\{F^+_{v_{j_\delta ,m_{j_\delta }}}\}$ of $m ^{(j_\delta )}=\sup(m_{j_\delta })+1$ real equivalent completions $F^+_{v_{j_\delta ,m_{j_\delta }}}$~.

Let 
\begin{align*}
F_\omega  &= \{F_{\omega _1},\dots,F_{\omega _{j,m_j}},\dots,F_{\omega _{r,m_r}}\}
\rresp{F_{\o \omega}  &= \{F_{\o\omega _1},\dots,F_{\o\omega _{j,m_j}},\dots,F_{\o\omega _{r,m_r}}\}}\end{align*}
denote the set of complex pseudo-ramified completions at the set of complex places $\omega $ \resp{$\o\omega $} and let
\begin{align*}
F^+_v  &= \{F^+_{v _1},\dots,F^+_{v _{j_\delta ,m_{j_\delta }}},\dots,F^+_{v _{r_\delta ,m_{r_\delta }}}\}
\rresp{F^+_{\o v}  &= \{F^+_{\o v _1},\dots,F^+_{\o v _{j_\delta ,m_{j_\delta }}},\dots,F^+_{\o v _{r_\delta ,m_{r_\delta }}}\}
}\end{align*}
be the corresponding set of real pseudo-ramified compeltions at the set of real places $v$ \resp{$\o v$}.

Then, the direct sum of the complex pseudo-ramified completions is given by:
\[ F_{\omega _\oplus}=\bigoplus_j\bigoplus_{m_j}F_{\omega _{j,m_j}}
\rresp{F_{\o\omega _\oplus}=\bigoplus_j\bigoplus_{m_j}F_{\o\omega _{j,m_j}}}\]
while the direct sum of the real pseudo-ramified completions is given by:
\[ F_{v _\oplus}=\bigoplus_{j_\delta }\bigoplus_{m_{j_\delta }}F_{v _{j_\delta ,m_{j_\delta }}}
\rresp{F_{\o v _\oplus}=\bigoplus_{j_\delta }\bigoplus_{m_{j_\delta }}F_{\o v _{j_\delta ,m_{j_\delta }}}}.\]
\vskip 11pt

And, a \lr pseudo-ramified adele semiring $\Aa_{F_\omega }$ (resp. $\Aa_{F_{\o\omega} }$~) can be introduced on the product of the basic  primary completions $F_{\omega _{j_p}}$ (resp. $F_{\o\omega _{j_p}}$~) and of their equivalent completions $F_{\omega _{{j_p},m_{j_p}}}$ (resp. $F_{\o\omega _{{j_p},m_{j_p}}}$~) over all primary complex places according to:
\begin{align*}
\Aa_{F_\omega } &= \prod_{j_p} F_{\omega _{j_p}}\prod_{m_{j_p}} F_{\omega _{{j_p},m_{j_p}}}\;, \qquad 1\le {j_p}\le r\le \infty \;, \quad m_{j_p}\ge 1\;, \\[11pt]
\text{(resp.} \quad
\Aa_{F_{\o \omega }} &= \prod_{j_p} F_{\o\omega _{j_p}}\prod_{m_{j_p}} F_{\o\omega _{{j_p},m_{j_p}}}\;).
\end{align*}

Similarly, a \lr adele semiring $\Aa_{F^+_v}$ (resp. $\Aa_{F^+_{\o v}}$~) can be introduced over all primary real places according to:
\begin{align*}
\Aa_{F^+_v} &= \prod_{j_{\delta_p} }F^+_{v_{j_{\delta_p} }} \prod_{m_{j_{\delta_p} }}F^+_{v_{j_{\delta_p} ,m_{j_{\delta_p} }}}\;, \qquad 1\le j_{\delta_p} \le r\le\infty \\[11pt]
\text{(resp.} \quad
\Aa_{F^+_{\o v}} &= \prod_{j_{\delta_p} }F^+_{\o v_{j_{\delta_p} }} \prod_{m_{j_{\delta_p} }}F^+_{\o v_{j_{\delta_p} ,m_{j_{\delta_p} }}}\;).
\end{align*}
\vskip 11pt

\subsection[The reductive bilinear algebraic semigroup $\Gn$]{\boldmath The reductive bilinear algebraic semigroup $\Gn$}

The set  $F_\omega$
(resp. $ F_{\o\omega}$~) of complex pseudo-ramified completions generates a tower of $r$ packets of completions following the complex places $\omega _j$ (resp. $\o\omega _j$~), $1\le j\le r$~: the \lr tower of $r$ packets of complex pseudo-ramified completions, restricted to the upper (resp. lower) half space, is a one-dimensional complex linear affine semigroup  noted $\SS^1_{\omega _L}$ (resp. $\SS^1_{\omega _R}$~).  For reasons developed in \cite{Pie1}, we are interested in the product, right by left, $\SS^1_{\omega _R}\times \SS^1_{\omega _L}$ of $\SS^1_{\omega _R}$
by its symmetric $\SS^1_{\omega _L}$ where $\SS^1_{\omega _R}\times \SS^1_{\omega _L}$ is a bilinear affine complex semigroup.

Similarly, a \lr tower of $r$ packets of real pseudo-ramified completions is a one-dimensional real affine semigroup  $\SS^1_{v_L}$ (resp. $\SS^1_{\o v_R}$~) and the product, right by left,
$\SS^1_{\o v_R}\times \SS^1_{v_L}$ of $\SS^1_{v_L}$ by $\SS^1_{\o v_R}$ is a bilinear affine real semigroup.
\vskip 11pt

The $n$-dimensional analog of $\SS^1_{\omega _R}\times \SS^1_{\omega _L}$ is a 
$n^2$-dimensional bilinear affine semigroup  which is a reductive bilinear algebraic 
semigroup $\Gn$ isomorphic to the bilinear algebraic semigroup of matrices $\GLn$ with 
entries in $\AaF$~.  Indeed,\linebreak $\GLn\equiv T^t_n({F_{\o \omega }}) \times T_n({F_\omega })$ is a condensed notation for the product of the group $T^t_n({F_{\o \omega }})$ of lower triangular matrices with entries in ${F_{\o \omega }}$ by the group $T_n({F_{\omega }})$ of upper triangular matrices with entries in ${F_{\omega }}$~.
\vskip 11pt

The bilinear algebraic semigroup $\Gn$ covers the corresponding linear algebraic group 
$G^{(n )}({F_{\o\omega }-F_ \omega })$~, where ${F_{\o\omega}-F_\omega }={F_{\o\omega }}\cup{F_\omega }$~, 
as it was justified in \cite{Pie1}, since the $n^2$-dimensional representation space 
$M_{{F_{\o\omega }}}\otimes M_{{F_\omega }}$ of $\GLn$ then coincides with the 
$n^2$-dimensional representation space $V$ of $\GL_n({F_{\o\omega }-F_\omega })$ under 
some conditions given in \cite{Pie1}.
\vskip 11pt

As the bilinear algebraic semigroup $\Gn$ is built over $\AaF$~, it is composed of $r$ conjugacy classes, $1\le j\le r$~, having multiplicities $m^{(r)}=\sup(m_r)+1$~, where $m^{(r)}$ denotes the number of equivalent representatives in the $r$-th conjugacy class.  Remark that the $r$ conjugacy classes of $\Gn$ correspond to the $r$ (bi)places of $\AaF$~.
\vskip 11pt

\subsection{Proposition}

{\em Let $\SMOD_{{F_\omega }}$ (resp. $\SMOD_{{F_{\o\omega} }}$~) denote the category of $T_n({F_\omega })$)-semimodules $M_{{F_\omega }}$
(resp. $T^t_n({F_{\o\omega }})$-semimodules $M_{{F_{\o\omega }}}$~) $\subset\GLn$-bisemimodules $M_{{F_{\o\omega }}}\otimes M_{{F_{\omega }}}$~.

Then, the $T_n( {F_\omega })$-semimodule $M_{ {F_{\omega }}}$ (resp. $T^t_n( {F_{\o\omega }})$-semimodule $M_{ {F_{\o\omega }}}$~) is a division $ {F_\omega }$-semialgebra (resp. a division $ {F_{\o\omega }}$-cosemialgebra).}
\vskip 11pt

\bpr Indeed, according to the appendix of \cite{Pie1}, the $ {F_\omega }$-semialgebra 
$M_{ {F_\omega }}$ (resp. $ {F_{\o\omega }}$-cosemialgebra $M_{ {F_{\o\omega }}}$~) 
over the  semiring $ {F_ \omega }$ (resp. $ {F_{\o\omega }}$~) is a semiring 
$M_{ {F_\omega }}$ (resp. $M_{ {F_{\o\omega }}}$~) such that:
\Bean
\item $(M_{ {F_\omega }},+)$ (resp. $(M_{ {F_{\o\omega }}},+)$~) is a unitary left $ {F_\omega }$-semimodule (resp. right $ {F_{\o\omega }}$-semimodule).

\item $\begin{aligned}[t]
a_L\ (m_L\ n_L)&=(a_L\ m_L)\ n_L = m_L\ (a_L\ n_L)\ , \; \forall\ a_L\in {F_\omega }\ , \; m_L,n_L\in M_{ {F_\omega }}\\ 
\text{(resp.\ }  (m_R\ n_R)\ a_R&=m_R\ (n_R\ a_R)=m_R\ (a_R\ n_R)\ , \;\forall\ a_R\in {F_{\o\omega }}\ , \; m_R,n_R\in M_{ {F_{\o\omega }}}\;).\end{aligned}$
\Ee

As $M_{ {F_\omega }}$ (resp. $M_{ {F_{\o\omega }}}$~) is a \lr division semiring,
$M_{ {F_\omega }}$ (resp. $M_{ {F_{\o\omega }}}$~) is a division $ {F_\omega }$-semialgebra (resp. $ {F_{\o\omega }}$-cosemialgebra).\epr

\subsection[Definition : $\AaF$-bisemialgebra]{\boldmath Definition : $\AaF$-bisemialgebra}

\Bena
\item The $ {F_\omega }$-semialgebra $M_{ {F_\omega }}$ is  a monoid $(M_{ {F_\omega }},\mu ,\eta )$ in  $\SMOD_{ {F_\omega }}$~, in the sense that:
\Bi
\item $M_{ {F_\omega }}$ is assumed to be a unitary $ {F_\omega }$-semimodule, i.e. a left vector semispace over $ {F_\omega }$ viewed as the center of $M_{F_{\omega }}$~;

\item $\mu :M_{ {F_{\o \omega}}}\otimes M_{ {F_\omega }}\to 
M_{ {F_\omega}} $ is a linear homomorphism;

\item $\eta : {F_\omega }\to M_{ {F_\omega }}$ is an injective homomorphism.
\Ei

\item The $ {F_{\o \omega }}$-cosemialgebra $M_{ {F_{\o \omega }}}$ is dually a comonoid $(M_{ {F_{\o\omega }}},\Delta ,\varepsilon )$ in $SMOD_{ {F_{\o\omega }}}$ in such a way that:
\Bi
\item $M_{ {F_{\o\omega }}}$ is assumed to be a unitary $M_{ {F_{\o\omega }}}$-semimodule, i.e. a right vector semispace over $ {F_{\o\omega }}$ so that $M_{ {F_{\o\omega }}}$ is the dual of $M_{ {F_\omega }}$~;

\item $\Delta : M_{ {F_{\o\omega }}}\to 
M_{ {F_{\o\omega }}}\times M_{ {F_{\omega }}}$ is a linear homomorphism called 
comultiplication;

\item $\varepsilon :M_{ {F_{\o\omega }}}\to  {F_{\o\omega }}$ is a linear form.
\Ei

\item $((M_{ {F_{\o\omega }}}\otimes 
M_{ {F_{\omega }}}),\mu ,\eta ,\Delta ,\varepsilon )$ 
is a division $( {F_{\o\omega }}\times  {F_{\omega }})$-bisemialgebra if
$(M_{ {F_{\omega }}},\mu ,\eta )$ is a division $ {F_\omega }$-semialgebra and if
$(M_{ {F_{\o\omega }}},\Delta ,\varepsilon )$ is a division 
$ {F_{\o\omega }}$-cosemialgebra.
\Ee
\vskip 11pt

\subsection{The bilinear parabolic subsemigroups}

As the \lr complex pseudo-ramified completions $F_{\omega _j}$ (resp. $F_{\o\omega _j}$~) and the \lr real pseudo-ramified completions $F^+_{v_{j_\delta }}$ (resp. $F^+_{\o v_{j_\delta }}$~) are assumed to be generated respectively from irreducible central complex completions $F_{\omega  ^1_j}$ (resp. $F_{\o\omega  ^1_j}$~) of rank $N\centerdot m^{(j_\delta )}$ and from irreducible central real completions
$F^+_{v^1_{j_\delta }}$ (resp. $F^+_{\o v^1_{j_\delta }}$~) of rank $N$~, a set of \lr irreducible complex pseudo-ramified completions
\begin{align*}
 {F_{\omega^1 } }
&= \{F_{\omega ^1_1},\dots,F_{\omega ^1_{j,m_j}},\dots,F_{\omega ^1_{r,m_r}}\}\\[11pt]
\text{(resp.} \quad 
  {F_{\o\omega^1 } }
&= \{F_{\o\omega ^1_1},\dots,F_{\o\omega ^1_{j,m_j}},\dots,F_{\o\omega ^1_{r,m_r}}\}
\;)\end{align*}
can be introduced, as well as a set of \lr irreducible real pseudo-ramified completions:
\begin{align*}
{F^+_{v^1} } 
&= \{F^+_{v^1_{1_\delta }} ,\dots,
F^+_{v^1_{j_\delta, m_{j_\delta }}},\dots
F^+_{v^1_{r_\delta, m_{r_\delta }}}\}\\[11pt]
\text{(resp.} \quad 
{F^+_{\o v^1} } 
&= \{F^+_{\o v^1_{1_\delta }} ,\dots,
F^+_{\o v^1_{j_\delta, m_{j_\delta }}},\dots
F^+_{\o v^1_{r_\delta, m_{r_\delta }}}\}
\;).\end{align*}

The $n$-dimensional smallest normal bilinear affine subsemigroup $\Gn$ (resp. $\Gnv$~) 
is the parabolic bilinear algebraic subsemigroup $\Pn$ (resp. $\Pnv$~) isomorphic to the 
bilinear algebraic semigroup of matrices
\begin{align*}
\GL_n( {F_{\o\omega^1 }}\times  {F_{\omega ^1}})
&= 
T_n^t ( {F _{\o\omega^1 }})\times T_n( {F_{\omega^1} })\\
\text{(resp.} \quad
\GLnv 
&= 
T_n^t ( {F^+_{\o v^1}})\times T_n( {F^+_{v^1}})
\;)\end{align*}
with entries in products of irreducible completions.

The parabolic bilinear subsemigroup $\Pn$ (resp. $\Pnv$~) can be considered as the unitary irreducible representation space of $\GLn$ (resp. $\GL_n( {F^+_{\o v}}\times  {F^+_v})$~), \cite{Pie1}.
\vskip 11pt

\subsection[The bialgebra of bifunctions on $\Gnv$]{\boldmath The bialgebra of bifunctions on $\Gnv$}

The smooth differentiable real-valued functions on the real bilinear algebraic semigroup $\Gnv$ are the tensor products (called bifunctions)
\[ \phi ^{(n)}_{G_R}(x_{g_R}) \otimes \phi ^{(n)}_{G_L}(x_{g_L}) \]
of smooth differentiable functions $\phi ^{(n)}_{G_L}(x_{g_L}) \in 
\widehat G_L^{(n)}( {F_v^+ })$~, $x_{g_L}\in G^{(n)}_L( {F^+_v})$~, of the algebra $\widehat G^{(n)}_L( {F^+_v})$ of these functions on the linear algebraic semigroup $G^{(n)}_L( {F^+_v})$~, localized in the upper half space, by the symmetric differentiable cofunctions
$\phi ^{(n)}_{G_R}(x_{g_R}) \in \widehat G_R^{(n)}( {F^+_{\o v}})$~, $x_{g_R}\in G^{(n)}_R( {F^+_{\o v}})$~, of the coalgebra $\widehat G^{(n)}_R( {F^+_{\o v}})$ of the cofunctions (or linear functionals) on the linear algebraic semigroup $G^{(n)}_R( {F^+_{\o v}})$~, localized in the lower half space.

As $\Gnv$ is partitioned into conjugacy classes on the biplaces $\o v_{j_\delta }\times v_{j_\delta }$~, $1\le j_\delta \le r$~, of $F^+_R\times F^+_L$~, we have to consider bifunctions
$\phi ^{(n)}_{G_{j_R}}(x_{j_\delta })\otimes \phi ^{(n)}_{G_{j_L}}(x_{j_\delta })$ on the conjugacy class representatives
$g^{(n)}\RL [j_\delta ,m_{j_\delta }]$~,
also noted $G^{(n)}(F^+_{\o v_{j_\delta ,m_{j_\delta }}} 
\times F^+_{v_{j_\delta ,m_{j_\delta }}})$~, of $\Gnv$~.

As the conjugacy classes $g^{(n)}\RL[j_\delta ]$ of $\Gnv$ form an increasing sequence
\[ g^{(n)}\RL[1] \subset \cdots \subset g^{(n)}\RL[j_\delta ] \subset \cdots\subset
g^{(n)}\RL[r] \;,\]
the bifunctions on the conjugacy class representatives $g^{(n)}\RL[j_\delta ,m_{j_\delta }]$ also form an increasing sequence:
\[ \phi ^{(n)}_{G_{1_R}}(x_{g_{1_R}}) \otimes \phi ^{(n)}_{G_{1_L}}(x_{g_{1_L}}) 
\subset \cdots \subset 
\phi ^{(n)}_{G_{j_R}}(x_{g_{j_R}}) \otimes \phi ^{(n)}_{G_{j_L}}(x_{g_{j_L}}) 
\subset \cdots \;.\]
The bialgebra of all differentiable real-valued measurable bifunctions
$ \phi ^{(n)}_{G_{R}}(x_{g_R })\otimes \phi ^{(n)}_{G_{R}}(x_{g_L })$ on $\Gnv$ satisfying:
\[\int_{\Gnv}\L| \phi ^{(n)}_{G_{R}}(x_{g_R} )\otimes \phi ^{(n)}_{G_{L}}(x_{g_L })
\R| dx_{g_R}\ dx_{g_L}<\infty \]
is noted $L^{1-1}\RL(\Gnv)$~.
\vskip 11pt

\subsection[(Bisemi)sheaf of differentiable (bi)functions on $\Gnv$]{\boldmath (Bisemi)sheaf of differentiable (bi)functions on \protect\linebreak $\Gnv$}

The set of real-valued differentiable functions $\phi ^{(n)}_{G_L}(x_{g_L})$
(resp. cofunctions $\phi ^{(n)}_{G_R}(x_{g_R})$~) on the \lr linear algebraic semigroup $G^{(n)}_L( {F^+_v})$ (resp. $G^{(n)}_R( {F^+_{\o v}})$~) is a \lr semisheaf of rings $\theta _{G^{(n)}_L}$ (resp. $\theta _{G^{(n)}_R}$~) because:
\Bean
\item for all basic conjugacy class representatives $g_L^{(n)}[j_\delta ,m_{j_\delta }=0]$
(resp. $g_R^{(n)}[j_\delta ,m_{j_\delta }=0]$~) of the topological semispace 
$G^{(n)}_L( {F^+_v})$ (resp. $G^{(n)}_R( {F^+_{\o v}})$~), we have a set 
$\theta _{G^{(n)}_L}(g^{(n)}_L)[j_\delta ,m_{j_\delta }=0]$
(resp. $\theta _{G^{(n)}_R}(g^{(n)}_R)[j_\delta ,m_{j_\delta }=0]$~);

\item for all pairs of basic conjugacy class representatives 
$g_L^{(n)}[j_\delta ,m_{j_\delta }=0] \subset g_L^{(n)}[j_\delta+1 ,m_{j_\delta +1}=0]$
(resp.
$g_R^{(n)}[j_\delta ,m_{j_\delta }=0] \subset g_R^{(n)}[j_\delta+1 ,m_{j_\delta +1}=0]$~) a restriction map
\begin{alignat*}{3}
\res_{g_L^{(n)}[j_\delta +1],g_L^{(n)}[j_\delta ]}
: \quad & \theta _{G^{(n)}_L}(g_L^{(n)}[j_\delta +1]) & \To & 
\theta _{G_L^{(n)}}(g^{(n)}_L[j_\delta ])\\
\text{(resp.} \quad 
\res_{g_R^{(n)}[j_\delta +1],g_R^{(n)}[j_\delta ]}
: \quad & \theta _{G^{(n)}_R}(g_R^{(n)}[j_\delta +1] ) & \To & 
\theta _{G_R^{(n)}}(g^{(n)}_R[j_\delta ])\;)
.\end{alignat*}
\Ee

a) and b) generate a presemisheaf of rings $\theta _{G_L^{(n)}}$ (resp. $\theta _{G_R^{(n)}}$~) because it is a sheaf of abelian semigroups for every \lr point $x_{g_L}$ (resp. $x_{g_R}$~) of $G^{(n)}_L( {F^+_{v}})$ (resp. $G^{(n)}_R( {F^+_{\o v}})$~) where 
$\theta _{G_L^{(n)}}(x_{g_L})$ (resp. $\theta _{G_R^{(n)}}(x_{g_R})$~) has the structure of a semiring \cite{Pie1}, \cite{Ser1}.

The presemisheaf $\theta _{G_L^{(n)}} $ (resp. $\theta _{G_R^{(n)}} $~) is a semisheaf of rings if for every collection $\{g^{(n)}_L[j_\delta ]\}^r_{j=1}$
(resp. $\{g^{(n)}_R[j_\delta ]\}^r_{j=1}$~) of basic conjugacy class representatives in
$G^{(n)}_L( {F^+_v})$ (resp. $G^{(n)}_R( {F^+_{\o v}})$~) with $g_L^{(n)}=\cup g^{(n)}_L[j_\delta ]$ (resp. $g_R^{(n)}=\cup g^{(n)}_R[j_\delta ]$~), the map
\begin{alignat*}{3}
\res_{g_L^{(n)}}
: \quad & \theta _{G^{(n)}_L}(g_L^{(n)}[j_\delta ]) & \To & 
\prod_{j_\delta }\theta _{G_L^{(n)}}(g^{(n)}_L[j_\delta ])\\
\text{(resp.} \quad 
\res_{g_R^{(n)}}
: \quad & \theta _{G^{(n)}_R}(g_R^{(n)}[j_\delta ]) & \To & 
\prod_{j_\delta }\theta _{G_R^{(n)}}(g^{(n)}_R[j_\delta ])\;)\end{alignat*}
is injective \cite{Mum}.

The set of real valued differentiable bifunctions $
\phi ^{(n)}_{G_{R}}(x_{g_R })\otimes \phi ^{(n)}_{G_{L}}(x_{g_L})$ on $\Gnv$ is a bisemisheaf of rings, noted $\theta _{G^{(n)}\RL}=\theta _{G^{(n)}_R}\otimes \theta _{G^{(n)}_L}$~, whose bisections are the differentiable bifunctions $\phi ^{(n)}_{G_{j_R}}(x_{g_{j_\delta }}) \otimes \phi ^{(n)}_{G_{j_L}}(x_{g_{j_\delta }}) $ on the conjugacy class representatives
$g^{(n)}\RL[j_\delta ,m_{j_\delta }]$ of $\Gnv$~.
\vskip 11pt

\subsection[(Bi)ideal (bi)semisheaf of differentiable (bi)functions on the bilinear parabolic semigroup $\Pnv$]{\boldmath (Bi)ideal (bi)semisheaf of differentiable (bi)functions on the bilinear parabolic semigroup $\Pnv$}

Let $F^+_{v^1_{j_\delta }}$ (resp. $F^+_{\o v^1_{j_\delta }}$~) and
$F^+_{v^1_{j_\delta,m_{j_\delta } }}$ (resp. $F^+_{\o v^1_{j_\delta,m_{j_\delta } }}$~) denote respectively the basic and the equivalent irreducible real completions of the $j_\delta $-th real place $v_{j_\delta }$ (resp. $\o v_{j_\delta }$~).

Let $P^{(n)}(F^+_{v^1_{j_\delta }})$ (resp. $P^{(n)}(F^+_{\o v^1_{j_\delta }})$~) be the $n$-dimensional \lr linear affine subsemigroup restricted to this $j_\delta $-th basic irreducible completion: it is then the $j_\delta $-th basic conjugacy class representative of $P^{(n)}(F^+_{v^1})$ (resp. $P^{(n)}(F^+_{\o v^1})$~).

On $P^{(n)}(F^+_{v^1_{j_\delta }})$ (resp. $P^{(n)}(F^+_{\o v^1_{j_\delta }})$~) we can introduce the ideal of complex-valued differentiable \lr functions
$\phi ^{(n)}_{P_{j_L}}(x_{p_{j_\delta }})$ (resp. $\phi ^{(n)}_{P_{j_R}}(x_{p_{j_\delta }})$~).

This ideal $\phi ^{(n)}_{P_{j_L}}(x_{p_{j_\delta }})$ (resp. $\phi ^{(n)}_{P_{j_R}}(x_{p_{j_\delta }})$~) is an equivalence class of all real valued differentiable functions 
$\phi ^{(n)}_{P_{j_L,m_j}}(x_{p_{j_\delta,m_{j_\delta } }})$ (resp. $\phi ^{(n)}_{P_{j_R}}(x_{p_{j_\delta,m_{j_\delta } }})$~) 
on the  conjugacy class equivalent representatives 
$P^{(n)}(F^+_{v^1_{j_\delta,m_{j_\delta } }})$ (resp. $P^{(n)}(F^+_{\o v^1_{j_\delta,m_{j_\delta } }})$~).

So, the set of \lr ideals $\phi ^{(n)}_{P_{j_L}}(x_{p_{j_\delta }})$ (resp. $\phi ^{(n)}_{P_{j_R}}(x_{p_{j_\delta }})$~) of real valued differentiable functions 
$\phi ^{(n)}_{G_{L}}(x_{g_L})$ (resp. cofunctions $\phi ^{(n)}_{G_{R}}(x_{g_R })$~) is a \lr ideal semisheaf of differentiable functions (resp. cofunctions) on the parabolic semigroup
$P^{(n)}( {F^+_{v^1}})$ (resp. $P^{(n)}( {F^+_{\o v^1}})$~).

And, the set of biideals 
$\{\phi ^{(n)}_{P_{j_R}}(x_{p_{j_\delta }})\otimes \phi ^{(n)}_{P_{j_L}}(x_{p_{j_\delta }})\}$ of differentiable bifunctions is a biideal bisemisheaf, noted $\theta _{P^{(n)}\RL}
= \theta _{P^{(n)}_R}\otimes \theta _{P^{(n)}_L}$~, whose bisections are the biideals on the conjugacy class representatives $P^{(n)}(F^+_{\o v^1_{j_\delta }}\times F^1_{v^1_{j_\delta }})$ of the bilinear parabolic semigroup $P^{(n)}( {F^+_{\o v^1}}\times  {F^+_{v^1}})$~.
\vskip 11pt

\subsection[(Bisemi)sheaf of complex valued differentiable (bi)functions on $\Gn$]{\boldmath (Bisemi)sheaf of complex valued differentiable (bi)functions on $\Gn$}

Considering the inclusion
\[ G^{(n)}( {F^+_{\o v}}\times  {F^+_v}) \hookrightarrow G^{(n)}( {F_{\o\omega }}\times  {F_\omega })\]
of the real bilinear algebraic semigroup $\Gnv$ into the corresponding complex bilinear 
algebraic semigroup $\Gn$ as envisaged in \cite{Pie1}, we can introduce the semisheaf
$\theta ^{(\cit)}_{G_L}$ (resp. $\theta ^{(\cit )}_{G_R}$~) of complex valued differentiable functions
$\phi ^{(\cit )}_{G^{(n)}_L}(x_{g_L})$ (resp. cofunctions $\phi ^{(\cit )}_{G^{(n)}_R}(x_{g_R})$~) on the \lr linear algebraic semigroup $G^{(n)}( {F_\omega })$ (resp. $G^{(n)}( {F_{\o\omega }})$~).  $\theta ^{(\cit )}_{G_L^{(n)}}$ (resp. $\theta ^{(\cit )}_{G_R^{(n)}}$~) can be defined similarly as it was done in section 1.7 and, in the following, the real case will be essentially considered.
\vskip 11pt

\section{singularization and versal deformation}

\subsection{The singularization}

\subsubsection{General statement of the singularization process}

The \lr semisheaf $\Theta _{G^{(n)}_L}$ (resp. $\Theta _{G^{(n)}_R}$~) is a semisheaf of smooth differentiable functions on the conjugacy class representatives $g^{(n)}_L[j_\delta ,m_{j_\delta }]$ (resp. $g^{(n)}_R[j_\delta ,m_{j_\delta }]$~) of the \lr real linear algebraic semigroup $G^{(n)}( {F^+_v}) \simeq T_n( {F^+_v})$
(resp. $G^{(n)}( {F^+_{\o v}}) \simeq T^t_n( {F^+_{\o v}})$~).

Under some external perturbation(s), singularities can be generated on the \lr semisheaf 
$\Theta _{G^{(n)}_L}$ (resp. $\Theta _{G^{(n)}_R}$~) in such a way that:
\Bean
\item these singularities are produced symmetrically on the functions and cofunctions respectively on the left and on the right conjugacy class representatives 
$g^{(n)}_L[j_\delta ,m_{j_\delta }]\in G^{(n)}( {F^+_v}) $ and
 $g^{(n)}_R[j_\delta ,m_{j_\delta }]\in G^{(n)}( {F^+_{\o v}}) $~): this results from the fact that $G^{(n)}( {F^+_{v}}) $ and $G^{(n)}( {F^+_{\o v}})$ are symmetrical algebraic semigroups localized in some small domains respectively in the upper and in the lower half space. So, an external perturbation affects in a similar way functions on the upper half space and cofunctions on the lower half space.

\item on each function $\phi ^{(n)}_{G_{j_L}}(x_{g_{j_\delta }})$ 
on $g^{(n)}_L[j_\delta ,m_{j_\delta }]$ (resp. cofunction 
$\phi ^{(n)}_{G_{j_R}}(x_{g_{j_\delta }})$\linebreak on $g^{(n)}_R[j_\delta ,m_{j_\delta }]$~), 
a same singularity (or a same set of singularities) is generated.
\Ee

Indeed, according to a), the external perturbation is assumed to affect similarly and symmetrically every function $\phi ^{(n)}_{G_{j_L}}(x_{g_{j_\delta }})$ and cofunction $\phi ^{(n)}_{G_{j_R}}(x_{g_{j_\delta }})$~.
\vskip 11pt

The process of generation of singularities will be called a singularization. It consists of a collapse of (a) normal crossings divisor(s) into a locus becoming singular: this is a contracting surjective morphism corresponding to the inverse of a resolution of singularities (see, for example, \cite{Abh}, \cite{Ber}, \cite{DeJ}, \cite{Hir1,Hir2,Hir3}, \cite{Zar1,Zar2,Zar3}).
\vskip 11pt

\subsubsection[Definition: singularization of regular $f$-schemes]{\boldmath Definition: singularization of regular $f$-schemes}

The \lr linear algebraic semigroup $G^{(n)}( {F^+_{v}}) $ (resp. $G^{(n)}( {F^+_{\o v}})$~) plus the \lr semisheaf $\Theta _{G^{(n)}_L}$ (resp. $\Theta _{G^{(n)}_R}$~) on it is a \lr affine semisheme $(G^{(n)}( {F^+_{v}}), \Theta _{G^{(n)}_L})$ (resp. $(G^{(n)}( {F^+_{\o v}}), \Theta _{G^{(n)}_R})$~) \cite{Mum}.

Every differentiable \lr function 
$\phi ^{(n)}_{G_{j_L}}(x_{g_{j_\delta }})$ (resp. $\phi ^{(n)}_{G_{j_R}}(x_{g_{j_\delta }})$~) of $\Theta _{G^{(n)}_L}$ (resp. $\Theta _{G^{(n)}_R}$~) being similarly affected by some external perturbation, will be considered as a ``prototype'' \lr scheme function, written in condensed form $\phi _L$ (resp. $\phi _R$~) and called a \lr $f$-scheme.

The singularization of the \lr regular $f$-scheme $\phi _L$ (resp. $\phi _R$~) is a contracting surjective morphism
\begin{align*}
\rho _L:\quad 
\phi _L &\to \phi _L^* \qquad (\ * \text{\ is for ``star'' symbolizing the singularities})\\
\text{(resp.} \quad
\rho _R:\quad   \phi _R&\to \phi _R^*\;)\end{align*}
yielding a \lr singular $f$-scheme $\phi _L^*$ (resp. $\phi _R^*$~) below $\phi _L$ (resp. $\phi _R$~) and verifying:
\Bena
\item $\phi _L$ (resp. $\phi_R$~) and $\phi _L^*$ (resp. $\phi _R^*$~) have the same dimension $n$~.
\item the singular $f$-scheme $\phi _L^*$ (resp. $\phi _R^*$~) is characterized by a singular 
locus $\Sigma _L$ (resp. $\Sigma _R$~) associated with the centre $Z_L$ (resp. $Z_R$~) 
of the singularization of $\phi ^*_L$ (resp. $\phi ^*_R$~).
\Ee
\vskip 11pt

\subsubsection{Proposition}

{\em The inverse morphism
\[ \rho _L^{-1}:\quad \phi _L^*\to \phi _L\;, \qquad \text{(resp.} \quad 
\rho _R^{-1}:\quad \phi _R^*\to \phi _R\;)\]
of the singularization of the \lr regular $f$-scheme $\phi _L$ (resp. $\phi _R$~) corresponds to the desingularization of $ \phi _L^*$ (resp. $\phi _R^*$~) if there exists an inclusion:
\[ h_L : \quad \phi _L\hookrightarrow \o\phi _L\qquad \text{(resp.} \quad 
h_R : \quad \phi _R\hookrightarrow \o\phi _R\;)\]
such that:
\Bena
\item $\o\phi _L$ (resp. $\o\phi _R$~) is a regular projective \lr $f$-scheme.
\item $\rho _L^{-1}(Z_L)\cup(\o\phi _L\smallsetminus \phi _L)$ (resp. 
$\rho _R^{-1}(Z_R)\cup(\o\phi _R\smallsetminus \phi _R)$~) is a closed $f$-subscheme of 
$\o\phi _L$ (resp. $\o\phi _R$~) identified with a normal crossings divisor $D_L$ 
(resp. $D_R$~).\Ee
}
\vskip 11pt

\bpr \quad The normal crossings divisor $D_L$ (resp. $D_R$~) is a regular closed $f$-subscheme of $\o\phi _L$ (resp. $\o\phi _R$~) and is the image under $\rho _L^{-1}\circ h_L$ (resp. $\rho _R^{-1}\circ h_R$~)  (of the centre $Z_L$ (resp. $Z_R$~)) of the singular locus $\Sigma _L$ (resp. $\Sigma _R$~) above which $\rho _L$ (resp. $\rho _R$~) is not an isomorphism.  The centre $Z_L$ (resp. $Z_R$~) is generally in the singular locus $\Sigma _L$ (resp. $\Sigma _R$~) but, up to now, no ad hoc definition of the centre working in any dimension has been discovered \cite{Hau}.\epr
\vskip 11pt

\subsubsection{Definition: Normal crossings divisor}

\Bi
\item A divisor has normal crossings if it can be defined locally by a monomial ideal.
\item A normal crossings divisor $D_L$ (resp. $D_R$~) will be assumed to be a closed 
$f$-subscheme function on one or on a set of real irreducible completions 
$F^+_{v^{j'_\delta }_{j_\delta }}$ (resp. $F^+_{\o v^{j'_\delta }_{j_\delta }}$~) of rank $N$ (see section 1.1).
\Ei
\vskip 11pt

\subsubsection{Proposition (singularization)}

{\em Let $D_L$ (resp. $D_R$~) be a normal crossings divisor included into the \lr regular $f$-scheme $\o\phi _L$ (resp. $\o\phi _R$~) in such a way that:
\[ \o\phi _L=\phi _L\cup D_L\qquad \text{(resp.} \quad \o\phi _R=\phi _R\cup D_R\;).\]

The singularization of $\o\phi _L$ (resp. $\o\phi _R$~) into a singular locus $\Sigma _L$ (resp. $\Sigma _R$~) is given by the contracting surjective morphism:
\[ \o\rho _L: \quad \o\phi _L\to  \phi _L^* \qquad \text{(resp.} \quad \o\rho _R:\quad\o\phi _R\to\phi ^*_R\;)\]
such that:
\Bean
\item $\Sigma _L\subset \phi _L^*$ (resp. $\Sigma _R\subset \phi _R^*$~) be the union of the 
(homotopic) image of $D_L\subset \o\phi _L$ (resp. $D_R\subset \o\phi _R$~) and of a closed 
singular sublocus $\Sigma ^S_L\subset \Sigma _L$ (resp. $\Sigma ^S_R\subset \Sigma _R$~) of 
$\phi _L^*$ (resp. $\phi ^*_R$~):
\begin{align*}
\Sigma _L &= \o\rho _L(D_L)\cup \Sigma ^S_L\\
\text{(resp.} \quad 
\Sigma _R &= \o\rho _R(D_R)\cup \Sigma ^S_R\;).\end{align*}

\item $\o\rho _L$ (resp. $\o\rho _R$~) restricted to:
\begin{align*}
\o\rho _L^{\rm is}: \quad \phi _L\smallsetminus \o\rho _L^{-1}(\Sigma ^S_L) &\To \phi _L^*\smallsetminus \Sigma _L\\
\text{(resp.} \quad 
\o\rho _R^{\rm is}: \quad \phi _R\smallsetminus \o\rho _R^{-1}(\Sigma ^S_R) &\To \phi _R^*\smallsetminus \Sigma _R\;)\end{align*}
be an isomorphism.\Ee}
\vskip 11pt

\bpr
\Bi
\item The centre $Z_L$ (resp. $Z_R$~) corresponds to the singular locus $\Sigma _L$ (resp. $\Sigma _R$~) in the case of curves because their singularities are isolated points.

For surfaces, the situation is more complicated because the singular locus consists of a finite number of isolated points and irreducible curves which may be singular \cite{Hau}.

But generally, $Z_L\subset \o\rho _L(D_L)$ (resp. $Z_R\subset \o\rho _R(D_R)$~).

\item On the other hand, the singular locus $\Sigma _L$ (resp. $\Sigma _R$~) may factorize into:
\[ \Sigma _L=E_L\centerdot I_L \qquad \text{(resp.} \quad \Sigma _R=E_R\centerdot I_R\;)\]
where:
\Bi
\item $E_L$ (resp. $E_R$~) is a power of the exceptional component $\o\rho _L(D_L)$ (resp. $\o\rho _R(D_R)$~), i.e. a monomial;
\item $I_L$ (resp. $I_R$~) is an ideal which has at each point of $\Sigma _L$ (resp. $\Sigma _R$~) order less than or equal to the order of $\o\rho _L(D_L)$ (resp. $\o\rho _R(D_R)$~) along $Z_L$ (resp. $Z_R$~).

$I_L$ (resp. $I_R$~) is the weak transform of $\o\rho _L(D_L)$ (resp. $\o\rho _R(D_R)$~).
\Ei

\item The singular sublocus $\Sigma ^S_L\subset \phi ^*_L$ (resp. $\Sigma ^S_R\subset \phi ^*_R$~) results from singularizations anterior to that of $\o\rho _L(D_L)$ (resp. $\o\rho _R(D_R)$~) and then becomes the singular locus of the following blowups of $\phi ^*_L\smallsetminus \o\phi _L(D_L)$ (resp. $\phi ^*_R\smallsetminus \o\phi _R(D_R)$~): this corresponds to a more general case than envisaged in proposition 2.1.3.

\item Let:
\begin{align*}
\o\rho _L^S: \quad D_L\To \Sigma _L\smallsetminus \Sigma _L^S\\
\text{(resp.} \quad 
\o\rho _R^S: \quad D_R\To \Sigma _R\smallsetminus \Sigma _R^S\;)\end{align*}
be the singularization map restricted to the singular locus $\Sigma _L\smallsetminus \Sigma _L^S$ (resp. $\Sigma _R\smallsetminus \Sigma _R^S$~).
\Ei

Then, $\Sigma _L\smallsetminus \Sigma _L^S$ (resp. $\Sigma _R\smallsetminus \Sigma _R^S$~) is the contracting homotopic image of $D_L$ (resp. $D_R$~) in the sense that the number $n_{D_L}$ (resp. $n_{D_R}$~) of irreducible real completions on which is defined $D_L$ (resp. $D_R$~) is superior or equal to the number $n_{(\Sigma _L\smallsetminus \Sigma _L^S)}$ (resp. $n_{(\Sigma _R\smallsetminus \Sigma _R^S)}$~) of irreducible real completions on which is defined $\Sigma _L\smallsetminus \Sigma _L^S$ (resp. $\Sigma _R\smallsetminus \Sigma _R^S$~):
\be n_{D_L} \le n_{(\Sigma _L\smallsetminus \Sigma _L^S)} \qquad \text{(resp.} \quad 
n_{D_R} \le n_{(\Sigma _R\smallsetminus \Sigma _R^S)} \;).\tag*{\eop}\ee
\vskip 11pt

\subsubsection{Proposition}

{\em The singularization $\o\rho _L:\o\phi _L\to\phi _L^*$ (resp. $\o\rho _R:\o\phi _R\to\phi _R^*$~) of $\o\phi _L$ (resp. $\o \phi _R$~) is assumed to be given by the following sequence of contracting surjective morphisms:
\begin{align*}
\o\phi _L &\equiv \phi _L^{*(0)}  \xrightarrow{\o\rho _L^{(1)}}  
\phi _L^{*(1)}   \xrightarrow{\o\rho _L^{(2)}}   \phi _L^{*(2)}   \xrightarrow{\quad}
\cdots   \xrightarrow{\o\rho _L^{(r-1)}}  
\phi _L^{*(r-1)}  \xrightarrow{\o\rho _L^{(r)}}  
\phi _L^{*(r)} \\
\text{(resp.} \quad
\o\phi _R &\equiv \phi _R^{*(0)}   \xrightarrow{\o\rho _R^{(1)}}  
\phi _R^{*(1)}   \xrightarrow{\o\rho _R^{(2)}} \phi _R^{*(2)} \xrightarrow{\quad}
\cdots   \xrightarrow{\o\rho _R^{(r-1)}}  
\phi _R^{*(r-1)}  \xrightarrow{\o\rho _R^{(r)}}  
\phi _R^{*(r)}\ ) \end{align*}
verifying
\Bena
\item the singular locus
\[ \Sigma _L\subset \phi ^*_L\equiv \phi _L^{*(r)} \qquad \text{(resp.} \quad
\Sigma _R\subset \phi ^*_R\equiv \phi _R^{*(r)}\ )\]
is given by
\begin{align*}
\Sigma _L
& \equiv \Sigma _L^{(r)}  =  \Sigma _L^{(1)}  \cup  \o\rho _L^{(2)}(D_L^{(1)})
 \cup  \cdots  \cup  \o\rho _L^{(\ell )}(D_L^{(\ell -1)})
  \cdots   \cup   \cdots   \cup   \o\rho _L^{(r)}(D_L^{(r-1)})\\
\text{(resp.} \quad 
\Sigma _R
& \equiv \Sigma _R^{(r)}  =  \Sigma _R^{(1)}  \cup   \o\rho _R^{(2)}(D_R^{(1)})
 \cup  \cdots   \cup  \o\rho _R^{(\ell )}(D_R^{(\ell -1)})
  \cdots   \cup   \cdots   \cup   \o\rho _R^{(r)}(D_R^{(r-1)})\ ).
\end{align*}

\item $\o\rho _L$ (resp. $\o\rho _R$~) restricted to:
\begin{align*}
\o\rho _L^{(\rm is)}: \quad &\o\phi _L\smallsetminus \o\rho _L^{-1}(\Sigma _L^{(r)})
\To \phi ^*_L\smallsetminus \Sigma _L^{(r)}\\
\text{(resp.} \quad 
\o\rho _R^{(\rm is)}: \quad &\o\phi _R\smallsetminus \o\rho _R^{-1}(\Sigma _R^{(r)})
\To \phi ^*_R\smallsetminus \Sigma _R^{(r)}\ )\end{align*}
is an isomorphism.

\item The orders of the singular subloci $\Sigma _L^{(\ell )}$ 
(resp. $\Sigma _R^{(\ell )}$~) form an increasing sequence (from left to right) 
parallely with the increase of $(\ell )$~, $1\le\ell \le r$~.
\Ee}
\vskip 11pt

\bpr
\Bi
\item Let $\o\rho _L^{(\rm sing)}$ (resp. $\o\rho _R^{(\rm sing)}$~) denote the complement of $\o\rho _L^{(\rm is)}$ (resp. $\o\rho _R^{(\rm is)}$~) in $\o\rho _L$ (resp. $\o\rho _R$~)~: it is defined by the increasing sequence of subloci:
\begin{align*}
\o\rho _L^{(\rm sing)} &:\quad D_L^{(0)}  \xrightarrow{\o\rho _L^{(1)}}  
\Sigma  _L^{(1)}\cup D_L^{(1)}   \xrightarrow{\o\rho _L^{(2)}}   
\Sigma  _L^{(2)}\cup D_L^{(2)}  \\
& \qquad \qquad \xrightarrow{\quad}   
\cdots   \xrightarrow{\o\rho _L^{(\ell)}}  
\Sigma  _L^{ (\ell )}\cup D_L^{(\ell )}  \xrightarrow{\quad}
\cdots \xrightarrow{\o\rho _L^{(r)}}  
\Sigma _L^{(r)} \\
\text{(resp.} \quad
\o\rho _R^{(\rm sing)} &:\quad D_R^{(0)}  \xrightarrow{\o\rho _R^{(1)}}  
\Sigma  _R^{(1)}\cup D_R^{(1)}   \xrightarrow{\o\rho _R^{(2)}}   
\Sigma  _R^{(2)}\cup D_R^{(2)}  \\
& \qquad \qquad  \xrightarrow{\quad}   
\cdots   \xrightarrow{\o\rho _R^{(\ell)}}  
\Sigma  _R^{ (\ell )}\cup D_R^{(\ell )}  \xrightarrow{\quad}
\cdots \xrightarrow{\o\rho _R^{(r)}}  
\Sigma _R^{(r)} \ ) \end{align*}
where
\Bi
\item $\Sigma _L^{(\ell )}=\Sigma _L^{(\ell -1)}\cup \o\rho _L^{(\ell )}(D_L^{(\ell -1)})$
 (resp. 
$\Sigma _R ^{(\ell )}=\Sigma _R^{(\ell -1)}\cup \o\rho _R^{(\ell )}(D_R^{(\ell -1)})$~);

\item $\Sigma _L^{(\ell )} $ (resp. $\Sigma _R^{(\ell )} $~) is the $\ell$-th singular sublocus generated by the composition of contracting surjective morphisms:
\[ \o\rho _L^{(1)} \circ \cdots \circ \o\rho _L^{(\ell-1 )}\circ \o\rho _L^{(\ell )} \qquad \text{(resp.} \quad
\o\rho _R^{(1)} \circ \cdots \circ \o\rho _R^{(\ell -1)}\circ \o\rho _R^{(\ell )} \ )\]
restricted to their action on the normal crossings divisors $D_L^{(0)}\cdots D_L^{(\ell -1)}$ (resp. $D_R^{(0)}\cdots D_R^{(\ell -1)}$~).
\Ei
\vskip 11pt

\item Let $\Sigma _L^{(\ell -1)}=E_L^{(\ell -1)}\centerdot I_L^{(\ell -1)}$ (resp.
$\Sigma _R^{(\ell -1)}=E_R^{(\ell -1)}\centerdot I_R^{(\ell -1)}$~) be the factorisation of the singular sublocus $\Sigma _L^{(\ell -1)}$ (resp. $\Sigma _R^{(\ell -1)}$~) as introduced in proposition 2.1.5.

And, let
\begin{align*}
\o\rho _L^{(\rm sing)(\ell )} : \quad \Sigma _L^{(\ell -1)}\cup D_L^{(\ell -1)}\To \Sigma _L^{(\ell )}\\
\text{(resp.} \quad
\o\rho _R^{(\rm sing)(\ell )} : \quad \Sigma _R^{(\ell -1)}\cup D_R^{(\ell -1)}\To \Sigma _R^{(\ell )}\ )\end{align*}
be the $\ell $-th surjective morphism restricted to the singular sublocus $\Sigma _L^{(\ell )}$ (resp. $\Sigma _R^{(\ell )}$~): $\o\rho _L^{(\rm sing)(\ell )} $ (resp. $\o\rho _R^{(\rm sing)(\ell )}$~) is a fibre bundle with fibre $D_L^{(\ell-1 )}$ (resp. $D_R^{(\ell -1)}$~).

Then, the order of $\Sigma _L^{(\ell )}$ (resp. $\Sigma _R^{(\ell )}$~) at each point $p_L$ (resp. $p_R$~), being the maximal power of the maximal ideal of $p_L$ (resp. $p_R$~), is superior to the order of $\Sigma _L^{(\ell-1 )}$ (resp. $\Sigma _R^{(\ell-1 )}$~).

Indeed, the factorisations of $\Sigma _L^{(\ell-1 )}$ (resp. $\Sigma _R^{(\ell-1 )}$~) and of $\Sigma _L^{(\ell )}$ (resp. $\Sigma _R^{(\ell )}$~) are respectively given by:
\[ \Sigma _L^{(\ell -1)} = E_L^{(\ell -1)}\centerdot I_L^{(\ell -1)}
\qquad \text{(resp.} \quad 
\Sigma _R^{(\ell -1)} = E_R^{(\ell -1)}\centerdot I_R^{(\ell -1)}\ )\]
and by
\[ \Sigma _L^{(\ell )} = E_L^{(\ell )}\centerdot I_L^{(\ell )}
\qquad \text{(resp.} \quad 
\Sigma _R^{(\ell )} = E_R^{(\ell )}\centerdot I_R^{(\ell )}\ )\]
in such a way that the order of $E_L^{(\ell )}$ (resp. $E_R^{(\ell )}$~) 
is superior or equal to the order of $E_L^{(\ell-1 )}$ (resp. $E_R^{(\ell-1 )}$~), 
taking into account \cite{Hau} that, at the beginning,
\[
\Sigma _L^{(0)}\equiv D^{(0)}_L = E_L^{(0)}\centerdot I_L^{(0)}
\qquad \text{(resp.} \quad 
\Sigma _R^{(0)}\equiv D^{(0)}_R = E_R^{(0)}\centerdot I_R^{(0)}\ )\]
is such that $E_L^{(0)}=1$ (resp. $E_R^{(0)}=1$~) and
$D_L^{(0)}=I_L^{(0)}$ (resp. $D_R^{(0)}=I_R^{(0)}$~).\epr
\Ei
\vskip 11pt

\subsubsection{Definition: Corank of the singular locus}

Let $P(x_L,y_L,z_L)$ (resp. $P(x_R,y_R,z_R)$~) be the polynomial characterizing the singular locus $\Sigma _L$ (resp. $\Sigma _R$~).  The number of variables of
$P(x_L,y_L,z_L)$ (resp. $P(x_R,y_R,z_R)$~) is the corank of $\Sigma _L$ (resp. $\Sigma _R$~).  

This corank is inferior or equal to 3 according to \cite{A-V-G1}.
\vskip 11pt

\subsubsection{Ideals of differentiable functions}

From the beginning of chapter 2, differentiable \lr functions of $\phi _{G_{j_L}}^{(n)}(x_{g_{j_\delta }})$ (resp. $\phi _{G_{j_R}}^{(n)}(x_{g_{j_\delta }})$~) on the conjugacy class representatives $g_L^{(n)}[j_\delta ,m_{j_\delta }]$ (resp. $g_R^{(n)}[j_\delta ,m_{j_\delta }]$~) were taken into account (see section 2.1.1).

If we consider the set $\{\phi _{G_{j_L}}^{(n)}(x_{g_{j_\delta }})\}_{m_{j_\delta }}$ 
(resp. $\{\phi _{G_{j_R}}^{(n)}(x_{g_{j_\delta }})\}_{m_{j_\delta }}$~) of differentiable functions on the $j_\delta $-th conjugacy class of $G^{(n)}( {F^+_v})$ (resp.
$G^{(n)}( {F^+_{\o v}})$~) restricted to the singular loci $\Sigma _L[j_\delta ,m_{j_\delta }]$ (resp.
$\Sigma _R[j_\delta ,m_{j_\delta }]$~), we introduce an ideal of differentiable functions characterized by the polynomial $P_{j_R}(x_L,y_L,z_L)$ (resp. $P_{j_L}(x_R,y_R,z_R)$~) (see section 1.8).
\vskip 11pt

\subsubsection{Simple germs of differentiable functions}

Assume that the singular locus $\Sigma _L[j_\delta ,m_{j_\delta }]$ (resp.
$\Sigma _R[j_\delta ,m_{j_\delta }]$~) is given by a singular point of finite codimension.  Then, the corresponding simple germs of differentiable functions are the following \cite{A-V-G1}:
\begin{alignat*}{5}
A_K : \quad &f(x)=x^{k+1}\;, \qquad  &&\quad k &\ge &1\;, \\
D_K : \quad & f(x,y) = x^2y+y^{k-1}\;,  &&\quad k&\ge & 4\;, \\
E_6 : \quad & f(x,y)=x^3+y^4\;, \\
E_7 : \quad & f(x,y) = x^3+xy^3\;, \\
E_8 : \quad & f(x,y)=x^3+y^5\;.\end{alignat*}
They are described by the classical Dynkin diagrams.

Applying proposition 2.1.6, we shall now see how it is possible to generate a singular point of corank 1 and myltiplicity $k$ by a singularization of type $A_k$~.
\vskip 11pt

\subsubsection{Proposition}

{\em  A singularization of type $A_k$~, given by the germ $y_L=x_L^{k+1}$ (resp. 
$y_R=x_R^{k+1}$~) of differentiable functions 
$\phi _{G_{j_L}}^{(n)}(x_{g_{j_\delta }})$ (resp. $\phi _{G_{j_R}}^{(n)}(x_{g_{j_\delta }})$~) on the $j_\delta $-th conjugacy class of $G^{(n)}( {F^+_v})$ (resp.
$G^{(n)}({F^+_{\o v}})$~), is generated by the following sequence of contracting surjective morphisms:
\begin{align*}
\o\rho _L^{(\rm sing)} &:\quad D_L^{(0)}  \xrightarrow{\o\rho _L^{(1)}}  
x_L \cup D_L^{(1)} \xrightarrow{\quad} \cdots  \xrightarrow{\o\rho _L^{(k)}}   
x  _L^{k}\cup D_L^{(k)} 
  \xrightarrow{\o\rho _L^{(k+1)}}  
x  _L^{ k+1} \\
\text{(resp.} \quad
\o\rho _R^{(\rm sing)} &:\quad D_R^{(0)}  \xrightarrow{\o\rho _R^{(1)}}  
x_R \cup D_R^{(1)} \xrightarrow{\quad} \cdots  \xrightarrow{\o\rho _R^{(k)}}   
x  _R^{k}\cup D_R^{(k)} 
  \xrightarrow{\o\rho _R^{(k+1)}}  
x  _R^{ k+1} \ ) \end{align*}
restricted to the singular subloci $\Sigma _L^{(k)}=x_L^k$ (resp. 
$\Sigma _R ^{(k)}=x_R^k$~), $1\le k\le k+1$~, in such a way that:
\Bena
\item the contracting surjective morphism:
\begin{align*}
\o \rho _L^{(k)} : &\quad x_L^{k-1}\cup D_L^{(k-1)}\To x_L^k\;, \qquad 1\le k\le k+1\;, \\
\text{(resp.} \quad 
\o \rho _R^{(k)} : &\quad x_R^{k-1}\cup D_R^{(k-1)}\To x_R^k\ ),\end{align*}
restricted to the singular sublocus $\Sigma _L^{(k)}$ (resp. $\Sigma _R^{(k)}$~), 
is a fibre bundle whose fibre $D_L^{(k-1)}$ (resp. $D_R^{(k-1)}$~), 
being a normal crossings divisor on a real irreducible completion 
$F^+_{v_{j_\delta }^{j_\delta '}}$ (resp. $F^+_{\o v_{j_\delta }^{j_\delta '}}$~) according to definition 2.1.4, collapses into the germ $x^k_L$ (resp. $x^k_R$~). The fibre $D_L^{(k-1)}$ (resp. $D_R^{(k-1)}$~) is thus a contracting fibre).
\vskip 11pt

\item $\begin{array}[t]{rrrrr}
& \o\rho_L : \quad &\o\phi _{G_{j_L}}^{(n)}(x_{g_{j_\delta }}) \cup (D_L^{(0)},\cdots,D_L^{(k)}) & \To & \phi _{G_{j_L}}^{*(n)}(x_{g_{j_\delta }})\\
\text{(resp.} \quad 
& \o\rho_R : \quad &\o\phi _{G_{j_R}}^{(n)}(x_{g_{j_\delta }}) \cup (D_R^{(0)},\cdots,D_R^{(k)}) & \To & \phi _{G_{j_R}}^{*(n)}(x_{g_{j_\delta }})\ )\end{array}$

is the singularization morphism of the differentiable function $\phi _{G_{j_L}}^{(n)}(x_{g_{j_\delta }})$\linebreak (resp. $\phi _{G_{j_R}}^{(n)}(x_{g_{j_\delta }})$~) generating a germ $y_L=x_L^{k+1}$ (resp. $y_R=x_R^{k+1}$~) on it.
\vskip 11pt

\item $\begin{array}[t]{rrrrr}
&\o\rho _L^{(\rm is)}: & \o \phi _{G_{j_L}}^{(n)}(x_{g_{j_\delta }})\smallsetminus
\o\rho _L^{-1}(\Sigma _L^{(k+1)}=x_L^{k+1}) & \to & 
 \phi _{G_{j_L}}^{*(n)}(x_{g_{j_\delta }})\smallsetminus
(\Sigma _L^{(k+1)}=x_L^{k+1}) \\
\text{(resp.} \quad 
&\o\rho _R^{(\rm is)}: & \o \phi _{G_{j_R}}^{(n)}(x_{g_{j_\delta }})\smallsetminus
\o\rho _R^{-1}(\Sigma _R^{(k+1)}=x_R^{k+1}) & \to & 
 \phi _{G_{j_R}}^{*(n)}(x_{g_{j_\delta }})\smallsetminus
(\Sigma _R^{(k+1)}=x_L^{k+1}) \ )\end{array}$

is an isomorphism.
\Ee
}
\vskip 11pt

\subsubsection{Corollary} Every contractive surjective morphism $\o\rho _L^{(k)}$ of the sequence $\o\rho _L^{(\rm sing)}$ of proposition 2.1.10 provides a germ  $y_L=x_L^{k+1}$ of type $A_k$ at a singular point of corank 1 and multiplicity $k$ in such a way that this sequence $\o\rho _L^{(\rm sing)}$ of singularizations generates the following sequence of simple germs:
\[ x^2_L \subset x^3_L \subset \cdots\subset x^k_L \subset x^{k+1}_L\subset \cdots\;, \qquad 2\le k<\infty \]
characterized by increasing finite multiplicities.
\vskip 11pt

\subsection{The versal deformation}

\subsubsection{Generalities on the versal deformation}

The versal deformation will be considered for germs of differentiable functions $\phi ^{(n)}_{G_{j_L}}(x_{g_{j_\delta }})$ (resp. $\phi ^{(n)}_{G_{j_R}}(x_{g_{j_\delta }})$~) of the \lr semisheaf $\theta _{G^{(n)}_L}$ (resp. $\theta _{G^{(n)}_L}$~).  These \lr differentiable functions will be written $\phi _{j_{\delta _L}}(x_L)$ (resp. $\phi _{j_{\delta _R}}(x_R)$~) where $x_L$ (resp. $x_R$~) is a $n$-tuple of numbers 
\begin{align*}
x_L & = (x_{1_L},x_{2_L},\cdots,x_{n_L})\in (F^+_L)^n\\
\text{(resp.}
\quad
x_R & = (x_{1_R},x_{2_R},\cdots,x_{n_R})\in (F^+_R)^n\ ).\end{align*}
So, on the $j_\delta $-th conjugacy class of the algebraic semigroup $G^{(n)}( {F^+_v})
\equiv T_n( {F^+_v})$ (resp. $G^{(n)}( {F^+_{\o v}})
\equiv T^t_n( {F^+_{\o v}})$~), there is:
\Bena
\item a set of \lr differentiable functions $\phi _{j_\delta ,m_{j_\delta }}(x_L)$ 
(resp. $\phi _{j_\delta ,m_{j_\delta }}(x_R)$~) which are the sections of the \lr semisheaf $\theta _{G^{(n)}_L}$ (resp. $\theta _{G^{(n)}_R}$~).

\item a germ $\phi _{j_\delta }(\omega _L)$ (resp. $\phi _{j_\delta }(\omega _R)$~) (or a set of germs) of differentiable functions, where $\omega _L$ \resp{$\omega _R$} denotes a $m$-tuple of numbers, $1\le m\le 3$~.
\Ee
This germ $\phi _{j_\delta }(\omega _L)$ (resp. $\phi _{j_\delta }(\omega _R)$~) is assumed to be:
\Bi
\item simple and of corank $\le 3$~;
\item associated with an isolated singularity of order $k$~.
\Ei

As the corank of the singularity is not superior to 3~, the $n$-tuple of  numbers of $(F^+_L)^n$ (resp. $(F^+_R)^n$~), restricted to a small domain centered on the singularity, will be rewritten according to:
\begin{align*}
x'_L &= (x_{1_L},\cdots,x_{n_L-m_L},\omega _{1_L},\cdots,\omega _{m_L})\;, \qquad 1\le m_L\le 3\\
\text{(resp.} \quad
x'_R &= (x_{1_R},\cdots,x_{n_R-m_R},\omega _{1_R},\cdots,\omega _{m_R})\ ).
\end{align*}

The finite determinacy was first envisaged for germs having an isolated singularity: this is the pioneer work of H. Grauert and H. Kerner \cite{G-K}, R. Thom, \cite{Tho1}, \cite{Tho2}, \cite{Lev}, J. Mather, \cite{Mat1}, \cite{Mat2}, \cite{Mat3}, V.I. Arnold \cite{Arn1}, J.C. Tougeron \cite{Tou}, B. Malgrange \cite{Mal}, and others.

Afterwards, this problem was generalized to functions having a fixed analytic set 
$\Sigma $ as critical set.  If $I$ denotes the ideal of functions on this critical set 
$\Sigma $~, the finite $I$-determinacy of these functions and their versal 
$I$-unfoldings were considered and proved in \cite{Sie} and \cite{Pel}.
\vskip 11pt

With this in view, the preparation theorem and the versal deformation will be recalled for germs $\phi _{j_\delta}(\omega _L)$ (resp. $\phi _{j_\delta}(\omega _R)$~) of differentiable functions having an isolated singularity of corank 1 and order $k$~.
\vskip 11pt

\subsubsection{The division theorem}

Let $x'_L = (x_{1_L},\cdots,x_{n_L-1},\omega _{L}) $ (resp. $x'_R = (x_{1_R},\cdots,x_{n_R-1},\omega _{R})$~) denote the  coordinates of $(F^+_L)^n$ (resp. $(F^+_R)^n$~).  

A germ   $\phi _{j_\delta}(\omega _L)$ (resp. $\phi _{j_\delta}(\omega _R)$~) has a singularity of corank 1 (then, $m=1$~) and order $k$ in $\omega _L$ (resp. $\omega _R$~) if $\phi _{j_\delta}(0,\omega _L) =\omega ^k_L\ e _{j_\delta}(\omega _L)$ (resp. $\phi _{j_\delta}(0,\omega _R) =\omega ^k_R\ e _{j_\delta}(\omega _R)$~), where
$e_{j_\delta }(\omega _L)$ (resp. $e_{j_\delta }(\omega _R)$~) is a differentiable unit, i.e. verifying $e_{j_{\delta _L}}(0)\neq 0$ (resp. $e_{j_{\delta _R}}(0)\neq 0$~).

Let $\theta [\omega _L]$ (resp. $\theta [\omega _R]$~) be the algebra of polynomials in $\omega _L$ (resp. $\omega _R$~) with coefficients $a_{ij_\delta }(x_L)$ (resp. $a_{ij_\delta }(x_R)$~), being ideals of functions defined on a domain $D_L\subset B_L$ (resp. $D_R\subset B_R$~) where:
\Bi
\item $B_L$ (resp. $B_R$~) is an upper (resp. lower) half open ball centered on $\omega _L$ (resp. $\omega _R$~) in $\phi _{j_\delta }(x_L)$ (resp. $\phi _{j_\delta }(x_R)$~);

\item $x_L=(x_{1_L},\cdots, x_{n_L-1})$ (resp. $x_R=(x_{1_R},\cdots, x_{n_R-1})$~) is the $(n-1)$-tuple of $x'_L$ (resp. $x'_R$~) in $(F^+_L)^{n-1}$ (resp. $(F^+_R)^{n-1}$~).
\Ei

If the germ $\phi _{j_\delta }(\omega _L)$ (resp. $\phi _{j_\delta }(\omega _R)$~) has order $k$ in $\omega _L$ (resp. $\omega _R$~), then, for every $n$-dimensional differentiable function (germ) $f_{j_{\delta _L}}$ (resp. $f_{j_{\delta _R}}$~), there exists a $(n-1)$-dimensional differentiable function (germ) $q_{j_{\delta _L}}$ (resp.
$q_{j_{\delta _R}}$~) and a polynomial:
\begin{align*}
R_{j_{\delta _L}}
&= \sum^s_{i=1}a_{ij_\delta }(x_L)\ \omega ^i_{j_{\delta _L}}\in \theta [\omega _L]\\
\text{(resp.} \quad
R_{j_{\delta _R}}
&= \sum^s_{i=1}a_{ij_\delta }(x_R)\ \omega ^i_{j_{\delta _R}}\in \theta [\omega _R]\ )
\end{align*}
with degree $s<k$ such that:
\begin{align*}
f_{j_{\delta _L}} &= \phi _{j_{\delta _L}}(\omega _L)\centerdot q_{j_{\delta _L}}+R_{j_{\delta _L}}\\
\text{(resp.} \quad
f_{j_{\delta _R}} &= \phi _{j_{\delta _R}}(\omega _R)\centerdot q_{j_{\delta _R}}+R_{j_{\delta _r}}\ )
\end{align*}
be the division theorem (adapted to the left and right cases) introduced by B. Malgrange in \cite{Mal}.  The Malgrange division theorem, closely related to the version of J. Mather \cite{Mat1}, \cite{Mat2}, \cite{Mat3}, is the differentiable version of the Weierstrass division theorem \cite{G-R}.
\vskip 11pt

\subsubsection{The division theorem (corank 2 case)}

The division theorem, recalled in section 2.2.2 for germs of differentiable functions having an isolated singularity of corank 1, can easily be generalized to germs $\phi _{j_\delta }(\omega _{1_L},\omega _{2_L})$ (resp. $\phi _{j_\delta }(\omega _{1_R},\omega _{2_R})$~) having an isolated singularity of corank 2.

Indeed, a germ $\phi _{j_\delta }(\omega _{1_L},\omega _{2_L})$ (resp. $\phi _{j_\delta }(\omega _{1_R},\omega _{2_R})$~) has a singularity of corank 2 and  order $k$ in $\omega _{1_L}$ and $\omega _{2_L}$ (resp. $\omega _{1_R}$ and $\omega _{2_R}$~) if:
\Bean
\item \bt[t]{ll}
& $\phi _{j_\delta }(\omega _{1_L},\omega _{2_L})=P_{j_\delta }(\omega _{1_L},\omega _{2_L})\ e_{j_\delta }(\omega _{1_L},\omega _{2_L})$\\
(resp. \quad 
& $\phi _{j_\delta }(\omega _{1_R},\omega _{2_R})= P_{j_\delta }(\omega _{1_R},\omega _{2_R})\ e_{j_\delta }(\omega _{1_R},\omega _{2_R})$\ )\te

where:
\Bi
\item $P_{j_\delta }(\omega _{1_L},\omega _{2_L})$ (resp.
$P_{j_\delta }(\omega _{1_R},\omega _{2_R})$~) is a polynomial of degree $k$ in $\omega _{1_L}$ or in $\omega _{2_L}$ (resp. $\omega _{1_R}$ or in $\omega _{2_R}$~);
\item $e_{j_\delta }(\omega _{1_L},\omega _{2_L})$ (resp. $e_{j_\delta }(\omega _{1_R},\omega _{2_R})$~) is a differentiable unit.
\Ei

\item the polynomial $R_{j_{\delta _L}}\in \theta [\omega _{1_L},\omega _{2_L}]$
(resp. $R_{j_{\delta _R}}\in \theta [\omega _{1_R},\omega _{2_R}]$~) of the quotient algebra has degree $\ell <k$~, in such a way that this quotient algebra be a finitely generated tensorial space of type $(0,2)$ and dimension $\ell\le k$~.
\Ee
\vskip 11pt

\subsubsection[singularization of the semisheaf $\theta _{G^{(n)}_L}$ (resp. $\theta _{G^{(n)}_R}$~)]{\boldmath singularization of the semisheaf $\theta _{G^{(n)}_L}$ (resp. $\theta _{G^{(n)}_R}$~)}

Let  $\theta _{G^{(n)}_L}$ (resp. $\theta _{G^{(n)}_R}$~) be the semisheaf of \lr smooth 
differentiable functions $\phi _{j_\delta }(x_L)$ (resp. $\phi _{j_\delta }(x_R)$~).

And, let
\begin{align*}
\o\rho _{G_L} &: \quad \theta _{G^{(n)}_L} \To \theta^* _{G^{(n)}_L} \\
\text{(resp.} \quad
\o\rho _{G_R} &: \quad \theta _{G^{(n)}_R} \To \theta^* _{G^{(n)}_R} \ )
\end{align*}
be the singularization of $ \theta _{G^{(n)}_L}$ (resp. $ \theta _{G^{(n)}_R}$~), in the sense of proposition 2.1.6, in such a way that
$ \theta^*_{G^{(n)}_L}$ (resp. $ \theta ^*_{G^{(n)}_R}$~) be the semisheaf whose sections
$\phi ^{*(n)}_{G_{j_{\delta_L}}}(x_{g_{j_\delta }})$ (resp. $\phi ^{*(n)}_{G_{j_{\delta _R}}}(x_{g_{j_\delta }})$~)
are the differentiable functions $\phi ^{(n)}_{G_{j_{\delta _L}}}(x_{g_{j_\delta }})$ (resp. $\phi ^{(n)}_{G_{j_{\delta _R}}}(x_{g_{j_\delta }})$~) endowed with germs $\phi _{j_\delta }(\omega _L )$ (resp. $\phi _{j_\delta }(\omega _R )$~) having degenerate singularities of corank 1~.
\vskip 11pt


\subsubsection{Proposition (Versal deformation)}

{\em The versal deformation of the semisheaf $\theta ^*_{G_L^{(n)}}$ (resp. $\theta ^*_{G_R^{(n)}}$~) is given by the contracting fibre bundle:
\begin{align*}
D_{S_L} : \quad & (\theta ^*_{G_L^{(n)}}\smallsetminus \theta _L(a))\times \theta [\omega _L] \To \theta ^*_{G_L^{(n)}}\\
\text{(resp.} \quad 
D_{S_R} : \quad & (\theta ^*_{G_R^{(n)}}\smallsetminus 
\theta _R(a))\times \theta [\omega _R] \To \theta ^*_{G_R^{(n)}}\ )\end{align*}
where:
\Bi
\item $\theta _L(a)$ (resp. $\theta _R(a)$~) is the (semi)sheaf of ideals 
$a_{ij_\delta }(x_L)$ (resp. $a_{ij_\delta }(x_R)$~) of differentiable functions as introduced in section 2.2.2.

\item $\theta [\omega _L]$ (resp. $\theta [\omega _R]$~) is the algebra of polynomials $R_{j_{\delta _L}}$ (resp. $R_{j_{\delta _R}}$~) introduced in section 2.2.2
\Ei
and whose fibre 
\begin{align*}
\theta _{S_L} &= \{\theta ^1(\omega ^1_L),\cdots,\theta ^1(\omega ^i_L),\cdots, \theta ^1(\omega ^s_L)\}\\
\text{(resp.} \quad 
\theta _{S_R} &= \{\theta ^1(\omega ^1_R),\cdots,\theta ^1(\omega ^i_R),\cdots, \theta ^1(\omega ^s_R)\}\ )\end{align*}
is the family of the (semi-)sheaves of the \lr base $S_L$ (resp. $S_R$~) of the versal deformation.
}
\vskip 11pt


\bpr \quad $\theta _{f_L}$ (resp. $\theta _{f_R}$~) is the (semi)sheaf of differentiable functions
\begin{align*}
f_{j_{\delta _L}} &= \phi _{j_\delta }(\omega _L)\centerdot q_{j_{\delta _L}}+R_{j_{\delta _L}}\\
\text{(resp.} \quad 
f_{j_{\delta _R}} &= \phi _{j_\delta }(\omega _R)\centerdot q_{j_{\delta _R}}+R_{j_{\delta _R}}\ )\end{align*}
introduced in section 2.2.2.

The polynomials $R_{j_{\delta _L}}\in \theta [\omega _L]$ (resp. $R_{j_{\delta _R}}\in \theta [\omega _R]$~) have as coefficients the ideals of functions $a_{ij_\delta }(x_L)$ (resp. $a_{ij_\delta }(x_R)$~) on the differentiable functions $\phi _{j_\delta }(x_L)$ (resp. $\phi _{j_\delta }(x_R)$~): so we have
\begin{align*}
a_{j_{\delta _L}}(x_L) & \subset \phi _{j_\delta} (x_L)\\
\text{(resp.} \quad
a_{j_{\delta _R}}(x_R) & \subset \phi _{j_\delta} (x_R)\ ).\end{align*}
Then, it appears that:
\begin{align*}
\theta ^*_{G^{(n)}_L}\times \theta _{S_L} &= (\theta ^*_{G^{(n)}_L}\smallsetminus \theta _L(a))\times \theta [\omega _L]\\
\text{(resp.} \quad
\theta ^*_{G^{(n)}_R}\times \theta _{S_R} &= (\theta ^*_{G^{(n)}_R}\smallsetminus \theta _R(a))\times \theta [\omega _R]\ ),\end{align*}
and, thus, that $\theta _{S_L}$ (resp. $\theta _{S_R}$~) is the fibre of the contracting 
fibre bundle $D_{S_L}$ (resp. $D_{S_R}$~) rewritten as follows \cite{Ste}:
\begin{align*}
D_{S_L} &: \quad \theta ^*_{G^{(n)}_L}\times \theta _{S_L}\To \theta ^*_{G^{(n)}_L}\\
\text{(resp.} \quad
D_{S_R} &: \quad \theta ^*_{G^{(n)}_R}\times \theta _{S_R}\To \theta ^*_{G^{(n)}_R}\ ).
\tag*{\eop}\end{align*}
\vskip 11pt

\subsubsection{Proposition}

{\em
Let $\theta ^{\rm vers}_{G^{(n)}_L} = \theta ^*_{G^{(n)}_L}\times \theta _{S_L}$ (resp.
$\theta ^{\rm vers}_{G^{(n)}_R} = \theta ^*_{G^{(n)}_R}\times \theta _{S_R}$~) denote the semisheaf unfolded from $\theta ^*_{G^{(n)}_L} $ (resp. $\theta ^*_{G^{(n)}_R} $~): it is the total space of the fibre bundle $D_{S_L}$ (resp. $D_{S_R}$~).

Let $\Sigma _{G^{(n)}_L}$ (resp. $\Sigma _{G^{(n)}_R}$~) be the singular locus of
$\theta ^*_{G^{(n)}_L} $ (resp. $\theta ^*_{G^{(n)}_R} $~): it is the sheaf $\theta _{\phi _{\omega _L}}$ (resp. $\theta _{\phi _{\omega _R}}$~) of germs $\phi _{j_\delta }(\omega _L)$ (resp. $\phi _{j_\delta }(\omega _R)$~) of differentiable functions.  

Then, we have that:
\Bean
\item the unfolded image $D^{-1}_{S_L}(\theta _{\phi _{\omega _L}})$ (resp.
$D^{-1}_{S_R}(\theta _{\phi _{\omega _R}})$~) of the singular locus $\Sigma _{G^{(n)}_L}$ (resp. $\Sigma _{G^{(n)}_R}$~) is the (semi)sheaf $\theta _{f_L}$ (resp. $\theta _{f_R}$~):
\begin{align*}
\theta _{f_L}=D^{-1}_{S_L}(\theta _{\phi _{\omega _L}})\\
\text{(resp.} \quad
\theta _{f_R}=D^{-1}_{S_R}(\theta _{\phi _{\omega _R}})\ );\end{align*}

\item $\begin{array}[t]{rl}
D^{-1}_{S_L}(\theta ^*_{G^{(n)}_L}\smallsetminus \theta _{\phi _{\omega _L}})
&= (\theta ^*_{G^{(n)}_L}\smallsetminus D^{-1}_{S_L}(\theta _{\phi _{\omega _L}}))\smallsetminus \theta _{S_L}\\
\text{(resp.} \quad
D^{-1}_{S_R}(\theta ^*_{G^{(n)}_R}\smallsetminus \theta _{\phi _{\omega _R}})
&= (\theta ^*_{G^{(n)}_R}\smallsetminus D^{-1}_{S_R}(\theta _{\phi _{\omega _R}}))\smallsetminus \theta _{S_R}
\ ).\end{array}$
\Ee
}
\vskip 11pt

\bpr \quad By versal deformation, $\theta _{S_L}$ (resp. $\theta _{S_R}$~) is the fibre of the sheaf $\theta _{\phi _{\omega _L}}$ (resp. $\theta _{\phi _{\omega _R}}$~) of germs of differentiable functions, outside of which $D_{S_L}$ (resp. $D_{S_R}$~) is an isomorphism: this is reflected by the equality b) of this proposition.
\vskip 11pt

\subsubsection{Definition}

The quotient algebra $\theta [\omega _L]$ (resp. $\theta [\omega _R]$~) of germs $\phi _{j_\delta }(\omega _L)$ (resp. $\phi _{j_\delta }(\omega _R)$~) having a singularity of corank 1 and multiplicity $(k-1)$ is the quotient of the algebra $\Es _{n_L}$ (resp. $\Es _{n_R}$~) of function germs (generally, it is the algebra of integer power series) by the graded ideal $I_{\phi _{\omega _L}}$ (resp. $I_{\phi _{\omega _R}}$~) of $\phi _{j_L}(\omega _L)$ (resp. $\phi _{j_R}(\omega _R)$~):
\begin{align*}
\theta [\omega _L] &= \Es _{n_L}\big/ I_{\phi _{\omega _L}}\;, \qquad n=1\\
\text{(resp.}\quad
\theta [\omega _R] &= \Es _{n_R}\big/ I_{\phi _{\omega _R}}\ )\end{align*}
where:
\begin{align*}
I_{\phi _{\omega _L}} &= \Es _{n_L} \ \langle \phi _L^{(1)},\cdots,\phi _L^{(k-1)}\rangle\\
\text{(resp.} \quad
I_{\phi _{\omega _R}} &= \Es _{n_R} \ \langle \phi _R^{(1)},\cdots,\phi _R^{(k-1)}\rangle\ )
\end{align*}
is generated by the partial derivatives $\phi _L^{(k-1)}$ (resp. $\phi _R^{(k-1)}$~) of $\phi _{j_\delta }(\omega _L)$ (resp. $\phi _{j_\delta }(\omega _R)$~), \cite{A-G-L-V}.

The quotient algebra is thus finitely generated: it is composed of the polynomials $R_{j_{\delta _L}}\in \theta [\omega _L]$ (resp. $R_{j_{\delta _R}}\in \theta [\omega _R]$~) (see section 2.2.2), which generate vector (semi)spaces of dimension $s<k$~, and proceeds from a set of contracting morphisms extended those considered in the singularization processes as developed in proposition 2.1.6.
\vskip 11pt

\subsubsection{Proposition}

{\em
\Bena
\item The versal deformations of the germ $y_L=\omega _L^{k+1}$ (resp. $y_R=\omega _R^{k+1}$~) of differentiable functions  $\phi ^{(n)}_{G_{j_{\delta _L}}}(x_{g_{j_\delta }})$ (resp. $\phi ^{(n)}_{G_{j_{\delta _R}}}(x_{g_{j_\delta }})$~) is generated by a sequence
\begin{align*}
D_{S_L^{(k+1)}} &= ( D^{(1)}_{S_L^{(k+1)}} ,\cdots,D^{(i)}_{S_L^{(k+1)}} ,\cdots,
D^{(k-1)}_{S_L^{(k+1)}}) \\
\text{(resp.} \quad
D_{S_R^{(k+1)}} &= ( D^{(1)}_{S_R^{(k+1)}} ,\cdots,D^{(i)}_{S_R^{(k+1)}} ,\cdots,
D^{(k-1)}_{S_R^{(k+1)}}) \ )\end{align*}
of $(k-2)$ contracting morphisms $ D^{(i)}_{S_L^{(k+1)}} $ (resp. $ D^{(i)}_{S_R^{(k+1)}} $~) extending the sequence of contracting surjective morphisms $\o\rho _L^{(\rm sing)}$ (resp. $\o\rho _R^{(\rm sing)}$~) of singularization according to the following diagram:
\begin{tiny}\begin{align*}
 \o\rho _L^{(\rm sing)}: \quad &D_L^{(0)} \to   \cdots\\
&  
\begin{array}{cccccccc}
 \xrightarrow{\o\rho _L^{(3)}} &\omega _L^{(3)}\cup D_L^{(3)}
& \to  \cdots \xrightarrow{\o\rho _L^{(i)}} &\omega _L^{(i)}\cup D_L^{(i)}
& \to  \cdots \xrightarrow{\o\rho _L^{(k)}} &\omega _L^{(k)}\cup D_L^{(k)}
&  \xrightarrow{\o\rho _L^{(k+1)}}& \omega _L^{(k+1)} \\
 &\cup && \cup && \cup && \cup\\
 & D^{'(1)}_L && D^{'(i-2)}_L && D^{'(k-2)}_L && D^{'(k-1)}_L \\
 &\begin{CD}@VV{D_{S_L^{(k+1)}}^{(1)}}V\end{CD} &&
  \begin{CD}@VV{D_{S_L^{(k+1)}}^{(i-2)}}V\end{CD} &&
  \begin{CD}@VV{D_{S_L^{(k+1)}}^{(k-2)}}V\end{CD} &&
  \begin{CD}@VV{D_{S_L^{(k+1)}}^{(k-1)}}V\end{CD} \\
 &\begin{array}[t]{l}
\omega _L^3\\ +a_1\omega ^1_L\end{array} & \To & 
 \begin{array}[t]{l}\omega ^i_L\\ +\sum\limits_{i=1}^{i-2}a_i\omega ^i_L \end{array}
& \To & 
 \begin{array}[t]{l}\omega ^k_L\\ +\sum\limits_{i=1}^{k-2}a_i\omega ^i_L \end{array}
& \To & 
 \begin{array}[t]{l}\omega ^{k+1}_L\\ +\sum\limits_{i=1}^{k-1}a_i\omega ^i_L \end{array}\end{array}\end{align*}\end{tiny}
(idem for the $R$-case)

where $D^{'(i-2)}_L$ is a normal crossings divisor on a real (or a set of) irreducible completion(s) $F^+_{v^{j'_\delta }_{j_\delta }}$ mapped onto the neighbourhood of the germ $y_L=\omega ^i_L$~.
\vskip 11pt

\item If $\phi ^{*(n,i)}_{G_{j_{\delta _L}}}(x_{g_{j_\delta }})$ denotes the function having a germ $y_L=\omega ^i_L$ of codimension $(i-2)$~, then, the $D^{(i-2)}_{S_L^{(k+1)}}$ contracting morphism corresponds to the contracting fibre bundle:
\[ D^{(i-2)}_{S_L^{(k+1)}}: \quad \phi ^{*(n,i)}_{G_{j_{\delta _L}}}(x_{g_{j_\delta }}) \times \omega _L^{i-2} \To \phi ^{*(n,i)}_{G_{j_{\delta _L}}}(x_{g_{j_\delta }})\]
in such a way that:
\Bi
\item $D^{(i-2)}_{S_L^{(k+1)}}\subset D_{S_L^{(k+1)}}$~;
\item $\omega _L^{i-2}$ is the contracting fibre, i.e. the divisor $D_L^{'(i-2)}$ in the neighbourhood of\linebreak $y_L=\omega _L^i$ on $a_{(i-2)j}(x_L)\subset \phi ^{*(n,i)}_{G_{j_{\delta _L}}}(x_{g_{j_{\delta_L} }})$~.
\Ei\Ee
}
\vskip 11pt

\bpr
\Bena
\item The versal deformation of a germ $y_L=\omega _L^{k+1}$ is an extension of its singularization as described in proposition 2.1.10.  Indeed, to the $i$-th contracting surjective morphism of singularization:
\[ \o\rho _L^{(i)}:\quad \omega _L^{i-1}\cup D_L^{(i-1)} \To \omega _L^i\;, 
\qquad 1\le i\le k+1\;, \]
introduced in proposition 2.1.10 as a contracting fibre bundle whose fibre 
$D_L^{(i-1)}$ (which is a normal crossings divisor) collapses into one point on the 
germ $y_L=\omega _L^{i-1}$~, corresponds the contracting fibre bundle 
$D^{(i-2)}_{S_L^{(k+1)}}$ of the versal deformation $D_{S_L^{(k+1)}}$~:
\[ D^{(i-2)}_{S_L^{(k+1)}} : \quad \phi ^{*(n,i)}_{G_{j_L}}(x_{g_{j_L}})\times \omega _L^{i-2} \To \phi ^{*(n,i)}_{G_{j_L}}(x_{g_{j_L}})\]
in such a way that the divisor $D_L^{'(i-2)}$ be projected in the neighbourhood of the singular germ $y_L=\omega _L^{k+1}$ where it is rewritten $\omega _L^{i-2}$~.

So, from the third contracting surjective morphism of singularization $\o\rho _L^{(3)}$~, where the singularity becomes degenerated, we can associate to each contracting surjective morphism of singularization $\o\rho _L^{(i)}$~, $3\le i\le k+1$~, a contracting fibre bundle $D^{(i-2)}_{S_L^{(k+1)}}$~, $i-2\le i\le k+1$~, of versal deformation.
\vskip 11pt

\item And, to the sequence of singularizations:
\[ \o\rho _L^{(3)}\;,\; \cdots\;,\;  \o\rho _L^{(i)}\;,\; \cdots\;,\; 
 \o\rho _L^{(k)} \;,  \o\rho _L^{(k+1)}\;, \]
corresponds the sequence of versal subdeformations:
\[ D^{(1)}_{S_L^{(k+1)}}\subset \cdots \subset 
D^{(i-2)}_{S_L^{(k+1)}}\subset \cdots \subset 
D^{(k-2)}_{S_L^{(k+1)}}\subset 
D^{(k-1)}_{S_L^{(k+1)}}\;,\]
which are embedded and which correspond to the sequence of $(k-1)$ embedded vector 
sub(semi)spaces generated by the polynomials:
\[R^{(1)}_{j_{\delta _L}} \subset \cdots \subset
R^{(i-2)}_{j_{\delta _L}} \subset \cdots \subset
R^{(k-2)}_{j_{\delta _L}} \subset 
R^{(k-1)}_{j_{\delta _L}} \;,\]
where $R^{(i-2)}_{j_{\delta _L}} $ is the truncated polynomial of the quotient algebra $\theta [\omega _L]$ introduced in section 2.2.2 and given by:
\[ R^{(i-2)}_{j_{\delta _L}}=\sum_{i=1}^{i-2} a_{ij_\delta }(x_L)\ \omega ^i_{j_{\delta _L}}\;.\]
\vskip 11pt

\item The order of the divisors $D_L^{'(i-2)}$~, projected in the neighbourhood of the singular germ $y_L=\omega _L^{k+1}$~, increases in function of the increase of the dimension of the generated vector sub(semi)spaces $R_j^{(i-2)}$~, $i-2\le i\le k+1$~, of the versal deformation, because the space around the singularity becomes ``over'' compact when the dimension of the versal deformation increases.  In fact, it will be proved in the following that the geometry of the space around the singularity deformed by versal unfolding is spherical.\epr
\Ee
\vskip 11pt

\subsubsection{Corollary}

{\em
The sequence $D_{S_L^{(k+1)}}$ (resp. $D_{S_R^{(k+1)}}$~) of contracting morphisms generating the versal deformations of singular germs of corank 1 and multiplicity $\ge 1$ explains why the quotient algebra $\theta [\omega _L]$ (resp. $\theta [\omega _R]$~) of the versal deformation is finitely determined \cite{Mat3}.
}\vskip 11pt


\subsection{The geometry of the versal deformation}

The geometry of the versal deformation will be envisaged for the sections $\phi ^{(n)}_{G_{j_L}}(x_{g_{j_\delta }})$ (resp. $\phi ^{(n)}_{G_{j_R}}(x_{g_{j_\delta }})$~) of the \lr semisheaf $\theta _{G^{(n)}_L}$ (resp. $\theta _{G^{(n)}_R}$~).  These sections are assumed to be differentiable functions having singular germs $\phi _{j_\delta }(\omega _L)$ (resp. $\phi _{j_\delta }(\omega _R)$~) of corank $m\le 3$ and multiplicity ``~$i$~'', $1\le i\le n$ (see section 2.2.1).

Let $\Sigma _{\phi ^{(n)}_{G_{j_L}}}$ (resp. $\Sigma _{\phi ^{(n)}_{G_{j_R}}}$~) denote the singular locus of a singular germ $\phi _{j_\delta }(\omega _L)$ (resp. $\phi _{j_\delta }(\omega _R)$~) of corank $m$ and multiplicity ``~$i$~'' and let $D_{\Sigma _{\phi _{G_{j_L}}}}$ (resp. $D_{\Sigma _{\phi _{G_{j_R}}}}$~) be the neighbourhood of this singular locus whose curvature is affected by the singularity.
\vskip 11pt

\subsubsection{Proposition}

{\em The geometry is hyperbolic on the neighbourhood $D_{\Sigma _{\phi _{G_{j_L}}}}$
(resp. $D_{\Sigma _{\phi _{G_{j_R}}}}$~) of the singular locus 
$\Sigma _{\phi ^{(n)}_{G_{j_L}}}$ (resp. $\Sigma _{\phi ^{(n)}_{G_{j_R}}}$~) of a not 
unfolded singular germ of corank $m\le 3$ and multiplicity ``~$i$~'', $1\le i\le n$~.
}
\vskip 11pt

\bpr \quad The main idea consists in showing that there is a deviation to euclidicity in
$D_{\Sigma _{\phi _{G_{j_L}}}}$
(resp. $D_{\Sigma _{\phi _{G_{j_R}}}}$~) of the differentiable function 
$\phi ^{(n)}_{G_{j_L}}(x_{g_{j_\delta }})$ (resp. $\phi ^{(n)}_{G_{j_R}}(x_{g_{j_\delta }})$~) of dimension $n$~.

This deviation to euclidicity can be evaluated by searching the conditions to which must satisfy the metric $ds^2=g_{ij}\ du^i\ du^j$ in the neighbourhood $D_{\Sigma _{\phi _{G_{j_L}}}}$
(resp. $D_{\Sigma _{\phi _{G_{j_R}}}}$~) of the singular locus.

General conditions in the Euclidean and non Euclidean cases were developed by E. Cartan in his classical book ``leçons sur la géométrie des espaces de Riemann'' \cite{Car} to which we refer.

The developments will be envisaged for the left and right cases without distinction: thus, we omit the ``~$L$'' and ``~$R$~'' and we simplify the notations:

\bt[t]{cccccc}
\textbullet & $D_{\Sigma _{\phi _{G_{j_L}}}}$&
and& $D_{\Sigma _{\phi _{G_{j_R}}}}$ & become & $D_\Sigma $~;\\
\textbullet & $\Sigma _{\phi ^{(n)}_{G_{j_L}}}$ &and&$\Sigma _{\phi ^{(n)}_{G_{j_R}}}$ &
become & $\Sigma $~;\\
\textbullet &  $\phi ^{(n)}_{G_{j_L}}(x_{g_{j_\delta }})$ &and&$\phi ^{(n)}_{G_{j_R}}(x_{g_{j_\delta }})$ & become & $\phi ^{(n)}$~.\te
\vskip 11pt

\Bena
\item Consider first the Euclidean case and remark that there does not exist in general a coordinate system giving to the Euclidean space a fix metric.  

Let $M$ be a point of coordinates $(u^1,\cdots,u^n)$ on $\phi ^{(n)}\smallsetminus D_\Sigma \smallsetminus \Sigma $ and let $(\vec e_1 ,\cdots,\vec e_n )$ be the basis vectors of the proper referential centred on the point $M$~, such that the components:
\Bi
\item $dM=du^i\ \vec e_i $~, $1\le i\le n$~,
\item $d\vec e_i
=\Gamma ^k_{ij}\ du^j\ \vec e_k$ where $\Gamma ^k_{ij}
=\half \L( \F{\partial g_{ik}}{\partial u^j} 
+ \F{\partial g_{ik}}{\partial u^i}
-\F{\partial g_{ij}}{\partial u^k}\R)$~,
\Ei
do not belong to a stratum on $(D_\Sigma \cup\Sigma )$~.

The integrability conditions of $d\vec e_i=\Gamma ^k_{ij}\ du^j\ \vec e_k$ are precisely the searched condition to which  the $g_{ij}$ must satisfy.

These integrability conditions can be obtained geometrically as follows.

To each point $M$ will correspond a point $P$ in the neighbourhood of $M$~.  This point $P$ is defined by its coordinates $(x^1,\cdots,x^n)$ with respect to the referential of the point $M$~.

The components of the differential of the point $P$ are given by:
\[Dx^i=dx^i+du^i+x^k\ \Gamma ^i_{kr}\ du^r\]
such that
\[ D_rx^i = \F{\partial x_i}{\partial u^r}+\delta ^i_r+x^k\ \Gamma ^i_{kr}\qquad
\text{with} \quad \delta ^i_r \begin{cases} =0 \;, \quad & \text{if}\ i\neq r\;,\\
=1 \;, \quad & \text{if}\ i=r\,.\end{cases}\]

Let then $M'$~, $M''$ and $M^{'''}$ be the points obtained as follows: the first 
$M'$ is obtained by increasing the coordinate $u^r$ by an infinitesimal parameter 
$\alpha $~, the second $M''$ by increasing the coordinate $u^s$ by an infinitesimal 
parameter $\beta $~, and the third $M^{'''}$ by increasing the coordinate $u^r$ by 
$\alpha$  and the coordinate $u^s$ by $\beta $~.

Let $P'$~, $P''$ and $P^{'''}$ be the points corresponding to the points $M'$~, $M''$ and $M^{'''}$ respectively.

The infinitesimal small vector $\overrightarrow{PP'}$ has the ``~$i$~'' contravariant components given by $\alpha D_rx^i$~.

On the other hand, to the elementary variations $\o{MM^{'''}}=\{\delta u^1=0,\cdots,\delta u^r=\alpha ,\cdots,\delta u^s=\beta ,\cdots,\delta u^n=0\}$ of the points $M$ will correspond the infinitesimal small vector $\overrightarrow{P^{''}P^{'''}}-\overrightarrow{PP'}$ whose components $\alpha \beta D_sD_rx^i$ are given by:
\begin{align*}
D_sD_rx^i
&= \F{\partial D_rx^i}{\partial u^s} + D_rx^k\Gamma^i_{ks}\\[11pt]
&= \F{\partial^2x^i}{\partial u^r\ \partial u^s}+\F{\partial x^k}{\partial u^s}\ \Gamma ^i_{kr}+ x^k\ \F{\partial\Gamma ^i_{kr}}{\partial u^s}+\F{\partial x^k}{\partial u^r}\ \Gamma ^i_{ks}+\Gamma ^i_{rs}+x^k\Gamma ^h_{kr}\Gamma ^i_{hs}\;.\end{align*}

Similarly, the infinitesimal small vector 
$\overrightarrow{P'P^{'''}}
-\overrightarrow{PP^{''}}$ has for components $\alpha \beta D_rD_sx^i$~.

An elementary calculus gives that $D_rD_sx^i- D_sD_rx^i=0$~, which leads to the searched integrability conditions of $d\vec e_i=\Gamma ^k_{ij}du^j\ \vec e_k$~:
\[\F{\partial \Gamma ^{kr}_i}{\partial u^s}-\F{\partial \Gamma ^{ks}_i}{\partial u^r}+
(\Gamma ^h_{ir}\Gamma ^k_{hs}-\Gamma ^h_{is}\Gamma ^k_{hr})=0\]
corresponding to the conditions to which the Euclidean metric $g_{ij}$ must satisfy.
\vskip 11pt

\item Consider the integrability conditions of $d\vec e_{i'}$ for the components 
$i'$~, $1\le i'\le m$~ , of a stratum $\Sigma $ of corank $m$~.

Each function $\phi ^{(n)}(u_1,\cdots,u_i,\cdots,u_{n-m},\cdots,v_1,\cdots,v_{i'},\cdots,v_m)$~, $1\le i\le n-m$~, $1\le i'\le m$~, having $M$ as critical point(s) \cite{Tho1} on $\Sigma $ must satisfy:
\[ \L.\F{\partial \phi ^{(n)}}{\partial v_1}\R|_M = \cdots = 
\L.\F{\partial \phi ^{(n)}}{\partial v_m}\R|_M =0\;, \]
i.e.
\begin{align*}
 &\lim_{\Delta v_{i'}\to0}
\big( \phi ^{(n)}(u_1,\cdots,u_i,\cdots,u_{n-m},v_1,\cdots,v_{i'}+\Delta v_{i'},\cdots,v_m)
\\
& \qquad \L.-
\phi ^{(n)}(u_1,\cdots,u_i,\cdots,u_{n-m},v_1,\cdots,v_{i'},\cdots,v_m)\big)\Big/{\Delta v_{i'}}
\R|_M=0\;,\quad 1\le i'\le m \;, \end{align*}
which implies that:
\begin{align*}
& \lim_{\Delta v_{i'}\to0}
 \phi ^{(n)}(u_1,\cdots,u_{n-m},v_1,\cdots,v_{i'}+\Delta v_{i'},\cdots,v_m)\\
&\qquad \qquad \qquad
=
\phi ^{(n)}(u_1,\cdots,u_{n-m},v_1,\cdots,v_{i'},\cdots,v_m)\end{align*}
in the neighbourhood $D_\Sigma $ of $\Sigma $~.

This means that, in the $1\le i'\le m$ dimensions of the singular locus $\Sigma $~, the differentials of the basic vectors $\vec e_{i'}$ must by given by:
\[ d\vec e_{i'} = \Gamma ^{k'}_{i'j'}dv^{j'}\vec e_{k'}-\kappa g_{i'k}dv^{k'}
\centerdot M(u^1,\cdots,u^{n-m},v^1,\cdots,v^m)\;,\quad
\text{with\ }\kappa \in\rit^+\;,\]
because $\|d\vec e_{i'}\|<\|d\vec e_i\|$ in the neighbourhood $D_\Sigma $ of the singular locus $\Sigma $~.

Similarly, the components of the differential of the point $P(x^{1'},\cdots,x^{n'})$ in the $i'$~, $1\le i'\le m$~, dimensions of $\Sigma $ will then be:
\[Dx^{i'}=dx^{i'}+dv^{i'}+x^{k'}\Gamma ^{i'}_{k'r'}dv^{r'}-\kappa g_{i'k'}dv^{k'}\]
such that
\[ D_{r'}x^{i'} = \F{\partial x^{i'}}{\partial v^{r'}} +\delta ^{i'}_{r'}+x^{k'} \Gamma ^{i'}_{k'r'}-\kappa g_{i'k'}\delta ^{k'}_{r'}\]
corresponds to the $r'$-th component of the infinitesimally small vector 
$\overrightarrow {PP'}/\alpha '$~.

Calculating $D_{r'}D_{s'}x^{i'}-D_{s'}D_{r'}x^{i'}$~, we find the integrability conditions of $d\vec e_{i'}$~:
\[\F{\partial \Gamma ^{k'}_{i'r'}}{\partial v^{s'}}- \F{\partial \Gamma ^{k'}_{i's'}}{\partial v^{r'}}
+ (\Gamma ^{h'}_{i'r'}\Gamma ^{k'}_{h's'} - \Gamma ^{h'}_{i's'}\Gamma ^{k'}_{h'r'})
= -\kappa \ (\delta ^{k'}_{s'}g_{i'r'}-\delta ^{k'}_{r'}g_{i's'})\]
which clearly do not correspond to the conditions given in 1) to which the coefficients $g_{ij}$ must satisfy in order that $ds^2$ be Euclidean.

We thus have a non-Euclidean hyperbolic metric of curvature ``~$-\kappa $~'' on each stratum $D_\Sigma $ in the neighbourhood of the singular locus $\Sigma $~.
\vskip 11pt

\item These developments correspond to those of Hironaka \cite{Hir3} who showed that there is a normal cone along the singular locus $\Sigma $~.\epr
\Ee
\vskip 11pt

\subsubsection{The limit set of the Kleinian group}

A Kleinian group $G$ of $\o{\rit}^n=\rit^n\cup \{\infty \}$ is the group of M\"obius transformations of $\o{\rit}^n$ if it acts discontinuously somewhere in $\rit^n$~.

The action of the Kleinian group $G$ of $\o{\rit}^n$ can be extended to $\o H^{n+1}=H^{n+1}\cup \o{\rit}^n$ where $H^{n+1}=\{(x_1,\cdots,x_{n+1})\in \rit^{n+1}:x_{n+1}>0\}$ is the $(n+1)$-dimensional hyperbolic space: $G$ thus acts as a group of isometrics of $H^{n+1}$ with the hyperbolic metric.

The orbit space $M_G$ associated with the Kleinian group $G$ is defined by:
\[ M_G=\o H^{n+1}\smallsetminus L(G))\big/G\]
where $L(G)$ denotes the limit set of a discrete Kleinian group $G$ \cite{Mil1}, \cite{Tuk1}, \cite{Tuk2}.

This limit set is the closure of the set of fixed points of non-elliptic elements of $G$ \cite{Abi}.  It is a nowhere dense set whose area measure is zero: this corresponds to the zero-measure problem of Ahlfors \cite{Ahl}.

A discrete Kleinian group $G$ is said to be elementary if $L(G)$ consists of at most two points.

An ordinary set $\Omega (G)$ of a Kleinian group $G$ is defined by $\Omega (G)=\o{\rit}^n\smallsetminus L(G)$~: it is the region of discontinuity of $G$~.

Recall that a M\"obius transformation $g$ of $\o{\rit}^n$ is loxodromic if it is a transformation of the form $g(x)=\lambda \alpha (x)$ where $x\in\rit^n$~, $\lambda >1$ and $\alpha \in O(n)$ is the orthogonal group of $\rit^n$~.  $g$ is hyperbolic if $\alpha =$ id., elliptic if $\lambda =1$ and parabolic if $g$ has the form $g(x)=\alpha (x)+a$ where $a\in\rit^n\smallsetminus \{0\}$ and $\alpha (a)=a$~.
\vskip 11pt

\subsubsection{Left and right actions of the Kleinian group}

Similarly, we can introduce \lr M\"obius transformations $g_L$ (resp. $g_R$~) acting discontinuously in $(F^+_L)^n$ (resp. $(F^+_R)^n$~) in the upper (resp. lower) half space and \lr actions of the Kleinian group on the upper (resp. lower) $n$-dimensional hyperbolic half space $H^n_L$ (resp. $H^n_R$~).

The \lr orbit space $M_{G_L}$ (resp. $M_{G_R}$~) associated with the \lr action of the Kleinian group is given by:
\[ M_{G_L}=H^n_L\smallsetminus L(G_L)\big/ G_L \qquad \text{(resp.} \quad
M_{G_R}=H^n_R\smallsetminus L(G_R)\big/ G_R \ )\]
where $L(G_L)$ (resp. $L(G_R)$~) denotes the limit set of the Kleinian group $G_L$ (resp. $G_R$~) acting on the upper (resp. lower) half space.

And, a \lr ordinary set $\Omega (G_L)$ (resp. $\Omega (G_R)$~) of $G_L$ (resp. $G_R$~) is defined by:
\[ \Omega _{G_L}= (\o F^+_L)^n\smallsetminus L(G_L) \qquad \text{(resp.} \quad
\Omega _{G_R}= (\o F^+_R)^n\smallsetminus L(G_R) \ )\]
where
\[ (\o F^+_L)^n = (F^+_L)^n\cup \{\infty \} \qquad \text{(resp.} \quad
(\o F^+_R)^n = (F^+_R)^n\cup \{\infty \} \ ).\]
\vskip 11pt

\subsubsection{Proposition}

{\em Let $\Sigma _{\phi ^{(n)}_{G_{j_L}}}$ (resp. $\Sigma _{\phi ^{(n)}_{G_{j_R}}}$~) denote the singular locus of a germ $\phi _{j_\delta }(\omega _L)$ (resp. $\phi _{j_\delta }(\omega _R)$~) of corank $m$~, $1\le m\le 3$~, on the differentiable function
$\phi ^{(n)}_{g_{j_L}}(x_{g_{j_\delta }})$ (resp. $\phi ^{(n)}_{g_{j_R}}(x_{g_{j_\delta }})$~) and let $D_{\Sigma _{\phi _{G_{j_L}}}}$ (resp. $D_{\Sigma _{\phi _{G_{j_R}}}}$~) be the neighbourhood of this singular locus.

Then, it can be asserted that:
\Bena
\item the limit set $L(G_L)$ (resp. $L(G_R)$~) of the Kleinian group $G_L$ (resp. $G_R$~) corresponds to the singular locus $\Sigma _{\phi ^{(n)}_{G_{j_L}}}$ (resp. $\Sigma _{\phi ^{(n)}_{G_{j_R}}}$~).

\item the ordinary set $\Omega (G_L)$ (resp. $\Omega (G_R)$~) of $G_L$ (resp. $G_R$~), characterized by a hyperbolic metric, corresponds to the neighbourhood $D_{\Sigma _{\phi _{g_{j_L}}}}$ (resp. $D_{\Sigma _{\phi _{g_{j_R}}}}$~) of the singular locus.
\Ee
}
\vskip 11pt

\bpr \Bena
\item The limit set $L(G_L)$ (resp. $L(G_R)$~) has a measure equal to zero and, thus, corresponds, by the Sard lemma, to the singular locus $\Sigma _{\phi ^{(n)}_{G_{j_L}}}$ (resp. $\Sigma _{\phi ^{(n)}_{G_{j_R}}}$~).

Furthermore, $L(G_L)$ (resp. $L(G_R)$~) is a nowhere dense set: so, we have that:
\[ L(G_L)=\Sigma _{\phi ^{(n)}_{G_{j_L}}}\qquad
\text{(resp.} \quad L(G_R)=\Sigma _{\phi ^{(n)}_{G_{j_R}}}\ ).\]

\item It results from section 2.3.3 that the ordinary set $\Omega (G_L)$ (resp. $\Omega (G_R)$~) of $ G_L$ (resp. $G_R$~) is characterized by a hyperbolic metric and, thus, that
$\Omega (G_L)$ (resp. $\Omega (G_R$~) corresponds to the neighbourhood 
$D_{\Sigma _{\phi _{g_{j_L}}}}$ (resp. $D_{\Sigma _{\phi _{g_{j_R}}}}$~) of the singular locus according to proposition 2.3.1.\epr
\Ee
\vskip 11pt

\subsubsection{Corollary}

{\em The hyperbolic geometry, characterizing the neighbourhood 
$D_{\Sigma _{\phi _{g_{j_L}}}}$ (resp. $D_{\Sigma _{\phi _{g_{j_R}}}}$~) of the singular locus, in such a way that:
\[ \Omega (G_L)=
D_{\Sigma _{\phi _{g_{j_L}}}} \qquad \text{(resp.} \quad 
\Omega (G_R)=
D_{\Sigma _{\phi _{g_{j_R}}}} \ ),\]
results from the sequence of contracting surjective morphisms of singularization as developed, for example, in proposition 2.1.10.
}
\vskip 11pt

\subsubsection[The unfolded stratum $f_{j_{\delta _L}}$ (resp. $f_{j_{\delta _R}}$~)]{\boldmath The unfolded stratum $f_{j_{\delta _L}}$ (resp. $f_{j_{\delta _R}}$~)}

Let $\phi _{j_\delta }(\omega _L)$ (resp. $\phi _{j_\delta }(\omega _R)$~) denote a singular germ of corank ``~$m$~'' and codimension ``~$s$~'' on a $n$-dimensional differentiable function $\phi ^{(n)}_{G_{j_L}}(x_{g_{j_\delta }})$ (resp.
$\phi ^{(n)}_{G_{j_R}}(x_{g_{j_\delta }})$~).

Let $f_{j_{\delta _L}}$ (resp. $f_{j_{\delta _R}}$~) denote the versal unfolding of the singular germ $\phi _{j_\delta }(\omega _L)$ (resp. $\phi _{j_\delta }(\omega _R)$~).

The dimension of $f_{j_{\delta _L}}$ (resp. $f_{j_{\delta _R}}$~) is in general equal to $d_{f_{j_\delta }}=s$~, where $m\le s\le n$~.

The unfolded function $f_{j_{\delta _L}}$ (resp. $f_{j_{\delta _R}}$~) is embedded in the function $\phi ^{(n)}_{G_{j_L}}(x_{g_{j_\delta }})$ (resp.
$\phi ^{(n)}_{G_{j_R}}(x_{g_{j_\delta }})$~) such that the complementary 
$f^\perp_{j_{\delta _L}}$ (resp. $f^\perp_{j_{\delta _R}}$~) of $f_{j_{\delta _L}}$ (resp. $f_{j_{\delta _R}}$~) on $\phi ^{(n)}_{G_{j_L}}(x_{g_{j_\delta }})$ (resp.
$\phi ^{(n)}_{G_{j_R}}(x_{g_{j_\delta }})$~) has dimension
$d_{f^\perp_{j_\delta }}=n-s$ where $s\le n$~.

The neighbourhood of $f_{j_{\delta _L}}$ (resp. $f_{j_{\delta _R}}$~) in $\phi ^{(n)}_{G_{j_L}}(x_{g_{j_\delta }})$ (resp.
$\phi ^{(n)}_{G_{j_R}}(x_{g_{j_\delta }})$~) affected by the versal deformation is denoted
$D_{f/\phi ^{(n)}_{G_{j_L}}}$ (resp. $D_{f/\phi ^{(n)}_{G_{j_L}}}$~).

Finally, let $\Sigma _{\phi _{g_{j_{L}}}}$ (resp. $\Sigma _{\phi _{g_{j_R}}}$~) denote the possible singular locus on the quotient algebra of the unfolded function 
$f_{j_{\delta _L}}$ (resp. $f_{j_{\delta _R}}$~) and let 
$D_{\Sigma _{\phi _{g_{j_L}}}}$ (resp. $D_{\Sigma _{\phi _{g_{j_R}}}} $~) be the neighbourhood of this singular locus.
\vskip 11pt

\subsubsection{Proposition}

{\em
The  neighbourhood $D_{f/\phi ^{(n)}_{G_{j_L}}}$ (resp. $D_{f/\phi ^{(n)}_{G_{j_L}}}$~) of 
the unfolded function $f_{j_{\delta _L}}$ \resp{$f_{j_{\delta _R}}$}
in $\phi ^{(n)}_{G_{j_L}}(x_{g_{j_\delta }})$ (resp.
$\phi ^{(n)}_{G_{j_R}}(x_{g_{j_\delta }})$~) is characterized by a spherical geometry except in the neighbourhood 
$D_{\Sigma _{f_{j_{\delta _L}}}}$ (resp.
$D_{\Sigma _{f_{j_{\delta _L}}}}$~) of the singular locus
$\Sigma _{\phi _{g_{j _L}}}$ (resp.
$\Sigma _{\phi _{g_{j _R}}}$~) where the geometry is of hyperbolic type.
}
\vskip 11pt

\bpr
\Bena
\item The left and right cases will not be distinguished as in proposition 2.3.1.

Let $M$ be a point of coordinates $(u^1,\cdots,u^n)$ on the differentiable function $\phi ^{(n)}$ in such a way that $M$ is not localized on the singular locus $\Sigma $ and its neighbourhood $D_\Sigma $ (see proposition 2.3.1 for the notations).

$(\vec e_1,\cdots,\vec e_n)$ will denote the basis vectors of the proper referential centred on $M$~.

As the stratum $\phi ^{(n)}\smallsetminus \Sigma \smallsetminus D_\Sigma $ of $\phi ^{(n)}$~, being not affected by the singular locus, is Euclidean, the differentials of $M$ and $\vec e_i$ are given by:
\[ dM=du^i\ \ e_i\;, \qquad d\vec e_i = \Gamma ^k_{ij}\ du^j\ \vec e_k\;.\]
The partial derivatives of the coordinates $x^i$ of a point $P$ in the neighbourhood of $M$ will be:
\[ D_rx^i = \F{\partial x^i}{\partial u^r}+\delta ^i_r+x^k\ \Gamma ^i_{kr}\;.\]

\item 	In consequence of the versal deformation of the singular locus $\Sigma $~, the basis vectors $(\vec e_1,\cdots,\vec e_n)$ will be increased by a small amount:

\[ \delta \ \vec \varepsilon _j = \kappa \ g_{kj}\ du^k\centerdot M(u^1,\cdots,u^n)\;, \qquad \text{with} \quad \kappa \in\rit^+\;.\]
In the dimensions of
\[ (f_{j_\delta }\cup D_{f/\phi ^{(n)}})\smallsetminus (\Sigma \cup D_\Sigma )\]
(as assumed in this proposition),

we then have that:
\[d\vec e_j=\Gamma ^k_{ij}\ du^i\ \vec e_k+\kappa \ g_{jk}\ du^k\centerdot M(u^1,\cdots,u^n)\]
leading to:
\[D_r\ x^j = \F{\partial x^j}{\partial u^r}+\delta ^j_r+x^k\ \Gamma ^j_{kr}+\kappa \ g_{jk}\ \delta ^k_r\;.\]
Proceeding as in proposition 2.3.1 to calculate the integrability conditions of $d\vec e_j$~, we find that the coefficients $g_{jk}$ correspond to a spherical metric of curvature $+\kappa >0$~.

The spherical geometry on 
$(f_{j_\delta} \cup D_{f/\phi ^{(n)}})\smallsetminus (\Sigma \cup D_\Sigma )$ 
is a consequence of the versal deformation of $\Sigma $ leading to an over compactness of these strata as resulting from proposition 2.2.8.
\vskip 11pt

\item In the neighbourhood $D_{\Sigma _{\phi _{g_{j_L}}}}$ of the singular locus $\Sigma _{\phi _{g_{j_L}}}$~, the geometry is hyperbolic, as developed in proposition 2.3.1.\epr
\Ee

\section{Spreading-out isomorphism and strange attractors}

\subsection{The spreading-out isomorphism}

The aim of this chapter is to introduce the blow-up of the versal deformation: it will be called spreading-out (isomorphism) and it is the analogue of the desingularization, also called a monoidal transformation.  So, in the prolongation of the singularization and of the versal deformation, the spreading-out isomorphism corresponds to the blow-up of a contracting morphism.
\vskip 11pt

\subsubsection{Characteristics of the versal deformation}

We refer to propositions 2.2.5 and 2.2.6 where the versal deformation of the semisheaf $\theta ^*_{G_L^{(n)}}$ (resp. $\theta ^*_{G_R^{(n)}}$~) of differentiable functions $\phi _{j_\delta }(x_{g_{j_{\delta _L}}})$ \resp{$\phi _{j_\delta }(x_{g_{j_{\delta _R}}})$} endowed with singular germs $\phi _{j_\delta }(\omega _L)$ \resp{$\phi _{j_\delta }(\omega _R)$} of corank 1 is given by the contracting fibre bundle:
\begin{align*}
 D_{S_L} : \quad \theta ^*_{G^{(n)}_L}\times \theta_{S_L} &\To \theta ^*_{G^{(n)}_L} \\[6pt]
\rresp{
D_{S_R} : \quad \theta ^*_{G^{(n)}_R}\times \theta_{S_R} &\To \theta ^*_{G^{(n)}_R}},\end{align*}
in such a way that the fibre 
\begin{align*}
\theta _{S_L}&=\{\theta ^1(\omega ^1_L),\cdots,\theta ^1(\omega ^i_L),\cdots,\theta ^1(\omega ^s_L)\}\\[6pt]
\rresp{\theta _{S_R}&=\{\theta ^1(\omega ^1_R),\cdots,\theta ^1(\omega ^i_R),\cdots,\theta ^1(\omega ^s_R)\}},
\end{align*}
 given by the family of sheaves of the base $S_L$ \resp{$S_R$} of the versal deformation, is projected onto the $(n-1)$-dimensional coefficient sheaf:
\begin{align*}
\theta _L(a) &= \{\theta ^{n-1}_L(a_1),\cdots,\theta ^{n-1}_L(a_i),\cdots,\theta ^{n-1}_L(a_s)\} \\[6pt]
\rresp{
\theta _R(a) &= \{\theta ^{n-1}_R(a_1),\cdots,\theta ^{n-1}_R(a_i),\cdots,\theta ^{n-1}_R(a_s)\}}
\end{align*}
whose sections $a_{ij_\delta }(x_L)\in\theta ^{n-1}_L(a_i)$
\resp{$a_{ij_\delta }(x_R)\in\theta ^{n-1}_R(a_i)$}, $1\le j_\delta \le r$~, 
are ideals of functions on $\phi _{j_\delta}(x_{g_{j_{\delta _L}}})$ \resp{$\phi _{j_\delta}(x_{g_{j_{\delta _L}}})$}.
\vskip 11pt

\subsubsection{Lemma}

{\em The semisheaves $\theta ^{n-1}_L(a_i)$ \resp{$\theta ^{n-1}_R(a_i)$}, $1\le i\le s$~, and $\theta ^{1}_L(\omega _L^i)$ \resp{$\theta ^{1}_R(\omega _R^i)$} are characterized by the same set of ranks.
}

\bpr
\Bena
\item The section $a_{ij_\delta }(x_L)\subset \phi _{j_\delta }(x_{g_{j_{\delta _L}}})$ being a subfunction of $\phi _{j_\delta }(x_{g_{j_{\delta _L}}})$ must be characterized by a rank
\[ n_{a_{ij_\delta }}=(h_{j_\delta }\cdot N)^{n-1}\]
where:
\Bi
\item the integer $h_{j_\delta }$ is a global residue degree verifying $h_{j_\delta }>j_\delta $~, with $j_\delta $ being the global residue degree of the conjugacy class $g^{(n)}_L[j_\delta ]$ (see section 1.6) on which $\phi _{j_\delta }(x_{g_{j_{\delta _L}}})$ is defined.

Note that the rank $r_{g^{(n)}_L}$ of $g^{(n)}_L[j_\delta ]$ is $r_{g^{(n)}_L}=(j_\delta \centerdot N)^n$~, $j_\delta \equiv j\in\nit$ \cite{Pie1}.

\item $N$ is the rank of a real irreducible central completion (see section 1.1).
\Ei

\item As $\theta _{S_L}$ \resp{$\theta _{S_R}$} is projected onto $\theta _L(a)$ \resp{$\theta _R(a)$} in such a way that the semisheaf $\theta ^1(\omega ^i_L)$ \resp{$\theta ^1(\omega ^i_R)$} be flat onto the semisheaf $\theta ^{n-1}_L(a_i)$ \resp{$\theta ^{n-1}_R(a_i)$}, the normal crossing divisor $\omega ^i_{j_{\delta _L}}\in \theta ^{1}_L(\omega ^i_L)$ \resp{$\omega ^i_{j_{\delta _R}}\in \theta ^{1}_R(\omega ^i_R)$} must have a rank $n_{\omega ^i_{j_\delta }}$ proportional or equal to the rank $n_{a_{ij_\delta }}$ of $a_{ij_\delta }(x_L)$~.  If $(n-1)\le 2$~, then we have that $n_{\omega ^i_{j_\delta }}=(h_{j_\delta } \centerdot N)^p$~, where $p\ge n-1$~.

Remark that we extend here the concept of rank of a (semi)module to the (semi)sheaf defined on this (semi)module.

\item Finally, let $n_{\omega ^i}=\{n_{\omega ^i_1},\cdots,n_{\omega ^i_{j_\delta }},\cdots,n_{\omega ^i_{r}}\}$ be the set of ranks of the base semisheaf $\theta ^1(\omega ^i_L)$ \resp{$\theta ^1(\omega ^i_R)$} and let $n_{a_i}=\{n_{a_{i_1}},\cdots,n_{a_{ij_\delta }},\cdots,n_{a_{i_r}}\}$ be the corresponding set of ranks of the semisheaf $\theta ^{n-1}_L(a_i)$ \resp{$\theta ^{n-1}_R(a_i)$} such that, if $n-1\le 2$~, $n_{\omega ^i_{j_\delta }}\ge n_{a_{ij_\delta }}$ in the sense of 1) and 2).\epr
\Ee
\vskip 11pt

\subsubsection{Galois antiautomorphic (semi)group}

\Be

\item Let $\Gal(\widetilde F_L^+\big/F^0)=\Aut_{F^0}\ \wt F^+_L$
\resp{$\Gal(\widetilde F_R^+\big/F^0)=\Aut_{F^0}\ \wt  F^+_R$} be the Galois automorphic group acting transitively on the set of ideals
\[
\wt  F^+_{v_1}  \subset \cdots   \subset   \wt  F^+_{v_{j_\delta }} 
 \subset   \cdots   \subset   \wt  F^+_{v_r}\qquad
\rresp{\wt  F^+_{\o v_1}  \subset \cdots   \subset   \wt  F^+_{\o v_{j_\delta }} 
 \subset   \cdots   \subset   \wt  F^+_{\o v_r}}
\]
forming an increasing sequence characterized by the extension degrees:
\be [\wt  F^+_{v_{j_\delta }}:F^0] \equiv [\wt  F^+_{\o v_{j_\delta }}:F^0]=*+j \centerdot N
\tag{see section 1.1}\ee
and, more particularly, by their global residue degrees forming the increasing sequence:
\[
f_{v_{1_\delta }}  \subset \cdots   \subset   f_{v_{j_\delta }} 
 \subset   \cdots  \subset  f_{v_{r_\delta }}\qquad
\rresp{f_{\o v_{1_\delta }}  \subset \cdots   \subset   f_{\o v_{j_\delta }} 
 \subset  \cdots  \subset   f_{\o v_{r_\delta }}}
\]
where: $f_{v_{j_\delta }}\equiv f_{\o v_{j_\delta }}=j$~, $j  \in\NN$~.

\item Inversely, we introduce the Galois antiautomorphic group 
$\Gal^{-1}(\widetilde F_L^+\big/F^0)=\Aut^{-1}_{F^0}\ \wt  F^+_L$
\resp{$\Gal^{-1}(\widetilde F_R^+\big/F^0)=\Aut^{-1}_{F^0}\ \wt  F^+_R$}
acting transitively on the set of ideals:
\[
\wt  F^+_{v_{s_\delta }}  \supset \cdots   \supset \wt  F^+_{v_{j_\delta }} 
 \supset  \cdots   \supset   \wt  F^+_{v_{1_\delta }}\qquad
\rresp{\wt  F^+_{\o v_{s_\delta }}  \supset \cdots   \supset   \wt  F^+_{\o v_{j_\delta }} 
 \supset   \cdots   \supset  \wt F^+_{\o v_{1_\delta }}}
\]
forming a  decreasing sequence characterized by a decreasing chain of global residue degrees:
\[
f_{v_{s_\delta }}  \supset \cdots   \supset   f_{v_{j_\delta }} 
 \supset   \cdots   \supset   f_{v_{1_\delta }}\qquad
\rresp{f_{\o v_{s_\delta }}  \supset \cdots   \supset   f_{\o v_{j_\delta }} 
 \supset   \cdots   \supset   f_{\o v_{1_\delta }}}
\]
with the condition that $s\le r$~.
\Ee

\subsubsection[$n$-dimensional representations of products, right by left, of Galois groups]{\bbf $n$-dimensional representations of products, right by left, of Galois groups}

\Be
\item Referring to \cite{Pie1}, we introduce the explicit $n$-dimensional representation:
\[ \Rep^{(n)}_{\Gal_{\wt F\RL}}: \quad \Gal(\wt F^+_L\big / F^0)\times \Gal (\wt F^+_L\big / F^0) \To G^{(n)}({F^+_{\o v}}\times {F^+_v})\]
of the product, right by left, of the Galois automorphic groups in such a way that the conjugacy class representatives of the bilinear algebraic semigroup $G^{(n)}({F^+_{\o v}}\times {F^+_v})$ form the increasing sequence:
\[ g^{(n)}\RL[1] \subset \cdots \subset g^{(n)}\RL[j_\delta ] \subset \cdots \subset g^{(n)}\RL[r]\]
characterized by the increasing sequence of their ranks
\be
(1\centerdot N)^{2n}\subset \cdots \subset (j \centerdot N)^{2n}
\subset \cdots \subset (r\centerdot N) ^{2n}\;, \quad j <r
\tag*{(see section 1.6).}\ee

\item Similarly, we can introduce the (inverse) $n$-dimensional representation:
\[ \Rep^{(n)}_{\Gal^{-1}_{\wt F\RL}}: \quad \Gal^{-1}(\wt F^+_L\big / F^0)\times \Gal^{-1} (\wt F^+_L\big / F^0) \To G^{-1(n)}({F^+_{\o v(s)}}\times {F^+_{v(s)}})\]
of the product, right by left, of the Galois antiautomorphic groups in such a way that the conjugacy class representatives of the inverse bilinear algebraic semigroup  form the decreasing sequence:
\[ g^{(n)}\RL[r] \supset \cdots \supset g^{(n)}\RL[j_\delta ] \supset \cdots \supset g^{(n)}\RL[s+1]\]
characterized by the decreasing chain of the  ranks
\[
(r\centerdot N)^{2n}\supset \cdots \supset (j \centerdot N)^{2n}
\supset \cdots \supset ((s+1)\centerdot N) ^{2n}\;, \]
with the condition that $s\le r$~.
\Ee
\vskip 11pt

\subsubsection{Proposition}

{\em Let $G^{(n)}({F^+_{\o v}}\times {F^+_{v}})$ be the bilinear algebraic semigroup 
over the product, right by left, of the sets of real completions $F^+_{\o v}$ and $F^+_v$~.

Then, the inverse bilinear algebraic semigroup 
$G^{-1(n)}({F^+_{\o v(s)}}\times {F^+_{v(s)}})$ with entries in the product, right by left, of 
the sets of real completions
over the first $s$ places, $s\le r$~, generates the following smooth endomorphism:
\[ E\big[ G^{(n)}({F^+_{\o v}}\times  {F^+_{v}})\big]
= G^{-1(n)}({F^+_{\o v(s)}}\times {F^+_{v(s)}})
\oplus G_I^{(n)}({F^+_{\o v {(r-s)}}}\times {F^+_{v{(r-s)}}})\]
where $G_I^{(n)}({F^+_{\o v {(r-s)}}}\times {F^+_{v {(t-s)}}})$ is a bilinear algebraic semigroup complementary of\linebreak
$G^{(-1)(n)}({F^+_{\o v(s)}}\times {F^+_{v(s)}})$ in the sense that it is defined with entries in the $(r-s)$ last places. 
}
\vskip 11pt

\bpr As $G^{-1(n)}({F^+_{\o v(s)}}\times {F^+_{v(s)}})$ is the 
$n$-dimensional representation of the product, right by left, of Galois antiautomorphic groups according to sections 3.1.3 and 3.1.4, its conjugacy classes form a decreasing sequence from the $r$-th 
biplace $\o v_r\times v_r$ 
until the $(s+1)$-th biplace 
$\o v_{s+1}\times v_{s+1}$ in such a way that 
$(r-s)$ conjugacy classes of
$G^{(n)}({F^+_{\o v}}\times {F^+_{v}})$ have been disconnected and generate the complementary bilinear algebraic semigroup
$G^{(n)}_I({F^+_{\o v {(r-s)}}}\times {F^+_{v {(r-s)}}})$~.

So, every smooth endomorphism $E[G^{(n)}({F^+_{\o v}}\times {F^+_{v}})]$ of $G^{(n)}({F^+_{\o v}}\times {F^+_{v}})$ is generated by the inverse bilinear semigroup
$G^{(-1)n}( {F^+_{\o v(s)}}\times {F^+_{v(s)}})$ in such a way that two non connected bilinear algebraic semigroups
$G^{(-1)(n)}({F^+_{\o v(s)}}\times  {F^+_{v(s)}})$ and
$G^{(n)}_I( {F^+_{\o v {(r-s)}}}\times  {F^+_{v {(r-s)}}})$ are produced from
$G^{(n)}( {F^+_{\o v}}\times  {F^+_{v}})$~.\epr
\vskip 11pt

\subsubsection{Corollary}

{\em Let $\theta _{G^{(n)}\RL}=\theta _{G^{(n)}_R}\otimes \theta _{G^{(n)}_R}$ be the bisemisheaf of rings on the bilinear algebraic semigroup 
$G^{(n)}( {F^+_{\o v}}\times {F^+_{v}})$ as introduced in section 1.7.

If $\theta _{G^{(n)}\RL}$ is noted
$\theta _{G^{(n)}( {F^+_{\o v}}\times  {F^+_{v}})}$~, every smooth endomorphism of it is given by:
\[ E[\theta _{G^{(n)}( {F^+_{\o v}}\times  {F^+_{v}})}]
= \theta _{G^{-1(n)}( {F^+_{\o v(s)}}\times  {F^+_{v(s)}})}
\oplus \theta_{G^{(n)}_I( {F^+_{\o v{(r-s)}}}\times  {F^+_{v {(r-s)}}})}\;.\]
}
\vskip 11pt

\bpr This is an adaptation of proposition 3.1.5 to the bisemisheaf
$\theta _{G^{(n)}( {F^+_{\o v}}\times {F^+_{v}})}$ on the bilinear algebraic semigroup $G^{(n)}( {F^+_{\o v}}\times  {F^+_{v}})$~.\epr
\vskip 11pt

\subsubsection[Decomposition of $G^{(n)}( {F^+_{\o v}}\times  {F^+_{v}})$ in irreducible completions]{\bbf Decomposition of $G^{(n)}( {F^+_{\o v}}\times  {F^+_{v}})$ in irreducible completions}

As the bilinear algebraic semigroup $G^{(n)}( {F^+_{\o v}}\times  {F^+_{v}})$ is constructed on products, right by left, of irreducible completions of rank $N$ and as its $j_\delta $-th, $m_{j_\delta }$-th conjugacy class representative
$g^{(n)}\RL[j_\delta ,m_{j_\delta }]$ counts $j_\delta $ products, right by left, of irreducible completions of rank $N$~, we have, in fact, in 
$G^{(n)}( {F^+_{\o v}}\times  {F^+_{v}})$~:
\[ r^{nr}_{\o v\times v}=\bigoplus^r_{j_\delta =1}(j_\delta \centerdot m^{(j_\delta )})^{n}\]
pairs of irreducible completions where $m^{(j_\delta )}=\sup m_{j_\delta }$ is the multiplicity of the $j_\delta $-th conjugacy class representative.

So, on $G^{(n)}( {F^+_{\o v}}\times  {F^+_{v}})$~, we have an increasing sequence:
\begin{align*}
G^{(n)}(F^+_{\o v_{1_\delta }}\times F^+_{v_{1_\delta }})_{\rm up}
\subset \cdots &\subset
G^{(n)}(F^+_{\o v_{j_\delta,m_{j_\delta } }}\times F^+_{v_{j_\delta,m_{j_\delta } }})_{\rm up}\\
&\subset \cdots \subset
G^{(n)}(F^+_{\o v^{r_\delta }_{r_\delta,m^{(r_\delta )} }}\times F^+_{v^{r_\delta }_{r_\delta,m^{(r_\delta )} }})_{\rm up}
\end{align*}
of 
sets of conjugacy class representatives 
where:
$G^{(n)}(F^+_{\o v_{j_\delta,m_{j_\delta } }}\times F^+_{v_{j_\delta,m_{j_\delta } }})_{\rm up}$ denotes a bilinear algebraic semigroup whose upper entry is the irreducible bicompletion
$(F^+_{\o v_{j_\delta,m_{j_\delta } }}\times F^+_{v_{j_\delta,m_{j_\delta } }})$~.
To simplify the notations, the increasing global residue degree of
$G^{(n)}({F^+_{\o v}}\times {F^+_v})$~, associated with its structure in products, right by left, of irreducible completions, will be noted $f$ running from 1 to $r^{nr}_{\o v\times v}$~.
\vskip 11pt

\subsubsection[The fibre $\theta _{S_L}$ \resp{$\theta _{S_R}$} of the versal deformation]{\bbf The fibre $\theta _{S_L}$ \resp{$\theta _{S_R}$} of the versal deformation}

These considerations can be transposed to the fibre
\begin{align*}
\theta _{S_L}&=\{\theta ^1(\omega ^1_L),\cdots,\theta ^1(\omega ^i_L),\cdots,\theta ^1(\omega ^s_L)\}\\
\rresp{\theta _{S_R}&=\{\theta ^1(\omega ^1_R),\cdots,\theta ^1(\omega ^i_R),\cdots,\theta ^1(\omega ^s_R)\}}\end{align*}
of the versal deformation in the following way:
\vskip 11pt
Let $\theta ^1(\omega ^i_L)$ \resp{$\theta ^1(\omega ^i_R)$} be the $i$-th sheaf of the base $S_L$ \resp{$S_R$} of the versal deformation.

This $i$-th (semi)sheaf $\theta ^1(\omega ^i_L)$ (and $\theta ^1(\omega ^i_R)$~) is characterized by the set of ranks $n_{\omega ^i} = \{n_{\omega ^i_{1_\delta }},\cdots,n_{\omega ^i_{j_\delta }},\cdots,n_{\omega ^i_{r_\delta }}\}$ according to lemma 3.1.2 where $n_{\omega ^i_{j_\delta }}$ refers to the rank of the normal crossings divisor $\omega ^i_{j_\delta }$~, which is the $i$-th generator of the versal  unfolding of a singular germ $\phi _{j_\delta }(\omega _L)$ \resp{$\phi _{j_\delta }(\omega _R)$} of corank 1 and codimension $s$ on the $j_\delta $-th differentiable function $\phi _{j_\delta }(x_{g_{j_{\delta _L}}})$
\resp{$\phi _{j_\delta }(x_{g_{j_{\delta _R}}})$} of the semisheaf
$\theta ^*_{G^{(n)}_L}$ \resp{$\theta ^*_{G^{(n)}_R}$} (see section 3.1.1 and proposition 2.2.8).

To this rank $n_{\omega ^i_{j_\delta }}=(h_{j_\delta }\centerdot N)^p$ (see lemma 3.1.2) corresponds the unramified rank or global residue degree $f_{\omega ^i_{j_\delta }}=(h_{j_\delta })^p=n_{\omega ^i_{j_\delta }}\big/ N^p$ which is the number of irreducible completions on the divisor $\omega ^i_{j_\delta }$~.

As in section 3.1.7, we shall label the set of irreducible completions in $\theta ^1(\omega ^i_L)$ (and on $\theta ^1(\omega ^i_R)$~) by a unique integer $f_i$ running over all the normal crossings divisor $\omega ^i_{j_\delta }$~, $1\le j_\delta \le r$~, in such a way that the maximal value of $f_i$ will be given by:
\[f^{\max}_i = \bigoplus^r_{j_\delta =1} \ \bigoplus_{m_{j_\delta }}
(h_{j_\delta ,m_{j_\delta }})^p\;, \qquad f_i^{\max}\in\NN\;.\]
So, we have $f^{\max}_i$ irreducible completions of rank $N$ on the $i$-th base semisheaf $\theta ^1(\omega ^i_L)$ \resp{$\theta ^1(\omega ^i_R)$} of the versal deformation.

Remark that  an integral irreducible closed subscheme of rank $N$ is defined on an irreducible completion of rank $N$~, the concept of rank being extended here from the topological (sub)space to the (sub)scheme on which it is defined.
\vskip 11pt

As a smooth endomorphism was introduced on the bilinear algebraic semigroup\linebreak
$G^{(n)}({F^+_{\o v}}\times {F^+_{v}})$ in proposition 3.1.5 and on the bisemisheaf $\theta _{G^{(n)}\RL}$ on it in corollary 3.1.6, a smooth endomorphism can be defined on the $i$-th base semisheaf $\theta ^1(\omega ^i_L)$ \resp{$\theta ^1(\omega ^i_R)$} as follows:
\vskip 11pt

\subsubsection{Proposition}

{\em Let $\theta ^1(\omega ^i_L)$ \resp{$\theta ^1(\omega ^i_R)$} be the $i$-th base semisheaf of the versal deformation of the semisheaf
$\theta ^*_{G_L^{(n)}}$ \resp{$\theta ^*_{G_R^{(n)}}$}.

Let $f_i^{\max}$ be the maximal value of its global residue degree counting the irreducible closed subschemes of rank $N$~.

Then, the following smooth endomorphism
\[ E_{\omega ^i_L}[\theta ^1(\omega ^i_L)_{f_i^{\max}}]
= \theta ^{*1}(\omega ^i_L)_{f^*_i}\oplus \theta ^1_I(\omega ^i_L)_{f^I_i}\;, \qquad \text{with} \quad f^I_i=f_i^{\max}-f^*_i\in\NN\;, \]
can be introduced on the semisheaf $\theta ^1(\omega ^i_L)_{f_i^{\max}}$ in such a way that it decomposes into two non connected complementary semisheaves
$\theta ^{*1}(\omega ^i_L)_{f_i^{*}}$ and 
$\theta ^{1}(\omega ^i_L)_{f_i^{I}}$ whose global residue degrees verify:
\[ f^{\max}_i=f^*_i+f^I_i\;.\]
}
\vskip 11pt

\bpr The semisheaf $\theta ^{*1}(\omega ^i_L)_{f_i^{*}}$ is a ``reduced'' semisheaf generated from the semisheaf
$\theta ^{1}(\omega ^i_L)_{f_i^{\max}}$ under the action of the Galois antiautomorphic group according to the endomorphism:
\[ E_{\omega ^i_L} : \quad
\theta ^{1}(\omega ^i_L)_{f_i^{\max}}\To
\theta ^{*1}(\omega ^i_L)_{f_i^{*}}\oplus 
\theta ^{1}_I(\omega ^i_L)_{f_i^{I}}\]
where:
\Bi
\item $\theta ^{i}_I(\omega ^i_L)_{f_i^{I}}$ is the semisheaf complementary of $\theta ^{*1}(\omega ^i_L)_{f_i^{*}}$ in the sense of proposition 3.1.5 and corollary 3.1.6.

\item $\theta ^{*1}(\omega ^i_L)_{f_i^{*}}$ is characterized by decreasing global residue degrees $f^*_i$ while
$\theta ^{1}_I(\omega ^i_L)_{f_i^{I}}$ is characterized by increasing global residue degrees $f^I_i$ in such a way that
\be f_i^{\max} = f^*_i+f^I_i\;, \qquad 0\le f^*_i\le f^{\max}_i\;, \quad 0\le f^I_i\le f_i^{\max}\;.\tag*{\eop}\ee
\Ei
\vskip 11pt

\subsubsection{Proposition}
 
{\em Every base semisheaf $\theta ^{1}(\omega ^i_L)$
\resp{$\theta ^{1}(\omega ^i_R)$} of the versal deformation of the semisheaf $\theta ^*_{G^{(n)}_L}$ \resp{$\theta ^*_{G^{(n)}_R}$}, $1\le i\le s$~, can generate under the smooth endomorphism $E_{\omega ^i_L}$ \resp{$E_{\omega ^i_R}$} the elements of the category $c(\theta ^1_{\omega ^i_L})$
\resp{$c(\theta ^1_{\omega ^i_R})$} of the $(f_i-1)$ pairs of semisheaves of rings:
\begin{multline*}
c(\theta ^1_{\omega ^i_L}) = \{ (\theta ^{*1}(\omega ^i_L)_{f^{\max}_i-1} \oplus \theta ^{1}_I(\omega ^i_L)_{1}),\cdots,
(\theta ^{*1}(\omega ^i_L)_{f^{*}_i} \oplus \theta ^{1}_I(\omega ^i_L)_{f^I_i}),\\
\cdots,
(\theta ^{*1}(\omega ^i_L)_{1} \oplus \theta ^{1}_I(\omega ^i_L)_{f_i^{\max}-1})\}\;, \quad 1\le f^*_i\le f_i^{\max}\;,
\end{multline*}
whose objects are two non connected semisheaves characterized by complementary global residue degrees verifying:
\[ f^{\max}_i=f^*_i+f^I_i\;.\]}
\vskip 11pt

\bpr This is a generalization of proposition 3.1.9 where $(f^{\max}_i-1)$ endomorphisms $E_{\omega ^i_L}$ are considered.\epr
\vskip 11pt

\subsubsection{Corollary}

{\em Let $f^*_i$ denote the global residue degree of the reduced semisheaf $\theta ^{*1}(\omega ^i_L)_{f^{*}_i} $~, $0\le f^*_i\le f_i^{*\max}$~.

Then, the smooth endomorphism $E_{\omega ^i_L}$ is maximal when $f^*_i=0$~.}
\vskip 11pt

\bpr If $f^*_i=0$~, then the semisheaf $\theta ^{1}(\omega ^i_L)_{f^{\max}_i} $ has been completely transformed under
$E_{\omega ^i_L}$ into the complementary base semisheaf
$\theta ^{1}_I(\omega ^i_L)_{f^{\max}_i} $~: this is equivalent to say that the base semisheaf $\theta ^{1}(\omega ^i_L)_{f^{\max}_i} $ has been totally disconnected from the semisheaf $\theta ^*_{G^{(n)}_L}$~.\epr
\vskip 11pt

\subsubsection{Proposition}

{\em Let $\theta ^{\rm vers}_{G^{(n)}_L}=\theta ^*_{G^{(n)}_L}\times \theta _{S_L}$ \resp{$\theta ^{\rm vers}_{G^{(n)}_R}=\theta ^*_{G^{(n)}_R}\times \theta _{S_R}$} be the versal deformation of the singular semisheaf $\theta ^*_{G^{(n)}_L}$ \resp{$\theta ^*_{G^{(n)}_R}$}, as introduced in proposition 2.2.6, where
\begin{align*}
\theta _{S_L} &=
\{\theta ^1(\omega ^1_L),\cdots,\theta ^1(\omega ^i_R),\cdots,\theta ^1(\omega ^s_L)\}\\[6pt]
\rresp{\theta _{S_R} &=
\{\theta ^1(\omega ^1_R),\cdots,\theta ^1(\omega ^i_R),\cdots,\theta ^1(\omega ^s_R)\}},\end{align*}
fibre of the contracting fibre bundle $D_{S_L}$ \resp{$D_{S_R}$}, is the family of the semisheaves of the base $S_L$ \resp{$S_R$} of the versal deformation.

Then, there exists a family of isomorphisms
\begin{align*}
&\Pi _{s_L} (f_1^I,\cdots,f_i^I,\cdots,f^I_s) :\\
& \qquad \theta ^*_{G^{(n)}_L}\times \theta _{S_L} 
 \To
\theta ^*_{G^{(n)}_L}\times \theta '_{S_L}
\cup \{ \theta ^1_I(\omega ^1_L)_{f^I_1},\cdots,
\theta ^1_I(\omega ^i_L)_{f^I_i},\cdots,\theta ^1_I(\omega ^s_L)_{f^I_s}\}\\[6pt]
\rresp{&\Pi _{s_R} (f_1^I,\cdots,f_i^I,\cdots,f^I_s) : \\
& \qquad \theta ^*_{G^{(n)}_R}\times \theta _{S_R} 
 \To
\theta ^*_{G^{(n)}_R}\times \theta '_{S_R}
\cup \{ \theta ^1_I(\omega ^1_R)_{f^I_1},\cdots,
\theta ^1_I(\omega ^i_R)_{f^I_i},\cdots,\theta ^1_I(\omega ^s_R)_{f^I_s}\}}
\end{align*}
disconnecting $f^I_1$ irreducible subsheaves of rank $N$ from the base semisheaf $\theta ^1(\omega ^1_L)_{f_i^{\max}}$ on $\theta ^*_{G^{(n)}_L}$~, \ldots, $f^I_i$ irreducible subsheaves of rank $N$ from the base semisheaf 
$\theta ^1(\omega ^i_L)_{f_i^{\max}}$ on $\theta ^*_{G^{(n)}_L}$~, \ldots, and so on, $1\le i\le s$~.

The set of complementary global residue degrees $(f^I_1,\cdots, f^I_i,\cdots,f^I_s)$ varies in such a way that $1\le f_1^I\le f_1^{\max},\cdots,1\le f_i^I\le f_i^{\max},\cdots,1\le f_s^I\le f_s^{\max}$ implying, for each set $(f^I_1,\cdots,f^I_i,\cdots,f^I_s)$ a family of isomorphisms $\Pi _s(f^I_1,\cdots,f^I_i,\cdots,f^I_s)$~.

The residue fibre $\theta '_{S_L}$ \resp{$\theta '_{S_R}$} is given by:
\begin{align*}
\theta '_{S_L}&= \theta _{S_L}\setminus \{
\theta ^1_I(\omega ^1_L)_{f^I_1},\cdots,
\theta ^1_I(\omega ^i_L)_{f^I_i},\cdots,\theta ^1_I(\omega ^s_L)_{f^I_s}\}\\[6pt]
\rresp{\theta '_{S_R}&= \theta _{S_R}\setminus \{
\theta ^1_I(\omega ^1_R)_{f^I_1},\cdots,
\theta ^1_I(\omega ^i_R)_{f^I_i},\cdots,\theta ^1_I(\omega ^s_R)_{f^I_s}\}}.\end{align*}
}

\bpr This proposition is a generalization of proposition 3.1.10 in such a way that the smooth endomorphism $E_{\omega ^i_L}$ \resp{$E_{\omega ^i_R}$}~, generating $(f_i^{\max}-1)$ pairs of semisheaves of the category $c(\theta ^1_{\omega ^i_L})$ \resp{$c(\theta ^1_{\omega ^i_R})$}, is extended to all the base semisheaves $\theta ^1(\omega ^i_L)$
\resp{$\theta ^1_{\omega ^i_R}$}, $1\le i\le s$~, of the considered versal deformation.\epr
\vskip 11pt

\subsubsection{Corollary}

{\em The family of isomorphisms 
$\Pi _{s_L} (f_1^I,\cdots,f_i^I,\cdots,f^I_s)$
\resp{$\Pi _{s_R} (f_1^I,\cdots,f_i^I,\cdots,f^I_s)$} is maximal in the $i$-th semisheaf $\theta ^1(\omega ^i_L)$
\resp{$\theta ^1(\omega ^i_R)$}, and, then, noted
$\Pi^{\max(i)} _{s_L} (f_1^I,\cdots,f_i^I,\linebreak \cdots, f^I_s)$
\resp{$\Pi^{\max (i)} _{s_R} (f_1^I,\cdots,f_i^I,\cdots,f^I_s)$}, if $f^*_i=0$~.
}
\vskip 11pt

\bpr If $f^*_i=0$~, then the base semisheaf $\theta ^1(\omega ^i_L)$ has been completely transformed, under $E_{\omega ^i_L}$~, into the disconnected complementary semisheaf $\theta ^1_I(\omega ^i_L)_{f^I_i}$
in such a way that $f^I_i=f_i^{\max}$~.\epr
\vskip 11pt

\subsubsection{Corollary}

{\em The family isomorphisms
$\Pi^{\max} _{s_L} (f_1^I,\cdots,f_i^I,\cdots,f^I_s)$
\resp{$\Pi^{\max} _{s_R} (f_1^I,\cdots,f_i^I,\cdots,f^I_s)$} is maximal if it is maximal in every semisheaf $\theta ^1(\omega ^i_L)$ \resp{$\theta ^1(\omega ^i_R)$ }, $1\le i\le s$~, of the base of the versal deformation.
}
\vskip 11pt

\bpr This is a generalization of corollary 3.1.3 to all the base semisheaves $\theta ^1(\omega ^i_L)$ \resp{$\theta ^1(\omega ^i_R)$} of the versal deformation, implying that:
\Bena
\item 
\textbullet \quad $f^*_1=0 \quad \and \quad f^I_1=f_1^{\max}$

\qquad $\vdots$

\textbullet \quad $f^*_i=0 \quad \and \quad f^I_i=f_i^{\max}$
\enlargethispage{\baselineskip}

\qquad $\vdots$

\textbullet \quad $f^*_s=0 \quad \and \quad f^I_s=f_s^{\max}$

\item \hspace{-1cm} $\begin{array}[t]{ll}
&\Pi^{\max}_{s_L} (f_1^I,\cdots,f_i^I,\cdots,f^I_s) :\\
& \qquad \theta ^*_{G^{(n)}_L}\times \theta _{S_L} 
 \to
\theta ^*_{G^{(n)}_L}
\cup \{ \theta ^1_I(\omega ^1_L)_{f^{\max}_1},\cdots,
\theta ^1_I(\omega ^i_L)_{f^{\max}_i},\cdots,\theta ^1_I(\omega ^s_L)_{f^{\max}_s}\}\\[6pt]
\rresp{&\Pi^{\max} _{s_R} (f_1^I,\cdots,f_i^I,\cdots,f^I_s) : \\
& \qquad \theta ^*_{G^{(n)}_R}\times \theta _{S_R} 
 \to
\theta ^*_{G^{(n)}_R}
\cup \{ \theta ^1_I(\omega ^1_R)_{f^{\max}_1},\cdots,
\theta ^1_I(\omega ^i_R)_{f^{\max}_i},\cdots,\theta ^1_I(\omega ^s_R)_{f^{\max}_s}\}}.\eop
\end{array}$
\Ee
\vskip 11pt

\subsubsection{Category of vertical tangent bundles}

Let $T_{V_{W_L}}=\{T_{V_{W^1_L}},\cdots,
T_{V_{W^i_L}},\cdots,
T_{V_{W^s_L}}\}$ \resp{$T_{V_{W_R}}=\{T_{V_{W^1_R}},\cdots,
T_{V_{W^i_R}},\cdots,
T_{V_{W^s_R}}\}$} denote the family of tangent vector bundles obtained by the projection of all the disconnected base semisheaves $\theta ^1_I(\omega ^i_L)$ \resp{$\theta ^1_I(\omega ^i_R)$}, $1\le i\le s$~, in the vertical tangent spaces
$T_{V_{W^i_L}}$ \resp{$T_{V_{W^i_R}}$} characterized by normal vector fields $\vec W_{i_L}$ \resp{$\vec W_{i_R}$}~.

The proper projective map of the tangent bundle
$T_{V_{W^i_L}}$ \resp{$T_{V_{W^i_R}}$} is given by:
\begin{align*}
\tau _{V_{W^i_L}}: \qquad 
T_{V_{W^i_L}}(\theta ^1_I(\omega ^i_L)_{f^I_i}) &\To \theta ^1_I(\omega ^i_L)_{f^I_i}\\[6pt]
\rresp{\tau _{V_{W^i_R}}: \qquad 
T_{V_{W^i_R}}(\theta ^1_I(\omega ^i_R)_{f^I_i}) &\To \theta ^1_I(\omega ^i_R)_{f^I_i}}
\end{align*}
so that
$\tau _{V_{W_L}}=\{\tau _{V_{W^i_L}}\}^s_{i=1}$
\resp{$\tau _{V_{W_R}}=\{\tau _{V_{W^i_R}}\}^s_{i=1}$}.

To the category $c(\theta ^1_I(\omega ^i_L))$ \resp{$c(\theta ^1_I(\omega ^i_R))$} of disconnected base semisheaves $\theta ^1_I(\omega ^i_L)$ \resp{$\theta ^1_I(\omega ^i_R)$}, $1\le i\le s$~, will then correspond  the category $c(T_{V_{W^i_L}}(\theta ^1_I(\omega ^i_L)$)
\resp{$c(T_{V_{W^i_R}}(\theta ^1_I(\omega ^i_R)))$} of sections of tangent vector bundles.
\vskip 11pt

\subsubsection{Proposition}

{\em The extension of the quotient algebra $\theta [{\omega _L}]$ \resp{$\theta ([{\omega _R}]$} of the versal deformation of the singular semisheaf $\theta ^*_{G^{(n)}_L}$ \resp{$\theta ^*_{G^{(n)}_R}$}, having an isolated degenerate singularity of corank 1 and codimension $s$ in each section of $\theta ^*_{G^{(n)}_L}$ \resp{$\theta ^*_{G^{(n)}_R}$}, is realized by the spreading-out isomorphism
\[SOT_L = (\tau _{V_{W_L}}\circ \Pi _{s_L})
\qquad \resp{\quad SOT_R = (\tau _{V_{W_R}}\circ \Pi _{s_R})}.\]
}
\vskip 11pt

\bpr Let $I_{\omega ^i_L}$ \resp{$I_{\omega ^i_R}$} be the kernel of the normal vector bundle $T_{V_{W^i_L}}$ \resp{$T_{V_{W^i_R}}$}.

Then, the exact sequence:
\begin{alignat*}{7}
0 &\To I_{\omega ^i_L} &\To & T_{V_{W^i_L}}(\theta ^1_I(\omega ^i_L)_{f^I_i}) &\ \xrightarrow{\tau _{V_{W^i_L}}}\ &
\theta ^1_I(\omega ^i_L)_{f^I_i} &\To & 0\\
\rresp{0 &\To I_{\omega ^i_R} &\To & T_{V_{W^i_R}}(\theta ^1_I(\omega ^i_R)_{f^I_i}) &\ \xrightarrow{\tau _{V_{W^i_R}}}\ &
\theta ^1_I(\omega ^i_R)_{f^I_i} &\To & 0}
\end{alignat*}
represents an extension of $\theta ^1_I(\omega ^i_L)_{f^I_i}$
\resp{$\theta ^1_I(\omega ^i_r)_{f^I_i}$} by the kernel $I_{\omega ^i_L}$ \resp{$I_{\omega ^i_R}$}.

And, the isomorphism $SOT_L=(\tau _{V_{W_L}}\circ \Pi _{s_L})$
\resp{$SOT_R=(\tau _{V_{W_R}}\circ \Pi _{s_R})$} constitutes an extension of the quotient algebra $\theta [{\omega _L}]$
\resp{$\theta [{\omega _R}]$} of the versal deformation of 
$\theta ^*_{G^{(n)}_L}$ \resp{$\theta ^*_{G^{(n)}_R}$} (see definition 2.2.7) since the base semisheaf $\theta _{S_L}$ \resp{$\theta _{S_R}$} has been partially or completely extracted from
$\theta ^*_{G^{(n)}_L}$ \resp{$\theta ^*_{G^{(n)}_R}$}: this also corresponds to an extension of the desingularization process as it will be described in the next chapter.\epr
\vskip 11pt

\subsubsection{Corollary}

{\em The extension of the quotient algebra of the versal deformation of the singular semisheaf $\theta ^*_{G^{(n)}_L}$ \resp{$\theta ^*_{G^{(n)}_R}$} is maximal if the spreading-out isomorphism is given by:
\[SOT^{\max}_L = (\tau _{V_{W_L}}\circ \Pi^{\max} _{s_L})
\qquad \resp{\quad SOT^{\max}_R = (\tau _{V_{W_R}}\circ \Pi ^{\max}_{s_R})}.\]
}
\vskip 11pt

\bpr
Indeed, in this case, the base semisheaf $\theta _{S_L}$ \resp{$\theta _{S_R}$} has been completely pulled out from the quotient algebra $\theta [{\omega _L}]$
\resp{$\theta [{\omega _R}]$} in the sense of corollary 3.1.14 and projected in the vertical tangent space $T_{V_{W_L}}$ \resp{$T_{V_{W_R}}$} according to the map $\tau _{V_{W_L}}$ \resp{$\tau _{V_{W_R}}$}.\epr
\vskip 11pt

\subsubsection{Proposition}

{\em  The spreading-out isomorphism $SOT_L$ \resp{$SOT_R$} is locally stable if the generated disconnected semisheaves
$T_{V_{W^i_L}}(\theta ^1_I(\omega ^i_L)_{f^I_i})$ \resp{$T_{V_{W^i_R}}(\theta ^1_I(\omega ^i_R)_{f^I_i})$} are locally free semisheaves.}
\vskip 11pt

\bpr That is to say that the semisheaves
$T_{V_{W^i_L}}(\theta ^1_I(\omega ^i_L)_{f^I_i})$ \resp{$T_{V_{W^i_R}}(\theta ^1_I(\omega ^i_R)_{f^I_i})$} are free of singularities.\epr
\vskip 11pt

\subsubsection{Proposition}

{\em The maximal number of complementary disconnected semisheaves
$T_{V_{W^i_L}}(\theta ^1_I(\omega ^i_L)_{f^I_i})$ \resp{$T_{V_{W^i_R}}(\theta ^1_I(\omega ^i_R)_{f^I_i})$} is equal to $s$~.
}
\vskip 11pt

\bpr Indeed, the integer ``~$s$~'' is the codimension of the versal deformation of the semisheaves $\theta ^*_{G_L^{(n)}}$ \resp{$\theta ^*_{G_R^{(n)}}$}.\epr
\vskip 11pt

\subsubsection{Gluing-up disconnected base semisheaves}

Let $\theta ^1_i(D\omega ^i_L)$ \resp{$\theta ^1_i(D\omega ^i_R)$} and
$\theta ^1_j(D\omega ^j_L)$ \resp{$\theta ^1_j(D\omega ^j_R)$} denote the $i$-th and $j$-th complementary semisheaves 
$T_{V_{W^i_L}}(\theta ^1_I(\omega ^i_L)_{f^I_i})$ \resp{$T_{V_{W^i_R}}(\theta ^1_I(\omega ^i_R)_{f^I_i})$} and
$T_{V_{W^j_L}}(\theta ^1_I(\omega ^j_L)_{f^I_j})$ \resp{$T_{V_{W^j_R}}(\theta ^1_I(\omega ^j_R)_{f^I_j})$} extracted from the base semisheaf $\theta _{S_L}$ \resp{$\theta _{S_R}$}.  These semisheaves can be glued together in a compact way according to:

For each pair $(i,j)$~, let $\Pi _{ij}$ be an isomorphism from
$\theta ^1_j(D(\omega ^i_L)\cap D(\omega ^j_L))$ to
$\theta ^1_i(D(\omega ^i_L)\cap D(\omega ^j_L))$ where
$D(\omega  ^i_L)$ and $D(\omega  ^j_L)$ denote the domains on which these semisheaves $\theta ^1_i(D\omega ^i_L)$ and $\theta ^1_j(D\omega ^i_L)$ are respectively defined.

Then, there exists a semisheaf $\theta ^1(D(\omega ^{i-j}_L))$~, defined on the connected domain $D(\omega ^{i-j}_L) = D(\omega ^i_L)\cup D(\omega ^j_L)$~, and an isomorphism $n_i$ from
$\theta ^1(D\omega ^i_L)$ to $\theta ^1_i(D\omega ^i_L)$ such that $\Pi _{ij}=n_i\circ n^{-1}_j$ in each point of 
$D(\omega ^i_L)\cap D(\omega ^j_L)$~, $\forall\ i,j$~, $1\le i,j\le s$~: this is an adapted version of a proposition of J.P. Serre \cite{Ser1}.

So, $\theta ^1(D(\omega ^{i-j}_L))$ is the semisheaf corresponding to the gluing-up of the semisheaves $\theta ^1_i(D\omega ^i_L)$and $\theta ^1_j(D\omega ^j_L)$~.

\vskip 11pt
Note that the right case ``~$R$~'' can be handled similarly and parallely.
\vskip 11pt

\subsubsection{Sequence of spreading-out isomorphims}

Let $\theta _{SOT(1)_L}$ \resp{$\theta _{SOT(1)_R}$} denote the family of disconnected base semisheaves\linebreak 
$\{ \theta ^1_i(D\omega ^{i}_L)\}^s_{i=1}$
\resp{$\{ \theta ^1_i(D\omega ^{i}_R)\}^s_{i=1}$} of the extension of the quotient algebra of the versal deformation $SOT(1)_L$ \resp{$SOT(1)_R$}.  This family of semisheaves, having been glued together according to section 3.1.20, covers partially the product
$\theta ^*_{G^{(n)}_L}\times \theta '_{S_L}$
\resp{$\theta ^*_{G^{(n)}_R}\times \theta '_{S_R}$} of the singular semisheaf
$\theta ^*_{G^{(n)}_L}$ \resp{$\theta ^*_{G^{(n)}_R}$} by the residue fibre $\theta '_{S_L}$ \resp{$\theta '_{S_R}$} of the versal deformation having not been disconnected from
$\theta ^*_{G^{(n)}_L}$ \resp{$\theta ^*_{G^{(n)}_R}$} by the isomorphisms
$\Pi _{S_L}(f^I_1,\cdots,f^I_i,\cdots,f^I_s)\in SOT(1)_L$
\resp{$\Pi _{S_R}(f^I_1,\cdots,f^I_i,\cdots,f^I_s)\in SOT(1)_R$}.

If the spreading-out isomorphism $SOT(1)_L$ \resp{$SOT(1)_R$} is not locally stable, as noted in proposition 3.1.18, then singular germs
$\omega ^i_{j_{\delta _L}}$ \resp{$\omega ^i_{j_{\delta _R}}$}, $1\le i\le s$~, $1\le j_\delta \le r$~, on the sections of the base semisheaves $\theta _{SOT(1)_L}$ \resp{$\theta _{SOT(1)_R}$} can be degenerated.

Consequently, a versal deformation of $\theta _{SOT(1)_L}$ \resp{$\theta _{SOT(1)_R}$} can be envisaged followed by a spreading-out isomorphism
$SOT(2)_L$ \resp{$SOT(2)_R$}.  The resulting family of disconnected base semisheaves of the extension of the quotient algebra of the versal deformation $SOT(2)_L$ \resp{$SOT(2)_R$} will be noted
$\theta _{SOT(2)_L}$ \resp{$\theta _{SOT(2)_R}$}.

So, a set of ``~$h$~'' versal deformations followed by ``~$h$~'' spreading-out isomorphisms can be envisaged until the disconnected base semisheaves $\theta _{SOT(h)_L}$ \resp{$\theta _{SOT(h)_R}$} are free or locally stable, $h\le (s-1)$~, $h\in\NN$~.
\vskip 11pt

\subsubsection{Remark}  The spreading-out isomorphism was envisaged for singular semisheaves $\theta ^*_{G^{(n)}_L}$ \resp{$\theta ^*_{G^{(n)}_R}$} having isolated singularities of corank 1~.  If it is referred to section 2.2.3 where the versal deformation of germs of corank 2 is considered, it is not difficult to develop the spreading-out isomorphism for unfolded germs of corank 2 (and corank 3) similarly as it was done for germs of corank 1~.
\vskip 11pt

\subsection{Inner dynamics and strange attractors}

This section envisages the spreading-out isomorphism from a differentiable dynamical point of view in such a way that strange attractors, related to the versal deformation of singular germs, blow up under the spreading-out isomorphism into new disconnected attractors.
\vskip 11pt

\subsubsection[Left and right diffeomorphisms of $G^{(n)}({F^+_{\o v}}\times  {F^+_v})$]{\bbf Left and right diffeomorphisms of $G^{(n)}( {F^+_{\o v}}\times  {F^+_v})$}

The generation of the bilinear algebraic semigroup $G^{(n)}( {F^+_{\o v}}\times  {F^+_v})$ can also be studied from a diffeomorphic point of view leading to an inner bidynamics.  In this respect, the differentiable biaction of the product, right by left,
$ {F^+_{\o v}}\times  {F^+_v}$ of 
the sets of completions 
of the bilinear algebraic semigroup
$G^{(n)}({F^+_{\o v}}\times {F^+_v})$ is considered.  This biaction is a homomorphism:
\[  {F^+_{\o v}}\times  {F^+_v}\To \Diff 
(G^{(n)}( {F^+_{\o v}}\times  {F^+_v}))\]
in such a way that
\[ ( {F^+_{\o v}}\times  {F^+_v}) \times
G^{(n)}( {F^+_{\o v}}\times  {F^+_v})\To
G^{(n)}( {F^+_{\o v}}\times  {F^+_v})\]
is differentiable.

$\Diff(G^{(n)}( {F^+_{\o v}}\times  {F^+_v}))$ denotes the group of all diffeomorphisms of
$G^{(n)}( {F^+_{\o v}}\times {F^+_v})$ splitting into:
\[\Diff(
G^{(n)}( {F^+_{\o v}}\times  {F^+_v}))
= \Diff_R(G_R^{(n)}( {F^+_{\o v}})
\times \Diff_L(G_L^{(n)}(  {F^+_v})\]
where $\Diff_L(G^{(n)}( {F^+_{v}}))$ \resp{$\Diff_R(G^{(n)}( {F^+_{\o v}}))$} denotes the semigroup of \lr diffeomorphisms.

If these diffeomorphisms are studied from the point of view of orbit structure, then a \lr generator $f_L\in \Diff_L(G_L^{(n)}( {F^+_{v}})$
\resp{$f_R\in \Diff_R(G_R^{(n)}( {F^+_{\o v}})$} must be taken into account as acting on an irreducible completion $F^+_{v^1_{1_\delta }}$
\resp{$F^+_{\o v^1_{1_\delta }}$} of rank $N$ in such a way that the orbits of
$F^+_{v^1_{1_\delta }}$ \resp{$F^+_{\o v^1_{1_\delta }}$} relative to $f_L$ \resp{$f_R$} are the \lr subsets 
$\{f_L^{j_\delta }(F^+_{v^1_{1_\delta }})\}^r_{j_\delta =1}$
\resp{$\{f_R^{j_\delta }(F^+_{\o v^1_{1_\delta }})\}^r_{j_\delta =1}$}
of the $j_\delta $ conjugacy classes of 
$G^{(n)}_L( {F^+_{v}})$ \resp{$G^{(n)}_R( {F^+_{\o v}})$}.
\vskip 11pt

\subsubsection{Proposition}
{\em A set of $(j_\delta )^{n}$ \lr orbits $f_L^{j_\delta }(F^+_{v^1_{1_\delta }})$ \resp{$f_R^{j_\delta }(F^+_{\o v^1_{1_\delta }})$} of $F^+_{v^1_{1 }}$ \resp{$F^+_{\o v^1_{1 }}$} relative to $f_L$ \resp{$f_R$} constitutes the structure of the $(j_\delta,m_{j_\delta })$-th conjugacy class representative 
$g^{(n)}_L[j_\delta ,m_{j_\delta }]$
\resp{$g^{(n)}_R[j_\delta ,m_{j_\delta }]$} of
$G^{(n)}_L( {F^+_{v}})$ \resp{$G^{(n)}_R({F^+_{\o v}})$}.

Each orbit $f_L^{j_\delta }(F^+_{v^1_{1_\delta }})$
\resp{$f_R^{j_\delta }(F^+_{\o v^1_{1_\delta }})$} is composed of $j_\delta $ irreducible completions 
$F^+_{v^{j'_\delta }_{j_\delta }}$
\resp{$F^+_{\o v^{j'_\delta }_{j_\delta }}$}, $1\le j'_\delta \le j_\delta $~, of rank $N$ and is associated with a Frobenius substitution:
\[ f_L(F^+_{v^1_{1_\delta }})\To f^{j_\delta }_L(F^+_{v^1_{1_\delta }})
\quad\rresp{f_R(F^+_{\o v^1_{1_\delta }})\To f^{j_\delta }_R(F^+_{\o v^1_{1_\delta }})}.\]
}
\vskip 11pt

\bpr \Be
\item As the \lr orbit
$ f^{j_\delta }_L(F^+_{v^1_{j_\delta }})$
\resp{$ f^{j_\delta }_R(F^+_{\o v^1_{j_\delta }})$} is the image of the map:
\[  f^{j_\delta }_L : \quad F^+_{v^1_{j_\delta }}\To
g^{(n)}_L[j_\delta ,m_{j_\delta }] \quad
\rresp{f^{j_\delta }_R : \quad F^+_{\o v^1_{j_\delta }}\To
g^{(n)}_R[j_\delta ,m_{j_\delta }]},\]
it must correspond to a divisor at $j_\delta$ irreducible completions
$F^+_{v^{j'_\delta }_{j_\delta }}$
\resp{$F^+_{\o v^{j'_\delta }_{j_\delta }}$}, $1\le j'_\delta \le j_\delta $~, of rank $N$ according to \cite{Pie1}.

As a result, the Frobenius substitution:
\[f_L\To f_L^{j_\delta }\qquad \rresp{f_R\To f_R^{j_\delta }}\]
on the generator $f_L\in \Diff_L(G_L^{(n)}( {F^+_{v}})$
\resp{$f_R\in \Diff_R(G_R^{(n)}({F^+_{\o v}})$} follows, expressing that we are dealing with one-dimensional components of representatives of the $j_\delta $-th conjugacy class of
$G_L^{(n)}( {F^+_{v}})$
\resp{$G_R^{(n)}({F^+_{\o v}})$}.

\item As the rank of the conjugacy class representative
$g^{(n)}_L[j_\delta ,m_{j_\delta }]$ is $r^{(n)}_{v_{j_\delta }}=(j_\delta \centerdot N)^n$ and its global residue is
$f^{(n)}_{v_{j_\delta }}=j^n_\delta $~, the number of \lr orbits
$f_L^{j_\delta }(F^+_{v^1_{j_\delta }})$
\resp{$f_R^{j_\delta }(F^+_{\o v^1_{j_\delta }})$} of
$g^{(n)}_L[j_\delta ,m_{j_\delta }]$
\resp{$g^{(n)}_R[j_\delta ,m_{j_\delta}]$} must be equal to
\be n_{O_{f^{j_\delta }_L}}\equiv n_{O_{f^{j_\delta }_R}}
=(j_\delta )^{n}\;.
\tag*{\eop}\ee
\Ee

\subsubsection{Proposition}

{\em Let $\Out(G^{(n)}({F^+_{\o v}}\times {F^+_v}))$ denote the group of outer automorphisms of
$G^{(n)}({F^+_{\o v}}\times {F^+_v})$~.  If the bilinear algebraic semigroup $G^{(n)}({F^+_{\o v}}\times {F^+_v})$ is assumed to be $(C^r)$-differentiable, then the isomorphism:
\[ I_{\Out\to\Diff}: \quad
\Out(G^{(n)}({F^+_{\o v}}\times {F^+_v}))
\To \Diff(G^{(n)}({F^+_{\o v}}\times {F^+_v})\]
follows.}
\vskip 11pt

\bpr According to \cite{Pie1}, the group of outer automorphisms
$\Out(G^{(n)}({F^+_{\o v}}\times {F^+_v}))$ of
$G^{(n)}({F^+_{\o v}}\times {F^+_v})$ is given by:

\[\Out(G^{(n)}({F^+_{\o v}}\times  {F^+_v})
=\Aut(G^{(n)}( {F^+_{\o v}}\times  {F^+_v})\big/
\Int(G^{(n)}( {F^+_{\o v}}\times  {F^+_v})\]
with $\Int(G^{(n)}( {F^+_{\o v}}\times  {F^+_v})=\Aut(
P^{(n)}( {F^+_{\o v^ 1}}\times  {F^+_{v^1}})$ where
$\Aut(P^{(n)}( {F^+_{\o v^1}}\times {F^+_{v^1}})$ is the group of automorphisms of the bilinear parabolic semigroup
$P^{(n)}( {F^+_{\o v^1}}\times  {F^+_{v^1}})$~.

As a result, $G^{(n)}( {F^+_{\o v}}\times  {F^+_v})$ acts on
$P^{(n)}( {F^+_{\o v^1}}\times {F^+_{v^1}})$ by conjugation generating by this way the conjugacy class representatives
$g^{(n)}\RL[j_\delta ,m_{j_\delta }]$ of
$G^{(n)}( {F^+_{\o v}}\times  {F^+_v})$~.

And, these conjugacy class representatives
$g^{(n)}\RL[j_\delta ,m_{j_\delta }]$ can be decomposed into normal crossing completions resulting from their compactification as developed in \cite{Pie1}.

So, we have that:
\[g^{(n)}\RL[j_\delta ,m_{j_\delta }]\simeq 
\prod_n(F^+_{\o v_{j_\delta },m_{j_\delta }} \times 
F^+_{v_{j_\delta },m_{j_\delta }})\] implying that
$g^{(n)}_L[j_\delta ,m_{j_\delta }]$
\resp{$g^{(n)}_R[j_\delta ,m_{j_\delta }]$} can be partitioned into $(j_\delta )^{n}$ completions from an algebraic point of view.  Indeed, it has been seen in \cite{Pie1} that
\[\Out(g^{(n)}\RL[j_\delta ,m_{j_\delta }])
\simeq \prod_n(\Gal(\widetilde F^+_{\o v_{j_\delta },m_{j_\delta }}\big/F^0) \times \Gal(\widetilde F^+_{v_{j_\delta },m_{j_\delta }}\big/F^0))\]
where $\Out(g^{(n)}\RL[j_\delta ,m_{j_\delta }])\subset
\Out(G^{(n)}( {F^+_{\o v}}\times  {F^+_v}))$~.

Thus, a one-to-one correspondence can be established between the elements of\linebreak $\Out(G^{(n)}( {F^+_{\o v}}\times  {F^+_v}))$ and the elements of 
$\Diff(G^{(n)}( {F^+_{\o v}}\times  {F^+_v}))$ leading to the generation of the conjugacy class representatives
$g^{(n)}_L[j_\delta ,m_{j_\delta }]$
\resp{$g^{(n)}_R[j_\delta ,m_{j_\delta }]$}
 which decompose into completions from an algebraic point of view or into orbits from a differentiable point of view in such a way that they correspond bijectively.\epr
 \vskip 11pt
 
\subsubsection{Inner bidynamics}

\Bena
\item By a one parameter  semigroup of \lr diffeomorphisms
$\Diff_L(G^{(n)}_L({F^+_{v}}))$ \resp{$\Diff_R(G^{(n)}_L({F^+_{\o v}}))$} of the algebraic semigroup
$G^{(n)}_L( {F^+_{v}})$ \resp{$G^{(n)}_L( {F^+_{\o v}})$}, we mean a continuous map:
\begin{align*}
f^{j_\delta }_L: \quad F^+_{v^1_{j_\delta }}\times g^{(n)}_L[j_\delta ,m_{j_\delta }] &\To g^{(n)}_L[j_\delta ,m_{j_\delta }]\\[11pt]
\rresp{f^{j_\delta }_R: \quad  g^{(n)}_R[j_\delta ,m_{j_\delta }]\times F^+_{\o v^1_{j_\delta }} &\To g^{(n)}_R[j_\delta ,m_{j_\delta }]}\end{align*}
for every conjugacy class representative of $G^{(n)}_L( {F^+_{v}})$
\resp{$G^{(n)}_R( {F^+_{\o v}})$} such that, for all 
$k_\delta ,\ell_\delta \in\nit$ verifying $j_\delta =k_\delta +\ell_\delta $~, we have that:
\begin{align*}
f^{j_\delta }_{L;k_\delta +\ell_\delta}(x_{j_{\delta ^{(1)}_L}})&= f^{j_\delta }_{L;k_\delta}\centerdot f^{j_\delta }_{L;\ell_\delta}(x_{j_{\delta ^{(1)}_L}}) \;, \qquad x_{j_{\delta ^{(1)}_L}}\in F^+_{v^1_{j_\delta }}\\[11pt]
\rresp{f^{j_\delta }_{R;k_\delta +\ell_\delta}(x_{j_{\delta ^{(1)}_R}})&= f^{j_\delta }_{R;k_\delta}\centerdot f^{j_\delta }_{R;\ell_\delta}(x_{j_{\delta ^{(1)}_R}}) \;, \qquad x_{j_{\delta ^{(1)}_L}}\in F^+_{\o v^1_{j_\delta }}},
\end{align*}
where $k_\delta $ and $\ell_\delta$ refer to the numbers of irreducible completions in $f^{j_\delta }_L$ \resp{$f^{j_\delta }_R$}~.

\item The tangent (semi)space to the one-parameter  semigroup of \lr diffeomorphisms of $G^{(n)}_L( {F^+_{v}})$
\resp{$G^{(n)}_R( {F^+_{\o v}})$} is the space of sections
$\Gamma (T(G^{(n)}_L( {F^+_{v}})))$
\resp{$\Gamma (T(G^{(n)}_R( {F^+_{\o v}})))$} of the tangent bundle
$T(G^{(n)}_L( {F^+_{v}}))$
\resp{$T(G^{(n)}_R( {F^+_{\o v}}))$} whose fibres in each point
$x_{j_{\delta _L}}\in g^{(n)}_L[j_\delta ,m_{j_\delta }]$
\resp{$x_{j_{\delta _R}}\in g^{(n)}_R[j_\delta ,m_{j_\delta }]$
} are given by the tangent vectors $\vec n_L(x_{j_{\delta _L}})$
\resp{$\vec n_R(x_{j_{\delta _R}})$} defined by:
\begin{align*}
\vec n_L(x_{j_{\delta _L}}) &= \F d{dt}\ \L(f_{L;t}^{j_\delta }(x_{j_{\delta ^{(1)}_L}})\R)_{t=0}\\[11pt]
\Big( \text{resp.} \quad \vec n_R(x_{j_{\delta _R}}) &= \F d{dt}\ \L(f_{R;t}^{j_\delta }(x_{j_{\delta ^{(1)}_R}})\R)_{t=0}\quad \Big).\end{align*}

This allows to generate a \lr internal dynamics of the algebraic semigroup $G^{(n)}_L( {F^+_{v}})$
\resp{$G^{(n)}_R( {F^+_{\o v}})$} taking into account that the tangent bundle $T(G^{(n)}_L( {F^+_{v}}))$
\resp{$T(G^{(n)}_R( {F^+_{\o v}}))$} has to be horizontal in order that the tangent vectors $\vec n_L(x_{j_{\delta _L}})$ \resp{$\vec n_R(x_{j_{\delta _R}})$} be rotational velocity vectors.

\item An internal bidynamics of the bilinear algebraic semigroup
$G^{(n)} ( {F^+_{\o v}}\times {F^+_{v}})$ can be reached by considering the horizontal tangent bibundle
\[ T(G^{(n)}( {F^+_{\o v}}\times {F^+_{v}}))
= T(G^{(n)}_R( {F^+_{\o v}}))\times T(G^{(n)}_L( {F^+_{v}}))\]
whose (bi)fibres in each bipoint $x_{j_{\delta _R}}\times x_{j_{\delta _L}}\in g^{(n)}\RL[j_\delta ,m_{j_\delta }]$ are the tangent bivectors $\vec n_R(x_{j_{\delta _R}})\times\vec n_L(x_{j_{\delta _L}})$~.
\Ee
\vskip 11pt

\subsubsection{Translated orbits in the neighbourhood of singular germs}

We are now interested by the dynamics around singularities on the sections $\phi^{\TAN}_{G_{j_{\delta _L}}}$ \resp{$\phi^{\TAN}_{G_{j_{\delta _R}}}$} of the space of sections 
$\Gamma (T(G^{(n)}_L( {F^+_{v}})))$
\resp{$\Gamma (T(G^{(n)}_R( {F^+_{\o v}})))$} of the tangent bundle
$T(G^{(n)}_L( {F^+_{v}}))$
\resp{$T(G^{(n)}_R( {F^+_{\o v}}))$} on the conjugacy class representatives $g^{(n)}_L[j_\delta ,m_{j_\delta }]$
\resp{$g^{(n)}_R[j_\delta ,m_{j_\delta }]$}.

Let then $\phi _{j_\delta }(\omega _L)$ \resp{$\phi _{j_\delta }(\omega _R)$} denote a singular germ of corank 1 and codimension $s$ on the $n$-dimensional real-valued differentiable function 
$\phi^{\TAN}_{G_{j_{\delta _L}}}(x^{\TAN}_{g_{j_\delta }})\in \Gamma (T(G^{(n)}_L( {F^+_{v}})))$ \resp{$\phi^{\TAN}_{G_{j_{\delta _R}}}(x^{\TAN}_{g_{j_\delta }})\in T(G^{(n)}_R( {F^+_{\o v}}))$} on $g^{(n)}_L[j_\delta ,m_{j_\delta }]$
\resp{$g^{(n)}_R[j_\delta ,m_{j_\delta }]$}.

As $g^{(n)}_L[j_\delta ,m_{j_\delta }]$
\resp{$g^{(n)}_R[j_\delta ,m_{j_\delta }]$} can be decomposed into a set of \lr orbits $f_L^{j_\delta }(F^+_{v^1_{1_\delta }})$
\resp{$f_R^{j_\delta }(F^+_{\o v^1_{1_\delta }})$} relative to $f_L\in\Diff(G^{(n)}_L( {F^+_v}))$ \resp{$f_R\in\Diff(G^{(n)}_R({ {F^+_{\o v}}}))$}, the differentiable function 
$\phi^{\TAN}_{G_{j_{\delta _L}}}$ \resp{$\phi^{\TAN}_{G_{j_{\delta _R}}}$} can also be decomposed into one-dimensional subfunctions $f^{j_\delta }_{L;\TAN}$ \resp{$f^{j_\delta }_{R;\TAN}$} corresponding to orbits translated from 
$f_L^{j_\delta }(F^+_{v^1_{1_\delta }})$
\resp{$f_R^{j_\delta }(F^+_{\o v^1_{1_\delta }})$} under the tangent bundle $T(G^{(n)}_L( {F^+_{v}}))$
\resp{$T(G^{(n)}_R( {F^+_{\o v}}))$}:
\[ T(G^{(n)}_L( {F^+_{v}})): \quad f^{j_\delta }_{L;\TAN} \To f^{j_\delta }_{L}\;, \quad \rresp{T(G^{(n)}_R( {F^+_{\o v}})): \quad f^{j_\delta }_{R;\TAN} \To f^{j_\delta }_{R}}.\]

As we are concerned with a singular germ $\phi _{j_\delta }(\omega _L)$ \resp{$\phi _{j_\delta }(\omega _R)$} on the differentiable function $\phi^{\TAN}_{G_{j_{\delta _L}}}(x^{\TAN}_{g_{j_\delta }})$ \resp{$\phi^{\TAN}_{G_{j_{\delta _R}}}(x^{\TAN}_{g_{j_\delta }})$} which is perturbed in the neighbourhood $D_{\phi^{\TAN}_{G_{j_{\delta _L}}}}$ \resp{$D_{\phi^{\TAN}_{G_{j_{\delta _R}}}}$} of $\phi _{j_\delta }(\omega _L)$ \resp{$\phi _{j_\delta }(\omega _R)$}, it shall be assumed that the ``translated'' orbits on $D_{\phi^{\TAN}_{G_{j_{\delta _L}}}}$ \resp{$D_{\phi^{\TAN}_{G_{j_{\delta _R}}}}$}
are perturbed one-dimensional subfunctions
$f^{a_{j_\delta }}_{L;\TAN}$ \resp{$f^{a_{j_\delta }}_{R;\TAN}$}
of $f^{j_\delta }_{L;\TAN}$ \resp{$f^{j_\delta }_{R;\TAN}$} characterized by a rank
$r_{f^{a_{j_\delta }}_{L;\TAN}}=r_{f^{a_{j_\delta }}_{R;\TAN}}=a\centerdot N$ where the global residue degree $a$ is inferior or equal to $j_\delta $~: $a\le j_\delta $~.

Before going closely into the study of the neighbourhoods of singular germs from a diffeomorphic point of view, the tangent space decomposition into contracting and expanding components will be recalled.
\vskip 11pt

\subsubsection{Splitting of the tangent space into stable, unstable and neutral subsets}

Let $f_{L;\TAN}\in \Diff_L(T(G^{(n)}_L( {F^+_{v}})))$
\resp{$f_{R;\TAN}\in \Diff_R(T(G^{(n)}_R( {F^+_{\o v}})))$} denote the \lr generator of the diffeomorphisms of the space of sections of the tangent bundle
$T(G^{(n)}_L( {F^+_{v}}))$
\resp{$T(G^{(n)}_R( {F^+_{\o v}}))$}.

\Bi
\item A point $x^{\TAN;N\omega }_{g_{j_{\delta _L}}}\in \phi ^{\TAN}_{G_{j_{\delta_L} }} $ \resp{$x^{\TAN;N\omega }_{g_{j_{\delta _R}}}\in \phi^{\TAN}_{G_{j_{\delta_R}} } $}
is non-wandering if, for every neighbourhood 
$ U_{\phi ^{\TAN}_{j_{\delta_L} }}$ \resp{$ U_{\phi ^{\TAN}_{j_{\delta _R}}}$} of 
$x^{\TAN;N\omega }_{g_{j_{\delta _L}}}$ \resp{$x^{\TAN;N\omega }_{g_{j_{\delta _R}}}$}, one has:
\begin{align*}
f_{L;\TAN}^{j_\delta }(U_{\phi ^{\TAN}_{j_{\delta_L} } })
\cap U_{\phi ^{\TAN}_{j_{\delta_L} }}&\neq \emptyset \\
\rresp{f_{R;\TAN}^{j_\delta }(U_{\phi ^{\TAN}_{j_{\delta_R} }})
\cap U_{\phi ^{\TAN}_{j_{\delta_R} } }&\neq \emptyset}.
\end{align*}

The set of non wandering points forms a closed invariant set noted
$\Omega _{\phi ^{\TAN}_{j_{\delta _L}}}$
\resp{$\Omega _{\phi ^{\TAN}_{j_{\delta _R}}}$}.  The other points of $\phi ^{\TAN}_{G_{j_{\delta L}}}$ \resp{$\phi ^{\TAN}_{G_{j_{\delta _R}}}$} are called wandering points and form invariant open subsets \cite{Sma}.

\item A linear automorphism $u_L$ \resp{$u_R$} of the tangent space
$\Gamma (T(G^{(n)}_L( {F^+_{v}})))$
\resp{$\Gamma (T(G^{(n)}_R( {F^+_{\o v}})))$} will be said to be flat, contracting or expanding if, under
\begin{align*}  u_L: \quad \Gamma (T(G^{(n)}_L( {F^+_{v}})))  &\To
\Gamma (T(G^{(n)}_L( {F^+_{v}}))) \\[11pt]
\rresp{u_R: \quad \Gamma (T(G^{(n)}_R( {F^+_{\o v}})))  &\To
\Gamma (T(G^{(n)}_R( {F^+_{\o v}})))},\end{align*}
the eigenvalues of
\begin{align*}
|u_L(f^{j_\delta }_{L;\TAN}(x^{\TAN}_{j_{\delta ^{(1)}_L}}))|
& =\lambda ^{j_\delta }
|f^{j_\delta }_{L;\TAN}(x^{\TAN}_{j_{\delta ^{(1)}_L}})|\;, \qquad \forall\ 1\le j_\delta \le r\;, \\[11pt]
\rresp{|u_R(f^{j_\delta }_{R;\TAN}(x^{\TAN}_{j_{\delta ^{(1)}_R}}))|
& =\lambda ^{j_\delta }
|f^{j_\delta }_{R;\TAN}(x^{\TAN}_{j_{\delta ^{(1)}_R}})|\;, \qquad \forall\ 1\le j_\delta \le r\;, }
\end{align*}
satisfy respectively  $|\lambda ^{j_\delta }|=1$~, $|\lambda ^{j_\delta }|<1$ or $|\lambda ^{j_\delta }|>1$~.

So, the tangent space $\Gamma (T(G^{(n)}_L( {F^+_{v}})))$
\resp{$\Gamma (T(G^{(n)}_R( {F^+_{\o v}})))$} exhibits a splitting into:
\Bi
\item stable subsets $E^s(f^{j_\delta }_{L;\TAN}(x^{\TAN}_{j_{\delta ^{(1)}_L}}))$ \resp{$E^s(f^{j_\delta }_{R;\TAN}(x^{\TAN}_{j_{\delta ^{(1)}_R}}))$} for which $u_L$ \resp{$u_R$} is contracting.
\item unstable subsets $E^u(f^{j_\delta }_{L;\TAN}(x^{\TAN}_{j_{\delta ^{(1)}_L}}))$ \resp{$E^u(f^{j_\delta }_{R;\TAN}(x^{\TAN}_{j_{\delta ^{(1)}_R}}))$} for which $u_L$ \resp{$u_R$} is expanding.
\item neutral subsets $E^n(f^{j_\delta }_{L;\TAN}(x^{\TAN}_{j_{\delta ^{(1)}_L}}))$ \resp{$E^n(f^{j_\delta }_{R;\TAN}(x^{\TAN}_{j_{\delta ^{(1)}_R}}))$} for which $u_L$ \resp{$u_R$} is flat.
\Ei

Classically, a subset being stable and unstable, i.e. exhibiting respectively  a volume contraction and a volume expansion under a linear automorphism of the tangent space, is said to be ``hyperbolic''.  This terminology will not be adopted here as it will appear in the following.
\Ei
\vskip 11pt

\subsubsection{Proposition}

{\em Let $\Gamma (T(G^{(n)}_L( {F^+_{v}})))$
\resp{$\Gamma (T(G^{(n)}_R( {F^+_{\o v}})))$} denote the space of sections of the tangent bundle to the algebraic semigroup
$G^{(n)}_L( {F^+_{v}})$
\resp{$G^{(n)}_R( {F^+_{\o v}})$}.  Then we have that:
\Bena
\item its stable subsets $E^s(f^{j_\delta }_{L;\TAN}(x^{\TAN}_{j_{\delta ^{(1)}_L}}))$ \resp{$E^s(f^{j_\delta }_{R;\TAN})(x^{\TAN}_{j_{\delta ^{(1)}_R}})$} are characterized by a hyperbolic geometry.

\item its unstable subsets $E^u(f^{j_\delta }_{L;\TAN}(x^{\TAN}_{j_{\delta ^{(1)}_L}}))$ \resp{$E^u(f^{j_\delta }_{R;\TAN})(x^{\TAN}_{j_{\delta ^{(1)}_R}})$} are characterized by a spherical geometry.

\item its neutral subsets $E^n(f^{j_\delta }_{L;\TAN}(x^{\TAN}_{j_{\delta ^{(1)}_L}}))$ \resp{$E^n(f^{j_\delta }_{R;\TAN}(x^{\TAN}_{j_{\delta ^{(1)}_R}}))$} are characterized by an euclidian geometry.
\Ee
}
\vskip 11pt

\bpr With reference to proposition 2.3.1, the developments will be envisaged for the left and right cases without distinction and the notations will be simplified as follows:

Let the point $x^{\TAN}_{j_{\delta ^{(1)}_L}}$  (and $x^{\TAN}_{j_{\delta ^{(1)}_R}}$~) be given by a point $M$ of coordinates $(u^1,\cdots,u^n)$ and let 
$f^{j_\delta }_{L;\TAN}(x^{\TAN}_{j_{\delta ^{(1)}_L}})$ (and {$f^{j_\delta }_{R;\TAN}(x^{\TAN}_{j_{\delta ^{(1)}_R}})$}~) define a point $P$ of coordinates $(x^1,\cdots,x^n)$ in the neighbourhood of $M$~.

The differential $d(\vec{MP})$ of the vector $\vec{MP}$ corresponds to a linear automorphism:
\begin{align*}
u: \quad \Gamma (T(G^{(n)}\LR( {F^+_{v}}))
&\To \Gamma (T(G^{(n)}\LR( {F^+_{v}})))\\
\vec{MP} \quad &\To \quad d(\vec{MP})\end{align*}
where $G^{(n)}\LR( {F^+_{v}})$ is a condensed notation for
$G^{(n)}_L( {F^+_{v}}) $ or $G^{(n)}_R( {F^+_{\o v}})$~.

$d(\vec{MP})$ can be expressed by the differential \cite{Car}
\[ dP=(Dx^1,\cdots,Dx^i,\cdots,Dx^n) \]
where $Dx^i=dx^i+du^i+x^k\ \Gamma ^i_{kr}\ du^r$~, $1\le i\le n$~, with
\[ \Gamma ^i_{kr} = \half\ \L(\F{\partial g_{ik}}{\partial u^r}
+ \F{\partial g_{ir}}{\partial u^k}
+ \F{\partial g_{kr}}{\partial u^i}\R)\;.\]
Let then $h$ be the covering of
$\Gamma (T(G^{(n)}\LR( {F^+_{v}})))$ by the euclidian space $\rit^n$~:
\[ h : \quad \Gamma (T(G^{(n)}\LR( {F^+_{v}})))\To \rit^n\]
given by:
\[ dP \To h(dP)\]
in such a way that:
\begin{align*}
h\circ u : \quad \Gamma (T(G^{(n)}\LR( {F^+_{v}}))& \To \rit^n\;, \\
\vec{MP} \quad &\To  h(d(\vec{MP}))\;.\end{align*}

Three possibilities occur:
\Bena
\item if \quad $\|\vec{MP}\|=\|h(d(\vec{MP}))\|$~, \quad the subsets of $\Gamma (T(G^{(n)}\LR( {F^+_{v}})))$ are locally Euclidian.

It follows that the norm of $\vec{MP}$ is conserved under the composition of maps $(h\circ u)$ and, thus, that the linear automorphism $u$ of $\Gamma (T(G^{(n)}\LR( {F^+_{v}})))$ is unitary.

Consequently, in:
\begin{align} \label{eq:*}
|u_L(f^{j_\delta }_{L;\TAN}(x^{\TAN}_{j_{\delta ^{(1)}_L}}))|
& =\lambda ^{j_\delta }
|f^{j_\delta }_{L;\TAN}(x^{\TAN}_{j_{\delta ^{(1)}_L}})|\;,  \tag{$*$}\\[11pt]
\rresp{|u_R(f^{j_\delta }_{R;\TAN}(x^{\TAN}_{j_{\delta ^{(1)}_R}}))|
& =\lambda ^{j_\delta }
|f^{j_\delta }_{R;\TAN}(x^{\TAN}_{j_{\delta ^{(1)}_R}})|},\notag
\end{align}
\Be
\item $|\lambda ^{j_\delta }|=1$~.
\item $u_L$ \resp{$ u_R$} is a flat automorphism characterized by an euclidian geometry on the subsets of $\Gamma (T(G^{(n)}\LR( {F^+_{v}})))$ which are neutral.
\Ee

\item if \quad $\|\vec{MP}\|>\|h(d(\vec{MP}))\|$~, \quad one must admit, according to proposition 2.3.1, that $d(\vec{MP})$ is given by the differential of $P$ whose components $Dx^i$ are:
\[ Dx^i=dx^i+du^i+x^k\ \Gamma ^i_{kr}\ du^r-\kappa \ g_{ik}\ du^k \qquad \text{with}\quad \kappa \in\rit\;.\]
The metric $g_{ik}$ is locally hyperbolic and the curvature of 
$\Gamma (T(G^{(n)}\LR( {F^+_{v}}))$ is locally negative or equal to $-\kappa$~.

The norm $\vec{MP}$ is not conserved under $(h\circ u)$~.

Thus, in \eqref{eq:*}, we have that:
\Be
\item $|\lambda ^{j_\delta }|<1$~.
\item $u_L$ \resp{$ u_R$} is a hyperbolic automorphism which is contracting.
\item the subsets of $\Gamma (T(G^{(n)}\LR( {F^+_{v}})))$ are stable and characterized by a hyperbolic geometry.
\Ee

\item if \quad $\|\vec{MP}\|<\|h(d(\vec{MP}))\|$~, \quad the components $Dx^i$ of $dP$ are given by:
\[ Dx^i=dx^i+du^i+x^k\ \Gamma ^i_{kr}\ du^r+\kappa \ g_{ik}\ du^k\;.\]
The metric $g_{ik}$ is localy spherical and the curvature of
$\Gamma (T(G^{(n)}\LR( {F^+_{v}})))$ is locally positive or equal to $+\kappa $~.

Then, in \eqref{eq:*}, we have that:
\Be
\item $|\lambda ^{j_\delta }|>1$~.
\item $u_L$ \resp{$ u_R$} is a ``spherical'' automorphism which is expanding.
\item the subsets of $\Gamma (T(G^{(n)}\LR( {F^+_{v}})))$ are locally unstable and characterized by a spherical geometry.\epr
\Ee
\Ee
\vskip 11pt

\subsubsection{Proposition: Singular hyperbolic attractor}

{\em Let $\phi _{j_\delta }(\omega _L)$ \resp{$\phi _{j_\delta }(\omega _R)$} be a singular germ (of corank 1 and codimension $s$~) on the $n$-dimensional real-valued differentiable function 
$\phi ^{\TAN}_{G_{j_{\delta _L}}}(x^{\TAN}_{g_{j_\delta }})\in\Gamma (T(G^{(n)}_L( {F^+_{v}}))) $ \resp{$\phi ^{\TAN}_{G_{j_{\delta _R}}} (x^{\TAN}_{g_{j_\delta }})\in \Gamma (T(G^{(n)}_R({F^+_{v}})))$}.
\vskip 11pt
The neighbourhood $D_{\phi ^{\TAN}_{G_{j_{\delta _L}}}}$ \resp{$D_{\phi ^{\TAN}_{G_{j_{\delta_R} }}}$} of the singular germ $\phi _{ {j_{\delta }}}(\omega _L)$ \resp{$\phi _{{j_{\delta }}}(\omega _R)$} on 
$\phi ^{\TAN}_{G_{j_{\delta_L} }}(x^{\TAN}_{g_{j_\delta }})$ \resp{$\phi ^{\TAN}_{G_{j_{\delta _R}}}(x^{\TAN}_{g_{j_\delta }})$} is a singular hyperbolic attractor $\Lambda ^{\TAN}_L$ \resp{$\Lambda ^{\TAN}_R$} with respect to $\Diff_L(T(G^{(n)}_L({F^+_{v}})))$ \resp{$\Diff_R(T(G^{(n)}_R( {F^+_{\o v}})))$} provided that:
\Bena
\item the orbits of $\Lambda ^{\TAN}_L$ \resp{$\Lambda ^{\TAN}_R$} are one-dimensional functions $f^{a_{j_\delta }}_{L;\TAN}$
\resp{$f^{a_{j_\delta }}_{R;\TAN}$} having a rank
$r^{a_{j_\delta }}_{f_L;\TAN}= r^{a_{j_\delta }}_{f_R;\TAN}=a\centerdot N$ and form the basin of attraction of 
$\Lambda ^{\TAN}_L$ \resp{$\Lambda ^{\TAN}_R$}.

\item the singularity is a non-wandering point.

\item the basin of attraction of $\Lambda ^{\TAN}_L$ \resp{$\Lambda ^{\TAN}_R$}  is a stable subset 
$E^s(f^{a_{j_\delta }}_{L;\TAN})$
\resp{$E^s(f^{a_{j_\delta }}_{R;\TAN})$} of
$\Gamma (T(G^{(n)}_L( {F^+_{v}})))$ \resp{$\Gamma (T(G^{(n)}_R( {F^+_{\o v}})))$} characterized by a hyperbolic geometry.
\Ee
}
\vskip 11pt

\bpr \Bena
\item According to section 3.2.5, the orbits on the neighbourhood
$D_{\phi ^{\TAN}_{G_{j_{\delta _L}}}}$ \resp{$D_{\phi ^{\TAN}_{G_{j_{\delta_R} }}}$} of the singularity
$\phi _{j_\delta }(\omega _L)$ \resp{$\phi _{j_\delta }(\omega _R)$} are one-dimensional subfunctions 
$f^{a_{j_\delta }}_{L;\TAN}$
\resp{$f^{a_{j_\delta }}_{R;\TAN}$} of 
$f^{j_\delta }_{L;\TAN}$
\resp{$f^{j_\delta }_{R;\TAN}$}.  Consequently, they constitute the basin of attraction of the singular hyperbolic attractor
$\Lambda ^{\TAN}_L$ \resp{$\Lambda ^{\TAN}_R$}.

\item Let $U_{\phi ^{\TAN}_{j_{\delta _L}}} \subset \Lambda ^{\TAN}_L$
\resp{$U_{\phi ^{\TAN}_{j_{\delta _R}}} \subset \Lambda ^{\TAN}_R$} be a neighbourhood of the singularity included into the singular hyperbolic attractor $\Lambda ^{\TAN}_L$ \resp{$\Lambda ^{\TAN}_R$}.

Then, the singularity is a non wandering point if we have, according to section 3.2.6:
\begin{align*}
f_{L;\TAN}^{a_{j_\delta }}(U_{\phi ^{\TAN}_{j_{\delta _L}}})
\cap U_{\phi ^{\TAN}_{j_{\delta_L}}}&\neq \emptyset \\
\rresp{f_{R;\TAN}^{a_{j_\delta} }(U_{\phi ^{\TAN}_{j_{\delta  _R)}}})
\cap U_{\phi ^{\TAN}_{j_{\delta _R}}}&\neq \emptyset}.
\end{align*}

\item If was proved in proposition 2.3.1 that the geometry is hyperbolic in the neighbourhood of the singularity.  Consequently, the basin of attraction of the singular hyperbolic attractor is a stable subset $E^s(f^{a_{j_\delta }}_{L;\TAN})$
\resp{$E^s(f^{a_{j_\delta }}_{R;\TAN})$} characterized by a hyperbolic geometry whose points are contracting under the linear automorphism $u_L$ \resp{$u_R$} of 
$\Gamma (T(G^{(n)}_L( {F^+_{v}})))$ \resp{$\Gamma (T(G^{(n)}_R( {F^+_{\o v}})))$}.\epr
\Ee\vskip 11pt

\subsubsection{The introduction of unfolded attractors}

\Bi
\item Let the singular germ $\phi _{j_\delta }(\omega _L)$ \resp{$\phi _{j_\delta }(\omega _R)$} of corank 1 and codimension $s$ be given by $y_L=\omega _L^{s+2}$ \resp{$y_R=\omega _R^{s+2}$}.

\item Its versal deformation corresponds to the development:
\begin{align*}
F(\omega _L,a_{ij}(x_L)) &= \omega _L^{s+2}+\sum^s_{i=1} a_{ij_\delta }(x_L)\ \omega ^i_{j_{\delta _L}}\\
\rresp{F(\omega _R,a_{ij}(x_R)) &= \omega _R^{s+2}+\sum^s_{i=1} a_{ij_\delta }(x_R)\ \omega ^i_{j_{\delta _R}}}\end{align*}
where:
\Bi
\item $a_{ij_\delta }(x_L)$ \resp{$a_{ij_\delta }(x_R)$} is a $(n-1)$-dimensional real valued differentiable function defined on the neighbourhood $D_{\phi ^{\TAN}_{G_{j_{\delta _L}}}}$ \resp{$D_{\phi ^{\TAN}_{G_{j_{\delta _R}}}}$} of the singularity $y_L=\omega _L^{s+2}$ \resp{$y_R=\omega _R^{s+2}$} on the $n$-dimensional differentiable function $\phi ^{\TAN}_{G_{j_{\delta _L}}}(x^{\TAN}_{g_{j_\delta }})$ \resp{$\phi ^{\TAN}_{G_{j_{\delta _R}}}(x^{\TAN}_{g_{j_\delta }})$}.

\item $\omega ^i_{j_{\delta _L}}$ \resp{$\omega ^i_{j_{\delta _R}}$} is a divisor, generator of the versal unfolding of $y_L=\omega _L^{s+2}$ \resp{$y_R=\omega _R^{s+2}$}.  It is localized on 
$D_{\phi ^{\TAN}_{G_{j_{\delta _L}}}}$ \resp{$D_{\phi ^{\TAN}_{G_{j_{\delta _R}}}}$} and projected on $a_{ij_\delta }(x_L)$
\resp{$a_{ij_\delta }(x_R)$} according to proposition 2.2.8.
\Ei

\item The set $\bigcup^s_{i=1} a_{ij_\delta }(x_L)$ \resp{$\bigcup^s_{i=1} a_{ij_\delta }(x_R)$} of functions on
$D_{\phi ^{\TAN}_{G_{j_{\delta _L}}}}$ \resp{$D_{\phi ^{\TAN}_{G_{j_{\delta _R}}}}$} can be partitioned into the one-dimensional functions $f^{a_{j_\delta }}_{L;\TAN}$
\resp{$f^{a_{j_\delta }}_{R;\TAN}$}, which are perturbed orbits of
$\Diff_L(T(G^{(n)}_L( {F^+_{v}})))$ \resp{$\Diff_R(T(G^{(n)}_L( {F^+_{\o v}})))$} (see proposition 3.2.8).

\item As the set $\bigcup^s_{i=1} \omega ^i_{j_{\delta _L}}$
\resp{$\bigcup^s_{i=1} \omega ^i_{j_{\delta _R}}$} of generators of the versal unfolding of the singularity $y_L=\omega _L^{s+2}$ \resp{$y_R=\omega _R^{s+2}$} is localized on $D_{\phi ^{\TAN}_{G_{j_{\delta _L}}}}$ \resp{$D_{\phi ^{\TAN}_{G_{j_{\delta _R}}}}$}, it constitutes the basin of an unfolded attractor $\Lambda ^{\TAN}_{unf_L}$ \resp{$\Lambda ^{\TAN}_{unf_R}$}:
\Bean
\item centred on the singularity $y_L=\omega _L^{s+2}$ \resp{$y_R=\omega _R^{s+2}$};
\item having as orbits the generators $\omega ^i_{j_{\delta _L}}$ \resp{$\omega ^i_{j_{\delta _R}}$}, $1\le i\le s$~, which can be rewritten according to:
\[
\omega ^i_{j_{\delta _L}}= f^{\omega ^i_{j_{\delta}}}_{L;\TAN}(\omega _{j_{\delta ^{(1)}_L}})\qquad
\rresp{\omega ^i_{j_{\delta _R}}= f^{\omega ^i_{j_{\delta}}}_{R;\TAN}(\omega _{j_{\delta ^{(1)}_R}})}\]
where:
\Bi 
\item $f^{\omega ^i_{j_{\delta}}}_{L;\TAN}$ \resp{$f^{\omega ^i_{j_{\delta}}}_{R;\TAN}$} is an orbital generator with respect to
$\Diff_L(\omega ^i_{j_{\delta _L}})$ \resp{$\Diff_R(\omega ^i_{j_{\delta _R}})$};
\item $\omega _{j_{\delta ^{(1)}_L}}$ \resp{$\omega _{j_{\delta ^{(1)}_R}}$} is a point of an irreducible completion of
$\omega^i _{j_{\delta _L}}$ \resp{$\omega ^i_{j_{\delta _R}}$}.
\Ei
\Ee
\Ei
\vskip 11pt

\subsubsection{Proposition}

{\em Let $\Lambda ^{\TAN}_{L}$ \resp{$\Lambda ^{\TAN}_{R}$} be a singular hyperbolic attractor centred on a singular germ 
$\phi _{j_\delta }(\omega _L)$ \resp{$\phi _{j_\delta }(\omega _R)$} of corank 1 and codimension $s$.

Then, the versal unfolding of the germ $\phi _{j_\delta }(\omega _L)$ \resp{$\phi _{j_\delta }(\omega _R)$}  involves the map:
\[VD_{\Lambda _L} : \quad \Lambda ^{\TAN}_L \To \Lambda ^{\TAN}_{str_L}\qquad 
\rresp{VD_{\Lambda _R} : \quad \Lambda ^{\TAN}_R \To \Lambda ^{\TAN}_{str_R}}\]
of the singular hyperbolic attractor 
 $\Lambda ^{\TAN}_L $ \resp{$\Lambda ^{\TAN}_{R}$} into the singular strange attractor
 \[ \Lambda ^{\TAN}_{str_L} = \Lambda ^{\TAN}_L\times \Lambda ^{\TAN}_{unf_L} \qquad \rresp{\Lambda ^{\TAN}_{str_R} = \Lambda ^{\TAN}_R\times \Lambda ^{\TAN}_{unf_R}}\]
 where $\Lambda ^{\TAN}_{unf_L} $ \resp{$\Lambda ^{\TAN}_{unf_R} $} is the unfolded attractor introduced in section 3.2.9.
 }
 \vskip 11pt
 
 \bpr  As the versal unfolding of the germ $\phi _{j_\delta }(\omega _L)$ \resp{$\phi _{j_\delta }(\omega _R)$}, given by $y_L=\omega _L^{s+2}$ \resp{$y_R=\omega _R^{s+2}$}, leads to the projection of $s$ divisors $\omega ^i_{j_{\delta _L}}$ \resp{$\omega ^i_{j_{\delta _R}}$} on the neighbourhood $D_{\phi^{\TAN}_{G_{j_{\delta _L}}}}$ \resp{$D_{\phi^{\TAN}_{G_{j_{\delta _R}}}}$} of the singularity and as these divisors $\omega ^i_{j_{\delta _L}}$ \resp{$\omega ^i_{j_{\delta _R}}$}  can be rewritten according to 
 $f^{\omega ^i_{j_\delta }}_{L;\TAN}(\omega _{j_{\delta ^{(1)}_L}})$
 \resp{$f^{\omega ^i_{j_\delta }}_{R;\TAN}(\omega _{j_{\delta ^{(1)}_R}})$} (see section 3.2.9), it is clear that the singular hyperbolic attractor $\Lambda ^{\TAN}_L$ \resp{$\Lambda ^{\TAN}_R$} will undergo an unfolding leading to the attractor
 \[ \Lambda ^{\TAN}_{str_L}= \Lambda ^{\TAN}_L\times \Lambda ^{\TAN}_{unf_L} \qquad 
 \rresp{\Lambda ^{\TAN}_{str_R}= \Lambda ^{\TAN}_R\times \Lambda ^{\TAN}_{unf_R} }\]
 which is a strange attractor.
 
 Indeed, according to proposition 2.3.7, the strata of the neighbourhood of the singularity perturbed by the versal deformation are characterized by a spherical geometry.  These strata are the generators (or the orbits) $f^{\omega ^i_{j_\delta }}_{L;\TAN}(\omega _{j_{\delta ^{(1)}_L}})$ \resp{$f^{\omega ^i_{j_\delta }}_{R;\TAN}(\omega _{j_{\delta ^{(1)}_R}})$} of the versal deformation projected on the orbits $f^{a_{j_\delta }}_{L;\TAN}$
 \resp{$f^{a_{j_\delta }}_{R;\TAN}$} of the singular hyperbolic attractor $\Lambda ^{\TAN}_L$ \resp{$\Lambda ^{\TAN}_R$} .
 
 Consequently, these strata are given by the functions:
\[ f^{str}_{L;\TAN} = f^{a_{j_\delta }\mid\omega ^i_{j_\delta }}_{L;\TAN} \times f^{\omega ^i_{j_\delta }}_{L;\TAN} \qquad 
\rresp{ f^{str}_{R;\TAN} = f^{a_{j_\delta }\mid\omega ^i_{j_\delta }}_{R;\TAN} \times f^{\omega ^i_{j_\delta }}_{R;\TAN}}\]
where  $ f^{a_{j_\delta }\mid\omega ^i_{j_\delta }}_{L;\TAN}$
\resp{$ f^{a_{j_\delta }\mid\omega ^i_{j_\delta }}_{R;\TAN}$} is a one-dimensional subfunction of
 $ f^{a_{j_\delta }}_{L;\TAN}$ \resp{ $ f^{a_{j_\delta }}_{R;\TAN}$} characterized by a rank  $ r_{f^{a_{j_\delta }\mid\omega ^i_{j_\delta }}_{L;\TAN}}=b\centerdot N$ inferior to the rank
$ r_{f^{a_{j_\delta }}_{L;\TAN}}=a\centerdot N$ of 
$ f^{a_{j_\delta }}_{L;\TAN}$ \resp{$ f^{a_{j_\delta }}_{R;\TAN}$} with $b<a$~, $b\in\nit$ (see proposition 3.2.8).
\vskip 11pt

These strata $ f^{str}_{L;\TAN}$ \resp{$ f^{str}_{R;\TAN}$} belong to the basin of the singular strange attractor
$\Lambda ^{\TAN}_{str_L}$ \resp{$\Lambda ^{\TAN}_{str_R}$} in such a way that:
\Bean
\item $ f^{a_{j_\delta }\mid\omega ^i_{j_\delta }}_{L;\TAN}\in \Lambda ^{\TAN}_{L}$ \resp{$ f^{a_{j_\delta }\mid\omega ^i_{j_\delta }}_{R;\TAN}\in \Lambda  ^{\TAN}_{R}$};

\item $ f^{\omega ^i_{j_\delta }}_{L;\TAN}\in \Lambda ^{\TAN}_{unf_L}$ \resp{$ f^{\omega ^i_{j_\delta }}_{R;\TAN}\in \Lambda ^{\TAN}_{unf_R}$}.
\Ee

These strata $ f^{str}_{L;\TAN}$ \resp{$ f^{str}_{R;\TAN}$} are characterized by a spherical geometry since their points are expanding under the automorphism $u_L$ \resp{$u_R$}.  Consequently, these strata $ f^{str}_{L;\TAN}$ \resp{$ f^{str}_{R;\TAN}$} constitue unstable subsets $E^u(f^{str}_{L;\TAN})$ \resp{$E^u( f^{str}_{R;\TAN})$} of $\Gamma (T(G^{(n)}_L( {F^+_{v}})))$ \resp{$\Gamma (T(G^{(n)}_R( {F^+_{\o v}})))$}.\epr
\vskip 11pt

\subsubsection{Proposition}

{\em The simplest singular strange attractor $\Lambda ^{\TAN}_{str_L}$ \resp{$\Lambda ^{\TAN}_{str_R}$} is composed of:
\Bena
\item a(n) (unfolded) singularity $y_L=\omega ^{s+2}_L$ \resp{$y_R=\omega ^{s+2}_R$};

\item a stable subset $E^s(f^{(a-b)_{j_\delta }}_{L;\TAN})$ \resp{$E^s(f^{(a-b)_{j_\delta }}_{R;\TAN})$} characterized by a hyperbolic geometry;

\item unstable subsets $E^u(f^{str}_{L;\TAN})$ \resp{$E^u(f^{str}_{R;\TAN})$} characterized by a spherical geometry.
\Ee}
\vskip 11pt

\bpr
\Bena
\item Let the unstable subsets $E^u(f^{str}_{L;\TAN})$ \resp{$E^u(f^{str}_{R;\TAN})$}  be given by the functions $f^{str}_{L;\TAN}$ \resp{$f^{str}_{R;\TAN}$} according to proposition 3.2.10.

Then, the stable subset $E^s(f^{(a-b)_{j_\delta }}_{L;\TAN})$ \resp{$E^s(f^{(a-b)_{j_\delta }}_{R;\TAN})$} characterized by a hyperbolic geometry, is:
\begin{align*}
E^s(f^{(a-b)_{j_\delta }}_{L;\TAN})
&= E^s(f^{a_{j_\delta }}_{L;\TAN})- \bigcup^s_{i=1} f^{a_{j_\delta }\mid\omega ^i_{j_\delta }}_{L;\TAN}\\[11pt]
\rresp{E^s(f^{(a-b)_{j_\delta }}_{R;\TAN})
&= E^s(f^{a_{j_\delta }}_{R;\TAN})- \bigcup^s_{i=1} f^{a_{j_\delta }\mid\omega ^i_{j_\delta }}_{R;\TAN}}\end{align*}
where:
\Bi
\item $E^s(f^{a_{j_\delta }}_{L;\TAN})$ \resp{$E^s(f^{a_{j_\delta }}_{R;\TAN})$} is the basin of attraction of the singular hyperbolic attractor $\Lambda ^{\TAN}_{L}$ \resp{$\Lambda ^{\TAN}_{R}$} given by the orbits $f^{a_{j_\delta }}_{L;\TAN}$ \resp{$f^{a_{j_\delta }}_{R;\TAN}$} having a rank 
$r_{f^{a_{j_\delta }}_{L;\TAN}}=a\centerdot N$ according to proposition 3.2.8.

\item $f^{a_{j_\delta }\mid\omega ^i_{j_\delta }}_{L;\TAN}$
\resp{$f^{a_{j_\delta }\mid\omega ^i_{j_\delta }}_{R;\TAN}$} is a subfunction of $f^{a_{j_\delta }}_{L;\TAN}$ \resp{$f^{a_{j_\delta }}_{R;\TAN}$} on which the generator  $\omega ^i_{j_{\delta _L}}$ \resp{$\omega ^i_{j_{\delta _R}}$} has been projected according to proposition 3.2.10.
\Ei

\item The concept of singular attractor was introduced in \cite{Pie3} and the notion of singular strange attractor was developed in \cite{Pie3} and independently in \cite{P-R}.

The fact that $\Lambda ^{\TAN}_{str_L}$ \resp{$\Lambda ^{\TAN}_{str_R}$} is a strange attractor corresponds to its description by H. Schuster \cite{Sch}:

\begin{quote}
``A strange attractor arises typically when the flow contracts the volume element in some directions, but stretches it along others.  To remain confined to a bounded domain, the volume element is folded at the same time''.
\end{quote}

An excellent literature on strange attractors can also be found in \cite{Rue1}, \cite{Rue2}; \cite{Mil2}, \cite{M-P}, \cite{E-R} and \cite{Wil}.\epr
\Ee
\vskip 11pt

\subsubsection{Structure of a general singular strange attractor}

According to proposition 3.2.10, a singular strange attractor is defined by:
\[\Lambda ^{\TAN}_{str_L}= \Lambda ^{\TAN}_L\times \Lambda ^{\TAN}_{unf_L}\qquad
\rresp{\Lambda ^{\TAN}_{str_R}= \Lambda ^{\TAN}_R\times \Lambda ^{\TAN}_{unf_R}}\]
where the unfolded attractor $\Lambda ^{\TAN}_{unf_L}$ \resp{$\Lambda ^{\TAN}_{unf_R}$} has as orbits the generators 
$\omega ^i_{j_{\delta _L}}$ \resp{$\omega ^i_{j_{\delta _R}}$}, $1\le i\le s$~, of the versal deformation of the singularity 
$y_L=\omega ^{s+2}_L$ \resp{$y_R=\omega ^{s+2}_R$} 
localized on the singular hyperbolic attractor $\Lambda ^{\TAN}_L$ \resp{$\Lambda ^{\TAN}_R$}.

But, these generators $\omega ^i_{j_{\delta _L}}$ \resp{$\omega ^i_{j_{\delta _R}}$} can carry singularities for $i\ge 2$ according to proposition 2.2.8 and section 3.1.21. As a result, a generator
$\omega ^i_{j_{\delta _L}}$ \resp{$\omega ^i_{j_{\delta _R}}$} carrying a singularity is a singular hyperbolic attractor, noted $\Lambda ^{\TAN}_{\omega ^i_{j_{\delta _L}}}$ \resp{$\Lambda ^{\TAN}_{\omega ^i_{j_{\delta _R}}}$}, according to proposition 3.2.8.

Consequently, the unfolded attractor can be rewritten as follows:
\[ \Lambda ^{\TAN}_{unf_L}=\bigcup^s_{i=1} \Lambda ^{\TAN}_{\omega ^i_{j_{\delta _L}}} \qquad
\rresp{\Lambda ^{\TAN}_{unf_R}=\bigcup^s_{i=1} \Lambda ^{\TAN}_{\omega ^i_{j_{\delta _R}}}}\]
where:
\Bi
\item $ \Lambda ^{\TAN}_{\omega ^i_{j_{\delta _L}}}$ \resp{$\Lambda ^{\TAN}_{\omega ^i_{j_{\delta _R}}}$} is a singular hyperbolic attractor for $i\ge 2$~.
\item $ \Lambda ^{\TAN}_{\omega ^1_{j_{\delta _L}}}$ \resp{$\Lambda ^{\TAN}_{\omega ^1_{j_{\delta _R}}}$} is a divisor.
\Ei

And, a general singular strange attractor can be decomposed according to:
\[ \Lambda ^{\TAN}_{str_L}= \Lambda ^{\TAN}_L\times \L( \bigcup^s_{i=1} \Lambda ^{\TAN}_{\omega ^i_{j_{\delta _L}}}\R) \qquad
\rresp{\Lambda ^{\TAN}_{str_R}= \Lambda ^{\TAN}_R\times \L( \bigcup^s_{i=1} \Lambda ^{\TAN}_{\omega ^i_{j_{\delta _R}}}\R)}.\]
\vskip 11pt

\subsubsection{Proposition}

{\em Let $SOT^{\max}_{j_{\delta _L}}=(\tau _{V_{\omega _{j_{\delta _L}}}}\circ \pi^{\max}_{s_{j_{\delta _L}}})$ \resp{$SOT^{\max}_{j_{\delta _R}}=(\tau _{V_{\omega _{j_{\delta _R}}}}\circ \pi^{\max}_{s_{j_{\delta _R}}})$} denote the maximal spreading-out isomorphism pulling out completely the generators $\omega ^i_{j_{\delta _L}}$ \resp{$\omega ^i_{j_{\delta _R}}$} of the versal deformation of the singularity $y_L=\omega _L^{s+2}$ \resp{$y_R=\omega _R^{s+2}$} from its neighbourhood $D_{\phi ^{\TAN}_{G_{j_{\delta _L}}}}$ \resp{$D_{\phi ^{\TAN}_{G_{j_{\delta _R}}}}$}.
\vskip 11pt

Then, the maximal spreading-out isomorphism
\begin{align*}
 SOT^{\max}_{j_{\delta _L}}: \quad \Lambda ^{\TAN}_{str_L} &\To \Lambda ^{\TAN}_L\bigoplus 
\bigcup^s_{i=1} \Lambda ^{\TAN}_{\omega ^i_{j_{\delta _L}}} \\[11pt]
\rresp{SOT^{\max}_{j_{\delta _R}}: \quad \Lambda ^{\TAN}_{str_R} &\To \Lambda ^{\TAN}_R\bigoplus 
\bigcup^s_{i=1}\Lambda ^{\TAN}_{\omega ^i_{j_{\delta _R}}}}\end{align*}
decomposes the general singular strange attractor $\Lambda ^{\TAN}_{str_L}$ \resp{$\Lambda ^{\TAN}_{str_R}$} into the original singular hyperbolic attractor $\Lambda ^{\TAN}_L$ \resp{$\Lambda ^{\TAN}_R$} and into a set of $(s-1)$ disconnected singular hyperbolic attractors and  a disconnected divisor.}
\vskip 11pt

\bpr The $SOT^{\max}_{j_{\delta _L}}$ \resp{$SOT^{\max}_{j_{\delta _R}}$} spreading-out isomorphism corresponds to an extension of the quotient algebra of the versal deformation of the singular semisheaf
$\theta ^*_{G_L^{(n)}}$ \resp{$\theta ^*_{G_R^{(n)}}$} developed in corollary 3.1.17 and restricted to the $(j_\delta ,m_{j_\delta })$ conjugacy class representative of $G^{(n)}_L( {F^+_{v}})$ \resp{$ G^{(n)}_R( {F^+_{\o v}})$}.

As the unfolded attractor $\Lambda ^{\TAN}_{unf_L}$ \resp{$\Lambda ^{\TAN}_{unf_R}$} is the union of hyperbolic attractors $\Lambda ^{\TAN}_{\omega ^i_{j_{\delta _L}}}$ \resp{$\Lambda ^{\TAN}_{\omega ^i_{j_{\delta _R}}}$} and of a divisor, which are the generators $\omega ^i_{j_{\delta _L}}$
\resp{$\omega ^i_{j_{\delta _R}}$} of the versal deformation of the singularity $y_L=\omega _L^{s+2}$ \resp{$y_R=\omega _R^{s+2}$}, the maximal spreading-out isomorphism $SOT^{\max}_{j_{\delta _L}}$ \resp{$SOT^{\max}_{j_{\delta _R}}$} is the inverse map
\[ SOT^{\max}_{j_{\delta _L}}=(VD_{\Lambda _L})^{-1} \qquad
\rresp{SOT^{\max}_{j_{\delta _R}}=(VD_{\Lambda _R})^{-1}}\]
of the versal unfolding
\[ VD_{\Lambda _L}: \quad \Lambda ^{\TAN}_L\To \Lambda ^{\TAN}_{str_L}\qquad 
\rresp{VD_{\Lambda _R}: \quad \Lambda ^{\TAN}_R\To \Lambda ^{\TAN}_{str_R}}\]
introduced in proposition 3.2.10.

Consequently, $SOT^{\max}_{j_{\delta _L}}$
\resp{$SOT^{\max}_{j_{\delta _R}}$} blows up the general singular strange attractor
$\Lambda ^{\TAN}_{str_L}$ \resp{$\Lambda ^{\TAN}_{str_R}$} into the original singular hyperbolic attractor $\Lambda ^{\TAN}_L$ \resp{$\Lambda ^{\TAN}_R$}, $(s-1)$ disconnected singular hyperbolic attractors $\Lambda ^{\TAN}_{\omega ^i_{j_{\delta _L}}}$ \resp{$\Lambda ^{\TAN}_{\omega ^i_{j_{\delta _R}}}$}, $i\ge 2$~, and a divisor $\Lambda ^{\TAN}_{\omega ^1_{j_{\delta _L}}}$ \resp{$\Lambda ^{\TAN}_{\omega ^1_{j_{\delta _R}}}$} corresponding to the generators of the versal deformation.\epr
\vskip 11pt

\subsubsection{Corollary}

{\em Let $VD_{\Lambda _L}$ \resp{$VD_{\Lambda _R}$} denote the versal deformation transforming a hyperbolic singular attractor $\Lambda ^{\TAN}_L$ \resp{$\Lambda ^{\TAN}_R$} into a general singular strange attractor $\Lambda ^{\TAN}_{str_L}$ \resp{$\Lambda ^{\TAN}_{str_R}$}.

Let $SOT^{\max}_{j_{\delta _L}}$ \resp{$SOT^{\max}_{j_{\delta _R}}$} be the maximal spreading-out isomorphism blowing up the general singular strange attractor.

Then, the following composition of maps:
\begin{alignat*}{3}
SOT^{\max}_{j_{\delta _L}} \circ VD_{\Lambda _L} : \quad \Lambda ^{\TAN}_L &\xrightarrow{VD_{\Lambda _L}} \Lambda ^{\TAN}_{str_L} &\xrightarrow{SOT^{\max}_{j_{\delta _L}}} & \Lambda ^{\TAN}_L\bigoplus \bigcup^s_{i=1} \Lambda ^{TAN}_{\omega ^i_{j_{\delta _L}}} \\
\rresp{SOT^{\max}_{j_{\delta _R}} \circ VD_{\Lambda _R} : \quad \Lambda ^{\TAN}_R &\xrightarrow{VD_{\Lambda _R}} \Lambda ^{\TAN}_{str_R} &\xrightarrow{SOT^{\max}_{j_{\delta _R}}} & \Lambda ^{\TAN}_R\bigoplus \bigcup^s_{i=1} \Lambda ^{TAN}_{\omega ^i_{j_{\delta _R}}}}
\end{alignat*}
is such that
\[ SOT^{\max}_{j_{\delta _L}}=(VD_{\Lambda _L})^{-1} \qquad
\rresp{SOT^{\max}_{j_{\delta _R}}=(VD_{\Lambda _R})^{-1}}.\]}

\bpr This is a consequence of proposition 3.2.13 and of the generation of the versal deformation given in proposition 2.2.8.

Note that the paper of \cite{B-L-M-P} introduces the explosion of singular cycles, phenomenon closed to the blowing up of strange attractors considered here and in \cite{Pie3}.\epr
\vskip 11pt

\section{Langlands global correspondences affected by degenerate singularities}

The aim of this chapter consists in showing in what extent it is possible to develop global correspondences of Langlands for a (bisemi)sheaf of rings on the real bilinear algebraic semigroup affected by degenerate singularities in the sense of chapters 2 and 3.
\vskip 11pt

\subsection[The transformation of the bisemisheaf of rings $\theta _{G^{(n)}_R}\times \theta _{G^{(n)}_L}$ on the real bilinear algebraic semigroup $\Gnv$ under degenerate  singularities]{\bbf The transformation of the bisemisheaf of rings $\theta _{G^{(n)}_R}\times \theta _{G^{(n)}_L}$ on the real bilinear algebraic semigroup $\Gnv$ under degenerate  singularities}

First, it will be recalled what is the $n$-dimensional real irreducible global correspondence of Langlands as developed in \cite{Pie1}: it consists in a bijection between the $n$-dimensional irreducible representation $\Irr\Rep^{(n)}_{W_{F^+\RL}}(W^{ab}_{F^+_R}\times W^{ab}_{F^+_L})$ of the product, right by left, of Weil groups and the irreducible cuspidal representation of $\Gnv$ given by\linebreak $\Irr\ELLIP(\GL_n( \Aa_{F^{+,T}_{\o v}} \times  \Aa_{F^{+,T}_v}))$~.  These concepts will thus be reviewed, and, among others, the definition of global Weil groups.
\vskip 11pt

\subsubsection{Global Weil groups}

Let $\Gal(\wF^+_{v_{j_\delta }}\big/ F^0)$ \resp{$\Gal(\wF^+_{\o v_{j_\delta }}\big/ F^0)$} denote the Galois subgroup of the extension $\wF^+_{v_{j_\delta }}$ \resp{$\wF^+_{\o v_{j_\delta }}$} in one-to-one correspondence with the  peudo-ramified completion $F^+_{v_{j_\delta }}$ \resp{$F^+_{\o v_{j_\delta }}$} having a rank given by:
\begin{align*}
[F^+_{v_{j_\delta }}:F^0] &=*+j_\delta \centerdot N\;, \qquad *\in\nit\;, \quad *<N\\
\rresp{[F^+_{\o v_{j_\delta }}:F^0] &=*+j_\delta \centerdot N}.
\end{align*}

And, let $\Gal(\dot{\wF}^+_{v_{j_\delta }}\big/ F^0)$ \resp{$\Gal(\dot{\wF}^+_{\o v_{j_\delta }}\big/ F^0)$} be the Galois subgroup of the peudo-ramified extension $\dot{\wF}^+_{v_{j_\delta }}$ \resp{$\dot{\wF}^+_{\o v_{j_\delta }}$} characterized by a degree:
\[ [\dot{\wF}^+_{v_{j_\delta }}:F^0] = j_\delta \centerdot N \qquad 
\rresp{[\dot{\wF}^+_{\o v_{j_\delta }}:F^0] = j_\delta \centerdot N}\]
(see section 1.1).

Then, according to \cite{Pie1}, the global Weil group $W^{ab}_{F^+_L}$ \resp{$W^{ab}_{F^+_R}$} is given by:
\[ W^{ab}_{F^+_L} = \bigoplus_{j_\delta } \ \bigoplus_{m_{j_\delta }}\Gal ( \dot{\wF}^+_{v_{j_\delta,m_{j_\delta } }}\big/F^0) \qquad
\rresp{W^{ab}_{F^+_R} = \bigoplus_{j_\delta } \ \bigoplus_{m_{j_\delta }}\Gal ( \dot{\wF}^+_{\o v_{j_\delta ,m_{j_\delta }}}\big/F^0)}\]
and, its product, right by left, is:
\[W^{ab}_{F^+_R}\times W^{ab}_{F^+_L}=\bigoplus_{j_\delta } \ \bigoplus_{m_{j_\delta }}
\Gal ( \dot{\wF}^+_{\o v_{j_\delta ,m_{j_\delta }}}\big/F^0)\times
\Gal ( \dot{\wF}^+_{v_{j_\delta ,m_{j_\delta }}}\big/F^0)\;.\]

The $n$-dimensional irreducible representation of the product, right by left, $W^{ab}_{F^+_R}\times W^{ab}_{F^+_L}$ of global Weil groups is given by:
\[\Irr\Rep^{(n)}_{W_{F^+\RL}}(W^{ab}_{F^+_R}\times W^{ab}_{F^+_L}) = 
G^{(n)}({F^+_{\o v_\oplus}}\times {F^+_{v_\oplus}})
\]
according to \cite{Pie1} (proposition 3.4.3).
\vskip 11pt

\subsubsection[Cuspidal representation of $\GLnV$]{\bbf Cuspidal representation of $\GLnV$}

Let $G^{(n)}( {F^{+,T}_{\o v}} \times  {F^{+,T}_{v}})$ be the reductive bilinear algebraic semigroup $\Gnv$ submitted to the toroidal compactification $\gamma ^c\RL$ as developed in section 3.4 of \cite{Pie1}.

Its conjugacy class representatives $g^{(n)}_{T\RL}[j_\delta ,m_{j_\delta }]=
g^{(n)}_{T_R}[j_\delta ,m_{j_\delta }] \times g^{(n)}_{T_L}[j_\delta ,m_{j_\delta }]$ are products, right by left, of $n$-dimensional real semitori $T^n_R[j_\delta ,m_{j_\delta }] \times T^n_L[j_\delta ,m_{j_\delta }]$~.

Each complex-valued bifunction $\phi ^{(n)}_{G^T_R}(x_{g^T_{j_{\delta _R}}}) \otimes
\phi ^{(n)}_{G^T_L}(x_{g^T_{j_{\delta _L}}})$ on the conjugacy class representative
$g^{(n)}_{T\RL}[j_\delta ,m_{j_\delta }]\in G^{(n)}( {F^{+,T}_{\o v}} \times  {F^{+,T}_{v}})$ is given by:
\begin{align*}
\phi ^{(n)}_{G^T_R}(x_{g^T_{j_{\delta _R}}}) \otimes
\phi ^{(n)}_{G^T_L}(x_{g^T_{j_{\delta _L}}})
&= T^n_R[j_\delta ,m_{j_\delta }] \times T^n_L[j_\delta ,m_{j_\delta }]\\
&= \lambda ^ {\half}(n,j_\delta ,m_{j_\delta })\ e^{-2\pi ij_\delta x} \otimes
\lambda ^ {\half}(n,j_\delta ,m_{j_\delta })\ e^{2\pi ij_\delta x}
\end{align*}
where:
\Bi
\item $\ds\lambda (n,j_\delta ,m_{j_\delta })=\prod^n_{c=1}\lambda _c(n,j_\delta ,m_{j_\delta })$ is a product of eigenvalues $\lambda _c(n,j_\delta ,m_{j_\delta })$ of the $j_\delta $-th coset representative $(U_{j_{\delta _R}}\times U_{j_{\delta _L}})$ of  the product $(T_R(n;t)\otimes T_L(n;t))$ of Hecke operators;

\item $\vec x=\ds\sum^n_{c=1}x_c\ \vec e_c$ is a vector of $(F^+_{v_{j_\delta }})^n$ and, more precisely, a point of $g^{(n)}_{L}[j_\delta ,m_{j_\delta }]$~.
\Ei

The sum of all bifunctions on the conjugacy class representatives of 
$G^{(n)}( {F^{+,T}_{\o v}} \times  {F^{+,T}_{v}})$ is:
\begin{align*}
\phi ^{(n)}_{G^T_R}(x_{g^T_{R}}) \otimes
\phi ^{(n)}_{G^T_L}(x_{g^T_{L}})
&= \bigoplus^r_{j_\delta =1} \ \bigoplus_{m_{j_\delta }}
(\phi ^{(n)}_{G^T_R}(x_{g^T_{j_{\delta _R}}}) \otimes
\phi ^{(n)}_{G^T_L}(x_{g^T_{j_{\delta _L}}}))\\[11pt]
&= \bigoplus_{j_\delta } \ \bigoplus_{m_{j_\delta }} (\lambda ^ {\half}(n,j_\delta ,m_{j_\delta })\ e^{-2\pi ij_\delta x} \otimes
\lambda ^ {\half}(n,j_\delta ,m_{j_\delta })\ e^{2\pi ij_\delta x})\\[11pt]
&=\ELLIP_R(n,j_\delta ,m_{j_\delta }) \otimes \ELLIP_L(n,j_\delta ,m_{j_\delta })
\end{align*}
where
\[\ELLIP\RL(n,j_\delta ,m_{j_\delta })=
\ELLIP_R(n,j_\delta ,m_{j_\delta }) \otimes \ELLIP_L(n,j_\delta ,m_{j_\delta })\]
constitutes the $n$-dimensional irreducible elliptic (and cuspidal) representation\linebreak
$\Irr\ELLIP(\GL_n(\Aa_ {F^{+,T}_{\o v}} \times \Aa_ {F^{+,T}_{v})})$ of
$\GL_n( {F^{+,T}_{\o v}} \times  {F^{+,T}_{v}})$~, where $\Aa_ {F^{+,T}_{\o v}} $ and $ \Aa_ {F^{+,T}_{v})}$ are toroidal adele semirings according to section 1.1.
\vskip 11pt

\subsubsection[The $n$-dimensional real irreducible global correspondence of Langlands ]{\bbf The $n$-dimensional real irreducible global correspondence of Langlands}

The global correspondence considered here is thus given by:
\[ \begin{array}{ccc}
\Irr\Rep^{(n)}_{W_{F^+\RL}}(W^{ab}_{F^+_R} \times W^{ab}_{F^+_L}) & \To &
\Irr\ELLIP\RL(\GL_n( \Aa_{F^{+,T}_{\o v}} \times \Aa_ {F^{+,T}_{v})})\\
\| & &\| \\
G^{(n)}( {F^{+,T}_{\o v_\oplus}} \times  {F^{+,T}_{v_\oplus})} & \To &
\ELLIP\RL(n,j_\delta ,m_{j_\delta })\end{array}\]
where:
\Bi
\item $\Irr\Rep^{(n)}_{W_{F^+\RL}}(W^{ab}_{F^+_R} \times W^{ab}_{F^+_L})$ is the sum of the products, right by left, of the equivalence classes of the irreducible $n$-dimensional Weil-Deligne representation of the bilinear global Weil group $(W^{ab}_{F^+_R} \times W^{ab}_{F^+_L})$ given by the algebraic bilinear real semigroup $G^{(n)}( {F^{+}_{\o v_\oplus}} \times  {F^{+}_{v_\oplus})})$~;
\item $\Irr\ELLIP\RL(\GL_n( \Aa_{F^{+,T}_{\o v}} \times  \Aa_{F^{+,T}_{v}})$ is the sum of the products, right by left, of the equivalence classes of the irreducible elliptic representation of
$\GL_n( {F^{+,T}_{\o v}} \times  {F^{+,T}_{v})})$ given by the $n$-dimensional global elliptic bisemimodule $\ELLIP\RL(n,j_\delta ,m_{j_\delta })$~.
\Ei
\vskip 11pt

\subsubsection[The singularization of the bisemisheaf on $\GL_n( {F^{+}_{\o v}} \times  {F^{+}_{v})}$]{\bbf The singularization of the bisemisheaf on $G^{(n)}( {F^{+}_{\o v}} \times  {F^{+}_{v})}$}

As we are concerned with the problem of (degenerate) singularities in the global program of Langlands, we shall take into account the set of real-valued differentiable bifunctions
$\{\phi ^{(n)}_{G_{j_R}}(x_{g_{j_{\delta _R}}}) \otimes
\phi ^{(n)}_{G_{j_L}}(x_{g_{j_{\delta _L}}})\}$ on the corresponding conjugacy class representatives
$\{ g ^{(n)}\RL[j_\delta ,m_{j_\delta }]\}$ of the bilinear real algebraic semigroup $\Gnv$~, as mentioned in section 1.6.  And, more precisely, we shall work with the bisemisheaf of rings
$\theta _{G^{(n)}\RL}= \theta _{G^{(n)}_R}\otimes \theta _{G^{(n)}_L}$~, whose (bi)sections are the differentiable bifunctions $\phi ^{(n)}_{G_{j_R}}(x_{g_{j_{\delta _R}}}) \otimes
\phi ^{(n)}_{G_{j_L}}(x_{g_{j_{\delta _L}}})$  (see section 1.7).

Let $ \theta _{G^{(n)}_L}$ \resp{$\theta _{G^{(n)}_R}$} be the semisheaf of \lr differentiable functions $\phi ^{(n)}_{G_{j_L}}(x_{g_{j_{L}}})$ \resp{$\phi ^{(n)}_{G_{j_R}}(x_{g_{j_{R}}})$}, rewritten in a condensed form according to $\phi _{j_\delta }(x_L)$ \resp{$\phi _{j_\delta }(x_R)$}.

The singularization of $\theta _{G^{(n)}_L}$ \resp{$\theta _{G^{(n)}_R}$} is given by the contracting surjective morphism (see section 2.1):
\[ \o\rho _{G_L}: \quad \theta _{G^{(n)}_L} \To \theta ^*_{G^{(n)}_L} \qquad 
\rresp{\o\rho _{G_R}: \quad \theta _{G^{(n)}_R} \To \theta ^*_{G^{(n)}_R}}\]
where $ \theta ^*_{G^{(n)}_L}$ \resp{$\theta ^*_{G^{(n)}_R}$} denotes the \lr singular semisheaf whose sections $\phi^* _{j_\delta }(x_L)$ \resp{$\phi^* _{j_\delta }(x_R)$} are endowed with germs 
$\phi _{j_\delta }(\omega _L)$ \resp{$\phi _{j_\delta }(\omega _R)$} having degenerate  singularities of corank 1~.

The corresponding singular bisemisheaf is $( \theta ^*_{G^{(n)}_R}\otimes\theta^* _{G^{(n)}_L})$ degenerate  from $( \theta _{G^{(n)}_R}\otimes\theta _{G^{(n)}_L})$ under the contracting surjective morphism:
\[ \o\rho _{G_R} \times \o\rho _{G_L}: \qquad \theta _{G^{(n)}_R}\times \theta _{G^{(n)}_L} \To \theta ^*_{G^{(n)}_R}\times \theta ^*_{G^{(n)}_L}\;.\]
\vskip 11pt

\subsubsection[The versal deformation of the singular bisemisheaf $(\theta ^*_{G^{(n)}_R}\times \theta ^*_{G^{(n)}_L})$]{\bbf The versal deformation of the singular bisemisheaf $(\theta ^*_{G^{(n)}_R}\times \theta ^*_{G^{(n)}_L})$}

Let 
\[ D_{S_L} : \quad \theta ^*_{G^{(n)}_L}\times \theta _{S_L} \To \theta ^*_{G^{(n)}_L} \qquad
\rresp{D_{S_R} : \quad \theta ^*_{G^{(n)}_R}\times \theta _{S_R} \To \theta ^*_{G^{(n)}_R}}\]
be the versal deformation of the singular semisheaf $\theta ^*_{G^{(n)}_L} $ \resp{$\theta ^*_{G^{(n)}_R} $}.  It is a contracting fibre bundle whose fibre $\theta _{S_L} = \{\theta ^1(\omega ^1_L),\cdots,\theta ^1(\omega ^i_L),\cdots,\theta ^1(\omega ^s_L)\}$
\resp{$\theta _{S_R} = \{\theta ^1(\omega ^1_R),\cdots,\theta ^1(\omega ^i_R),\cdots,\theta ^1(\omega ^s_R)\}$} is the family of (semi)sheaves of the \lr base $S_L$ \resp{$S_R$} of the versal deformation.

As developed in section  2.1, the versal deformation of a degenerate  singularity of corank 1 and codimension $s$ on a section $\phi^* _{j_\delta }(x_L)$ \resp{$\phi^* _{j_\delta }(x_R)$} of 
$ \theta ^*_{G^{(n)}_L}$ \resp{$\theta ^*_{G^{(n)}_R}$} is yielded by a sequence of $s$ contracting morphisms extending the corresponding sequence of contracting surjective morphisms of singularization.
\vskip 11pt

\subsubsection{The spreading-out isomorphism}

It consists in a blow-up of the versal deformation $ \theta ^*_{G^{(n)}_L}\times \theta _{S_L}$ \resp{$ \theta ^*_{G^{(n)}_R}\times \theta _{S_R}$} of the singular semisheaf 
$ \theta ^*_{G^{(n)}_L}$ \resp{$ \theta ^*_{G^{(n)}_R}$}.

This blow-up is maximal if the spreading-out isomorphism $SOT^{\max}_L=(\tau _{V_{\omega _L}}\circ \pi ^{\max}_{s_L})$ \resp{$SOT^{\max}_R=(\tau _{V_{\omega _R}}\circ \pi ^{\max}_{s_R})$} is the inverse of the versal deformation:
\[ SOT^{\max}_L=(D_{S_L})^{-1} \qquad \resp{SOT^{\max}_R=(D_{S_R})^{-1}}.\]
It is then given by the map:
\begin{align*}
SOT^{\max}_L: \quad \theta ^*_{G^{(n)}_L}\times \theta _{S_L} &\To \theta ^*_{G^{(n)}_L}\cup \theta _{S_L}\\
\rresp{SOT^{\max}_R: \quad \theta ^*_{G^{(n)}_R}\times \theta _{S_R} &\To \theta ^*_{G^{(n)}_R}\cup \theta _{S_R}}\end{align*}
projecting the family of sheaves of the \lr base $S_L$ \resp{$S_R$} of the versal deformation in the vertical tangent space $T_{V_{\omega _L}}$ \resp{$T_{V_{\omega _R}}$} (see proposition 3.1.16 and corollary 3.1.17).

Let $\theta ^*_{SOT(1)_L}$ \resp{$\theta ^*_{SOT(1)_R}$} be the family $\theta _{S_L}$
\resp{$\theta _{S_R}$} of disconnected base semisheaves having been glued together according to section 3.1.20:

The semisheaf $\theta ^*_{SOT(1)_L} \simeq \theta _{S_L}$ \resp{$\theta ^*_{SOT(1)_R} \simeq \theta _{S_R}$} covers partially the singular semisheaf $\theta ^*_{G^{(n)}_L}$ \resp{$\theta ^*_{G^{(n)}_R}$} but can be affected by singularities in its sections\ldots\ which can also undergo a versal deformation.

For the simplicity, we shall consider that we are only confronted with the singular semisheaf $\theta ^*_{G^{(n)}_L}$ \resp{$\theta ^*_{G^{(n)}_R}$} and its first cover 
$\theta ^*_{SOT(1)_L} $ \resp{$\theta ^*_{SOT(1)_R} $}, having possible singularities.
\vskip 11pt

\subsubsection{Contracting morphisms of singularization}

The bisemisheaf $\theta _{G^{(n)}_R}\times\theta _{G^{(n)}_L}$~, being defined on the bilinear algebraic semigroup $\Gnv$~, constitutes an $n$-dimensional irreducible real representation 
$\Irr\Rep^{(n)}_{W_{F^+\RL}}(W^{ab}_{F^+_R} \times W^{ab}_{F^+_L})$ of the bilinear global Weil group.

The fact of considering contracting morphisms of singularization leads to the transformation of
$\theta _{G^{(n)}_R}\times\theta _{G^{(n)}_L}$ into:
\begin{multline*}
\theta _{G^{(n)}_R}\times\theta _{G^{(n)}_L} 
\xrightarrow{\o\rho _{G_R}\times\o\rho_{G_L}} 
\theta ^*_{G^{(n)}_R}\times\theta ^*_{G^{(n)}_L} \\
\xrightarrow{D_{S_R}\times D_{S_L}} 
 (\theta ^*_{G^{(n)}_R}\times\theta _{S_R}) \otimes(\theta ^*_{G^{(n)}_L}\times\theta _{S_l})\\
 \xrightarrow{SOT^{\max}_R \times SOT^{\max}_L} 
 (\theta ^*_{G^{(n)}_R}\cup\theta^*_{SOT(1)_R}) \otimes (\theta ^*_{G^{(n)}_L}\cup\theta^*_{SOT(1)_L})
 \end{multline*}
 where:
 \Bi
 \item $\o\rho _{G_R}\times\o\rho_{G_L}$ is the contracting morphism of singularization;
 \item $D_{S_R}\times D_{S_L}$ is the contracting morphism of versal deformation;
 \item $SOT^{\max}_R \times SOT^{\max}_L$ is the (contracting) blow-up of the versal deformation.
 \Ei
 
But, the bisemisheaves $(\theta ^*_{G^{(n)}_R}\otimes\theta ^*_{G^{(n)}_L})$ and 
$(\theta^*_{SOT(1)_R} \otimes \theta^*_{SOT(1)_L}) $~, being affected by singularities, cannot be endowed with a cuspidal representation.

To reach this objective, it is necessary to:
\Bena
\item desingularize these bisemisheaves;
\item submit them to a toroidal compactification.
\Ee

This is the aim of the next section.
\vskip 11pt

\subsection{Langlands global correspondences despite of degenerate  singularities}

\subsubsection[Desingularizing the bisemisheaf $(\theta ^*_{G^{(n)}_R}\otimes\theta ^*_{G^{(n)}_L})$]{\bbf Desingularizing the bisemisheaf $(\theta ^*_{G^{(n)}_R}\otimes\theta ^*_{G^{(n)}_L})$}

The desingularization of the semisheaf $ \theta ^*_{G^{(n)}_L}$ \resp{$ \theta ^*_{G^{(n)}_R}$} corresponds to the classical monoidal transformations applied to the singularities on the sections 
$\phi^* _{j_\delta }(x_L)$ \resp{$\phi^* _{j_\delta }(x_R)$} of 
$ \theta ^*_{G^{(n)}_L}$ \resp{$\theta ^*_{G^{(n)}_R}$}.

A desingularization of $\phi^* _{j_\delta }(x_L)$ \resp{$\phi^* _{j_\delta }(x_R)$} is described succinctly in proposition 2.1.3 and corresponds to the inverse morphism of a singularization developed in section 2.1.

More concretely, if we want to desingularize a germ $\phi _{j_\delta }(\omega _L)$ \resp{$\phi _{j_\delta }(\omega _R)$} on $\phi^* _{j_\delta }(x_L)$ \resp{$\phi^* _{j_\delta }(x_R)$}, given by the degenerate  singularity $y_L=\omega _L^{s+2}$ \resp{$y_R=\omega _R^{s+2}$} of corank 1 and codimensions $s$~, we have to consider the following sequence of expanding morphisms:
\begin{alignat*}{3}
\o\rho\des_L &: \quad \omega ^{s+2}_L 
\ \xrightarrow{\o\rho_L^{(-1)(s+2)}}\
\omega ^{s+1}_L\cup D_L^{(s+1)}
&\ \xrightarrow{\o\rho_L^{(-1)(s+1)}} \ &
\omega ^{s}_L\cup D_L^{s}
\To \ldots \qquad \\[11pt]
&&&  \qquad \xrightarrow{\o\rho_L^{(-1)(1)}} 
\ \omega _L\cup D_L^{(1)}\\[15pt]
\rresp{\o\rho\des_R &: \quad \omega ^{s+2}_R
\ \xrightarrow{\o\rho_R^{(-1)(s+2)}}\
\omega ^{s+1}_R\cup D_R^{(s+1)}
 &\ \xrightarrow{\o\rho_R^{(-1)(s+1)}} \ &
\omega ^{s}_R\cup D_R^{s}
\To \ldots \qquad \\[11pt]
&&&\qquad \xrightarrow{\o\rho_R^{(-1)(1)}} \
\omega _R\cup D_R^{(1)}}\end{alignat*}
where:
\Bena
\item the expanding morphism of desingularization
\[ \o\rho_L^{(-1)(s+1)} : \quad \omega _L^{s+1} \To \omega _L^s\cup D s_L \qquad 
\rresp{\o\rho_R^{(-1)(s+1)} : \quad \omega _R^{s+1} \To \omega _R^s\cup D s_R }\]
is a projective morphism, disconnecting the divisor $D^s_L$ \resp{$D^s_R$} from the singular sublocus
\[ \Sigma ^{(s+1)}_L = \omega _L^{s+1}\qquad 
\rresp{\Sigma ^{(s+1)}_R = \omega _R^{s+1}}.\]

\item $\begin{array}[t]{rcl}
\o\rho\des_L : \quad \phi^* _{j_\delta }(x_L) & \To & \phi _{j_\delta }(x_L)\cup (D_L^{(s+1)}, D_L^{(s)}, \cdots, D_L^{(1)})\\[11pt]
\rresp{\o\rho\des_R : \quad \phi^* _{j_\delta }(x_R) & \To & \phi _{j_\delta }(x_R)\cup (D_R^{(s+1)}, D_R^{(s)}, \cdots, D_R^{(1)})}\end{array}$

is the desingularization of $y_L=\omega _L^{s+2}$ \resp{$y_R=\omega _R^{s+2}$} on 
$\phi^* _{j_\delta }(x_L)$ \resp{$\phi^* _{j_\delta }(x_R)$} disconnecting (by projection) the set of divisors $(D_L^{(s+1)}, D_L^{(s)}, \cdots, D_L^{(1)})$ \resp{$(D_R^{(s+1)},\linebreak D_R^{(s)}, \cdots, D_R^{(1)})$} which generate a real projective subscheme of dimension $(s-1)$~.

As a result, $\o\phi _{j_\delta }(x_L)$ \resp{$\o\phi _{j_\delta }(x_R)$} becomes a smooth section of $ \theta _{G^{(n)}_L}$ \resp{$\theta _{G^{(n)}_R}$}.

\item the expanding morphism of desingularization $\o\rho\des_L$ \resp{$\o\rho\des_R$} is an isomorphism outside the singular locus $\Sigma _L=\omega _L^{s+2}$ \resp{$\Sigma _R=\omega _R^{s+2}$} according to:
\begin{align*}
\o\rho_L^{(-1){\rm is}} : \quad \phi _{j_\delta }(x_L)\smallsetminus (\Sigma _L) &\To
\o\phi  _{j_\delta }(x_L) \smallsetminus (\o\rho_L^{(-1)}(\Sigma _L))\\[11pt]
\rresp{\o\rho_R^{(-1){\rm is}} : \quad \phi _{j_\delta }(x_R)\smallsetminus (\Sigma _R) &\To
\o\phi  _{j_\delta }(x_R) \smallsetminus (\o\rho_R^{(-1)}(\Sigma _R))}.
\end{align*}
\Ee

This desingularization process is exactly the inverse of the singularization developed in proposition 2.1.10.

The desingularization or resolution of singularities on all the sections 
$\phi^* _{j_\delta }(x_L)$ \resp{$\phi^* _{j_\delta }(x_R)$} of the singular semisheaf
$ \theta ^*_{G^{(n)}_L}$ \resp{$\theta ^*_{G^{(n)}_R}$} is given by the set of expanding morphisms:
\[ \o\rho\des_{G_L} : \quad   \theta ^*_{G^{(n)}_L}   \To   \theta  _{G^{(n)}_L}\qquad 
\rresp{\o\rho\des_{G_R} : \quad   \theta ^*_{G^{(n)}_R}   \To   \theta  _{G^{(n)}_R} }\]
in such a way that:
\[ \o\rho\des_{G_L}=\o\rho^{-1}_{G_L} \qquad \rresp{\o\rho\des_{G_R}=\o\rho^{-1}_{G_R}}\]
where $\o\rho_{G_L}$ \resp{$\o\rho_{G_R}$} denotes the set of contracting morphisms of singularization of $ \theta _{G^{(n)}_L}$ \resp{$\theta _{G^{(n)}_R}$}.

And, the resolution of singularities of the singular semisheaf $ \theta ^*_{G^{(n)}_R}\otimes \theta ^*_{G^{(n)}_L}$ is given by:
\[ \o\rho\des_{G_R} \times \o\rho\des_{G_L} : \quad \theta ^*_{G^{(n)}_R}\otimes\theta ^*_{G^{(n)}_L} \To \theta  _{G^{(n)}_R}\otimes\theta  _{G^{(n)}_L}\;.\]
\vskip 11pt

\subsubsection[Resolution of singularities of the covering bisemisheaf $\theta^*_{SOT(1)_R} \otimes \theta^*_{SOT(1)_L}$]{\bbf Resolution of singularities of the covering bisemisheaf $\theta^*_{SOT(1)_R} \otimes \theta^*_{SOT(1)_L}$}

As the sections of the semisheaf $ \theta ^*_{G^{(n)}_L}$ \resp{$\theta ^*_{G^{(n)}_R}$} are endowed with singularities of corank 1 and codimension $s$~, the semisheaf 
$\theta^*_{SOT(1)_L} $ \resp{$\theta^*_{SOT(1)_R}$} (being the family of ``~$s$~'' base semisheaves of the versal deformation having been glued together) can be affected on its sections by singularities of corank 1 and maximal codimensions equal to $(s-2)$ according to section 3.1 and proposition 2.2.8.

So, a resolution of singularities of this covering bisemisheaf $\theta^*_{SOT(1)_R} \otimes \theta^*_{SOT(1)_L}$ must be envisaged as it was done for the bisemisheaf $ \theta ^*_{G^{(n)}_R}\otimes \theta ^*_{G^{(n)}_L}$~.

The resolution of singularities of $\theta^*_{SOT(1)_R} \otimes \theta^*_{SOT(1)_L}$ will then be given by the morphism:
\[ \o\rho\des_{SOT(1)_R} \times \o\rho\des_{SOT(1)_L} : \quad \theta ^*_{SOT(1)_R}\otimes\theta ^*_{SOT(1)_L} \To \theta  _{SOT(1)_R}\otimes\theta  _{SOT(1)_L}\]
where $\theta  _{SOT(1)_R}\otimes\theta  _{SOT(1)_L}$ denotes the free corresponding bisemisheaf.
\vskip 11pt

\subsubsection[Global holomorphic representation of $ \theta _{G^{(n)}_R}\otimes \theta _{G^{(n)}_L}$]{\bbf Global holomorphic representation of $ \theta _{G^{(n)}_R}\otimes \theta _{G^{(n)}_L}$}

As the bisemisheaf $ \theta _{G^{(n)}_R}\otimes \theta _{G^{(n)}_L}$ was desingularized from its corresponding singular equivalent $ \theta ^*_{G^{(n)}_R}\otimes \theta ^*_{G^{(n)}_L}$~, a $n$-dimensional irreducible global holomorphic representation can be worked out for it, as it was developed in section 3.1 of \cite{Pie1}.  We shall recall it briefly.

The sections of the semisheaf $ \theta _{G^{(n)}_L}$ \resp{$\theta _{G^{(n)}_R}$} are the complex-valued differentiable functions (resp. cofunctions):
\begin{align*}
 f_{v_{j_\delta },m_{j_\delta }}(z^{j_\delta }) : \quad g_L^{(n)}[j_\delta ,m_{j_\delta }] 
 &\To F_{\omega _j}\\
\rresp{f_{\o v_{j_\delta },m_{j_\delta }}(z^{*j_\delta }) : \quad g_L^{(n)}[j_\delta ,m_{j_\delta }] 
 &\To F_{\o \omega _j}}\end{align*}
on the conjugacy class representatives $g_L^{(n)}[j_\delta ,m_{j_\delta }]$ \resp{$g_R^{(n)}[j_\delta ,m_{j_\delta }]$} of $G^{(n)}( {F^+_v})$ \resp{$G^{(n)}( {F^+_{\o v}})$} on which $ \theta _{G^{(n)}_L}$ \resp{$\theta _{G^{(n)}_R}$} is defined.

$z^{j_\delta }$ \resp{$z^{*j_\delta }$} are the coordinate functions on 
$g_L^{(n)}[j_\delta ,m_{j_\delta }]$ \resp{$g_R^{(n)}[j_\delta ,m_{j_\delta }]$} with respect to the charts:
\begin{align*}
 c_{j_\delta ,m_{j_\delta }}z^{j_\delta } : \quad 
g_L^{(n)}[j_\delta ,m_{j_\delta }]  &\To g_L^{(n)}[j  ,m_{j }] \\
\rresp{c^*_{j_\delta ,m_{j_\delta }}z^{*j_\delta } : \quad 
g_R^{(n)}[j_\delta ,m_{j_\delta }]  &\To g_R^{(n)}[j  ,m_{j }] }\end{align*}
where $g_L^{(n)}[j  ,m_{j }]$ \resp{$g_R^{(n)}[j  ,m_{j }]$} are the corresponding complex conjugacy class representatives.

If the conjugacy class representatives $g_L^{(n)}[j _\delta  ,m_{j _\delta }]$ \resp{$g_R^{(n)}[j_\delta   ,m_{j_\delta  }]$} are glued together, then, a Laurent polynomial corresponding to the mapping:
\[ f_v(z) : \quad G^{(n)}(F^+_v) \To  {F_\omega }\qquad 
\rresp{f_{\o v}(z^*) : \quad G^{(n)}(F^+_{\o v}) \To  {F_{\o \omega} }}\]
is given, on $G^{(n)}(F^+_{v_\oplus})$ \resp{$G^{(n)}(F^+_{\o v_\oplus})$}~, by:
\begin{align*}
f_v(z) &= \sum^r_{j_\delta =1} \ \sum_{m_{j_\delta }} c_{j_\delta ,m_{j_\delta }}\ z^{j_\delta }\;, \qquad 1\le j_\delta \le r\le \infty \\[11pt]
\rresp{f_{\o v}(z^*) &= \sum^r_{j_\delta =1} \ \sum_{m_{j_\delta} } c^*_{j_\delta ,m_{j_\delta }}\ z^{*j_\delta }\;.}
\end{align*}
where:
\quad $F^+_{v_\oplus}=\ds\bigoplus^r_{j_\delta =1} \ \bigoplus_{m_{j_\delta }} F^+_{v_{j_\delta ,m_{j_\delta }}}$ \quad \resp{$F^+_{\o v}=\ds\bigoplus^r_{j_\delta =1}\ \bigoplus_{m_{j_\delta }} F^+_{\o v_{j_\delta ,m_{j_\delta }}}$}.



Note that the Laurent polynomial $f_v(z)$ \resp{$f_{\o v}(z^*)$} is the sum of the Laurent monomials $f_{v_{j_\delta }}(z^{j_\delta })$ \resp{$f_{\o v_{j_\delta }}(z^{*j_\delta })$} on $g_L^{(n)}[j_\delta ,m_{j_\delta }]$ \resp{$g_R^{(n)}[j_\delta ,m_{j_\delta }]$}:
\[ f_v(z) = \sum_{j_\delta } \ \sum_{m_{j_\delta }} f_{v_{j_\delta ,m_{j_\delta }}}(z^{j_\delta })
\quad 
\rresp{f_{\o v}(z^*) = \sum_{j_\delta } \ \sum_{m_{j_\delta }} f_{\o v_{j_\delta ,m_{j_\delta }}}(z^{*j_\delta })}.\]
So, on $G^{(n)}(F^+_{v_\oplus})$ \resp{$G^{(n)}(F^+_{\o v_\oplus})$}, the function $f_v(z)$ \resp{$f_{\o v}(z^*)$}, defined in a neighbourhood of a point $z_0$ \resp{$z_0^*$} of $\cit^n$~, is holomorphic at $z_0$ \resp{$z_0^*$} if we have the multiple power series development:
\begin{align*}
f_v(z) &= \sum^\infty _{j_\delta =1} \ \sum_{m_{j_\delta }} c_{j_\delta ,m_{j_\delta }}\ (z_1-z_{01})^{j_\delta }\ \cdots\ (z_n-z_{0n})^{j_\delta }\\[11pt]
\rresp{f_{\o v}(z^*) &= \sum^\infty _{j_\delta =1} \ \sum_{m_{j_\delta }} c^*_{j_\delta ,m_{j_\delta }}\ (z^*_1-z^*_{01})^{j_\delta }\ \cdots\ (z^*_n-z^*_{0n})^{j_\delta }}
\end{align*}
where:
\Bi
\item $z_1,z_{01},\cdots,z_n,z_{0n}$ are complex functions of one real variable;
\item $z_i: F^+_{v^1_{i_\sigma }}\To F_{\omega ^1_i}$

\item $c_{j_\delta ,m_{j_\delta }}$ \resp{$c^*_{j_\delta ,m_{j_\delta }}$} is in one-to-one correspondence with the product of the square roots of the eigenvalues of the $(j_\delta ,m_{j_\delta })$-th coset representative
$U_{j_\delta ,m_{j_{\delta_R} }} \times U_{j_\delta ,m_{j_{\delta_L} }}$ of the product $T_R(n;t)\otimes T_L(n;t)$ of the Hecke operators.
\Ei

And, the global holomorphic representation $\Irr\hol^{(n)}(\theta _{G^{(n)}_R} \otimes\theta _{G^{(n)}_L})$ of the bisemisheaf $\theta _{G^{(n)}_R} \otimes\theta _{G^{(n)}_L}$ is given by the morphism:
\[\Irr\hol^{(n)}_{\theta _{G\RL}}: \quad \theta _{G^{(n)}_R} \otimes\theta _{G^{(n)}_L}\To f_{\o v}(z^*) \otimes f_v(z)\]
where $f_{\o v}(z^*)\otimes f_v(z)$ is the holomorphic bifunction obtained by gluing together and adding the bisections of the bisemisheaf $\theta _{G^{(n)}_R} \otimes\theta _{G^{(n)}_L}$~.
\vskip 11pt

\subsubsection[Holomorphic representation of the covering bisemisheaf $\theta  _{SOT(1)_R}\otimes\theta  _{SOT(1)_L}$]{\bbf Holomorphic representation of the covering bisemisheaf\linebreak $\theta  _{SOT(1)_R}\otimes\theta  _{SOT(1)_L}$}

The singularities of the covering bisemisheaf $\theta ^* _{SOT(1)_R}\otimes\theta ^* _{SOT(1)_L}$ having been resolved, the free bisemisheaf $\theta  _{SOT(1)_R}\otimes\theta  _{SOT(1)_L}$ can be endowed with a holomorphic representation as it was done in section 4.2.3 for $\theta _{G^{(n)}_R} \otimes\theta _{G^{(n)}_L}$~.

Let $g^{(n)}_{SOT(1)_L}[j_{\delta_{(\rm cov})},m_{j_{\delta_{(\rm cov})}}]$ 
\resp{$g^{(n)}_{SOT(1)_R}[j_{\delta ({\rm cov})},m_{j_{\delta ({\rm cov})}}]$ } denote the conjugacy class representative of the algebraic semigroup $G^{(n)}( {F^+_{v_{\rm cov}}})$
\resp{$G^{(n)}( {F^+_{\o v_{\rm cov}}})$} covering $G^{(n)}( {F^+_{v}})$ \resp{$G^{(n)}( {F^+_{\o v}})$}  where
${F^+_{\o v_{\rm cov}}}$ denotes the set of covering completions in such a way that the
completions $F^+_{v_{j_{\delta_{(\rm cov})},m_{j_{\delta _{\rm cov}}}}}$ are characterized by ranks
 $r_{j_{\delta_{(\rm cov})}}=j_{\delta _{\rm cov}}\centerdot N$ inferior or equal to corresponding ranks $r_{j_\delta }=j_\delta \centerdot N$ of $F^+_{v_{j_\delta ,m_{j_\delta }}}$~: so, $j_{\delta _{(\rm cov)}}\le j_\delta $~.
 
 The sections of the semisheaf $\theta ^* _{SOT(1)_L}$ \resp{$\theta ^* _{SOT(1)_R}$} are the complex-valued differentiable functions (resp. cofunctions):
\begin{align*}
f_{v_{j_{\delta _{\rm cov)}},m_{j_{\delta _{(\rm cov)}}}}}(z^{j_{\delta _{(\rm cov)}}}): \quad
g^{(n)}_{SOT(1)_L}[j_{\delta _{(\rm cov)}},m_{j_{\delta _{(\rm cov)}}}]&\To F_{\omega _{j_{(\rm cov)}}}\\
\rresp{f_{\o v_{j_{\delta _{\rm cov)}},m_{j_{\delta _{(\rm cov)}}}}}(z^{*j_{\delta _{(\rm cov)}}}): \quad
g^{(n)}_{SOT(1)_R}[j_{\delta _{(\rm cov)}},m_{j_{\delta _{(\rm cov)}}}]&\To F_{\o\omega _{j_{(\rm cov)}}}}\;\end{align*}
If the conjugacy class representatives $g^{(n)}_{SOT(1)_L}[j_{\delta _{(\rm cov)}},m_{j_{\delta _{(\rm cov)}}}]$ \resp{$g^{(n)}_{SOT(1)_R}[j_{\delta _{(\rm cov)}},m_{j_{\delta _{(\rm cov)}}}]$} are glued together, then the Laurent polynomial
\begin{align*}
f_{v_{\rm cov}}(z_{\rm cov}) &= \sum^r_{j_{\delta_{(\rm cov)}}=1} \ \sum_{m_{j_{\delta_{(\rm cov)}}}}
c_{j_{\delta _{(\rm cov)}},m_{j_{\delta _{(\rm cov)}}}}\ z^{j_{\delta _{(\rm cov)}}}\;, \qquad 1\le j_{\delta _{(\rm cov)}}\le r\le \infty \\
\rresp{f_{\o v_{\rm cov}}(z^*_{\rm cov}) &= \sum^r_{j_{\delta_{(\rm cov)}}=1} \ \sum_{m_{j_{\delta_{(\rm cov)}}}}
c^*_{j_{\delta _{(\rm cov)}},m_{j_{\delta _{(\rm cov)}}}}\ z^{*j_{\delta _{(\rm cov)}}}}
\end{align*}
can be defined on $G^{(n)}( {F^+_{v_{(\rm cov)}}})$ \resp{$G^{(n)}( {F^+_{\o v_{(\rm cov)}}})$} where $c_{j_{\delta _{(\rm cov)}},m_{j_{\delta _{(\rm cov)}}}}$ is  in one-to-one correspondence with the product of the square roots of the eigenvalues of the $(j_{\delta _{(\rm cov)}},m_{j_{\delta _{(\rm cov)}}})$-th coset representative of the product of the Hecke operators.

And, the corresponding holomorphic function (resp. cofunction) will be given by the multiple power series development:
\begin{align*}
f_{v_{\rm cov}}(z_{\rm cov}) &= \sum^r_{j_{\delta_{(\rm cov)}}=1} \ \sum_{m_{j_{\delta_{(\rm cov)}}}}
c_{j_{\delta _{(\rm cov)}},m_{j_{\delta _{(\rm cov)}}}}\ (z_1-z_{01})^{j_{\delta _{(\rm cov)}}} \ \cdots \ (z_n-z_{0n})^{j_{\delta _{(\rm cov)}}}\\[11pt]
\rresp{f_{\o v_{\rm cov}}(z^*_{\rm cov}) &= \sum^r_{j_{\delta_{(\rm cov)}}=1} \sum_{m_{j_{\delta_{(\rm cov)}}}}
c^*_{j_{\delta _{(\rm cov)}},m_{j_{\delta _{(\rm cov)}}}}\ (z^*_1-z^*_{01})^{j_{\delta _{(\rm cov)}}} \ \cdots \ (z^*_n-z^*_{0n})^{j_{\delta _{(\rm cov)}}}}.\end{align*}
The global homomorphic representation  $\Irr\hol^{(n)}(\theta _{SOT(1)_R}\otimes \theta _{SOT(1)_L})$ of the covering bisemisheaf $\theta _{SOT(1)_R}\otimes \theta _{SOT(1)_L}$ is given by the morphism:
\[\Irr\hol^{(n)}_{\theta_{SOT(1)\RL}}: \quad \theta _{SOT(1)_R}\otimes \theta _{SOT(1)_L})
\To f_{\o v_{\rm cov}}(z^*_{\rm cov}) \otimes f_{v_{\rm cov}}(z_{\rm cov})\;.\]
\vskip 11pt

\subsubsection[Covering $n$-dimensional representation of Weil groups]{\bbf Covering $n$-dimensional representation of Weil groups}

As the covering bisemisheaf $\theta _{SOT(1)_R}\otimes \theta _{SOT(1)_L})$ is defined on the covering algebraic bilinear semigroup  $G^{(n)}( {F^+_{\o v_{(\rm cov)}}}\times  {F^+_{v_{(\rm cov)}}})$, a $n$-dimensional irreducible real representation $\Irr\Rep^{(n)}_{F^{\rm cov}\RL}(W^{ab}_{F^{\rm cov}_R} \times W^{ab}_{F^{\rm cov}_L} )$ of the bilinear global Weil group $(W^{ab}_{F^{\rm cov}_R} \times W^{ab}_{F^{\rm cov}_L})$ must correspond to it.

The global Weil groups $W^{ab}_{F^{\rm cov}_R}$ and $ W^{ab}_{F^{\rm cov}_L}$ may be defined as in section 4.1.1 by:
\begin{align*}
W^{ab}_{F^{\rm cov}_R }&= \bigoplus_{j_{\delta _{(\rm cov)}}} \
\bigoplus_{m_{j_{\delta _{(\rm cov)}}}} 
\Gal(\dot{\widetilde F}^+_{\o v_{j_{\delta _{(\rm cov)}},m_{j_{\delta _{(\rm cov)}}}}}\big/F^0) 
\\[11pt]
W^{ab}_{F^{\rm cov}_L} &= 
\bigoplus_{j_{\delta _{(\rm cov)}}} \
\bigoplus_{m_{j_{\delta _{(\rm cov)}}}} 
\Gal(\dot{\widetilde F}^+_{v_{j_{\delta _{(\rm cov)}},m_{j_{\delta _{(\rm cov)}}}}}\big/F^0) 
\end{align*}
where $\dot{\widetilde F}^+_{v_{j_{\delta _{(\rm cov)}},m_{j_{\delta _{(\rm cov)}}}}}$ is the ramified Galois extension corresponding to the completion\linebreak
$ F^+_{v_{j_{\delta _{(\rm cov)}},m_{j_{\delta _{(\rm cov)}}}}}$~.
\vskip 11pt

\subsubsection{Proposition}


{\em On the bisemisheaf $\theta _{G^{(n)}_R} \otimes\theta _{G^{(n)}_L}$ affected by degenerate  singularities, the following global holomorphic correspondences exist:
\[\begin{psmatrix}[colsep=0cm,rowsep=1.5cm]
\Irr\Rep^{(n)}_{W_{F^+\RL}}(W^{ab}_{F^+_R}\times W^{ab}_{F^+_L}) &&&
 \Irr\hol^{(n)}(\theta _{G^{(n)}_R}\otimes \theta _{G^{(n)}_L})\\
\theta _{G^{(n)}_R}\otimes \theta _{G^{(n)}_L} &&&
 f_{\o v}(z^*) \otimes f_v(z)\\
&&\theta^*_{G^{(n)}_R}\otimes \theta^*_{G^{(n)}_L}&  \\
 & &(\theta ^*_{G^{(n)}_R}\times \theta _{S_R} )\otimes (\theta ^*_{G^{(n)}_L}\times \theta _{S_L}) &\\
 \qquad\qquad(\theta ^*_{G^{(n)}_R}\otimes \theta ^*_{G^{(n)}_L}) &\cup & (\theta ^*_{SOT(1)_R}\otimes \theta ^*_{SOT(1)_L}) & \\
 \theta _{SOT(1)_R}\otimes \theta _{SOT(1)_L} & & &
f_{\o v_{\rm cov}}(z^*_{\rm cov}) \otimes f_{v_{\rm cov}}(z_{\rm cov}) \\
 \Irr\Rep ^{(n)}_{W_{F^{\rm cov}\RL}}(W^{ab}_{F^{\rm cov}_R} \times W^{ab}_{F^{\rm cov}_L})
&&&\Irr\hol^{(n)}(\theta _{SOT(1)_R}\otimes \theta _{SOT(1)_L})
\psset{nodesep=.5cm}
\everypsbox{\scriptstyle}
\ncline[arrows=->]{1,1}{1,4}
\ncline[arrows=->]{2,1}{2,4}
\ncline[doubleline=true]{1,1}{2,1}
\ncline[doubleline=true]{1,4}{2,4}
\ncline[arrows=->]{2,1}{3,3}>{\o\rho_{G_R}\times \o\rho_{G_L}}
\ncline[arrows=->]{3,3}{4,3}>{D_{S_R}\times D_{S_L}}
\ncline[arrows=->]{4,3}{5,1}>{SOT^{\max}_R\times SOT^{\max}_L}
\ncline[arrows=->]{5,1}{2,1}>{\o\rho^{(\rm desing)}_{G_R} \times \o\rho^{(\rm desing)}_{G_L}}
\ncline[arrows=->]{5,3}{6,1}>{\o\rho^{(\rm desing)}_{SOT(1)_R} \times \o\rho^{(\rm desing)}_{SOT(1)_L}}
\ncline[arrows=->]{6,1}{6,4}
\ncline[arrows=->]{7,1}{7,4}
\ncline[doubleline=true]{6,1}{7,1}
\ncline[doubleline=true]{6,4}{7,4}
\end{psmatrix}\]
where:
\Bena
\item $\Irr\Rep^{(n)}_{W_{F^+\RL}}(W^{ab}_{F^+_R}\times W^{ab}_{F^+_L}) \to
 \Irr\hol^{(n)}(\theta _{G^{(n)}_R}\otimes \theta _{G^{(n)}_L})$ is the global holomorphic correspondence on the bisemisheaf $\theta _{G^{(n)}_R}\otimes \theta _{G^{(n)}_L}$ submitted consecutively to:
 \Be
 \item the singularization morphism $\o\rho_{G_R}\times \o\rho_{G_L}$~;
 \item the versal deformation $D_{S_R}\times D_{S_L}$~;
 \item the blow- up $SOT^{\max}_R\times SOT^{\max}_L$ of the versal deformation;
 \item the desingularization $\o\rho^{(\rm desing)}_{G_R} \times \o\rho^{(\rm desing)}_{G_L}$~.
 \Ee
 
 \item $\Irr\Rep ^{(n)}_{W_{F^{\rm cov}\RL}}(W^{ab}_{F^{\rm cov}_R} \times W^{ab}_{F^{\rm cov}_L})
 \to \Irr\hol^{(n)}(\theta _{SOT(1)_R}\otimes \theta _{SOT(1)_L})$ is the global holomorphic correspondence on the covering bisemisheaf $\theta _{SOT(1)_R}\otimes \theta _{SOT()_L}$ generated by the versal deformation $D_{S_R}\times D_{S_L}$ followed by the spreading-out isomorphism 
 $SOT^{\max}_R\times SOT^{\max}_L$ of the singular bisemisheaf $\theta ^*_{G^{(n)}_R}\otimes \theta ^*_{G^{(n)}_L}$~.
 \Ee
}
\vskip 11pt

\bpr This diagram proceeds from the preceding developments.  It then results that the blow-up of the versal deformation of the singular bisemisheaf $\theta ^*_{G^{(n)}_R}\otimes \theta ^*_{G^{(n)}_L}$ generates the covering bisemisheaf $\theta _{SOT(1)_R}\otimes \theta _{SOT()_L}$ from which a new global holomorphic correspondence can be established.\epr
\vskip 11pt

\subsubsection{Toroidal compactification}

In order to get a possible automorphic representation of the bisemisheaf $\theta _{G^{(n)}_R}\otimes \theta _{G^{(n)}_L}$ on the bilinear algebraic semigroup $\Gnv$ and of the covering bisemisheaf $\theta _{SOT(1)_R}\otimes \theta _{SOT(1)_L}$ on the covering bilinear algebraic semigroup $G^{(n)}( {F^+_{\o v_{\rm cov}}}\times {F^+_{v_{\rm cov}}})$~, a toroidal compactification of these bilinear algebraic semigroups $\Gnv$ and 
$G^{(n)}({F^+_{\o v_{\rm cov}}}\times {F^+_{v_{\rm cov}}})$ must be realized.

According to propositions 3.2.2 and 3.2.3, the conjugacy class representatives\linebreak $g^{(n)}_L[j_\delta ,m_{j_\delta }]\in G^{(n)}( {F^+_v})$ \resp{$g^{(n)}_R[j_\delta ,m_{j_\delta }]\in G^{(n)}( {F^+_{\o  v}})$} decompose into $(j_\delta )^{n}$ completions of rank $N$~.

The toroidal compactification of the linear algebraic semigroup $G^{(n)}( {F^+_{v}})$\linebreak \resp{$G^{(n)}( {F^+_{\o  v}})$} can be carried out by considering the horizontal rotational tangent bundle (see section 3.2.4):
\begin{align*}
\tau _{G_L^{(n)}(\ {F^+_{v}})} : \quad T(G_L^{(n)}( {F^{+,T}_{v}})) &\To G_L^{(n)}( {F^{+,T}_{v}})\\[11pt]
\rresp{\tau _{G_R^{(n)}(  {F^+_{\o v}})} : \quad T(G_R^{(n)}( {F^{+,T}_{\o v}})) &\To G_R^{(n)}( {F^{+,T}_{\o  v}})}
\end{align*}
whose total space
$T(G_L^{(n)}( {F^{+,T}_{v}}))$
\resp{$T(G_R^{(n)}( {F^{+,T}_{\o v}}))}$
is a projective linear algebraic semigroup $PG_L^{(n)}( {F^{+,T}_{v}})$ \resp{$PG_R^{(n)}( {F^{+,T}_{\o v}})$}.

The sections $g^{(n)}_{T_L}[j_\delta ,m_{j_\delta }]$ \resp{$g^{(n)}_{T_R}[j_\delta ,m_{j_\delta }]$} of the total space $T(G_L^{(n)}( {F^{+,T}_{v}}))$ \resp{$T(G_R^{(n)}( {F^{+,T}_{\o v}}))$} are $n$-dimensional real semitori whose one-dimensional fibres are semicircles obtained by toroidal deformation of the completions of rank $N\centerdot j_\delta $ of $g^{(n)}_L[j_\delta ,m_{j_\delta }]$ \resp{$g^{(n)}_R[j_\delta ,m_{j_\delta }]$} under the action of the horizontal tangent bundle
 $\tau _{G_L^{(n)}( {F^{+,T}_{v}})} $ \resp{$\tau _{G_R^{(n)}( {F^{+,T}_{\o v}})} $}.
 
 Remark that this kind of toroidal compactification is in one-to-one correspondence with the toroidal compactification of the conjugacy class representatives $g^{(n)}_L[j_\delta ,m_{j_\delta }]$ \resp{$g^{(n)}_R[j_\delta ,m_{j_\delta }]$} given in terms of the projective emergent isomorphism $\gamma ^c_L$ \resp{$\gamma ^c_R$} introduced in sections 3.2 and 3.3 of \cite{Pie1}.
 \vskip 11pt
 
 The toroical compactification of the covering algebraic semigroup $G_L^{(n)}(  {F^{+}_{v_{\rm cov}}})$ \resp{$G_R^{(n)}(  {F^{+}_{\o v_{\rm cov}}})$} can be performed similarly, i.e. by considering the horizontal rotational tangent bundle:
\begin{align*}
\tau _{G_L^{(n)}( {F^{+,T}_{v_{\rm cov}}})} : \quad
T(G_L^{(n)}({F^{+,T}_{v_{\rm cov}}})) 
&\To G_L^{(n)}(  {F^{+,T}_{v_{\rm cov}}})\\[11pt]
\rresp{\tau _{G_R^{(n)}( {F^{+,T}_{\o v_{\rm cov}}})} : \quad 
T(G_R^{(n)}(  {F^{+,T}_{\o v_{\rm cov}}})) 
&\To G_R^{(n)}(  {F^{+,T}_{\o  v_{\rm cov}}})
}
\end{align*}
in such a way that:
\Bena
\item the total space of 
$\tau _{G_L^{(n)}( {F^{+,T}_{v_{\rm cov}}})} $
\resp{$\tau _{G_R^{(n)}( {F^{+,T}_{\o v_{\rm cov}}})}$} is given by
 $T(G_L^{(n)}( {F^{+,T}_{v_{\rm cov}}}))$ 
 \resp{$T(G_R^{(n)}( {F^{+,T}_{\o v_{\rm cov}}}))$}
~.

\item the sections of the total space $T(G_L^{(n)}( {F^{+,T}_{v_{\rm cov}}}))$
\resp{$T(G_R^{(n)}( {F^{+,T}_{v_{\rm cov}}}))$} covering the sections 
$g^{(n)}_{T_L}[j_\delta ,m_{j_\delta }]$ \resp{$g^{(n)}_{T_R}[j_\delta ,m_{j_\delta }]$} are not necessarily complete $n$-dimensional real semitori because their one-dimensional fibres are in one-to-one correspondence with the completions of 
$G_L^{(n)}( {F^{+}_{v_{\rm cov}}})$
\resp{$G_R^{(n)}( {F^{+}_{v_{\rm cov}}})$} having a rank 
$r_{j_{\delta _{\rm cov}}}=j_{\delta _{\rm cov}}\centerdot N$ inferior or equal to the rank $r_{j_\delta }=j_{\delta }\centerdot N$ of the completions of $g^{(n)}_L[j_\delta ,m_{j_\delta }]$ \resp{$g^{(n)}_R[j_\delta ,m_{j_\delta }]$} (see section 4.2.4).

\Ee
\vskip 11pt

The bisemisheaf on the toroidal bilinear algebraic semigroup 
$G^{(n)}( {F^{+,T}_{\o v}}\times {F^{+,T}_{v}})$ will be noted $\theta _{G^{(n)}_{T_R}}\otimes 
\theta _{G^{(n)}_{T_L}}$ and the bisemisheaf on the covering toroidal bilinear algebraic semigroup 
$G^{(n)}( {F^{+,T}_{\o v_{\rm cov}}}\times  {F^{+,T}_{v_{\rm cov}}})$ will be written
$\theta ^{\rm cov}_{G^{(n)}_{T_R}}\otimes 
\theta ^{\rm cov}_{G^{(n)}_{T_L}}$~.
\vskip 11pt

\subsubsection{Proposition}

{\em Let
\[ \Irr\hol^{(n)}(\theta _{G^{(n)}_R} \otimes \theta _{G^{(n)}_L}: \quad
\theta _{G^{(n)}_R} \otimes \theta _{G^{(n)}_L} \To f_{\o v}(z^*) \otimes f_v(z)\]
be the global holomorphic representation of the bisemisheaf $\theta _{G^{(n)}_R} \otimes \theta _{G^{(n)}_L}$ given by the holomorphic bifunction $f_{\o v}(z^*) \otimes f_v(z)$ as introduced in section 4.2.3.

Then, $\tau ^{\rm tor}(\Irr\hol^{(n)}_{\theta_{G\RL}})$~, denoting the toroidal compactification of the global holomorphic representation of $\theta _{G^{(n)}_R} \otimes \theta _{G^{(n)}_L}$~, generates the corresponding elliptic representation of the bisemisheaf $\theta _{G^{(n)}_{T_R}} \otimes \theta _{G^{(n)}_{T_L}}$ according to:
\[\begin{psmatrix}[colsep=0.5cm,rowsep=1.5cm]
\Irr\hol^{(n)}(\theta _{G^{(n)}_R} \otimes \theta _{G^{(n)}_L}): &
 \theta _{G^{(n)}_R} \otimes \theta _{G^{(n)}_L} && f_{\o v}(z^*) \otimes f_v(z)\\
\Irr\ELLIP(\theta _{G^{(n)}_{T_R}} \otimes \theta _{G^{(n)}_{T_L}}): &
 \theta _{G^{(n)}_{T_R}} \otimes \theta _{G^{(n)}_{T_L}} && \ELLIP\RL(n,j_\delta ,m_{j_\delta })
\psset{arrows=->,nodesep=.3cm}
\everypsbox{\scriptstyle}
\ncline{1,2}{1,4}
\ncline{2,2}{2,4}
\ncline{1,1}{2,1}>{\tau ^{\rm tor}\Irr\hol^{(n)}_{\theta _{G\RL}}}
\ncline{1,2}{2,2}
\ncline{1,4}{2,4}
\end{psmatrix}\]
where: 
\[ \ELLIP\RL(n,j_\delta ,m_{j_\delta })
= \ELLIP_R(n,j_\delta ,m_{j_\delta })\otimes \ELLIP_L(n,j_\delta ,m_{j_\delta })\]
being the global elliptic representation of $\theta _{G^{(n)}_{T_R}} \otimes \theta _{G^{(n)}_{T_L}}$~, is the product, right by left, of $n$-dimensional real global elliptic semimodules given by:
\begin{align*}
\ELLIP_L &= \bigoplus^r_{j_\delta =1}\ \bigoplus_{m_{j_\delta }}\lambda ^{\half}(n,j_\delta ,m_{j_\delta })\ e^{2\pi ij_\delta x}\;, \\[11pt]
\ELLIP_R &= \bigoplus^r_{j_\delta =1}\ \bigoplus_{m_{j_\delta }}\lambda ^{\half}(n,j_\delta ,m_{j_\delta })\ e^{-2\pi ij_\delta x}\;,
\end{align*}
with:
\Bi
\item $\ds\vec x = \sum^n_{c=1} x_c\ \vec e_c$ a vector of $(F^+_{v^1})^n$~;
\item $\ds\lambda (n,j_\delta ,m_{j_\delta })=\prod^n_{c=1}\lambda _c(n,j_\delta ,m_{j_\delta })$ according to section 4.1.2.
\Ei
}

\bpr $ \ELLIP\RL(n,j_\delta ,m_{j_\delta })$  is also the $n$-dimensional irreducible elliptic representation $ \Irr\ELLIP\RL(\GL_n(\Aa_ {F^{+,T}_{\o v}} \times  \Aa_{F^{+,T}_{v}} ))$ of 
$\GL_n( {F^{+,T}_{\o v}} \times  {F^{+,T}_{v}} )$ as developed in section 4.1.2.

In section 4.2.7, the toroidal compactification of the bisemisheaf $\theta _{G^{(n)}_R} \otimes \theta _{G^{(n)}_L}$ into the bisemisheaf $\theta _{G^{(n)}_{T_R}} \otimes \theta _{G^{(n)}_{T_L}}$ was realized.

It remains to prove that the holomorphic representation of $\theta _{G^{(n)}_R} \otimes \theta _{G^{(n)}_L}$~, given by the holomorphic bifunction $f_{\o v}(z^*) \otimes f_v(z)$~, is in one-to-one correspondence with the global elliptic bisemimodule 
$ \ELLIP\RL(n,j_\delta ,m_{j_\delta })$ under the toroidal compactification of  $f_{\o v}(z^*) \otimes f_v(z)$~.

This is evident, since $f_v(z)$ \resp{$f_{\o v}(z^*)$} decomposes into:
\begin{align*}
f_v(z) &= 
\sum^r_{j_\delta =1} \ \sum_{m_{j_\delta }} c_{j_\delta ,m_{j_\delta }}\ z^{j_\delta }\;, \\[15pt]
\rresp{f_{\o v}(z^*) &= \sum^r_{j_\delta =1} \ \sum_{m_{j_\delta }} c^*_{j_\delta ,m_{j_\delta }}\ z^{*j_\delta }}\;, \end{align*}
according to section 4.2.3, in such a way that:
\Bi
\item each term $c_{j_\delta ,m_{j_\delta }}\ z^{j_\delta }$ \resp{$c^*_{j_\delta ,m_{j_\delta }}\ z^{*j_\delta }$} is a complex-valued differentiable function on the conjugacy class representative 
$g^{(n)}_L[j_\delta ,m_{j_\delta }]$ \resp{$g^{(n)}_R[j_\delta ,m_{j_\delta }]$};
\item the coefficient $c_{j_\delta ,m_{j_\delta }}$ of $f_v(z)$ corresponds to the coefficient 
$\lambda ^{\half}(n,j_\delta ,m_{j_\delta })$ of\linebreak  $ \ELLIP_L(n,j_\delta ,m_{j_\delta })$ according to sections 4.1.2 and 4.2.3.
\Ei

The toroidal deformation of the conjugacy class representatives $g^{(n)}_L[j_\delta ,m_{j_\delta }]$ \resp{$g^{(n)}_R[j_\delta ,m_{j_\delta }]$} into their toroidal equivalents
$g^{(n)}_{T_L}[j_\delta ,m_{j_\delta }]$ \resp{$g^{(n)}_{T_R}[j_\delta ,m_{j_\delta }]$} is such that the one-dimensional fibres of $g^{(n)}_L[j_\delta ,m_{j_\delta }]$ \resp{$g^{(n)}_R[j_\delta ,m_{j_\delta }]$}, which are completions of rank $j_\delta \centerdot N$~, are transformed into semicircles, which are the fibres of $g^{(n)}_{T_L}[j_\delta ,m_{j_\delta }]$ \resp{$g^{(n)}_{T_R}[j_\delta ,m_{j_\delta }]$} (see section 4.2.7).

So, the toroidal deformation of $g^{(n)}_{L}[j_\delta ,m_{j_\delta }]$ \resp{$g^{(n)}_{R}[j_\delta ,m_{j_\delta }]$} corresponds to the mapping:
\begin{align*}
\tau ^{\rm tor}[j_\delta ,m_{j_\delta }]: \quad
g^{(n)}_{L}[j_\delta ,m_{j_\delta }] &\To 
g^{(n)}_{T_L}[j_\delta ,m_{j_\delta }]\\[11pt]
c_{j_\delta ,m_{j_\delta }}\ z^{j_\delta } &\To 
\lambda ^{\half}(n,j_\delta ,m_{j_\delta })\ e^{+2\pi ij_\delta x}\;, \end{align*}
on the $[j_\delta ,m_{j_\delta }]$-th conjugacy class representative of $G^{(n)}( {F^+_v})$~, sending $c_{j_\delta ,m_{j_\delta }}\ z^{j_\delta }$ into the $n$-dimensional real semitorus 
$\lambda ^{\half}(n,j_\delta ,m_{j_\delta })\ e^{2\pi ij_\delta x}$~. 

By adding the toroidal deformations on all conjugacy class representatives $g^{(n)}_{R}[j_\delta ,m_{j_\delta }]\otimes g^{(n)}_{L}[j_\delta ,m_{j_\delta }]$ of $\Gnv$~, we get the searched one-to-one correspondence:
\begin{align*}
\tau ^{\rm tor}[F^+_{\o v}\times F^+_v] : \quad
\bigoplus_{j_\delta ,m_{j_\delta }} (g^{(n)}_{R}[j_\delta ,m_{j_\delta }]\otimes g^{(n)}_{L}[j_\delta ,m_{j_\delta }]) &\To 
\bigoplus_{j_\delta ,m_{j_\delta }} (g^{(n)}_{T_R}[j_\delta ,m_{j_\delta }]\otimes g^{(n)}_{T_L}[j_\delta ,m_{j_\delta }])\\[15pt]
f_{\o v}(z^*)\otimes f_v(z) &\To \ELLIP\RL(n,j_\delta ,m_{j_\delta })\;.\tag*{\eop}\end{align*}
\vskip 11pt

\subsubsection{Proposition}

{\em Let
\[ \Irr\hol^{(n)}(\theta _{SOT(1)_R} \otimes \theta _{SOT(1)_L}): \quad
\theta _{SOT(1)_R} \otimes \theta _{SOT(1)_L} \To f_{\o v_{\rm cov}}(z^*_{\rm cov}) \otimes f_{v_{\rm cov}}(z_{\rm cov})\]
denote the global holomorphic representation of the covering bisemisheaf $\theta _{SOT(1)_R} \otimes \theta _{SOT(1)_L}$ given by the holomorphic bifunction $f_{\o v_{\rm cov}}(z^*_{\rm cov}) \otimes f_{v_{\rm cov}}(z_{\rm cov})$   introduced in section 4.2.4.

Then, $\tau ^{\rm tor}(\Irr\hol^{(n)}_{\theta_{SOT(1)\RL}})$~, denoting the toroidal compactification of the global holomorphic representation of $(\theta _{SOT(1)_R} \otimes \theta _{SOT(1)_L})$~, generates the corresponding partial elliptic representation of the covering bisemisheaf $(\theta _{SOT(1)_{T_R}} \otimes \theta _{SOT(1)_{T_L}})$ according to:
\[\begin{psmatrix}[colsep=0.5cm,rowsep=1.5cm]
\Irr\hol^{(n)}\ba{l}\theta _{SOT(1)_R} \\[-8pt] \ \ \otimes\; \theta _{SOT(1)_L}\ea: &
 \begin{array}{l}\theta _{SOT(1)_R}\\[-8pt] \ \ \otimes\; \theta _{SOT(1)_L}\end{array} && f_{\o v_{\rm cov}}(z^*_{\rm cov}) \otimes f_{v_{\rm cov}}(z_{\rm cov})\\
\Irr\ELLIP^{\rm part}\ba{l}\theta _{SOT(1)_{T_R}}\\[-8pt] \ \ \otimes\; \theta _{SOT(1)_{T_L}}\ea: &
 \begin{array}{l}\theta _{SOT(1)_{T_R}} \\[-8pt]\ \ \otimes\; \theta _{SOT(1)_{T_L}}\end{array} && \ELLIP^{\rm part}\RL(n,j_{\delta_{(\rm cov)}} ,m_{j_{\delta_{(\rm cov)}}})
\psset{arrows=->,nodesep=.3cm}
\everypsbox{\scriptstyle}
\ncline{1,2}{1,4}
\ncline{2,2}{2,4}
\ncline{1,1}{2,1}>{\tau ^{\rm tor}\Irr\hol^{(n)}_{\theta _{SOT(1)\RL}}}
\ncline{1,2}{2,2}
\ncline{1,4}{2,4}
\end{psmatrix}\]
where: 
\[ \ELLIP^{\rm part}\RL(n,j_{\delta_{(\rm cov)}} ,m_{j_{\delta_{(\rm cov)}}})
= \ELLIP^{\rm part}_R(n,j_{\delta_{(\rm cov)}} ,m_{j_{\delta_{(\rm cov)}}})\otimes \ELLIP^{\rm part}_L(n,j_{\delta_{(\rm cov)}} ,m_{j_{\delta_{(\rm cov)}}})\]
 is the product, right by left, of $n$-dimensional real partial global elliptic semimodules given by:
\begin{align*}
\ELLIP^{\rm part}_L(n,j_{\delta_{(\rm cov)}} ,m_{j_{\delta_{(\rm cov)}}}) &= 
\bigoplus^r_{j_{\delta_{(\rm cov)}} =1}\bigoplus_{m_{j_{\delta_{(\rm cov)}} }}\lambda ^{\half}(n,j_{\delta_{(\rm cov)}} ,m_{j_{\delta_{(\rm cov)}}})\ e^{i\cdot j_{\delta_{(\rm cov)}} \cdot x}\;, 
\tag*{$i\cdot j_{\delta _{\rm cov}}\cdot x\le \pi$}\\[15pt]
\ELLIP^{\rm part}_R(n,j_{\delta_{(\rm cov)}} ,m_{j_{\delta_{(\rm cov)}}}) &= \bigoplus^r_{j_{\delta_{(\rm cov)}} =1}\bigoplus_{m_{j_{\delta_{(\rm cov)}} }}\lambda ^{\half}(n,j_{\delta_{(\rm cov)}} ,m_{j_{\delta_{(\rm cov)}}})\ e^{-i\cdot j_{\delta_{(\rm cov)}} \cdot x}\;,
\end{align*}
with:
\Bi
\item $x\in (F^+_{v^1_{(\rm cov)}})^n$~;
\item $\lambda (n,j_{\delta_{(\rm cov)}} ,m_{j_{\delta_{(\rm cov)}} })$  defined as in section 4.2.8.
\Ei
}
\vskip 11pt

\bpr This can be proved similarly as it was done in proposition 4.2.8: that is to say that the holomorphic representation of $(\theta _{SOT(1)_R} \otimes \theta _{SOT(1)_L}$)~, given by the covering holomorphic bifunction $f_{\o v_{\rm cov}}(z^*_{\rm cov}) \otimes f_{v_{\rm cov}}(z_{\rm cov})$~, is in one-to-one correspondence with the global partial elliptic bisemimodule $\ELLIP^{\rm part}\RL(n,j_{\delta_{(\rm cov)}} ,m_{j_{\delta_{(\rm cov)}}})$~.

The procedure is exactly the same as given in the proof of proposition 4.2.8.  The only difference lies in the fact that 
$\ELLIP^{\rm part}\RL(n,j_{\delta_{(\rm cov)}} ,m_{j_{\delta_{(\rm cov)}}})$ is a global ``{\bf partial}'' elliptic bisemimodule.  
This results from the fact that, under the toroidal deformation of the conjugacy class representatives covering 
$g^{(n)}_{L}[j_\delta ,m_{j_\delta }]$ \resp{$g^{(n)}_{R}[j_\delta ,m_{j_\delta }]$}, the one-dimensional fibres covering $g^{(n)}_{L}[j_\delta ,m_{j_\delta }]$ \resp{$g^{(n)}_{R}[j_\delta ,m_{j_\delta }]$} are completions of rank 
$j_{\delta _{(\rm cov)}}\centerdot N$ transformed into incomplete semicircles covering the semicircles of 
$g^{(n)}_{T_L}[j_\delta ,m_{j_\delta }]$ \resp{$g^{(n)}_{T_R}[j_\delta ,m_{j_\delta }]$} having ranks $j_\delta\centerdot N\ge j_{\delta _{(\rm cov)}}\centerdot N$~.\epr
\vskip 11pt

\subsubsection{Proposition}

{\em On the bisemisheaf $\theta _{G^{(n)}_R} \otimes \theta _{G^{(n)}_L}$~, affected by degenerate  singularities, the following global correspondences of Langlands, prolonging the global holomorphic correspondences of proposition 4.2.6, are:

\scalebox{.73}{$\begin{psmatrix}[colsep=0cm,rowsep=2cm]
\Irr\Rep^{(n)}_{W_{F^+\RL}}(W^{ab}_{F^+_R}\times W^{ab}_{F^+_L})  &&
 \Irr\hol^{(n)}(\theta _{G^{(n)}_R}\otimes \theta _{G^{(n)}_L}) & \xrightarrow{\tau ^{\rm tor}\Irr\hol^{(n)}_{\theta _{G\RL}}}\ & \Irr\ELLIP(\theta _{G^{(n)}_{T_R}}\otimes \theta _{G^{(n)}_{T_L}})
\\
\theta _{G^{(n)}_R}\otimes \theta _{G^{(n)}_L} &&
 f_{\o v}(z^*) \otimes f_v(z) &  & \ELLIP\RL(n,j_\delta ,m_{j_\delta })\\
&&\theta^*_{G^{(n)}_R}\otimes \theta^*_{G^{(n)}_L}&  \\
 & &(\theta ^*_{G^{(n)}_R}\times \theta _{S_R} )\otimes (\theta ^*_{G^{(n)}_L}\times \theta _{S_L}) &\\
 \qquad\qquad(\theta ^*_{G^{(n)}_R}\otimes \theta ^*_{G^{(n)}_L}) &\cup & (\theta ^*_{SOT(1)_R}\otimes \theta^* _{SOT(1)_L}) & \\
 \theta ^*_{SOT(1)_R}\otimes \theta ^*_{SOT(1)_L} & & 
f_{\o v_{\rm cov}}(z^*_{\rm cov}) \otimes f_{v_{\rm cov}}(z_{\rm cov})  && \ELLIP^{\rm part}\RL(n,j_{\delta _{(\rm cov)}},m_{j_{\delta _{(\rm cov)}}})\\
 \Irr\Rep ^{(n)}_{W_{F^{\rm cov}\RL}}(W^{ab}_{F^{\rm cov}_R} \times W^{ab}_{F^{\rm cov}_L})
\hspace*{5mm} &\To&\hspace*{5mm}  \Irr\hol^{(n)}(\theta _{SOT(1)_R}\otimes \theta _{SOT()_L})
&  & \Irr\ELLIP^{\rm part}(\theta _{SOT(1)_{T_R}}\otimes\theta _{SOT(1)_{T_L}})
\psset{nodesep=.5cm}
\everypsbox{\scriptstyle}
\ncline[arrows=->]{2,1}{2,3}
\ncline[doubleline=true]{1,1}{2,1}
\ncline[doubleline=true]{1,3}{2,3}
\ncline[arrows=->]{2,1}{3,3}>{\o\rho_{G_R}\times \o\rho_{G_L}}
\ncline[arrows=->]{3,3}{4,3}>{D_{S_R}\times D_{S_L}}
\ncline[arrows=->]{4,3}{5,1}>{SOT^{\max}_R\times SOT^{\max}_L}
\ncline[arrows=->]{5,1}{2,1}>{\o\rho^{(\rm desing)}_{G_R} \times \o\rho^{(\rm desing)}_{G_L}}
\ncline[arrows=->]{5,3}{6,1}>{\o\rho^{(\rm desing)}_{SOT(1)_R} \times \o\rho^{(\rm desing)}_{SOT(1)_L}}
\ncline[arrows=->]{6,3}{6,5}
\ncline[arrows=->]{7,3}{7,5}
\ncline[arrows=->]{6,1}{6,3}
\ncline[arrows=->]{2,3}{2,5}
\ncline[arrows=->]{1,1}{1,3}
\ncline[doubleline=true]{6,1}{7,1}
\ncline[doubleline=true]{6,3}{7,3}
\ncline[doubleline=true]{1,5}{2,5}
\ncline[doubleline=true]{6,5}{7,5}
\end{psmatrix}$}
\vskip 33pt
}

\bpr This diagram shows how the Langlands global correspondence of section 3.4 of \cite{Pie1} can be extended by considering degenerate  singularities on the sections of the bisemisheaf 
$\theta _{G^{(n)}_R} \otimes \theta _{G^{(n)}_L}$ on the bilinear algebraic semigroup
$\Gnv$ constituting the irreducible $n$-dimensional Weil-Deligne representation
$\Irr\Rep^{(n)}_{W_{F^+\RL}}(W^{ab}_{F^+_R}\times W^{ab}_{F^+_L})$ of the bilinear global Weil group
$(W^{ab}_{F^+_{\o v}}\times W^{ab}_{F^+_v})$~.

So, let:
\[\begin{psmatrix}[colsep=1.5cm,rowsep=1.5cm]
\Irr\Rep^{(n)}_{W_{F^+\RL}}(W^{ab}_{F^+_R}\times W^{ab}_{F^+_L}) &&
 \Irr\ELLIP(\theta _{G^{(n)}_{T_R}} \otimes \theta _{G^{(n)}_{T_L}})
\\
\theta _{G^{(n)}_R} \otimes \theta _{G^{(n)}_L} && \ELLIP\RL(n,j_\delta ,m_{j_\delta })
\psset{nodesep=.3cm}
\everypsbox{\scriptstyle}
\ncline[doubleline=true]{1,1}{2,1}
\ncline[doubleline=true]{1,3}{2,3}
\ncline[arrows=->]{1,1}{1,3}
\ncline[arrows=->]{2,1}{2,3}
\end{psmatrix}\]
be the irreducible Langlands global correspondence of section 3.4 of \cite{Pie1} applied to the free bisemisheaf $\theta _{G^{(n)}_R} \otimes \theta _{G^{(n)}_L}$~.

The fact of considering:
\Bena
\item a singularization  $(\o\rho_{G_R}\times \o\rho_{G_L})$ of $(\theta _{G^{(n)}_R} \otimes \theta _{G^{(n)}_L})$~;
\item a versal deformation $(D_{S_R}\times D_{S_L})$ of the singular bisemisheaf
$(\theta^* _{G^{(n)}_R} \otimes \theta^* _{G^{(n)}_L})$~;
\item a spreading-out $(SOT^{\max}_R\times SOT^{\max}_L)$ of the unfolded bisemisheaf
$(\theta^* _{G^{(n)}_R} \times \theta^* _{S_R}) \otimes (\theta^* _{G^{(n)}_L} \otimes \theta^* _{S_L})$ generating the singular bisemisheaf $(\theta^* _{G^{(n)}_R} \otimes \theta^* _{G^{(n)}_L})$ and the covering singular bisemisheaf
$(\theta^* _{SOT(1)_R} \otimes \theta^* _{SOT(1)_L})$~;
\item a desingularization $(\o\rho^{(\rm desing)}_{G_R} \times \o\rho^{(\rm desing)}_{G_L})$ of the bisemisheaf $(\theta^* _{G^{(n)}_R} \otimes \theta^* _{G^{(n)}_L})$ and a desingularization
$(\o\rho^{(\rm desing)}_{G_R} \times \o\rho^{(\rm desing)}_{G_L})$ of the covering bisemisheaf 
$(\theta^* _{SOT(1)_R} \otimes \theta^* _{SOT(1)_L})$ 
\Ee
allows to recover the desingularized original bisemisheaf $\theta _{G^{(n)}_R} \otimes \theta _{G^{(n)}_L}$~, to which the above mentioned Langlands global correspondence can be reformulated, and an additional desingularized covering bisemisheaf $(\theta^* _{SOT(1)_R} \otimes \theta^* _{SOT(1)_L})$ to which the following Langlands global correspondence can be stated:
\[\begin{psmatrix}[colsep=1.5cm,rowsep=1.5cm]
\Irr\Rep^{(n)}_{W_{F^{\rm cov}\RL}}(W^{ab}_{F^{\rm cov}_R}\times W^{ab}_{F^{\rm cov}_L}) &&
 \Irr\ELLIP^{\rm part}(\theta _{SOT(1)_{R}} \otimes \theta _{SOT(1)_{L}})
\\
\theta _{SOT(1)_R} \otimes \theta _{SOT(1)_L} && \ELLIP^{\rm part}\RL(n,j_{\delta_{(\rm cov)}} ,m_{j_{\delta_{(\rm cov)}}})
\psset{nodesep=.3cm}
\everypsbox{\scriptstyle}
\ncline[doubleline=true]{1,1}{2,1}
\ncline[doubleline=true]{1,3}{2,3}
\ncline[arrows=->]{1,1}{1,3}
\ncline[arrows=->]{2,1}{2,3}
\end{psmatrix}\]
where:
\Bi
\item $\Irr\Rep^{(n)}_{W_{F^{\rm cov}\RL}}(W^{ab}_{F^{\rm cov}_R}\times W^{ab}_{F^{\rm cov}_L})$ is the sum of products, right by left, of the equivalence classes of the irreducible $n$-dimensional Weil-Deligne representation of the covering bilinear global Weil group 
$(W^{ab}_{F^{\rm cov}_R}\times W^{ab}_{F^{\rm cov}_L})$ given by the covering bisemisheaf 
$\theta _{SOT(1)_{R}} \otimes \theta _{SOT(1)_{L}}$~;
\item $\Irr\ELLIP(\theta _{SOT(1)_{R}} \otimes \theta _{SOT(1)_{L}})$ is the sum of the products, right by left, of the equivalence classes of the irreducible elliptic representation of the covering bisemisheaf 
$(\theta _{SOT(1)_{R}} \otimes \theta _{SOT(1)_{L}})$~, given by the $n$-dimensional ``partial'' global elliptic bisemimodule $\ELLIP^{\rm part}\RL(n,j_{\delta_{(\rm cov)}} ,m_{j_{\delta_{(\rm cov)}}})$~.\epr
\Ei
\vskip 11pt

\subsubsection[Langlands {\sc reducible} global correspondences on bisemisheaves over reducible bilinear algebraic semigroups affected by degenerate  singularities]{Langlands {\bfseries\sc reducible} global correspondences on bisemisheaves over reducible bilinear algebraic semigroups affected by degenerate  singularities}

The correspondences considered here can be worked out similarly as it was done in chapter 4 fo \cite{Pie1}, and, more particularly, in proposition 4.2.14.

So, on the basis of this proposition 4.2.14 of \cite{Pie1}, Langlands reducible global correspondences can be developed as it was done in the diagram of the preceding proposition 4.2.10 by taking into account that:
\Bena
\item To each irreducible representation of a reducible bilinear algebraic semigroup, a diagram of Langlands global correspondence can be established as in proposition 4.2.10.

\item The bisemisheaves on the irreducible representations, decomposing the reducible representations of a bilinear algebraic semigroup, generate each one:
\Be
\item a Langlands global correspondence on the original desingularized bisemisheaf;
\item a Langlands global correspondence on the covering desingularized bisemisheaf, generated from the original bisemisheaf submitted to versal deformation and spreading-out morphism.
\Ee
\Ee

\section{Langlands global correspondences over monodromy}

In this chapter, the monodromy \cite{Ebe1}, \cite{Ebe2}, \cite{Gro}, \cite{Gri} of (isolated) singularities on the (bisemi)sheaf of differentiable (bi)functions on the complex bilinear algebraic semigroup $\gn$ is analysed and the Langlands global correspondences on the non singular fibres generated by monodromy are developed in the irreducible  and reducible cases.

\subsection[The monodromy of isolated singularities on irreducible complex semisheaves $\theta^{\cit}_{G^{(n)}\LR}$]{\bbf The monodromy of isolated singularities on irreducible complex semisheaves $\theta^{\cit}_{G^{(n)}\LR}$}

\subsubsection{Complex semisheaves on complex algebraic semigroups}

As the Picard-Lefschetz theory is the complex analogue of Morse theory \cite{Del3}, \cite{H-Z}, our attention will be focused on singularities on the complex-valued differentiable functions $\phi^{(n)}_{G^{(\cit)}_L}(z_{g_L})$ \resp{cofunctions $\phi^{(n)}_{G^{(\cit)}_R}(z_{g_R})$} on the \lr complex linear algebraic semigroup $G^{(n)}(F_\omega )$ \resp{$G^{(n)}(F_{\o \omega })$}, taking into account the inclusion $G^{(n)}(F^+_{\o v}\times F^+_{v})\hookrightarrow \gn$ of the real bilinear algebraic semigroup $G^{(n)}(F^+_{\o v}\times F^+_{v})$ into the corresponding complex equivalent
$\gn$ according to section 1.9.

These complex-valued differentiable functions $\phi^{(n)}_{G^{(\cit)}_L}(z_{g_L})$ \resp{cofunctions $\phi^{(n)}_{G^{(\cit)}_R}(z_{g_R})$} are the sections $\phi^{(n)}_{G^{(\cit)}_{j_L}}(z_{g_{j_L}})$ \resp{$\phi^{(n)}_{G^{(\cit)}_{j_R}}(z_{g_{j_R}})$}, $1\le j\le r\le\infty $~,  of a \lr semisheaf $\theta^{(\cit)}_{G^{(n)}_L}$ \resp{$\theta^{(\cit)}_{G^{(n)}_R}$} on the corresponding conjugacy class representatives $g^{(n)}_L[j,m_j]$ \resp{$g^{(n)}_R[j,m_j]$} of the complex linear algebraic semigroup $G^{(n)}(F_\omega )\equiv T_n(F_\omega)$ \resp{$G^{(n)}(F_{\o\omega })\equiv T^t_n(F_{\o\omega)}$}.
\vskip 11pt

\subsubsection{Monodromy in expanding phase}

Remark that the singularisations and the versal deformations, developed in chapter 2, were envisaged in a contracting phase in the sense that:
\Bena
\item the singularisation of a regular $f$-scheme is a contracting surjective morphism (see proposition 2.1.6).
\item the versal deformation of a semisheaf can be described as a contracting fibre bundle according to propositions 2.2.5 and 2.2.6.
\Ee

And, the spreading-out, introduced as the blow-up of the versal deformation in section 3.1, also occurs naturally in a contracting phase, but the projective map of the tangent bundle on the disconnected base semisheaves has to be viewed as an expanding morphism\ldots\ of blow-up (see section 3.1.15).

The contracting phase of a manifold reflects the fact that its submanifolds become closer and closer with respect to a fixed measure.

The monodromy, studied in this chapter, arises in an expanding phase as it will be justified in the following.

This expanding phase, reflecting the expansion of the submanifolds of a given manifold with respect to a fixed measure, is assumed to generate locally contracting surjective morphisms of singularisations as introduced in section 2.1.
\vskip 11pt

\subsubsection{Types of singularities}

\Bi
\item Let $\phi^{(n)}_{G^{(\cit)}_{j_L}}(z_{g_{j_L}})$ \resp{$\phi^{(n)}_{G^{(\cit)}_{j_R}}(z_{g_{j_R}})$} be, as complex-valued differentiable function, a $j$-th section representative of the \lr semisheaf $\theta^{(\cit)}_{G^{(n)}_L}$ \resp{$\theta^{(\cit)}_{G^{(n)}_R}$}.

\item $\phi^{(n)}_{G^{(\cit)}_{j_L}}(z_{g_{j_L}})$ \resp{$\phi^{(n)}_{G^{(\cit)}_{j_R}}(z_{g_{j_R}})$}, submitted locally to contracting surjective morphism(s) of singularisation(s), can be:
\Bean
\item either a Morse function, i.e. a function having a non-degenerate singular point at zero where a local coordinate system $(z_1,\cdots,z_n)$ in $(\cit^n,0)$ exists for which:
\[ \phi^{(n)}_{G^{(\cit)}_{j_L}}(z_1,\cdots,z_n) = \sum^n_{i=1} z^2_i \qquad 
\rresp{\phi^{(n)}_{G^{(\cit)}_{j_R}}(z^*_1,\cdots,z^*_n) = \sum^n_{i=1} (z^*_i)^2};\]
\vskip 11pt

\item or a function having a degenerate singular point \cite{Cam1}, \cite{Cam2} of type $A_k$~, $D_k$~, $E_6$~, $E_7$ or $E_8$~, as mentioned in section 2.19.  In $(\cit^n,0)\approx (\rit^{2n},0)$~, a local real coordinate system $(x_1,\cdots,x_{2n})$ can be found for which (case $A_k$~):
\begin{align*}
\phi^{(2n)}_{G^{(\rit)}_{j_L}}(x_{1_{L_j}},\cdots,x_{2n_{L_j}}) &= x^{k+1}_{1_{L_j}}+\sum^{2n}_{i=2} x^2_{i_{L_j}}\\[11pt]
\rresp{\phi^{(2n)}_{G^{(\rit)}_{j_R}}(x_{1_{R_j}},\cdots,x_{2n_{R_j}}) &= x^{k+1}_{1_{R_j}}+\sum^{2n}_{i=2} x^2_{i_{R_j}}}, \quad x_{i_R}=-x_{i_L}\;.\end{align*}
\Ee

A small deformation (for example, of versal type) allows to split up a compound singular point generally into simpler ones.

For instance, a small deformation of the function
\[ \phi^{(2n)}_{G^{(\rit)}_{j_L}}(x_{1_{L_j}},\cdots,x_{2n_{L_j}})=
x^{k+1}_{1_{L_j}}+\sum^{2n}_{i=2} x^2_{i_{L_j}}\]
allows to transform it into:
\[\wt \phi^{(2n)}_{G^{(\rit)}_{j_L}}(x_{1_{L_j}},\cdots,x_{2n_{L_j}})=
\phi^{(2n)}_{G^{(\rit)}_{j_L}}(x_{1_{L_j}},\cdots,x_{2n_{L_j}})-\varepsilon  \ x_{1_{L_j}}\;;\]
this is a function having $k$ singular points \cite{H-Z}, generally simpler than the original one.

\vskip 11pt
\item All the sections of the \lr semisheaf $\theta ^{(\cit)}_{G_L^{(n)}}$ \resp{$\theta ^{(\cit)}_{G_R^{(n)}}$}, localized in some open ball, are assumed to be affected by the same kind of external perturbations, and, thus, by the same kind of singularities.
\vskip 11pt

\item Assume, on the other hand, that, instead of considering only the semisheaf $\theta ^{(\cit)}_{G_L^{(n)}}$ \resp{$\theta ^{(\cit)}_{G_R^{(n)}}$}, we have this semisheaf
$\theta ^{(\cit)}_{G_L^{(n)}}$ \resp{$\theta ^{(\cit)}_{G_R^{(n)}}$} covered by one or several semisheaves $\{\theta ^*_{SOT(1)_L},\theta ^*_{SOT(2)_L}\}$ \resp{$\{\theta ^*_{SOT(1)_R},\theta ^*_{SOT(2)_R}\}$} due to the versal deformations and spreading-out isomorphisms as described precedingly and especially in section 4.1.5.  Then, we can state the following lemma.
\Ei
\vskip 11pt

\subsubsection{Lemma}

{\em The monodromy on the sections of the semisheaf $\theta ^{(\cit)}_{G_L^{(n)}}$ \resp{$\theta ^{(\cit)}_{G_R^{(n)}}$}, covered partially by semisheaves generated by spreading-out isomorphisms, can not occur before the blow-up of these covering semisheaves.}
\vskip 11pt

\bpr The monodromy arises only in an expanding phase, i.e. in a phase whose distance between objects increases. So, as the semisheaves $\{\theta ^*_{SOT(1)_L},\theta ^*_{SOT(2)_L}\}$ \resp{$\{\theta ^*_{SOT(1)_R},\theta ^*_{SOT(2)_R}\}$}, generated by blow-up of the versal deformations, cover only partially by patches the semisheaf $\theta ^{(\cit)}_{G_L^{(n)}}$ \resp{$\theta ^{(\cit)}_{G_R^{(n)}}$}, the expanding phase gives rise to:
\Bean
\item a blow-up of the covering semisheaves $\{\theta ^*_{SOT(1)_L},\theta ^*_{SOT(2)_L}\}$ \resp{$\{\theta ^*_{SOT(1)_R},\theta ^*_{SOT(2)_R}\}$} disconnecting them completely from $\theta ^{(\cit)}_{G_L^{(n)}}$ \resp{$\theta ^{(\cit)}_{G_R^{(n)}}$}.
\item contracting surjective morphisms of singularisations due to the strong perturbation of the expanding phase.
\item monodromy groups on the sections $\phi^{*(n)}_{G^{(\cit)}_{j_L}}$ \resp{$\phi^{*(n)}_{G^{(\cit)}_{j_R}}$} of the semisheaf
$\theta ^{(\cit)}_{G_L^{(n)}}$ \resp{$\theta ^{(\cit)}_{G_R^{(n)}}$}.\epr
\Ee\vskip 11pt

\subsubsection[Monodromy for non degenerate singularities of corank $2n$]{\bbf Monodromy for non degenerate singularities of corank $2n$}

\Bi
\item Assume that each section $\phi^{(2n)}_{G^{(\rit)}_{j_L}}(x_{1_{L_j}},\cdots,x_{2n_{L_j}})$ \resp{$\phi^{(2n)}_{G^{(\rit)}_{j_R}}(x_{1_{R_j}},\cdots,x_{2n_{R_j}})$}  of the semisheaf $\theta ^{(\rit)}_{G_L^{(n)}}$ \resp{$\theta ^{(\rit)}_{G_R^{(n)}}$} is a Morse function, i.e. a function affected by an isolated non degenerate singularity on a domain $U_{j_L}$ \resp{$U_{j_R}$} included into the conjugacy class representative $g^{(n)}_L[j,m_j]$ \resp{$g^{(n)}_R[j,m_j]$}.  Then, on the domain(s) $U_{j_L}\subset g^{(n)}_L[j,m_j]$ \resp{$U_{j_R}\subset g^{(n)}_R[j,m_j]$}, 
$\phi^{(2n)}_{G^{(\rit)}_{j_L}}(U_{j_L})$ \resp{$\phi^{(2n)}_{G^{(\rit)}_{j_R}}(U_{j_R})$} is described locally by:
\[ \phi^{(2n)}_{G^{(\rit)}_{j_L}}(U_{j_L})=\sum^{2n}_{i=1}x^2_{i_{L_j}} \qquad 
\rresp{\phi^{(2n)}_{G^{(\rit)}_{j_R}}(U_{j_R})=\sum^{2n}_{i=1}x^2_{i_{R_j}}}.\]
\vskip 11pt

\item The critical level set of $\phi^{(2n)}_{G^{(\rit)}_{j_L}}(U_{j_L})$ \resp{$\phi^{(2n)}_{G^{(\rit)}_{j_R}}(U_{j_R})$} is the singular fibre $F^{(2n-1)}_{0_{j_L}}$ \resp{$F^{(2n-1)}_{0_{j_R}}$} given by
\[ \phi^{(2n)}_{G^{(\rit)}_{j_L}}(U_{j_L})=\sum^{2n}_{i=1}x^2_{i_{L_j}}=0 \qquad 
\rresp{\phi^{(2n)}_{G^{(\rit)}_{j_R}}(U_{j_R})=\sum^{2n}_{i=1}x^2_{i_{R_j}}=0}.\]
\vskip 11pt

\item The non-singular fibres $F^{(2n-1)}_{\lambda _{j_L}}$ \resp{$F^{(2n-1)}_{\lambda _{j_R}}$} of $\phi^{(2n)}_{G^{(\rit)}_{j_L}}$ \resp{$\phi^{(2n)}_{G^{(\rit)}_{j_R}}$} are diffeomorphic to the space $TS^{2n-1}_{L_j}$ \resp{$TS^{2n-1}_{R_j}$} of the tangent bundle to a sphere $S^{2n-1}_{L_j}$ \resp{$S^{2n-1}_{R_j}$} of real dimension $(2n-1)$ and radius $r_{L_j}=1$ \resp{$r_{R_j}=1$}.

The non-singular fibres $F^{(2n-1)}_{\lambda _{j_L}}$ \resp{$F^{(2n-1)}_{\lambda _{j_R}}$} are thus diffeomorphic to:
\begin{align*}
TS^{2n-1}_{L_j}&=\{(x_{1_{L_j}},\cdots,x_{2n_{L_j}}) \mid x^2_{1_{L_j}}+\cdots+x^2_{2n_{L_j}}=\lambda _{L_j}\}\\
\rresp{TS^{2n-1}_{R_j}&=\{(x_{1_{R_j}},\cdots,x_{2n_{R_j}}) \mid x^2_{1_{R_j}}+\cdots+x^2_{2n_{R_j}}=\lambda _{R_j}\}}\end{align*}
while the sphere
\begin{align*}
S^{2n-1}_{L_j}&=\{(x_{1_{L_j}},\cdots,x_{2n_{L_j}}) \mid x^2_{1_{L_j}}+\cdots+x^2_{2n_{L_j}}=1\}\\
\rresp{S^{2n-1}_{R_j}&=\{(x_{1_{R_j}},\cdots,x_{2n_{R_j}}) \mid x^2_{1_{R_j}}+\cdots+x^2_{2n_{R_j}}=1\}}\end{align*}
is diffeomorphic to the vanishing cycle $\Delta ^{(2n-1)}_{L_j} \subset F^{(2n-1)}_{\lambda _{L_j}}$ \resp{$\Delta ^{(2n-1)}_{R_j} \subset F^{(2n-1)}_{\lambda _{R_j}}$}.
\vskip 11pt

\item As the vanishing cycle $\Delta ^{(2n-1)}_{L_j}$ \resp{$\Delta ^{(2n-1)}_{R_j}$} is diffeomorphic to the unit sphere $S^{2n-1}_{L_j}$ \resp{$S^{2n-1}_{R_j}$}, it must correspond to a function $\phi ^{(2n-1)}_{P_{L_j}}(x_{p_{j_L}})$\linebreak \resp{$\phi ^{(2n-1)}_{P_{R_j}}(x_{p_{j_R}})$} on the \lr real conjugacy class representative $P^{(2n-1)}(F^+_{v^1_{j,m_j}})$ \resp{$P^{(2n-1)}(F^+_{\o v^1_{j,m_j}})$} of the \lr linear parabolic subgroup $P^{(2n-1)}(F^+_{v^1})$ \resp{$P^{(2n-1)}(F^+_{\o v^1})$}.  Indeed, according to \cite{Pie1}, the bilinear parabolic affine subsemigroup $P^{(2n-1)}(F^+_{\o v^1}\times F^+_{v^1})$ can be considered as the unitary irreducible representation (space) of the algebraic bilinear semigroup  $\GL_{2n-1}(F^+_{\o v}\times F^+_v)$~.
\vskip 11pt

\item Let $\gamma _{j_L}:[0,1]\to \Delta ^{(2n-1)}_{L_j}$ \resp{$\gamma _{j_R}:[0,1]\to \Delta ^{(2n-1)}_{R_j}$} be a closed loop on the vanishing cycle $\Delta ^{(2n-1)}_{L_j}$ \resp{$\Delta ^{(2n-1)}_{R_j}$}.
\Ei
\vskip 11pt

\subsubsection{Proposition}

{\em Let $(\phi^{(2n)}_{G^{(\rit)}_{j_L}}(U_{j_L}), F^{(2n-1)}_{0 _{j_L}}, F^{(2n-1)}_{\lambda _{j_L}}, \Delta ^{(2n-1)}_{L_j},\gamma _{j_L})$ \resp{$(\phi^{(2n)}_{G^{(\rit)}_{j_R}}(U_{j_R}), F^{(2n-1)}_{0 _{j_R}}, F^{(2n-1)}_{\lambda _{j_R}}, \Delta ^{(2n-1)}_{R_j},\linebreak \gamma _{j_R})$} be the 5-th tuple introduced in section 5.1.5.

Then, the mapping
\[ h_{\gamma _{j_L}} : \quad F^{(2n-1)}_{\lambda _{j_L}} \To F^{(2n-1)}_{\lambda _{j_L}} \qquad 
\rresp{h_{\gamma _{j_R}} : \quad F^{(2n-1)}_{\lambda _{j_R}} \To F^{(2n-1)}_{\lambda _{j_R}}}\]
of the non-singular fibre  $ F^{(2n-1)}_{\lambda _{j_L}}$ \resp{$ F^{(2n-1)}_{\lambda _{j_R}}$} into itself is the monodromy of the closed loop $\gamma _{j_L}$ \resp{$\gamma _{j_R}$} realized by the conjugation action of the $j$-th conjugacy class representative of the restricted linear algebraic semigroup $G^{(2n-1)}(F^{+{\rm res}}_{v_j})$ \resp{$G^{(2n-1)}(F^{+{\rm res}}_{\o v_j})$} on the $j$-th conjugacy class representative of the linear parabolic subsemigroup 
$P^{(2n-1)}(F^{+}_{v^1_j})$ \resp{$P^{(2n-1)}(F^{+}_{\o v^1_j})$} where $F^{+{\rm res}}_{v_j}$ \resp{$F^{+{\rm res}}_{\o v_j}$} is the $j$-th real completion restricted to the domain $U_{j_L}$ \resp{$U_{j_R}$}.}
\vskip 11pt

\bpr 
\Bena 
\item The monodromy $h_{\gamma _{j_L}}$ \resp{$h_{\gamma _{j_R}}$} is associated with an injective mapping
\begin{align*} I_{\Delta _{L_j}\to F_{\lambda _{j_L}}} : \quad \Delta ^{(2n-1)}_{L_j} &\To F^{(2n-1)}_{\lambda _{j_L}} \\  
\rresp{I_{\Delta _{R_j}\to F_{\lambda _{j_R}}} : \quad \Delta ^{(2n-1)}_{R_j} &\To F^{(2n-1)}_{\lambda _{j_R}}}\end{align*}
inflating the vanishing cycle $\Delta ^{(2n-1)}_{L_j}$ \resp{$\Delta ^{(2n-1)}_{R_j}$} into the non-singular fibre $F^{(2n-1)}_{\lambda _{j_L}}$ \resp{$F^{(2n-1)}_{\lambda _{j_R}}$}.  This injective mapping $I_{\Delta _{L_j}\to F_{\lambda _{j_L}}}$ \resp{$I_{\Delta _{R_j}\to F_{\lambda _{j_R}}}$} is in one-to-one correspondence with the injective mapping
\begin{align*}
I_{S^{2n-1}_{L_j}\to TS^{2n-1}_{L_j}} : \quad S^{2n-1}_{L_j} &\To TS^{2n-1}_{L_j} \\
\rresp{I_{S^{2n-1}_{R_j}\to TS^{2n-1}_{R_j}} : \quad S^{2n-1}_{R_j} &\To TS^{2n-1}_{R_j}},\end{align*}
i.e. the inverse of the projective mapping of the tangent bundle introduced in section 5.1.5.

$I_{S^{2n-1}_{L_j}\to TS^{2n-1}_{L_j}}$ \resp{$I_{S^{2n-1}_{R_j}\to TS^{2n-1}_{R_j}}$} then inflates the sphere $S^{2n-1}_{L_j}$ \resp{$S^{2n-1}_{R_j}$}, diffeomorphic to the vanishing cycle 
$\Delta ^{(2n-1)}_{L_j}$ \resp{$\Delta ^{(2n-1)}_{R_j}$}, into $TS^{2n-1}_{L_j}$ \resp{$TS^{2n-1}_{R_j}$}, diffeomorphic to the non-singular fibre $ F^{(2n-1)}_{\lambda _{j_L}}$ \resp{$ F^{(2n-1)}_{\lambda _{j_R}}$}, in such a way that the following diagram:
\[ \begin{CD} \Delta ^{2n-1)}_{L_j} @>>{I_{\Delta _{L_j}\to F_{\lambda _{j_L}}}}> F^{(2n-1)}_{\lambda _{j_L}} \\
@V{\wr}VV @V{\wr}VV\\
S^{2n-1}_{L_j} @>>{I_{S^{2n-1}_{L_j}}\to TS^{2n-1}_{L_j}}> TS^{2n-1}_{L_j} \end{CD}\]
be commutative.
\vskip 11pt

\item This injective mapping $I_{\Delta _{L_j}\to F_{\lambda _{j_L}}}$ \resp{$I_{\Delta _{R_j}\to F_{\lambda _{j_R}}}$} corresponds to the conjugation action of the $j$-th conjugacy class representative $G^{(2n-1)}(F^{+(\rm res)}_{v_j})$ \resp{$G^{(2n-1)}(F^{+(\rm res)}_{\o v_j})$} of
$G^{(2n-1)}( {F^{+(\rm res)}_{v}})$ \resp{$G^{(2n-1)}( {F^{+(\rm res)}_{\o v}})$} 
on the $j$-th conjugacy class representative $P^{(2n-1)}(F^{+}_{v^1_j})$ \resp{$P^{(2n-1)}(F^{+}_{\o v^1_j})$} of the linear parabolic subsemigroup, as developed in \cite{Pie1}, 
since the linear parabolic subsemigroup  
$P^{(2n-1)}( {F^{+(\rm res)}_{v^1}})$ \resp{$P^{(2n-1)}( {F^{+(\rm res)}_{\o v^1}})$} 
can be considered as the unitary irreducible representation space of $\GL_{2n-1}( {F^{+(\rm res)}_{v}})$ \resp{$\GL_{2n-1}( {F^{+(\rm res)}_{\o v}})$}.
\vskip 11pt

\item The set of non-singular fibres $F^{(2n-1)}_{\lambda _{j_L}}(t)$ \resp{$F^{(2n-1)}_{\lambda _{j_R}}(t)$}, $t\in[0,1]$ of $\gamma _{j_L}$ \resp{$\gamma _{j_R}$},
generates a sheaf $\Fs_{F^{(2n-1)}_{\lambda _{j_L}}}$ \resp{$\Fs_{F^{(2n-1)}_{\lambda _{j_R}}}$} on the etale sites above $U_{j_L}$ \resp{$U_{j_R}$} \cite{Del1}, \cite{Del2}.\epr
\Ee
\vskip 11pt

\subsubsection{Proposition}

{\em The injective mapping
\begin{align*} I_{\Delta _{L_j}\to F_{\lambda _{j_L}}} : \quad \Delta ^{(2n-1)}_{L_j} &\To F^{(2n-1)}_{\lambda _{j_L}} \\  
\rresp{I_{\Delta _{R_j}\to F_{\lambda _{j_R}}} : \quad \Delta ^{(2n-1)}_{R_j} &\To F^{(2n-1)}_{\lambda _{j_R}}},\end{align*}
being the inverse of the projective mapping of the tangent bundle $TB_{j_L}(\Delta ^{(2n-1)}_{L_j},F^{(2n-1)}_{\lambda _{j_L}},\linebreak I^{-1}_{\Delta _{L_j}\to F_{\lambda _{j_L}}})$
\resp{$TB_{j_R}(\Delta ^{(2n-1)}_{R_j},F^{(2n-1)}_{\lambda _{j_R}}, I^{-1}_{\Delta _{R_j}\to F_{\lambda _{j_R}}})$}, is such that:
\Bena
\item the vanishing cycle $\Delta ^{(2n-1)}_{L_j}$ \resp{$\Delta ^{(2n-1)}_{R_j}$} is characterized by a rank $r_{\Delta ^{(2n-1)}_{j}}=N^{2n-1}$~.

\item the non-singular fibre $F^{(2n-1)}_{\lambda _{j_L}}$ \resp{$F^{(2n-1)}_{\lambda _{j_R}}$} is characterized by a rank $r_{F^{(2n-1)}_{\lambda _{j}}}\le (j\centerdot N)^{2n-1}$~.

\item the fibre $F_{I_{\Delta _{L_j}\to F_{\lambda _{j_L}}}}$ \resp{$F_{I_{\Delta _{R_j}\to F_{\lambda _{j_R}}}}$} of the tangent bungle $TB_{j_L}$ \resp{$TB_{j_R}$} has a rank $r_{F_{I_{\Delta _{L_j}\to F_{\lambda _{j_L}}}}}$ verifying: 
\[r_{F_{I_{\Delta _{L_j}\to F_{\lambda _{j_L}}}}}\le (j\centerdot N)^{2n-1}\;.\]
\Ee}
\vskip 11pt

\bpr \Bena
\item According to section 5.1.5, the vanishing cycle $\Delta ^{(2n-1)}_{L_j}$ \resp{$\Delta ^{(2n-1)}_{R_j}$} is a function on the $j$-th conjugacy class representative of the parabolic subgroup $P^{(2n-1)}(F^{+}_{v^1})$ \resp{$P^{(2n-1)}(F^{+}_{\o v^1})$}.  So we have that its rank is given by \cite{Pie1}:
\[ r_{\Delta ^{(2n-1)}_{L_j}}=N^{2n-1}\;.\]

\item As $F^{(2n-1)}_{\lambda _{j_L}}$ \resp{$F^{(2n-1)}_{\lambda _{j_R}}$} is a non-singular fibre above the $j$-th conjugacy class representative of the linear algebraic semigroup 
$G^{(n)}( {F^{+}_{v}})$ \resp{$G^{(n)}( {F^{+}_{\o v}})$} characterized by a rank:
\[ r_{g^{(n)}[j,m_j]} = (j\centerdot N)^{2n}\]
and as $F^{(2n-1)}_{\lambda _{j_L}}$ \resp{$F^{(2n-1)}_{\lambda _{j_R}}$} results from an inflation mapping from the vanishing cycle $\Delta ^{(2n-1)}_{L_j}$ \resp{$\Delta ^{(2n-1)}_{R_j}$} in such a way that the inflation of $F^{(2n-1)}_{\lambda _{j_L}}$ \resp{$F^{(2n-1)}_{\lambda _{j_R}}$} from $\Delta ^{(2n-1)}_{L_j}$ \resp{$\Delta ^{(2n-1)}_{R_j}$} is proportional to the conjugacy action of the $j$-th conjugacy class representative $g_L^{(n)}[j,m_j]\in G^{(n)}(F^+_v)$ \resp{$g_R^{(n)}[j,m_j]\in G^{(n)}(F^+_{\o v})$} with respect to the $j$-th conjugacy class representative of the parabolic subgroup $P^{(2n)}(F^{+}_{v^1})$ \resp{$P^{(2n)}(F^{+}_{\o v^1})$}, we have that:
\[ r_{F^{(2n-1)}_{\lambda _j}}\le (j\centerdot N)^{2n-1}\]
since $(j\centerdot N)^{2n}$ is the rank of $g_L^{(n)}[j,m_j]$ \resp{$g_R^{(n)}[j,m_j]$}.\epr
\Ee
\vskip 11pt

\subsubsection{Definition: Monodromy operator \cite{Ber1}, \cite{Ber2}, \cite{Chm}}

If we have that
\begin{align*}
h_{\gamma ^*_{j_L}}: \quad H_{2n-1}(\Fs_{F^{(2n-1)}_{\lambda _{j_L}}};\ZZ) &\To 
H_{2n-1}(\Fs_{F^{(2n-1)}_{\lambda _{j_L}}};\ZZ)\\
\rresp{h_{\gamma ^*_{j_R}}: \quad H_{2n-1}(\Fs_{F^{(2n-1)}_{\lambda _{j_R}}};\ZZ) &\To 
H_{2n-1}(\Fs_{F^{(2n-1)}_{\lambda _{j_R}}};\ZZ)},\end{align*}
the action $h_{\gamma ^*_{j_L}}$ \resp{$h_{\gamma ^*_{j_R}}$} of $h_{\gamma _{j_L}}$ \resp{$h_{\gamma _{j_R}}$} in the homology group $H_{2n-1}(\Fs_{F^{(2n-1)}_{\lambda _{j_L}}};\ZZ)$ \resp{$H_{2n-1}(\Fs_{F^{(2n-1)}_{\lambda _{j_R}}};\ZZ)$} of the sheaf $\Fs_{F^{(2n-1)}_{\lambda _{j_L}}}$ \resp{$\Fs_{F^{(2n-1)}_{\lambda _{j_R}}}$} of non-singular fibres 
$F^{(2n-1)}_{\lambda _{j_L}}(t)$ \resp{$F^{(2n-1)}_{\lambda _{j_R}}(t)$}, is the monodromy operator of the closed loop $\gamma _{j_L}$ \resp{$\gamma _{j_L}$} \cite{H-Z}.

Indeed, the homology group $H_{2n-1}(\Fs_{F^{(2n-1)}_{\lambda _{j_L}}};\ZZ)$ \resp{$H_{2n-1}(\Fs_{F^{(2n-1)}_{\lambda _{j_R}}};\ZZ)$} is generated by the homology class of the vanishing cycle $\Delta ^{(2n-1)}_{L_j}$ \resp{$\Delta ^{(2n-1)}_{R_j}$} \cite{H-Z}.
\vskip 11pt

\subsubsection{Definition: The surjective mapping}

The surjective mapping
\begin{align*}
r_{F_{\lambda _{j_L}}\to F_{0_{j_L}}} : \quad F^{(2n-1)}_{\lambda _{j_L}}(t) &\To F^{(2n-1)}_{0_{j_L}}\;, \qquad 0\le t\le 1\\
\rresp{r_{F_{\lambda _{j_R}}\to F_{0_{j_R}}} : \quad F^{(2n-1)}_{\lambda _{j_R}}(t) &\To F^{(2n-1)}_{0 _{j_R}}}
\end{align*}
of the non-singular fibre(s) $F^{(2n-1)}_{\lambda _{j_L}}(t)$ \resp{$F^{(2n-1)}_{\lambda _{j_R}}(t)$} into the singular fibre $F^{(2n-1)}_{0_{j_L}}$ \resp{$F^{(2n-1)}_{0_{j_R}}$} is the retraction of the monodromy.
\vskip 11pt

\subsubsection{Monodromy for degenerate singularities}

If a degenerate singularity (for example of type $A_k$ (see section 5.1.3)) decomposes by deformation into $k$ elementary non degenerate singular points, the single monodromy envisaged in the case of a unique non degenerate singularity, becomes a monodromy group where the loop $\gamma _{j_L}$ \resp{$\gamma _{j_R}$} runs through the fundamental group $\Pi _1(V_{j_L}-\{\omega _{i_L}\},x_{0_L})$ \resp{$\Pi _1(V_{j_R}-\{\omega _{i_R}\},x _{0_R})$} of the complementary of the set of critical values $\omega_{i_L}$
 \resp{$\omega_{i_R}$}
, $1\le i\le k$~, where:
\Bi
\item $V_{j_L}$ \resp{$V_{j_R}$} is a compact domain in $\cit$ included into $U_{j_L}$ \resp{$U_{j_R}$};
\item $\gamma _{j_L}(0) =\gamma _{j_L}(1)=x _{0_L}$ \resp{$\gamma _{j_R}(0) =\gamma _{j_R}(1)=x _{0_R}$}.
\Ei

The complementary of the set of critical values in $V_{j_L}$ \resp{$V_{j_R}$}  is a loop beginning and ending at $x_{0_L}$ \resp{$x _{0_R}$} and passing round the critical values $\omega _{i_L}$ \resp{$\omega _{i_R}$}.

The domain $V_{j_L}$ \resp{$V_{j_R}$} without the $k$ critical values $\{\omega _{i_L}\mid i=1,\cdots,k\}$ \resp{$\{\omega _{i_R}\mid i=1,\cdots,k\}$} of $\phi ^{(2n)}_{G^{(\rit)}_{j_L}}(x_{i_L})$
\resp{$\phi ^{(2n)}_{G^{(\rit)}_{j_R}}(x_{i_R})$} is homotopically equivalent to a bunch of $k$ circles.  So, the fundamental group $\Pi _1(V_{j_L}-\{\omega _{i_L}\},x _{0_L})$ \resp{$\Pi _1(V_{j_R}-\{\omega _{i_R}\},x _{0_R})$} is a free group at $k$ generators.

If $\{v _{i_L}\mid i=1,\cdots,k\}$ \resp{$\{v _{i_R}\mid i=1,\cdots,k\}$} is a set of  paths defining a set of vanishing cycles $\Delta ^{(2n-1)}_{i_{L_j}}\in H_{2n-1}(\Fs_{F^{(2n-1)}_{i_{\lambda _{j_L}}}})$ \resp{$\Delta ^{(2n-1)}_{i_{R_j}}\in H_{2n-1}(\Fs_{F^{(2n-1)}_{i_{\lambda _{j_R}}}})$}, $1\le i\le k$~, then the fundamental group $\Pi _1(V_{j_L}-\{\omega _{i_L}\},x _{0_L})$ \resp{$\Pi _1(V_{j_R}-\{\omega _{i_R}\},x _{0_R})$} is generated by the simple loops $\gamma _{1_{j_L}}, \cdots,\gamma _{k_{j_L}}$ \resp{$\gamma _{1_{j_R}}, \cdots,\gamma _{k_{j_R}}$} associated with the paths $v _{1_L},\cdots,v _{k_L}$ \resp{$v _{1_R},\cdots,v _{k_R}$}.

The monodromy group of $\phi^{(2n)}_{G^{(\rit)}_{j_L}}(U_{j_L})$ \resp{$\phi^{(2n)}_{G^{(\rit)}_{j_R}}(U_{j_R})$} is the image of the homomorphism of the fundamental group $\Pi _1(V_{j_L}-\{\omega _{i_L}\},x _{0_L})$ \resp{$\Pi _1(V_{j_R}-\{\omega _{i_R}\},x _{0_R})$} into the group $\Aut( H_{2n-1}(\Fs_{F^{(2n-1)}_{i_{\lambda _{j_L}}}}))$
\resp{$\Aut( H_{2n-1}(\Fs_{F^{(2n-1)}_{i_{\lambda _{j_R}}}}))$} of automorphisms of
$H_{2n-1}(\Fs_{F^{(2n-1)}_{i_{\lambda _{j_L}}}})$ \resp{$H_{2n-1}(\Fs_{F^{(2n-1)}_{i_{\lambda _{j_R}}}})$} which associate with the loop $\gamma _{i_{j_L}}$ \resp{$\gamma _{i_{j_R}}$} the monodromy operator:
\begin{align*}
h_{\gamma ^*_{i_{j_L}}}: \quad H_{2n-1}(\Fs_{F^{(2n-1)}_{i_{\lambda _{j_L}}}};\ZZ) &\To 
H_{2n-1}(\Fs_{F^{(2n-1)}_{i_{\lambda _{j_L}}}};\ZZ)\\
\rresp{h_{\gamma ^*_{i_{j_R}}}: \quad H_{2n-1}(\Fs_{F^{(2n-1)}_{i_{\lambda _{j_R}}}};\ZZ) &\To 
H_{2n-1}(\Fs_{F^{(2n-1)}_{i_{\lambda _{j_R}}}};\ZZ)},\end{align*}
where $\Fs_{F^{(2n-1)}_{i_{\lambda _{j_L}}}}$ \resp{$\Fs_{F^{(2n-1)}_{i_{\lambda _{j_R}}}}$} is the $i$-th sheaf of non-singular fibres $F^{(2n-1)}_{i_{\lambda _{j_L}}}(t)$ \resp{$F^{(2n-1)}_{i_{\lambda _{j_R}}}(t)$}.  

As developed in \cite{Pie1} and in proposition 5.1.6, 
\begin{align*}
\Aut( H_{2n-1}(\Fs_{F^{(2n-1)}_{i_{\lambda _{j_L}}}})) &= \Aut (\Delta ^{(2n-1)}_{i_{L_j}}) \simeq \Aut (P^{(2n-1)}(F^+_{v^1_j}))\\
\rresp{\Aut( H_{2n-1}(\Fs_{F^{(2n-1)}_{i_{\lambda _{j_R}}}})) &= \Aut (\Delta ^{(2n-1)}_{i_{R_j}}) \simeq \Aut (P^{(2n-1)}(F^+_{\o v^1_j}))}
\end{align*}
where $\Aut (P^{(2n-1)}(F^+_{v^1_j}))$ \resp{$\Aut (P^{(2n-1)}(F^+_{\o v^1_j}))$} is the group of Galois automorphisms of the linear parabolic subsemigroup $P^{(2n-1)}(F^+_{v^1_j})$ \resp{$P^{(2n-1)}(F^+_{\o v^1_j})$} on the $j$-th irreducible completion $F^+_{v^1_j}$ \resp{$F^+_{\o v^1_j}$} of rank $N$~.
\vskip 11pt

\subsubsection{Monodromy for a set of non degenerate singularities}

The monodromy, envisaged in section 5.1.9 for a degenerate singularity splitting up into elementary non degenerate singularities, is also valid for a set of critical points $\omega _{i_L}$ \resp{$\omega _{i_R}$}, $1\le i\le k$~, of $\phi^{(2n)}_{G^{(\rit)}_{j_L}}(\omega _{i_L})$ \resp{$\phi^{(2n)}_{G^{(\rit)}_{j_R}}(\omega _{i_R})$} which are not degenerated and such that their critical values $\omega _{i_L}$ \resp{$\omega _{i_R}$} are distinct.

Then, we can state the following proposition.
\vskip 11pt

\subsubsection{Proposition}

{\em Assume that every section $\phi^{(2n)}_{G^{(\rit)}_{j_L}}(U_{j_L})$ \resp{$\phi^{(2n)}_{G^{(\rit)}_{j_R}}(U_{j_R})$} of the semisheaf $\theta ^{(\rit)}_{G^{(2n)}_L}$ \resp{$\theta ^{(\rit)}_{G^{(2n)}_R}$}, $1\le j\le r\le \infty $~, is endowed with a set of $k$ non degenerate singularities $\omega _{i_L}$ \resp{$\omega _{i_R}$}, $1\le i\le k$~, on the domain $U_{j_L}$ \resp{$U_{j_R}$}.

Let $(\phi^{(2n)}_{G^{(\rit)}_{j_L}}(U_{j_L}), F^{(2n-1)}_{i_{0_{j_L}}}, \Fs_{F^{(2n-1)}_{i_{\lambda _{j_L}}}} , \Delta ^{(2n-1)}_{i_{L_j}} , \gamma _{i_{j_L}})$
\resp{$(\phi^{(2n)}_{G^{(\rit)}_{j_R}}(U_{j_R}), F^{(2n-1)}_{i_{0_{j_R}}}, \Fs_{F^{(2n-1)}_{i_{\lambda _{j_R}}}} ,\linebreak \Delta ^{(2n-1)}_{i_{R_j}} , \gamma _{i_{j_R}})$} be the $i$-th 5-tuple associated with the $i$-th singularity on the $j$-th section of 
$\theta ^{(\rit)}_{G^{(2n)}_L}$ \resp{$\theta ^{(\rit)}_{G^{(2n)}_R}$}, as introduced in proposition 5.1.6.

Then, the $i$ mappings:
\begin{align*}
h_{\gamma _{i_{j_L}}} : \quad F^{(2n-1)}_{i_{\lambda _{j_L}}} &\To F^{(2n-1)}_{i_{\lambda _{j_L}}}\;, \qquad 1\le i\le k\;, \\
\rresp{h_{\gamma _{i_{j_R}}} : \quad F^{(2n-1)}_{i_{\lambda _{j_R}}} &\To F^{(2n-1)}_{i_{\lambda _{j_r}}}}\end{align*}
of the non-singular fibres $F^{(2n-1)}_{i_{\lambda _{j_L}}}$ \resp{$F^{(2n-1)}_{i_{\lambda _{j_R}}}$} into themselves, associated with the ``~$i$~'' non degenerate singularities 
$\omega _{i_L}$ \resp{$\omega _{i_R}$}:
\Bena
\item are the monodromies of the closed loops $\gamma _{i_{j_L}}$ \resp{$\gamma _{i_{j_R}}$};
\item generate the sheaves $\Fs_{F^{(2n-1)}_{i_{\lambda _{j_L}}}}$ \resp{$\Fs_{F^{(2n-1)}_{i_{\lambda _{j_R}}}}$} of non singular fibres $F^{(2n-1)}_{i_{\lambda _{j_L}}}(t)$ \resp{$F^{(2n-1)}_{i_{\lambda _{j_R}}}(t)$} on the etale sites $U_{j_L}$ \resp{$U_{j_R}$};
\item are associated with the injective mappings:
\begin{align*}
I_{\Delta _{i_{L_j}}\to F_{i_{\lambda _{j_L}}}} : \quad \Delta ^{(2n-1)}_{i_{L_j}} &\To F^{(2n-1)}_{i_{\lambda _{j_L}}}(t)\\
\rresp{I_{\Delta _{i_{R_j}}\to F_{i_{\lambda _{j_R}}}} : \quad \Delta ^{(2n-1)}_{i_{R_j}} &\To F^{(2n-1)}_{i_{\lambda _{j_R}}}(t)}
\end{align*}
inflating the vanishing cycles $\Delta ^{(2n-1)}_{i_{L_j}}$ \resp{$\Delta ^{(2n-1)}_{i_{R_j}}$} into the non-singular fibres $F^{(2n-1)}_{i_{\lambda _{j_L}}}(t)$ \resp{$F^{(2n-1)}_{i_{\lambda _{j_R}}}(t)$} and being in one-to-one correspondence with the corresponding injective mappings:
\begin{align*}
I_{S^{2n-1}_{i_{L_j}}\to TS^{2n-1}_{i_{L_j}}} : \quad S^{2n-1}_{i_{L_j}} &\To TS^{2n-1}_{i_{L_j}}\\
\rresp{I_{S^{2n-1}_{i_{R_j}}\to TS^{2n-1}_{i_{R_j}}} : \quad S^{2n-1}_{i_{R_j}} &\To TS^{2n-1}_{i_{R_j}}};\end{align*}
\item result from the monodromy operators:
\begin{align*}
h_{\gamma ^*_{i_{j_L}}}: \quad H_{2n-1}(\Fs_{F^{(2n-1)}_{i_{\lambda _{j_L}}}}) &\To 
H_{2n-1}(\Fs_{F^{(2n-1)}_{i_{\lambda _{j_L}}}}) \\
\resp{h_{\gamma ^*_{i_{j_R}}}: \quad H_{2n-1}(\Fs_{F^{(2n-1)}_{i_{\lambda _{j_R}}}}) &\To 
H_{2n-1}(\Fs_{F^{(2n-1)}_{i_{\lambda _{j_R}}}})}
\end{align*}
on the vanishing cycles $H_{2n-1}(\Fs_{F^{(2n-1)}_{i_{\lambda _{j_L}}}}) \equiv
\Delta ^{(2n-1)}_{i_{L_j}}$ \resp{$H_{2n-1}(\Fs_{F^{(2n-1)}_{i_{\lambda _{j_R}}}}) \equiv
\Delta ^{(2n-1)}_{i_{R_j}}$}.
\Ee}
\vskip 11pt

\subsubsection[Monodromy sheaves above $\theta ^{(\rit)}_{G^{(2n)}_L}$ \resp{$\theta ^{(\rit)}_{G^{(2n)}_R}$}]{\bbf Monodromy sheaves above $\theta ^{(\rit)}_{G^{(2n)}_L}$ \resp{$\theta ^{(\rit)}_{G^{(2n)}_R}$}}

Let $\theta ^{(\rit)}_{G^{(2n)}_R}\otimes \theta ^{(\rit)}_{G^{(2n)}_L}$ be the bisemisheaf on the bilinear affine semigroup  $G^{(2n)}( {F^+_{\o v}}\times  {F^+_v})$~.

If every section of $\theta ^{(\rit)}_{G^{(2n)}_L}$ \resp{$\theta ^{(\rit)}_{G^{(2n)}_R}$} is endowed with $k$ isolated non degenerate singularities, then the set of bisheaves
$\{\Fs_{F^{(2n-1)}_{i_{\lambda _{j_R}}}}\otimes\Fs_{F^{(2n-1)}_{i_{\lambda _{j_L}}}}\}^k_{i=1}$ of non singular bifibres 
$F^{(2n-1)}_{i_{\lambda _{j_R}}}(t)\otimes F^{(2n-1)}_{i_{\lambda _{j_L}}}(t)$ 
are generated by monodromy above every bisection $\phi ^{(2n)}_{G^{(\rit)}_{j_R}}(U_{j_R}) \otimes
\phi ^{(2n)}_{G^{(\rit)}_{j_L}}(U_{j_L})$ 
of $(\theta  ^{(\rit)}_{G^{(2n)}_{R}}\otimes\theta  ^{(\rit)}_{G^{(2n)}_{L}})$~.

Let $\Fs_{F^{(2n-1)}_{i_{\lambda _{j_L}}}}$ \resp{$\Fs_{F^{(2n-1)}_{i_{\lambda _{j_R}}}}$} be the sheaf of non-singular fibres of the monodromy of the closed loop $\gamma _{i_{j_L}}$ \resp{$\gamma _{i_{j_R}}$} associated with the $i$-th singularity on the $j$-th section of $\theta  ^{(\rit)}_{G^{(2n)}_{L}}$ \resp{$\theta  ^{(\rit)}_{G^{(2n)}_{R}}$}.

Consider that there are $b_i$~, $b_i\in\nit$~, non-singular fibres in the sheaf $\Fs_{F^{(2n-1)}_{i_{\lambda _{j_L}}}}$ \resp{$\Fs_{F^{(2n-1)}_{i_{\lambda _{j_R}}}}$}.

As it was assumed that all the sections of $\theta  ^{(\rit)}_{G^{(2n)}_{L}}$ \resp{$\theta  ^{(\rit)}_{G^{(2n)}_{R}}$} are endowed with the same kind of $k$ isolated non degenerate singularities, $b_i$ sheaves $\Fs_{F^{(2n-1)}_{i_{\lambda _{L}}}}(\beta _i)$ \resp{$\Fs_{F^{(2n-1)}_{i_{\lambda _{R}}}}(\beta _i)$}, $1\le \beta _i\le b_i$~, whose sections are the non-singular fibres above the sections of $\theta  ^{(\rit)}_{G^{(2n)}_{L}}$ \resp{$\theta  ^{(\rit)}_{G^{(2n)}_{R}}$}, can be envisaged as generated by monodromy from the $i$-th singularities on all the sections of $\theta  ^{(\rit)}_{G^{(2n)}_{L}}$ \resp{$\theta  ^{(\rit)}_{G^{(2n)}_{R}}$}.

So, above the $i$-th singularity on all the sections of $\theta  ^{(\rit)}_{G^{(2n)}_{L}}$ \resp{$\theta  ^{(\rit)}_{G^{(2n)}_{R}}$}, a set\linebreak $\{\Fs_{F^{(2n-1)}_{i_{\lambda _{L}}}}(\beta _i)\}^{b_i}_{\beta _i=1}$ \resp{$\{\Fs_{F^{(2n-1)}_{i_{\lambda _{R}}}}(\beta _i)\}^{b_i}_{\beta _i=1}$} of $b_i$ monodromy sheaves are generated and, above the set of $i$ singularities on all the sections of $\theta  ^{(\rit)}_{G^{(2n)}_{L}}$ \resp{$\theta  ^{(\rit)}_{G^{(2n)}_{R}}$}, a set of $i\times b_i$ monodromy sheaves can be constructed above $\theta  ^{(\rit)}_{G^{(2n)}_{L}}$ \resp{$\theta  ^{(\rit)}_{G^{(2n)}_{R}}$}.

Let $\Fs_{F^{(2n-1)}_{i_{\lambda _{L}}}}(\beta _i)$ \resp{$\Fs_{F^{(2n-1)}_{i_{\lambda _{R}}}}(\beta _i)$} be the $\beta _i$-th monodromy sheaf, generated from the $i$-th singularity on all the sections of $\theta  ^{(\rit)}_{G^{(2n)}_{L}}$ \resp{$\theta  ^{(\rit)}_{G^{(2n)}_{R}}$}.

The sections of $\Fs_{F^{(2n-1)}_{i_{\lambda _{L}}}}(\beta _i)$ \resp{$\Fs_{F^{(2n-1)}_{i_{\lambda _{R}}}}(\beta _i)$} are in one-to-one correspondence with the conjugacy class representatives of $G^{(2n)}_L(F^+_v)$ \resp{$G^{(2n)}_L(F^+_{\o v})$}: they are thus labelled by the pairs of integers $(j,m_j)$~, $1\le j\le r\le\infty $~, and their ranks verify: 
\[ r_{F^{(2n-1)}_{i_{\lambda _j}}}\le (j\centerdot N)^{2n-1}\;, \]
according to proposition 5.1.7, $F^{(2n-1)}_{i_{\lambda _{j_L}}}(\beta _i)\in 
\Fs_{F^{(2n-1)}_{i_{\lambda _{L}}}}(\beta _i)$ \resp{$F^{(2n-1)}_{i_{\lambda _{j_R}}}(\beta _i)\in 
\Fs_{F^{(2n-1)}_{i_{\lambda _{R}}}}(\beta _i)$} being the $j$-th section.
\vskip 11pt

\subsubsection{Proposition}

{\em Let $\theta  ^{(\rit)}_{G^{(2n)}_{R}}\times \theta  ^{(\rit)}_{G^{(2n)}_{L}}$ be the bisemisheaf on $G^{(2n)}(F^+_{\o v}\times F^+_v)$ such that its linear sections are endowed with $k$ isolated non-degenerate singularities of the same type.

Let $\{\Fs_{F^{(2n-1)}_{i_{\lambda _{R}}}}(\beta _i)\otimes \Fs_{F^{(2n-1)}_{i_{\lambda _{L}}}}(\beta _i)\}_{i,\beta _i}$, $1\le i\le k$~, $1\le \beta _i\le b_i$~, be the set of $k\times b_i$ monodromy bi(semi)sheaves above $\theta  ^{(\rit)}_{G^{(2n)}_{R}}\otimes \theta  ^{(\rit)}_{G^{(2n)}_{L}}$ according to section 5.1.13.

Then, it results that:
\Bena
\item a global holomorphic representation $\Irr \hol^{(2n)} (\theta  ^{(\rit)}_{G^{(2n)}_{R}}\otimes
\theta  ^{(\rit)}_{G^{(2n)}_{L}})$ of the bisemisheaf $\theta  ^{(\rit)}_{G^{(2n)}_{R}} \otimes
\theta  ^{(\rit)}_{G^{(2n)}_{L}}$ is given by the morphism::
\[ \Irr \hol^{(2n)} _{\theta  ^{(\rit)}_{G\RL}} : \quad \theta  ^{(\rit)}_{G^{(2n)}_{R}}\otimes
\theta  ^{(\rit)}_{G^{(2n)}_{L}} \To f_{\o v}(z^*)\otimes f_v(z)\]
where $f_{\o v}(z^*)\otimes f_v(z)$ is the holomorphic bifunction getting by gluing together and adding the bisections of the desingularized bisemisheaf $\theta  ^{(\rit)}_{G^{(2n)}_{R}}\otimes
\theta  ^{(\rit)}_{G^{(2n)}_{L}}$.

\item $k\times b_i$ global holomorphic representations
\begin{multline*}
 \Irr \hol^{(2n-1)} _{\Fs_{F^{(2n-1)}_{i_{\lambda \RL}}}} : \quad
\Fs_{F^{(2n-1)}_{i_{\lambda _R}}}(\beta_i) \otimes \Fs_{F^{(2n-1)}_{i_{\lambda _L}}}(\beta_i)
\To f_{\o v_{\rm mon}}(z^*_{\beta _i}) \otimes f_{v_{\rm mon}}(z_{\beta _i})\;, \\ \forall\ 1\le i\le k\ , \; 1\le \beta _i\le b_i\end{multline*}
of the monodromy bisemisheaves 
$\Fs_{F^{(2n-1)}_{i_{\lambda _R}}}(\beta_i) \otimes \Fs_{F^{(2n-1)}_{i_{\lambda _L}}}(\beta _i)$ can be stated

where $f_{\o v_{\rm mon}}(z^*_{\beta _i}) \otimes f_{v_{\rm mon}}(z_{\beta _i})$ is the $(i,\beta _i)$-th holomorphic bifunction obtained by gluing and adding the sections of 
$\Fs_{F^{(2n-1)}_{i_{\lambda _R}}}(\beta_i) \otimes \Fs_{F^{(2n-1)}_{i_{\lambda _L}}}(\beta _i)$~.
\Ee
}
\vskip 11pt

\bpr
\Bena
\item If the bisemisheaf $\theta  ^{(\rit)}_{G^{(2n)}_{R}}\otimes
\theta  ^{(\rit)}_{G^{(2n)}_{L}}$ is desingularized, then a $2n$-dimensional irreducible global holomorphic representation can be envisaged for it, as it was done in section 4.2.3, in the sense that a multiple power series development $f_{\o v}(z^*) \otimes f_{v}(z)$ can be associated to it where
\begin{align*}
f_v(z) &= \sum_{j,m_j} c_{j,m_j}(z_1-z_{01})^j\ \cdots\ (z_{2n}-z_{02n})^j\\
\rresp{f_{\o v}(z^*) &= \sum_{j,m_j} c^*_{j,m_j}(z^*_1- z^*_{01})^j\ \cdots\ (z^*_{2n} -z^*_{02n})^j}
\end{align*}
with:
\Bi
\item $z_1,z_{0n},\cdots,z_{2n},z_{02n}$ are complex functions of one real variable;
\item $c_{j,m_j}$ is in one-to-one correspondence with the product of the square roots of the eigenvalues of the $(j,m_j)$-th coset representative $U_{j,m_{j_R}}\times U_{j,m_{j_L}}$ of the product of Hecke operators;
\item the sum $\sum\limits_{j,m_j}$ runs over the conjugacy class representatives of $G^{(2n)}_L$ \resp{$G^{(2n)}_R$} which are glued together.
\Ei

\item Assume that the sections $F^{(2n-1)}_{i_{\lambda _{j_L}}}(\beta_i)$ \resp{$F^{(2n-1)}_{i_{\lambda _{j_R}}}(\beta_i)$} of the monodromy sheaf 
$\Fs_{F^{(2n-1)}_{i_{\lambda _{L}}}}(\beta_i)$ \resp{$\Fs_{F^{(2n-1)}_{i_{\lambda _{R}}}}(\beta_i)$} have ranks given by $r_{F^{(2n-1)}_{i_{\lambda _{j}}}}=(j\centerdot N)^{2n-1}$ according to section 5.1.13.

If these sections, which are complex-valued differentiable functions, are glued together, a holomorphic function (resp. cofunction) given by the multiple power series development:
\begin{align*}
f_{v_{\rm mon}}(z_{\beta _i}) &= \sum_{j,m_j} c_{j_{\beta _i},m_{j_{\beta _i}}}(z_{1_{\beta _i}}-z_{01_{\beta _i}})^j \ \cdots\ (z_{2n-1)_{\beta _i}}-z_{0(2n-1)_{\beta _i}})^j\\
\rresp{f_{\o v_{\rm mon}}(z^*_{\beta _i}) &= \sum_{j,m_j} c^*_{j_{\beta _i},m_{j_{\beta _i}}}(z^*_{1_{\beta _i}}-z^*_{01_{\beta _i}})^j \ \cdots\ (z^*_{2n-1)_{\beta _i}}-z^*_{0(2n-1)_{\beta _i}})^j}\end{align*}
can be associated with them (see, for instance, section 4.2.4).
And, a global holomorphic representation $\Irr \hol^{(2n-1)} _{\Fs_{F^{(2n-1)}_{i_{\lambda \RL}}}}$ of the monodromy bisemisheaf $\Fs_{F^{(2n-1)}_{i_{\lambda _R}}}(\beta _i) \otimes
\Fs_{F^{(2n-1)}_{i_{\lambda _L}}}(\beta _i)$~, as given in this proposition, can be envisaged.\epr
\Ee
\vskip 11pt

\subsubsection[Monodromy $n$-dimensional representations of global Weil groups]{\bbf Monodromy $n$-dimensional representations of global Weil groups}

As the monodromy bisemisheaves $\Fs_{F^{(2n-1)}_{i_{\lambda _R}}}(\beta _i) \otimes
\Fs_{F^{(2n-1)}_{i_{\lambda _L}}}(\beta _i)$ are of algebraic type and are defined above the algebraic bilinear semigroup $\Gnv$~, $(2n-1)$-dimensional irreducible real representations 
$\Irr\Rep^{(2n-1)}_{W^{\rm mon}_{F^+\RL}} (W^{ab}_{F^+_{R_{\rm mon}}}(\beta _i) \times W^{ab}_{F^+_{L_{\rm mon}}}(\beta _i))$ of the bilinear global Weil groups
$(W^{ab}_{F^+_{R_{\rm mon}}}(\beta _i) \times W^{ab}_{F^+_{L_{\rm mon}}}(\beta _i))$ can be introduced for them, similarly as it was done in section 4.2.5.

And, as it was developed in section 4.1.1, the desingularized bisemisheaf $\theta  ^{(\rit)}_{G^{(2n)}_{R}}\otimes
\theta  ^{(\rit)}_{G^{(2n)}_{L}}$ on the algebraic bilinear semigroup $G^{(2n)}(F^+_{\o v}\times F^+_v)$ constitutes a $2n$-dimensional irreducible representation 
$\Irr\Rep^{(2n)}_{W_{F^+\RL}} (W^{ab}_{F^+_{R}}\times W^{ab}_{F^+_{L}})$ of the product, right by left, $W^{ab}_{F^+_{R}}\times W^{ab}_{F^+_{L}}$ of global Weil groups.
\vskip 11pt

\subsubsection{Proposition}

{\em Let $\{\Fs_{F^{(2n-1)}_{i_{\lambda _R}}}(\beta _i) \otimes
\Fs_{F^{(2n-1)}_{i_{\lambda _L}}}(\beta _i)\}_{i,\beta _i}$~, $1\le i\le k$~, $1\le \beta_i\le b_i$~, be the set of $k\times \beta_i$ monodromy bisemisheaves above the desingularized bisemisheaf $\theta  ^{(\rit)}_{G^{(2n)}_{R}}\otimes
\theta  ^{(\rit)}_{G^{(2n)}_{L}}$~.

Then, the following global holomorphic correspondences are:
\Bean
\item $\begin{array}[t]{ccc}
\Irr\Rep^{(2n)}_{W_{F^+\RL}} (W^{ab}_{F^+_{R}}\times W^{ab}_{F^+_{L}}) & \To & \Irr \hol^{(2n)} (\theta  ^{(\rit)}_{G^{(2n)}_{R}}\otimes
\theta  ^{(\rit)}_{G^{(2n)}_{L}})\\[11pt]
 \mbox{\Huge{$\|$}} &&\mbox{\Huge{$\|$}} \\[11pt]
\theta  ^{(\rit)}_{G^{(2n)}_{R}}\otimes
\theta  ^{(\rit)}_{G^{(2n)}_{L}} &\To& f_{\o v}(z^*)\otimes f_v(z)
\end{array}$

for the desingularized bisemisheaf $\theta  ^{(\rit)}_{G^{(2n)}_{R}}\otimes
\theta  ^{(\rit)}_{G^{(2n)}_{L}}$~;

\item \scalebox{.9}{$\begin{array}[t]{ccc}
\Irr\Rep^{(2n-1)}_{W^{\rm mon}_{F^+\RL}} (W^{ab}_{F^+_{R_{\rm mon}}}(\beta _i)\times W^{ab}_{F^+_{L_{\rm mon}}}(\beta i)) &\to&  \Irr \hol^{(2n-1)} (\Fs_{F^{(2n-1)}_{i_{\lambda _R}}}(\beta _i) \otimes
\Fs_{F^{(2n-1)}_{i_{\lambda _L}}}(\beta _i))\\[11pt]
 \mbox{\Huge{$\|$}}&& \mbox{\Huge{$\|$}} \\[11pt]
\Fs_{F^{(2n-1)}_{i_{\lambda _R}}}(\beta _i) \otimes
\Fs_{F^{(2n-1)}_{i_{\lambda _L}}}(\beta _i) &\to&  f_{\o v_{\rm mon}}(z^*_{\beta _i})\otimes f_{v_{\rm mon}}(z_{\beta _i})
\end{array}$}

\mbox{}\hfill $\forall\ i,\beta _i$~,

for the monodromy bisemisheaves.
\Ee
}
\vskip 11pt

\subsubsection[Toroidal compactification of $G^{(2n)}(F^+_{\o v}\times F^+_v)$]{\bbf Toroidal compactification of $G^{(2n)}(F^+_{\o v}\times F^+_v)$}

As in section 4.2.7, a toroidal compactification of the bilinear algebraic semigroup
$G^{(2n)}(F^+_{\o v}\times F^+_v)$ can be envisaged in such a way that its linear conjugacy class representatives $g^{(2n)}_L[j,m_j]$ \resp{$g^{(2n)}_R[j,m_j]$} be transformed into $2n$-dimensional real semitori $g^{(2n)}_{T_L}[j,m_j]$ \resp{$g^{(2n)}_{T_R} [j,m_j]$}.

Let
\begin{align*}
\tau ^{\rm tor}[j,m_j] : \quad g^{(2n)}_L[j,m_j] &\To g^{(2n)}_{T_L}[j,m_j]\\
c_{j,m_j}\ z_j &\To \lambda ^{\half}(2n,j,m_j)\ e^{2\pi ijx}\;, \quad x\in\rit^{2n}\;, \\[11pt]
\rresp{\tau ^{\rm tor}[j,m_j] : \quad g^{(2n)}_R[j,m_j] &\To g^{(2n)}_{T_R}[j,m_j]\\
c^*_{j,m_j}\ z^*_j &\To \lambda ^{\half}(2n,j,m_j)\ e^{-2\pi ijx}}
\end{align*}
be the toroidal deformation of $g^{(2n)}_L[j,m_j]$ \resp{$g^{(2n)}_R[j,m_j]$} transforming
$G^{(2n)}(F^+_{\o v}\times F^+_v)$ into its toroidal equivalent 
$G^{(2n)}(F^{+,T}_{\o v}\times F^{+,T}_v)$~.
\vskip 11pt

\subsubsection{Proposition}

{\em Let
\[ \Irr \hol^{(2n-1)}_{\theta  ^{(\rit)}_{G\RL}} : \quad 
\theta  ^{(\rit)}_{G^{(2n)}_{R}} \otimes \theta  ^{(\rit)}_{G^{(2n)}_{L}} \To f_{\o v}(z^*)\otimes f_v(z)\]
be the global holomorphic representation of the bisemisheaf 
$\theta  ^{(\rit)}_{G^{(2n)}_{R}} \otimes \theta  ^{(\rit)}_{G^{(2n)}_{L}}$ on $G^{(2n)}( {F^+_{\o v}}\times  {F^+_v})$~.

Then, the toroidal compactification $\tau ^{\rm tor}(\Irr \hol^{(2n)}_{\theta  ^{(\rit)}_{G\RL}})$ of the global holomorphic representation of the bisemisheaf 
$\theta  ^{(\rit)}_{G^{(2n)}_{R}} \otimes \theta  ^{(\rit)}_{G^{(2n)}_{L}}$ generates the corresponding elliptic representation according to:
\[\begin{psmatrix}[colsep=0cm,rowsep=2cm]
\Irr \hol^{(2n)}_{\theta  ^{(\rit)}_{G\RL}} : &\hspace*{1cm} & \qquad &
\theta  ^{(\rit)}_{G^{(2n)}_{R}} \otimes \theta  ^{(\rit)}_{G^{(2n)}_{L}} & \To & f_{\o v}(z^*)\otimes f_v(z)\\
\Irr\ELLIP^{(2n)}_{\theta ^{(\rit)}_{G\RL}} : & \qquad & \qquad &
\theta  ^{(\rit)}_{G^{(2n)}_{T_R}} \otimes \theta  ^{(\rit)}_{G^{(2n)}_{T_L}} & \To &
\ELLIP\RL(2n,j,m_j)
\psset{nodesep=.5cm}
\everypsbox{\scriptstyle}
\ncline[arrows=->]{1,1}{2,1}>{\tau ^{\rm tor}(\Irr \hol^{(2n)}_{\theta  ^{(\rit)}_{G\RL}})}
\ncline[arrows=->]{1,4}{2,4}
\ncline[arrows=->]{1,6}{2,6}
\end{psmatrix}\]
where 
\[\ELLIP\RL(2n,j,m_j)=\ELLIP_R(2n,j,m_j) \otimes_{(D)} \ELLIP_L(2n,j,m_j)\;,\]
being the global elliptic representation of the bisemisheaf 
$\theta  ^{(\rit)}_{G^{(2n)}_{T_R}} \otimes \theta  ^{(\rit)}_{G^{(2n)}_{T_L}}$~,  is the product, right by left, of $2n$-dimensional real global elliptic semimodules given by:
\begin{align*}
\ELLIP_L(2n,j,m_j) &= \bigoplus_{j,m_j} \lambda ^{\half}(2n,j,m_j)\ e^{2\pi ijx}\;, \quad x\in (F^+_{v^1})^n\;, \\[15pt]
\text{and} \quad \ELLIP_R(2n,j,m_j) &= \bigoplus_{j,m_j} \lambda ^{\half}(2n,j,m_j)\ e^{-2\pi ijx}.
\end{align*}
}
\vskip 11pt

\bpr This proposition is an adaptation of proposition 4.2.8.\epr
\vskip 11pt

\subsubsection{Proposition}

{\em \Bena
\item A cuspidal representation, given by the elliptic representation $\ELLIP\RL(2n,j,m_j)$~, corresponds to the desingularized bisemisheaf $\theta  ^{(\rit)}_{G^{(2n)}_{R}} \otimes \theta  ^{(\rit)}_{G^{(2n)}_{L}}$~.
\vskip 11pt
\item On the monodromy bisemisheaves $\{\Fs_{F^{(2n-1)}_{i_{\lambda _R}}}(\beta _i) \otimes
\Fs_{F^{(2n-1)}_{i_{\lambda _L}}}(\beta _i)\}_{i,\beta _i}$~, $1\le i\le k$~,\linebreak $1\le \beta_i\le b_i$~, above $\theta  ^{(\rit)}_{G^{(2n)}_{R}} \otimes \theta  ^{(\rit)}_{G^{(2n)}_{L}}$~, no cuspidal representation of the elliptic type can be found, except if surgeries are performed.
\Ee}
\vskip 11pt

\bpr According to section 5.1.5, the non-singular fibres (or sections) $F^{(2n-1)}_{i_{\lambda _{j_L}}}$ \resp{$F^{(2n-1)}_{i_{\lambda _{j_R}}}}$ of the monodromy sheaf $\Fs_{F^{(2n-1)}_{i_{\lambda _L}}}(\beta _i)$ \resp{$\Fs_{F^{(2n-1)}_{i_{\lambda _R}}}(\beta _i)$} are diffeomorphic to the spheres $TS^{(2n-1)}_{i_{L_j}}$ \resp{$TS^{(2n-1)}_{i_{R_j}}$}. So, no elliptic representation can be found for the monodromy sheaf $\Fs_{F^{(2n-1)}_{i_{\lambda _L}}}(\beta _i)$ \resp{$\Fs_{F^{(2n-1)}_{i_{\lambda _R}}}(\beta _i)$} since $TS^{(2n-1)}_{i_{L_j}}$ \resp{$TS^{(2n-1)}_{i_{R_j}}$} cannot be transformed bijectively into a $(2n-1)$-dimensional real semitorus (being able to constitute an equivalence class representative of a global elliptic representation according to proposition 5.1.18).\epr
\vskip 11pt

\subsubsection{Proposition}

{\em On the desingularized bisemisheaf $\theta  ^{(\rit)}_{G^{(2n)}_{R}} \otimes \theta  ^{(\rit)}_{G^{(2n)}_{L}}$ and its monodromy bisemisheaves\linebreak $\{\Fs_{F^{(2n-1)}_{i_{\lambda _R}}}(\beta _i)\otimes \Fs_{F^{(2n-1)}_{i_{\lambda _L}}}(\beta _i)\}_{i,\beta _i}$~, the only irreducible  global correspondences of\linebreak Langlands:
\[ \begin{array}{ccc}
\Irr\Rep^{(2n)}_{W_{F^+\RL}} (W^{ab}_{F^+_{R}}\times W^{ab}_{F^+_{L}}) &\To &
\Irr\ELLIP (\theta ^{(\rit)}_{G^{(2n)}_{T_R}} \otimes \theta ^{(\rit)}_{G^{(2n)}_{T_L}})\\[11pt]
 \mbox{\Huge{$\|$}} &&  \mbox{\Huge{$\|$}}\\[11pt]
\theta ^{(\rit)}_{G^{(2n)}_R} \otimes \theta ^{(\rit)}_{G^{(2n)}_L} &\To & \ELLIP\RL(2n,j,m_j)
\end{array}\]
exists.
}
\vskip 11pt

\bpr \Bena
\item The above mentioned Langlands correspondence results from the toroidal compactification of the irreducible holomorphic correspondence  $\Irr\hol^{(2n)} (\theta ^{(\rit)}_{G^{(2n)}_R} \otimes \theta ^{(\rit)}_{G^{(2n)}_L})$~, introduced in proposition 5.1.16:

\[\begin{psmatrix}[colsep=1.5cm,rowsep=1.5cm]
\Irr\Rep^{(2n)}_{W_{F^+\RL}} (W^{ab}_{F^+_{R}}\times W^{ab}_{F^+_{L}})
& \Irr\hol^{(2n)} (\theta ^{(\rit)}_{G^{(2n)}_R} \otimes \theta ^{(\rit)}_{G^{(2n)}_L})\\
\theta ^{(\rit)}_{G^{(2n)}_R} \otimes \theta ^{(\rit)}_{G^{(2n)}_L} & 
f_{\o v}(z^*) \otimes f_{v}(z)\\
& \Irr\ELLIP (\theta ^{(\rit)}_{G^{(2n)}_R} \otimes \theta ^{(\rit)}_{G^{(2n)}_L}) \\ 
& \ELLIP\RL(2n,j,m_j)
\psset{nodesep=.3cm}
\everypsbox{\scriptstyle}
\ncline[arrows=->]{1,1}{1,2}
\ncline[arrows=->]{1,1}{3,2}
\ncline[arrows=->]{2,2}{3,2}>{\tau ^{\rm tor}(\Irr \hol^{(2n-1)}_{\theta  ^{(\rit)}_{G\RL}})}
\ncline[doubleline=true]{1,1}{2,1}
\ncline[doubleline=true]{1,2}{2,2}
\ncline[doubleline=true]{3,2}{4,2}
\end{psmatrix}\]
 \item According to proposition 5.1.19, no Langlands correspondence exists for the monodromy bisemisheaves because no bijection can be found between the non-singular fibres, diffeomorphic to $TS^{2n-1}$~, and $(2n-1)$-dimensional real semitori.
 \epr
 \Ee
 \vskip 11pt

\subsection{The monodromy of isolated singularities on reducible complex bisemisheaves}

\subsubsection{Singular reducible bisemisheaves}

In chapter 4 of \cite{Pie1}, the possible reducibilities of the representation $\Rep (\GL_{2n}(F_{\o\omega }\times F_\omega ))$ of the bilinear algebraic semigroup $\GL_{2n}(F_{\o\omega }\times F_\omega )$ were introduced.  They are of three types:
\Bean
\item partially reducible if:
\[\Rep (\GL_{2n=2n_1+\cdots+2n_s}(F_{\o\omega }\times F_\omega ))
= \mathop{\boxplus}\limits^{2n_s}_{2n_\ell =2n_1} \Rep (\GL_{2n_\ell }(F_{\o\omega }\times F_\omega ))\]
for any partition $2n=2n_1+\cdots+ 2n_\ell +\cdots+2n_s$ of $2n$~;

\item orthogonally completely reducible if:
\[\Rep (\GL_{2n=2_1+\cdots+ 2_n}(F_{\o\omega }\times F_\omega ))
= \mathop{\boxplus}\limits^{n}_{\ell =1} \Rep (\GL_{2_\ell }(F_{\o\omega }\times F_\omega ))\]

\item non orthogonally completely reducible if:
\begin{align*}
\Rep (\GL_{2n\RL}(F_{\o\omega }\times F_\omega ))
&= \mathop{\boxplus}\limits^{2n}_{2\ell_R=2\ell _L =1} \Rep (\GL_{2_{\ell\RL} }(F_{\o\omega }\times F_\omega ))\\
&= \mathop{\boxplus}\limits_{2k_R\neq \ell _L} \Rep (T^t_{2_{k_R}}(F_{\o\omega })\times T_{2_{\ell _L}}(F_\omega ))\;.\end{align*}
\Ee

The  analytic representation spaces over these reducible bilinear algebraic semigroups are respectively the following bisemisheaves of differentiable bifunctions:
\Bean
\item the partially reducible bisemisheaf $\theta ^{(\cit)}_{\GL_{2n=2n_1+\cdots+ 2n_\ell +\cdots+2n_s}(F_{\o\omega }\times F_\omega )}$ over the bilinear algebraic semigroup $\GL_{2n=2n_1+\cdots+ 2n_\ell +\cdots+2n_s}(F_{\o\omega }\times F_\omega )$~;

\item the orthogonal completely reducible bisemisheaf $\theta^{(\cit)}_ {\GL_{2n=2_1+\cdots+ 2_n}(F_{\o\omega }\times F_\omega )}$ over the bilinear algebraic semigroup
$\GL_{2n=2_1+\cdots+ 2_\ell +\cdots+2_r}(F_{\o\omega }\times F_\omega )$~;

\item the non orthogonal completely reducible bisemisheaf
$\theta^{(\cit)}_{\GL_{2n\RL}(F_{\o\omega }\times F_\omega )}$ over the bilinear algebraic semigroup $\GL_{2n\RL}(F_{\o\omega }\times F_\omega )$~.
\Ee

Being concerned by the monodromy group action on these (``singular'') bisemisheaves, the only relevant reducible bisemisheaf is the orthogonal completely reducible bisemisheaf.  Indeed,  the monodromy group action on the partially reducible bisemisheaf decomposing into:
\[ \theta ^{(\cit)}_{\GL_{2n=2n_1+\cdots+ 2n_s}(F_{\o\omega }\times F_\omega )}
= \mathop{\boxplus}\limits^{n_s}_{n_\ell =n_1} \theta^{(\cit)}_{\GL_{2n_\ell }(F_{\o\omega }\times F_\omega )}\]
can amount to the general irreducible case treated in section 5.1 for $n_\ell \ge 2$~, $n_1\le n_\ell \le n_s$~.

On the other hand, as the non orthogonal completely reducible bisemisheaves
\[ \theta ^{(\cit)}_{T^t_{2_{k_R}}}(F_{\o\omega }) \otimes \theta ^{(\cit)}_{T^t_{2_{\ell_L}}}(F_{\omega })\in \theta ^{(\cit)}_{\GL_{2n\RL}}(F_{\o\omega }\times F_\omega )\;, \quad 1\le k,\ell\le n\;, \]
on the off-diagonal algebraic linear semigroups $T^t_{2_{k_R}}(F_{\o\omega })$  and 
$T^t_{2_{\ell_L}}(F_{\omega })$ are generated from the orthogonal ones 
$\theta^{(\cit)}_{T^t_{2_{\ell_R}}}(F_{\o\omega }) \otimes \theta^{(\cit)}_{T^t_{2_{\ell_L}}}(F_{\omega })$~, their monodromy group actions are not really pertinent.

So, the only relevant reducible bisemisheaf from the monodromy point of view is the orthogonal completely reducible bisemisheaf 
\[\theta^{(\cit)}_{\GL_{2n=2_1+\cdots+2_n}}(F_{\o\omega }\times F_\omega )
=\mathop{\boxplus}\limits^n_{\ell=1} \theta^{(\cit)}_{\GL_{2_\ell}}(F_{\o\omega }\times F_\omega )\]
 whose irreducible elements $\theta^{(\cit)}_{\GL_{2_\ell}}(F_{\o\omega }\times F_\omega ))$~, $1\le \ell\le n$~, being able to generate monodromy groups, will be studied in this section.
 \vskip 11pt
 
\subsubsection[Critical level sets of the bisemisheaf $\theta^{(\cit)}_{\GL_{2}(F_{\o\omega }\times F_\omega )}$]{\bbf Critical level sets of the bisemisheaf $\theta^{(\cit)}_{\GL_{2}(F_{\o\omega }\times F_\omega )}$}

Let
\[
\theta^{(\cit)}_{\GL_{2}(F_{\o\omega }\times F_\omega )} \; \equiv \;
\theta^{(\cit)}_{\GL_{2}(F_{\o\omega })} \otimes \theta^{(\cit)}_{\GL_{2}( F_\omega )}
\; (\ \equiv \;
\theta^{(\cit)}_{T^t_{2}(F_{\o\omega })} \otimes \theta^{(\cit)}_{T_{2}( F_\omega )}\ )\]
be a bisemisheaf of complex-valued differentiable bifunctions $\phi ^{(2)}_{G^{(\cit)}_{g_{j_R}}}(z_{g_{j_R}})\otimes \phi ^{(2)}_{G^{(\cit)}_{g_{j_L}}}(z_{g_{j_L}})$~, $1\le j\le r\le \infty $~, over the conjugacy class representatives $g^{(2)}\RL(j,m_j]$ of the complex algebraic bilinear semigroup $G^{(2)}(F_{\o\omega }\times F_\omega )$~.

Assume that, on a domain $U^{(2)}_{j_R} \times U^{(2)}_{j_L} \subset g^{(2)}\RL[j,m_j]$~, each bifunction $\phi ^{(2)}_{G^{(\cit)}_{g_{j_R}}}(z_{g_{j_R}})\otimes \phi ^{(2)}_{G^{(\cit)}_{g_{j_L}}}(z_{g_{j_L}})$ is locally a Morse bifunction described by:
\[\phi ^{(2)}_{G^{(\cit)}_{g_{j_R}}}(U^{(2)}_{j_R}) \otimes \phi ^{(2)}_{G^{(\cit)}_{g_{j_L}}}(U^{(2)}_{j_L}) = z^2_{j_1}+z^2_{j_2}\;, \quad (z_1,z_2)\in\cit^2\;.\]

At $z_{j_1}=z_{j_2}=0$~, the bifunction $\phi ^{(2)}_{G^{(\cit)}_{g_{j_R}}}(U^{(2)}_{j_R}) \otimes \phi ^{(2)}_{G^{(\cit)}_{g_{j_L}}}(U^{(2)}_{j_L})$ has a non-degenerate singularity.

The critical level set of $\phi ^{(2)}_{G^{(\cit)}_{g_{j_R}}}(U^{(2)}_{j_R}) \otimes \phi ^{(2)}_{G^{(\cit)}_{g_{j_L}}}(U^{(2)}_{j_L})$ is the singular fibre $F^{(1)}_{0_{j_R}} \times F^{(1)}_{0_{j_L}}$ given by:
\[\phi ^{(2)}_{G^{(\cit)}_{g_{j_R}}}(U^{(2)}_{j_R}) \otimes \phi ^{(2)}_{G^{(\cit)}_{g_{j_L}}}(U^{(2)}_{j_L}) = z^2_{j_1}+z^2_{j_2}=0\;;\]
it consists in two complex lines intersecting at 0 \cite{H-Z}.
\vskip 11pt

\subsubsection{Proposition}

{\em Let 
\[ F^{(1)}_{0_{j_R}} \times F^{(1)}_{0_{j_L}}=z^2_{j_1}+z^2_{j_2}=0\;, \quad \forall\ j\;, \quad 1\le j\le r\le \infty \;, \]
be the singular bifibre of the $j$-th bisection $\phi ^{(2)}_{G^{(\cit)}_{g_{j_R}}}(U^{(2)}_{j_R}) \otimes \phi ^{(2)}_{G^{(\cit)}_{g_{j_L}}}(U^{(2)}_{j_L})$ of the bisemisheaf $\theta^{(\cit)}_{\GL_2( {F_{\o \omega }}\times  {F_\omega })}$~.

Then, we have that:
\Bena
\item the corresponding non singular bifibres $F^{(1)}_{\lambda _{j_R}}(t) \times F^{(1)}_{\lambda _{j_L}}(t)$ are diffeomorphic to the  product, right by left, $T^2_{\lambda _{j_R}}(t) \times T^2_{\lambda _{j_L}}(t) $ of two semitori.

\item the homology group $H_1(F^{(1)}_{\lambda _{j_L}};\ZZ)\simeq \ZZ$ \resp{$H_1(F^{(1)}_{\lambda _{j_R}};\ZZ)\simeq \ZZ$} of the semitorus $T^2_{\lambda _{j_L}}$ \resp{$T^2_{\lambda _{j_R}}$} is generated by the upper (resp. lower) semicircle $\Delta ^{(1)}_{L_j}$ \resp{$\Delta ^{(1)}_{R_j}$} on $T^2_{\lambda _{j_L}}$ \resp{$T^2_{\lambda _{j_R}}$} in such a way that when its radius tends to the unity, the semicircle shrinks to the singularity and is then called the vanishing semicycle characterized by a rank $r_{\Delta ^{(1)}_{j}}=N$~.

\item The covanishing semicycle 
$\nabla  ^{(1)}_{L_j}$ \resp{$\nabla  ^{(1)}_{R_j}$} is a line on $T^2_{\lambda _{j_L}}$ \resp{$T^2_{\lambda _{j_R}}$} perpendicular to  the vanishing semicycle $\Delta ^{(1)}_{L_j}$ \resp{$\Delta ^{(1)}_{R_j}$}  and is characterized by a rank $r_{F^{(1)}_{\lambda j}}\approx j\centerdot N$~.
\Ee}
\newpage

\bpr
\Bena
\item The $j$-th critical level sets are the non-singular bifibres $F^{(1)}_{\lambda _{j_R}}(t) \times F^{(1)}_{\lambda _{j_L}}(t)$ described by the equations:
\[ z^{2}_{j_1}(t) + z^2_{j_2}(t)=\lambda_j (t)\;, \quad \lambda \neq 0\;; \]
they are diffeomorphic to the cylinder $S^1\times \rit^1$~.

Indeed, the Riemann surface of the function $z_{j_2}(t)=\sqrt{\lambda _j(t)-z^2_{j_1}}$ is formed from two copies of the complex $z_{j_1}$-plane glued together along the cut $(-\lambda _j(t),+\lambda _j(t))$ \cite{H-Z}. Each copy of the cut plane is homeomorphic to half a cylinder (and a 2-dimensional semitorus).  The line of the cut is a circle on the cylinder encircling the critical value and given by the equation:
\[\lambda _j(t) = r_{\lambda _j(t)}\ e^{2\pi ijt}\;, \quad 0\le t\le 1\;.\]
As $t$ increases, both branch points $z_{j_1}=\pm\sqrt{\lambda _j(t)}=(\pm\sqrt{r_{\lambda _j(t)}\ e^{\pi ijt}})$ move around $z_{j_1}=0$ in the positive direction. As $t$ varies from 0 to 1~, each of these points performs a revolution and changes place with the other.

Thus, as $\lambda _j(t)$ encircles the singularity, a corresponding series of pairs $\{T^2_{\lambda _{j_R}}(t),T^2_{\lambda _{j_L}}(t)\}_t$ of two-dimensional semitori are generated.

\item The circle on the cylinder encircling the critical value is the vanishing cycle given by the equation:
\[\lambda ^{(V)}_j(t)\simeq e^{2\pi ijt}\;, \]
i.e. when the radius $r_{\lambda _j(t)}\simeq 1$~.

Indeed, in this case, the \lr vanishing semicycle $\Delta ^{(1)}_{L_j}$ \resp{$\Delta ^{(1)}_{R_j}$}, given by the equation:
\[\lambda _j(t)\simeq e^{2\pi ijt}\;, \]
restricted to the upper (resp. lower) half plane, corresponds to a function $\phi ^{(1)}_{P_{j_L}}(x_{p_{j_L}})$ \resp{$\phi ^{(1)}_{P_{j_R}}(x_{p_{j_R}})$} on the \lr conjugacy class representative $P^{(2)}(F^+_{v^1_{j,m_j}})$ \resp{$P^{(2)}(F^+_{\o v^1_{j,m_j}})$} of the linear parabolic subgroup $P^{(2)}(F^+_{v^1})$ \resp{$P^{(2)}(F^+_{\o v^1})$} according to section 5.1.5.

And, thus, the rank of $\Delta ^{(1)}_{L_j}$ and of $\Delta ^{(1)}_{R_j}$ is given by $r_{\Delta ^{(1)}_{j}}=N$~.

\item As the covanishing semicycle $\nabla  ^{(1)}_{L_j}$ \resp{$\nabla  ^{(1)}_{R_j}$} is a semicircle on $T^2_{\lambda _{j_L}}$ \resp{$T^2_{\lambda _{j_R}}$} perpendicular to 
$\Delta ^{(1)}_{L_j}$ \resp{$\Delta ^{(1)}_{R_j}$}, it will be assumed to have a rank 
$r_{\nabla _{j^{(1)}}}\simeq j\centerdot N$ since it is defined on the $j$-th conjugacy class representative of the linear algebraic semigroup $G^{(2)}(F_\omega )$ \resp{$G^{(2)}(F_{\o\omega} )$} (see, for example, proposition 3.2.2. to illustrate this point).\epr
\Ee
\vskip 11pt

\subsubsection{Proposition}

{\em Let $(\phi ^{(2)}_{G^{(\cit)}_{g_{j_R}}}(U^{(2)}_{j_R}) \otimes \phi ^{(2)}_{G^{(\cit)}_{g_{j_L}}}(U^{(2)}_{j_L}), F^{(1)}_{0_{j_R}} \times F^{(1)}_{0_{j_L}} ,
F^{(1)}_{\lambda _{j_R}} (t)\times F^{(1)}_{\lambda _{j_L}} (t),\Delta ^{(1)}_{R_j} \times
\Delta ^{(1)}_{L_j})$ be the 4-th bituple introduced in proposition 5.2.3.

Then the mapping
\[ h^{(1)}_{\gamma _{j\RL}}: \quad 
F^{(1)}_{\lambda _{j_R}} (t)\times F^{(1)}_{\lambda _{j_L}} (t) \To
F^{(1)}_{\lambda _{j_R}} (t)\times F^{(1)}_{\lambda _{j_L}} (t)\]
of the non singular bifibre into itself is the monodromy of the product, right by left, 
$\Delta ^{(1)}_{R_j} \times
\Delta ^{(1)}_{L_j}$ of the vanishing semicycles realized by the conjugacy action of the $j$-th conjugacy class representative of the bilinear algebraic semigroup $G^{(2)}( {F^+_{\o v}} \times  {F^+_{v }})$ on the corresponding conjugacy class representative of the bilinear parabolic subsemigroup $P^{(2)}( {F_{\o v ^1}} \times  {F_{v^1}})$~.
}
\vskip 11pt

\bpr \Bena
\item This proposition is a particular case of the one which was treated in proposition 5.1.6 in the sense that the monodromy $h^{(1)}_{\gamma _{\RL}}$ is associated with the injective mapping:
\[ I_{\Delta ^{(1)}_{\lambda _{j\RL}}\to F^{(1)}_{j\RL}}: \quad
\Delta ^{(1)}_{R_j} \times \Delta ^{(1)}_{L_j} \To
F^{(1)}_{\lambda _{j_R}}(t) \times F^{(1)}_{\lambda _{j_L}}(t) \]
inflating the \lr vanishing semicycle $\Delta ^{(1)}_{L_j}$ \resp{$\Delta ^{(1)}_{R_j}$}, characterized by a rank $r_{\Delta ^{(1)}_j}=N$~, into the \lr non singular fibre $F^{(1)}_{\lambda _{j_L}}(t)$ \resp{$F^{(1)}_{\lambda _{j_R}}(t)$}, characterized by a rank $r_{F^{(1)}_{\lambda _j}}=m_j\ (j\centerdot N)$ where $m_j$ is the number of $1D$-fibres perpendicular to $\Delta ^{(1)}_{L_j}$ \resp{$\Delta ^{(1)}_{R_j}$} into $F^{(1)}_{\lambda _{j_L}}(t)$ \resp{$F^{(1)}_{\lambda _{j_R}}(t)$} diffeomorphic to the $2D$-semitorus 
$T^2_{\lambda _{j_L}}(t)$ \resp{$T^2_{\lambda _{j_R}}(t)$}.

\item The inflation action of $I_{\Delta ^{(1)}_{j\RL}\to F^{(1)}_{\lambda _{j\RL}}}$ on $\Delta ^{(1)}_{R_j} \times \Delta ^{(1)}_{L_j}$ corresponds to the conjugation of the $j$-th conjugacy class representative of $G^{(2)}(F^+_{\o v } \times F^+_{v })$ on the corresponding conjugacy class representative of $P^{(2)}(F^+_{\o v ^1} \times F^+_{v^1 })$~.\epr
\Ee
\vskip 11pt

\subsubsection{Number of non singular bifibres}

Assume that each bisection $\phi ^{(2)}_{G^{(\cit)}_{g_{j_R}}}(U^{(2)}_{j_R}) \otimes \phi ^{(2)}_{G^{(\cit)}_{g_{j_L}}}(U^{(2)}_{j_L})$ of $\theta^{(\cit)}_{\GL_2( {F_{\o \omega }}\times  {F_\omega })}$ is endowed with the same singular bifibre $F^{(1)}_{0 _{j_R}} \times F^{(1)}_{0 _{j_L}}=z^2_{j_1}+z^2_{j_2}=0$~.

The number of non singular bifibres $(F^{(1)}_{\lambda _{j_R}} (t)\times F^{(1)}_{\lambda _{j_L}}(t))$~, corresponding to $F^{(1)}_{0 _{j_R}} \times F^{(1)}_{0 _{j_L}}$~, depends on the  expanding phase responsible for the monodromy $h^{(1)}_{\gamma _{j\RL}}$ (see section 5.1.2).  Indeed, it corresponds to this expanding phase the inverse mapping
\[ r^{-1}_{F^{(1)}_{\lambda _{j\RL}}\to F^{(1)}_{0 _{j\RL}}} :
\quad F^{(1)}_{0 _{j_R}}\times F^{(1)}_{0 _{j_L}}\To
F^{(1)}_{\lambda  _{j_R}}(t) \times F^{(1)}_{\lambda  _{j_L}}(t)\] of the retraction of the monodromy (see section 5.1.9):
\[ r_{F^{(1)}_{\lambda _{j\RL}}\to F^{(1)}_{0 _{j\RL}}} :
\quad  F^{(1)}_{\lambda  _{j_R}}(t) \times F^{(1)}_{\lambda  _{j_L}}(t)\To F^{(1)}_{0 _{j_R}}\times F^{(1)}_{0 _{j_L}}\;.\]
Let $\beta $ be the number of non singular bifibres $(F^{(1)}_{\lambda  _{j_R}}(t) \times F^{(1)}_{\lambda  _{j_L}}(t))$ above each bisection of $\theta ^{(\cit)}_{\GL_2(F_{\o\omega }\times F_\omega )}$~.
\vskip 11pt

\subsubsection{Proposition}

{\em If each bisection of the bisemisheaf $\theta ^{(\cit)}_{\GL_2(F_{\o\omega }\times F_\omega )}$ is endowed with the same singular bifibre $F^{(1)}_{0 _{j_R}}\times F^{(1)}_{0 _{j_L}}= z^2_{j_1}+z^2_{j_2}=0$~, then a set of $\beta $ bisemisheaves
$\{\theta ^{(\cit)\rm mon}_{\GL_2(F_{\o\omega }\times F_\omega )}(b)\}^\beta _{b=1}$~, isomorphic to the desingularized bisemisheaf $\theta ^{(\cit)}_{\GL_2(F_{\o\omega }\times F_\omega )}$~, can be generated by monodromy if $\beta $ is the number of non singular bifibres above each bisection of $\theta ^{(\cit)}_{\GL_2(F_{\o\omega }\times F_\omega )}$~.
}

\bpr As there are $\beta $ identical non singular bifibres $\{F^{(1)}_{\lambda  _{j_R}}(b) \times F^{(1)}_{\lambda  _{j_L}}(b)\}^\beta _{b=1}$ (assuming the one-to-one correspondence $t\leftrightarrow b$~) above each bisection of $\theta ^{(\cit)}_{\GL_2(F_{\o\omega }\times F_\omega )}$~, $\beta $ bisemisheaves $\theta ^{(\cit)\rm mon}_{\GL_2(F_{\o\omega }\times F_\omega )}(b)$~, $1\le b\le \beta $~, can be built from these nonsingular bifibres in such a way that the sections of $\theta ^{(\cit)\rm mon}_{\GL_2(F_{\o\omega }\times F_\omega )}(b)$ are in one-to-one correspondence with the sections of $\theta ^{(\cit)}_{\GL_2(F_{\o\omega }\times F_\omega )}$~.  Furthermore, the sections of $\theta ^{(\cit)}_{\GL_2(F_{\o\omega }\times F_\omega )}$~, and the corresponding sections of the monodromy bisemisheaves $\theta ^{(\cit)\rm mon}_{\GL_2(F_{\o\omega }\times F_\omega )}(b)$ are isomorphic according to proposition 5.2.4.

So, the monodromy bisemisheaves $\theta ^{(\cit)\rm mon}_{\GL_2(F_{\o\omega }\times F_\omega )}(b)$~, $1\le b\le \beta $~, are isomorphic to (or ``copies of'') the original desingularized bisemisheaf
$\theta ^{(\cit)}_{\GL_2(F_{\o\omega }\times F_\omega )}$~.\epr
\vskip 11pt

\subsubsection{Proposition}

{\em Let $\theta ^{(\cit)}_{\GL_2(F_{\o\omega }\times F_\omega )}$ be a complex bisemisheaf, whose bisections are endowed with the same singular bifibres $F^{(1)}_{0 _{j_R}}\times F^{(1)}_{0 _{j_L}}= z^2_{j_1}+z^2_{j_2}=0$~, $\forall\ j$~, $1\le j\le r\le \infty $~, and let 
$\{\theta ^{(\cit)\rm mon}_{\GL_2(F_{\o\omega }\times F_\omega )}(b)\}^\beta _{b=1}$ be the set of $\beta $ monodromy bisemisheaves.

Then, it results that:
\Bena
\item a global holomorphic representation:
\[ \Irr \hol^{(1)}_{\theta  ^{(\cit)}_{G\RL}}: \quad
\theta ^{(\cit)}_{\GL_2(F_{\o\omega }\times F_\omega )} \To f_{\o\omega }(z^*_m)\otimes f_\omega (z_m)\]
corresponds to the desingularized ground bisemisheaf $\theta ^{(\cit)}_{\GL_2(F_{\o\omega }\times F_\omega )}$~;

\item $\beta $ global holomorphic representations:
\[ \Irr \hol^{(1)}_{\theta  ^{(\cit)\rm mon}_{G\RL}}: \quad
\theta ^{(\cit)\rm mon}_{\GL_2(F_{\o\omega }\times F_\omega )}(b) \To f_{\o\omega }(z^*_{m_b})\otimes f_\omega (z_{m_b})\;,\quad 1\le b\le\beta \;,\]
can be associated with the monodromy bisemisheaves $\theta ^{(\cit)\rm mon}_{\GL_2(F_{\o\omega }\times F_\omega )}(b)$~.
\Ee}
\vskip 11pt

\bpr\Bena
\item If $\theta ^{(\cit)}_{\GL_2(F_{\o\omega }\times F_\omega )}$ is desingularized and if its \lr linear sections are glued together, a global holomorphic representation can be given to it by the holomorphic bifunction $f_{\o\omega }(z^*_m)\otimes f_\omega (z_m)$ where:
\begin{align*}
f_{\omega }(z_m) &= \sum_{j,m_j} c_{j,m_j} (z_m-z_{m_0})^j\\
\rresp{f_{\o\omega }(z^*_m) &= \sum_{j,m_j} c^*_{j,m_j} (z^*_m-z^*{m_0})^j}\end{align*}
with:
\Bi
\item $z_m,z_{m_0}$ \resp{$z^*_m,z^*_{m_0}$} complex (resp. conjugate complex) variables;
\item $c_{j,m_j}$ \resp{$c^*_{j,m_j}$} coefficients (see proposition 5.1.14).
\Ei

\item $\beta$ global holomorphic representations of monodromy type are given by the holomorphic bifunctions $f_{\o\omega }(z^*_{m_b})\otimes f_\omega (z_{_b}m)$~, $1\le b\le\beta $~, in such a way that they are equivalent to the ground holomorphic bifunction $f_{\o\omega }(z^*_m)\otimes f_\omega (z_m)$ taking into account the proposition 5.2.6.\epr
\Ee
\vskip 11pt

\subsubsection{Proposition}

{\em Let $\{\theta ^{(\cit)\rm mon}_{\GL_2(F_{\o\omega }\times F_\omega )}(b)\}^\beta _{b=1}$ be the $\beta $ monodromy bisemisheaves above the desingularized ground bisemisheaf 
$\theta ^{(\cit)}_{\GL_2(F_{\o\omega }\times F_\omega )}$~.

Then, the global holomorphic correspondences are the following:
\Be
\item $ \begin{array}[t]{ccc}
\Irr\Rep^{(1)}_{W_{F\RL}} (W^{ab}_{F_{R}}\times W^{ab}_{F_{L}}) &\To &
\Irr\hol^{(1)} (\theta ^{(\cit)}_{\GL_2(F_{\o\omega }\times F_\omega )})\\[11pt]
\mbox{\Huge{$\|$}} && \mbox{\Huge{$\|$}}\\[11pt]
\theta ^{(\cit)}_{\GL_2(F_{\o\omega }\times F_\omega )} &\To & f_{\o\omega }(z^*_m)\otimes f_\omega (z_m)
\end{array}$

where $\Irr\Rep^{(1)}_{W_{F\RL}} (W^{ab}_{F_{R}}\times W^{ab}_{F_{L}})$ is the irreducible complex representation of the bilinear global Weil group $(W^{ab}_{F_{R}}\times W^{ab}_{F_{L}})$ given  by the ground bisemisheaf $\theta ^{(\cit)}_{\GL_2(F_{\o\omega }\times F_\omega )}$~.

\item $ \begin{array}[t]{ccc}
\Irr\Rep^{(1)}_{W^{\rm mon}_{F\RL}} (W^{ab}_{F_{R_{\rm mon}}}(b)\times W^{ab}_{F_{L_{\rm mon}}}(b)) &\To &
\Irr\hol^{(1)} (\theta ^{(\cit)\rm mon}_{\GL_2( {F_{\o \omega }}\times  {F_\omega })}(b))\\[11pt]
\mbox{\Huge{$\|$}} && \mbox{\Huge{$\|$}}\\[11pt]
\theta ^{(\cit)\rm mon}_{F_{\o\omega }\times F_\omega }(b) &\To & f_{\o\omega }(z^*_{m_b})\otimes f_\omega (z_{m_b})
\end{array}$

\mbox{}\hfill $1\le b\le \beta \;,$

for the $\beta $ monodromy bisemisheaves $\theta ^{(\cit)\rm mon}_{\GL_2(F_{\o\omega }\times F_\omega )}(b)$~.
\Ee}
\vskip 11pt

\bpr This proposition is an adaptation of proposition 5.1.16 to the products, right by left, of $1D$-complex (semi)sheaves.\epr
\vskip 11pt

\subsubsection{Toroidal compactification}

The ground bisemisheaf $\theta ^{(\cit)}_{\GL_2(F_{\o\omega }\times F_\omega )}$ and the $\beta $ monodromy bisemisheaves $\theta ^{(\cit)\rm mon}_{\GL_2(F_{\o\omega }\times F_\omega )}(b)$ are defined over the bilinear algebraic semigroup $\GL_2(F_{\o\omega }\times F_\omega )$~.  So, a toroidal compactification of the linear conjugacy class representatives $g^{(2)}_L[j,m_j]$ of $\GL_2(F_\omega )$ \resp{$g^{(2)}_R[j,m_j]$ of $\GL_2(F_{\o\omega })$} can be realized by the mappings:
\begin{align*}
\tau_{\cit} ^{\rm tor}[j,m_j] : \quad g^{(2)}_L[j,m_j] &\To g^{(2)}_{T_L}[j,m_j]\\
c_{j,m_j}\ z_m^j &\To \lambda ^{\half}(2,j,m_j)\ e^{2\pi ijz_m}\;, \quad z_m\in\cit\;, \\[11pt]
\rresp{\tau_{\cit} ^{\rm tor}[j,m_j] : \quad g^{(2)}_R[j,m_j] &\To g^{(2)}_{T_R}[j,m_j]\\
c^*_{j,m_j}\ z^{*j}_m &\To \lambda ^{\half}(2,j,m_j)\ e^{-2\pi ijz_m}}\quad \forall\ 1\le j\le r\le\infty \;,
\end{align*}
where $g^{(2)}_{T_L}[j,m_j]=\lambda ^{\half}(2,j,m_j)\ e^{2\pi ijz_m}$ \resp{$g^{(2)}_{T_R}[j,m_j]=\lambda ^{\half}(2,j,m_j)\ e^{-2\pi ijz_m}$} is a two-dimensional real semitorus localized in the upper (resp. lower) half space.
\vskip 11pt

\subsubsection{Proposition}

{\em
\Bena
\item A cuspidal representation, given right by the product $\EIS\RL(1,j,m_j)=\EIS_R(1,j,m_j)\times  \EIS_R(1,j,m_j)$~, of the (truncated) Fourier development of a normalized right cusp form by its left equivalent, can be associated with the desingularized ground bisemisheaf $\theta ^{(\cit)}_{\GL_2(F_{\o\omega }\times F_\omega )}$~.

\item Similarly, a cuspidal representation given by $\EIS^{\rm mon}\RL(1,j,m_j)$ corresponds to each monodromy bisemisheaf $\theta ^{(\cit)\rm mon}_{\GL_2(F_{\o\omega }\times F_\omega )}$ on the bilinear algebraic semigroup $\GL_2(F_{\o\omega }\times F_\omega )$ compactified toroidally.
\Ee}
\vskip 11pt

\bpr
\Bena
\item The toroidal compactification $\tau_{\cit} ^{\rm tor}(\Irr\hol^{(1)} (\theta ^{(\cit)}_{\GL_2(F_{\o\omega }\times F_\omega )}))$ of the global holomorphic representation of the ground bisemisheaf $\theta ^{(\cit)}_{\GL_2(F_{\o\omega }\times F_\omega )}$ generates the corresponding cuspidal representation $\Irr\cusp^{(1)} _{\theta ^{(\cit)}_{\GL_{2\RL}}}$ according to:
\[\begin{psmatrix}[colsep=0cm,rowsep=2cm]
\Irr \hol^{(1)}_{\theta  ^{(\cit)}_{G\RL}} : &\hspace*{2cm} & \qquad &
\theta  ^{(\cit)}_{\GL_2(F_{\o\omega }\times F_\omega )} & \To & f_{\o \omega }(z^*_m)\otimes f_\omega (z_m)\\
\Irr\cusp^{(1)}_{\theta ^{(\cit)}_{\GL_{2\RL}}} : & \qquad & \qquad &
\theta  ^{(\cit)}_{\GL_2(F^T_{\o\omega }\times F^T_\omega )} & \To &
\EIS\RL(1,j,m_j)
\psset{nodesep=.5cm}
\everypsbox{\scriptstyle}
\ncline[arrows=->]{1,1}{2,1}>{\tau ^{\rm tor}_{\cit}(\Irr \hol^{(1)}(\theta  ^{(\cit)}_{\GL_2( {F_{\o\omega }}\times  {F_\omega }})}
\ncline[arrows=->]{1,4}{2,4}
\ncline[arrows=->]{1,6}{2,6}
\end{psmatrix}\]
where:
\Bi
\item $\EIS\RL(1,j,m_j)=\EIS_R(1,j,m_j)\times_{(D)}\EIS_R(1,j,m_j)$~, being the global cuspidal representation of the ground bisemisheaf $\theta ^{(\cit)}_{\GL_2(F^T_{\o\omega }\times F^T_\omega )}$ is the product, right by left, of the (truncated) Fourier development of the cusp forms \cite{Pie1}:
\begin{align*}
\EIS_L(1,j,m_j) &= \bigoplus_{j,m_j} \lambda ^{\half}(1,j,m_j)\ e^{2\pi ijz_m}\;,\\
\EIS_R(1,j,m_j) &= \bigoplus_{j,m_j} \lambda ^{\half}(1,j,m_j)\ e^{-2\pi ijz_m}\;,\end{align*}
with $\lambda ^{\half}(1,j,m_j)$ being the square root of the product of the eigenvalues of the coset representative $U_{j,m_{j_R}}\times U_{j,m_{j_L}}$ of the product, right by left, of Hecke operators;

\item $\GL_2(F^T_{\o\omega }\times F^T_\omega )$ is the bilinear algebraic semigroup whose conjugacy class representatives $g^{(2)}_{T_R}[j,m_j] \times g^{(2)}_{T_L}[j,m_j]$ have undergone a toroidal compactification.
\Ei

\item Each monodromy bisemisheaf $\theta ^{(\cit)\rm mon}_{\GL_2(F^T_{\o\omega }\times F^T_\omega )}$~, having been compactified toroidally, gives rise to a similar cuspidal representation 
$\Irr\cusp^{(1)}(\theta ^{(\cit)\rm mon}_{\GL_2(F^T_{\o\omega }\times F^T_\omega )})=\EIS\RL(1,j,m_j)$ since the monodromy bisemisheaves are isomorphic to (or copies of) the ground bisemisheaf according to proposition 5.2.6.\epr
\Ee
\vskip 11pt

\subsubsection{Proposition}

{\em On the desingularized ground bisemisheaf $\theta ^{(\cit)}_{\GL_2(F_{\o\omega }\times F_\omega )}$ and its monodromy bisemisheaves\linebreak $\{\theta ^{(\cit)\rm mon}_{\GL_2(F_{\o\omega }\times F_\omega )}(b)\}^\beta _{b=1}$~, we have the following irreducible Langlands global correspondences:
\Be
\item $ \begin{array}[t]{ccc}
\Irr\Rep^{(1)}_{W_{F\RL}} (W^{ab}_{F_{R}}\times W^{ab}_{F_{L}}) &\To &
\Irr\cusp (\theta ^{(\cit)}_{\GL_2(F^T_{\o\omega }\times F^T_\omega )})\\[11pt]
\mbox{\Huge{$\|$}} && \mbox{\Huge{$\|$}}\\[11pt]
\theta ^{(\cit)}_{\GL_2(F_{\o\omega }\times F_\omega )} &\To & \EIS\RL(1,j,m_j)
\end{array}$

\item $ \begin{array}[t]{ccc}
\Irr\Rep^{(1)}_{W^{\rm mon}_{F\RL}} (W^{ab}_{F_{R_{\rm mon}}}(b)\times W^{ab}_{F_{L_{\rm mon}}}(b)) &\To &
\Irr\cusp (\theta ^{(\cit)\rm mon}_{\GL_2(F^T_{\o\omega }\times F^T_\omega )}(b))\\[11pt]
\mbox{\Huge{$\|$}} && \mbox{\Huge{$\|$}}\\[11pt]
\theta ^{(\cit)\rm mon}_{\GL_2(F_{\o\omega }\times F_\omega )}(b) &\To & \EIS^{\rm mon}\RL(1,j,m_j)_b
\end{array}$
\Ee
}
\vskip 11pt

\bpr These Langlands global correspondences result from the preceding sections summarized in the two following diagrams:
\Be
\item \scalebox{.92}{$ \begin{array}[t]{ccccc}
\Irr\Rep^{(1)}_{W_{F\RL}} (W^{ab}_{F_{R}}\times W^{ab}_{F_{L}}) &\To &
\Irr\hol^{(1)} (\theta ^{(\cit)}_{\GL_2(F_{\o\omega }\times F_\omega )})
& \xrightarrow{\tau ^{\rm tor}_{\cit}} &
\Irr\cusp (\theta ^{(\cit)}_{\GL_2(F^T_{\o\omega }\times F^T_\omega )})\\[11pt]
\mbox{\Huge{$\|$}} && \mbox{\Huge{$\|$}} && \mbox{\Huge{$\|$}}\\[11pt]
\theta ^{(\cit)}_{\GL_2(F_{\o\omega }\times F_\omega )} &\To & 
f_{\o\omega }(z^*_m)\otimes f_\omega (z_m) 
& \xrightarrow{\tau ^{\rm tor}_{\cit}} &
\EIS\RL(1,j,m_j)_b
\end{array}$}

\item \scalebox{.8}{$ \begin{array}[t]{ccccc}
\Irr\Rep^{(1)}_{W^{\rm mon}_{F\RL}} (W^{ab}_{F_{R_{\rm mon}}}(b)\times W^{ab}_{F_{L_{\rm mon}}}(b)) &\To &
\Irr\hol^{(1)} (\theta ^{(\cit)\rm mon}_{\GL_2(F_{\o\omega }\times F_\omega )}(b))
& \xrightarrow{\tau ^{\rm tor}_{\cit}} &
\Irr\cusp (\theta ^{(\cit)\rm mon}_{\GL_2(F^+_{\o\omega }\times F^+_\omega )}(b)\\[11pt]
\mbox{\Huge{$\|$}} && \mbox{\Huge{$\|$}} && \mbox{\Huge{$\|$}}\\[11pt]
\theta ^{(\cit)\rm mon}_{\GL_2(F_{\o\omega }\times F_\omega )}(b) &\To & 
f_{\o\omega }(z^*_{m_b})\otimes f_\omega (z_{m_b}) 
& \xrightarrow{\tau ^{\rm tor}_{\cit}} &
\EIS^{\rm mon}\RL(1,j,m_j)_b
\end{array}$}

\mbox{}\epr
\Ee

\end{document}